\documentclass{amsart}
\usepackage[utf8]{inputenc}
\usepackage{amssymb,amsmath,amsfonts}
\usepackage[cm]{fullpage}
\usepackage{wrapfig}
\usepackage{graphicx}
\usepackage{courier}
\usepackage{caption}
\usepackage{subcaption}
\usepackage{float}
\usepackage{yfonts}
\usepackage{mathrsfs}
\usepackage{mathtools}
\usepackage{color}
\usepackage{colortbl}
\usepackage{tikz}
\usepackage{makecell}
\usepackage{multirow}
\usepackage{multicol}
\usepackage{enumitem} \setlist[enumerate]{leftmargin=.5in} \setlist[itemize]{leftmargin=.5in}
\usepackage{hyperref}
\usepackage{subfiles} 
\usepackage[T1]{fontenc}
\usepackage{scalerel,stackengine}
\usepackage[style=numeric,backend=biber]{biblatex}
\newtheorem{Lemma}{Lemma}[section]
\newtheorem{Theorem}{Theorem}
\newtheorem{Proposition}[Lemma]{Proposition}

\newtheorem{Remark}[Lemma]{Remark}

\newtheorem{Criteria}[Lemma]{Criteria}
\newtheorem{Criterion}[Lemma]{Criterion}

\usepackage{todonotes}
 
\newcommand{\revised}[1]{\textcolor{black!100}{#1}}
\newcommand{\uvec}[1]{\boldsymbol{\hat{\textbf{#1}}}}

\def\uVecOne{\uvec{e}^{\scriptscriptstyle 1}}
\def\uVecTwo{\uvec{e}^{\scriptscriptstyle 2}}
\def\uVecThree{\uvec{e}^{\scriptscriptstyle 3}}
\def\uVecJ{\uvec{e}^{\scriptscriptstyle j}}
\def\vorticity{\boldsymbol{\omega}}

\definecolor{Gray}{gray}{0.9}

\def\Ray{\mathsf{R}}
\def\Pra{\mathsf{P}}
\def\Rot{\mathsf{S}}
\def\Sha{\mathsf{k}_1}
\def\parameters{\boldsymbol{\mathfrak{P}}}

\graphicspath{ {./pics/} } 

\title{Rotating Rayleigh-B\'enard convection: Attractors, bifurcations and heat transport via a Galerkin hierarchy}

\author[F1.RW]{Roland Welter$^1$}
\address{$^1$Universität Hamburg\\
	Fachbereich Mathematik \\
	Bundesstraße 55 \\
	20146 Hamburg, DE}
\email{roland.welter@uni-hamburg.de}


\begin{document}

	\maketitle
	
	\begin{abstract}
		Motivated by the need for energetically consistent climate models, the Boussinessq-Coriolis (BC) equations are studied with a focus \revised{on the averaged vertical heat transport, ie the Nusselt number.  A set of formulae are derived by which arbitrary Fourier truncations of the BC model can be explicitly generated, and Criteria are given which precisely guarantee that such truncated models obey energy relations consistent with the PDE.  The Howard-Krishnamurti-Coriolis (HKC) hierarchy of such energetically consistent ODE models is then implemented in MATLAB, with code available on GitHub.  Several theoretical results are proven to support a numerical analysis.  Well-posedness and convergence of the HKC hierarchy toward the BC model are proven, as well as the existence of an attractor for the BC model.  Since the rate of convergence is unknown, explicit upper and lower bounds on the attractor dimension are proven so as to provide guidance for the required spatial resolution for an accurate approximation of the Nusselt number.  Finally, a series of numerical studies are performed using MATLAB, which investigate the required spatial resolution and indicate the presence of multiple stable values of the Nusselt number,} setting the stage for an energetically consistent analysis of convective heat transport. 
	\end{abstract}
	
	\section{Introduction}
	
	In the study of climate, general circulation models aim to accurately represent phenomena in the atmosphere and ocean using equations from fluid dynamics.   In order to be computationally tractable, these models must make approximations and, in so doing, can violate energy conservation and other physical principles \cite{eden2019energy}.  The main goal of this paper is to develop a mathematically rigorous framework for studying vertical heat transport in atmospheric convection via a combination of theoretical and numerical results, with a particular focus on energetic consistency.
	
	The models considered in this paper are the three-dimensional, non-dimensionalized Boussinesq Oberbeck equations with a Coriolis force
	\begin{equation}
		\label{BoussinesqCoriolis}  \begin{split}
			\frac{1}{\Pra} \Big [ \partial_{t} \textbf{u} + \textbf{u} \cdot \nabla \textbf{u} \Big ] + \nabla p & = \Delta \textbf{u} + \pi \Ray T \uVecThree - \Rot \hspace{.5 mm} \uVecThree \times \textbf{u} , \\
			\nabla \cdot \textbf{u} & = 0 , \\
			\partial_t T + \textbf{u}\cdot \nabla T & =  \Delta T ,
	\end{split} \end{equation}
	and its Galerkin truncations.  Here $\textbf{u} = (u_1,u_2,u_3)$, $T$, $p$ are the non-dimensionalized velocity, temperature and pressure of a fluid lying in a box $\textbf{x} = (x_1,x_2,x_3) \in [0,\frac{2\pi}{\Sha}]\times[0,\frac{2\pi}{\mathsf{k}_2}] \times [0,\pi]$ for shape parameters $\mathsf{k}_i$, and $\uvec{e}^{\scriptscriptstyle i}$ denotes the unit vector in the $i^{th}$ direction.  The fluid is assumed to be periodic in $x_1,x_2$ and it is assumed that the fluid configuration is horizontally aligned, ie it is independent of the variable $x_2$.  Due to this condition it suffices to consider the behavior on the plane $x_2 =0$, hence the domain is taken to be $\Omega = [0,\frac{2\pi}{\Sha}] \times [0,\pi]$.  The fluid is assumed to be heated from below and cooled from above, hence the non-dimensionalized temperature $T$ must satisfy  the boundary conditions
	\begin{equation} \label{BC_Temp} T(x_1,0,t) = 1 \hspace{1 cm} \text{ , } \hspace{1 cm} T(x_1,\pi,t) = 0. \end{equation}
	Furthermore, the fluid is assumed to satisfy the impenetrability and free slip conditions at the top and bottom boundaries
	\begin{equation} \label{BC_Velocity} u_3 = 0 \hspace{1 cm} \text{ , }  \hspace{1 cm} \partial_{x_3} u_1 = \partial_{x_3} u_2 = 0  \hspace{2 cm} \text{ for } \hspace{.5 cm} x_3 = 0,1. \end{equation}
	The parameters $\Pra, \Ray$ are the Prandtl and Rayleigh numbers, whereas the parameter $\Rot$, referred to as the "rotation number", is proportional to the angular speed of the rotating frame.  
	
	This model was chosen because it is a simple model which still captures many features of atmospheric convection.  Vertical heat transport is of interest, so no hydrostatic assumption is made.  Rotation plays an important role, so the velocity field must be 3d.  The fully 3d case still presents intractable theoretical problems, so the horizontal alignment assumption is made to provide a nice analytic setting.  This assumption is invariant under the full 3d evolution, though such solutions will most often be dynamically unstable.  By inserting approximate values relevant for atmospheric convection as in \ref{app:PhysParam} one finds that the parameter values should be roughly as follows:
	\begin{equation} \label{ParameterVals} \Pra \approx 1 \hspace{1 cm} \text{ , } \hspace{1 cm}  \Ray \approx 10^{16} \hspace{1 cm} \text{ , } \hspace{1 cm} \Rot \approx \left \{ \begin{array}{ll}
			0  & \text{ near the equator},  \\
			10^{15} & \text{ near the poles} .
		\end{array} \right . \end{equation}
	Note that the rotation number $\Rot$ rather than the Rossby number is used in this paper, since the case $\Rot = 0$ is a relevant parameter value for geophysical flows.  Finally, the free-slip condition implies that the fluid can slip along the boundary with zero viscosity.  Despite the small viscosity this is somewhat unphysical, but it simplifies the analysis considerably.    
	
	As mentioned, the main quantity of interest is the average vertical heat transport, as measured by the Nusselt number.  Using the notation $\langle \cdot \rangle$ to denote the volume integral and an overbar to denote an infinite time average
	\[ \langle f \rangle := \int_{\Omega} f(\textbf{x}) d\textbf{x} \hspace{.5 cm} \text{ , } \hspace{.5 cm} \bar{f} := \lim_{t \to \infty} \frac{1}{t} \int_0^t f(s) ds , \]
	the Nusselt number is then defined as the vertical heat flux averaged over time and space:
	\begin{equation} \label{NusseltDef} \mathsf{Nu} := 1 + \frac{\Sha}{2\pi^2} \overline{\langle u_3 T \rangle} .  \end{equation}
	The Nusselt number depends in a complicated way on the parameters, $\Ray,\Pra,\Rot, \Sha$ and on the initial states $(\textbf{u}_0,T_0)$, and this dependence is analyzed here.  An important aspect for this analysis is energetic consistency.  Energy conservation is not to be expected in a model such as \eqref{BoussinesqCoriolis}, but rather the kinetic and potential energy should obey the following definite rules: 
	\begin{equation} \label{EnergyEvolution} \frac{d}{dt} \langle \frac{1}{2}|\textbf{u}|^2 \rangle = - \Pra \langle |\nabla \textbf{u} |^2 \rangle + \pi \Pra \Ray \langle u_3 T \rangle \hspace{.25 cm} \text{ , } \hspace{.25 cm} \frac{d}{dt} \langle (1-\frac{x_3}{\pi}) T \rangle = - \frac{1}{\pi} \langle u_3 T \rangle + \langle (1-\frac{x_3}{\pi}) \partial_{x_3}^2 T \rangle . \end{equation}
	These will be referred to as the balance equations for kinetic and potential energy, and a numerical approximation of \eqref{BoussinesqCoriolis} will be said to be energetically consistent if these are satisfied.  While there are abundant numerical studies regarding the Nusselt number \cite{wen_goluskin_doering_2022, WenDingChiniKerswell2022, iyer2020classical, StevensClercxLohse_2010}, the question of energetic consistency, and its influence on the approximate Nusselt number, remains relatively unaddressed.  Another important aspect for the present analysis concerns how much spatial and temporal resolution are required for an accurate representation of the flow, which also influence the approximate Nusselt number.  Finally, "stable" Nusselt number values, i.e. those realized by solutions with dynamic stability, are of particular interest, since these are most relevant for climate applications. 
	
\begin{figure}[H]
		\begin{center}
			\begin{tabular}{ccc}
				\includegraphics[height=50mm]{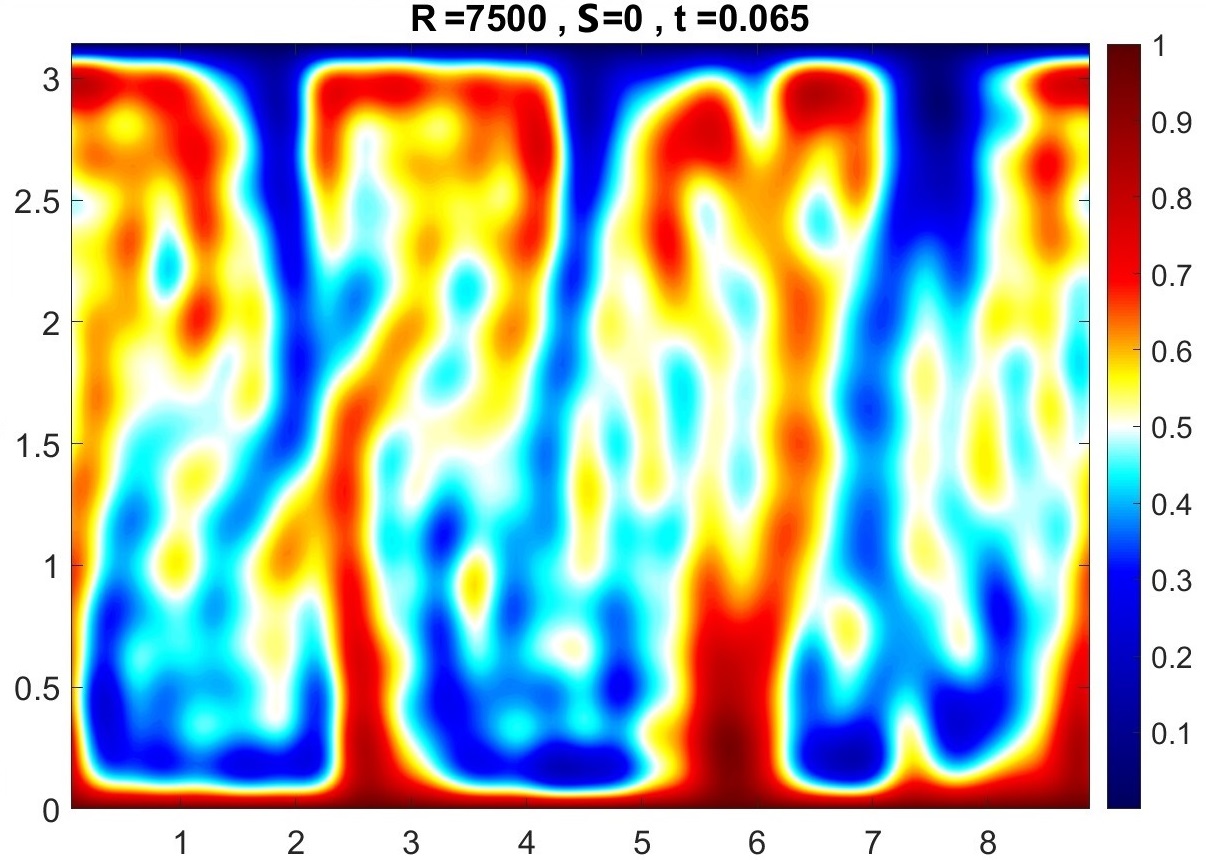} & \hspace{.05 cm} & \includegraphics[height=50mm]{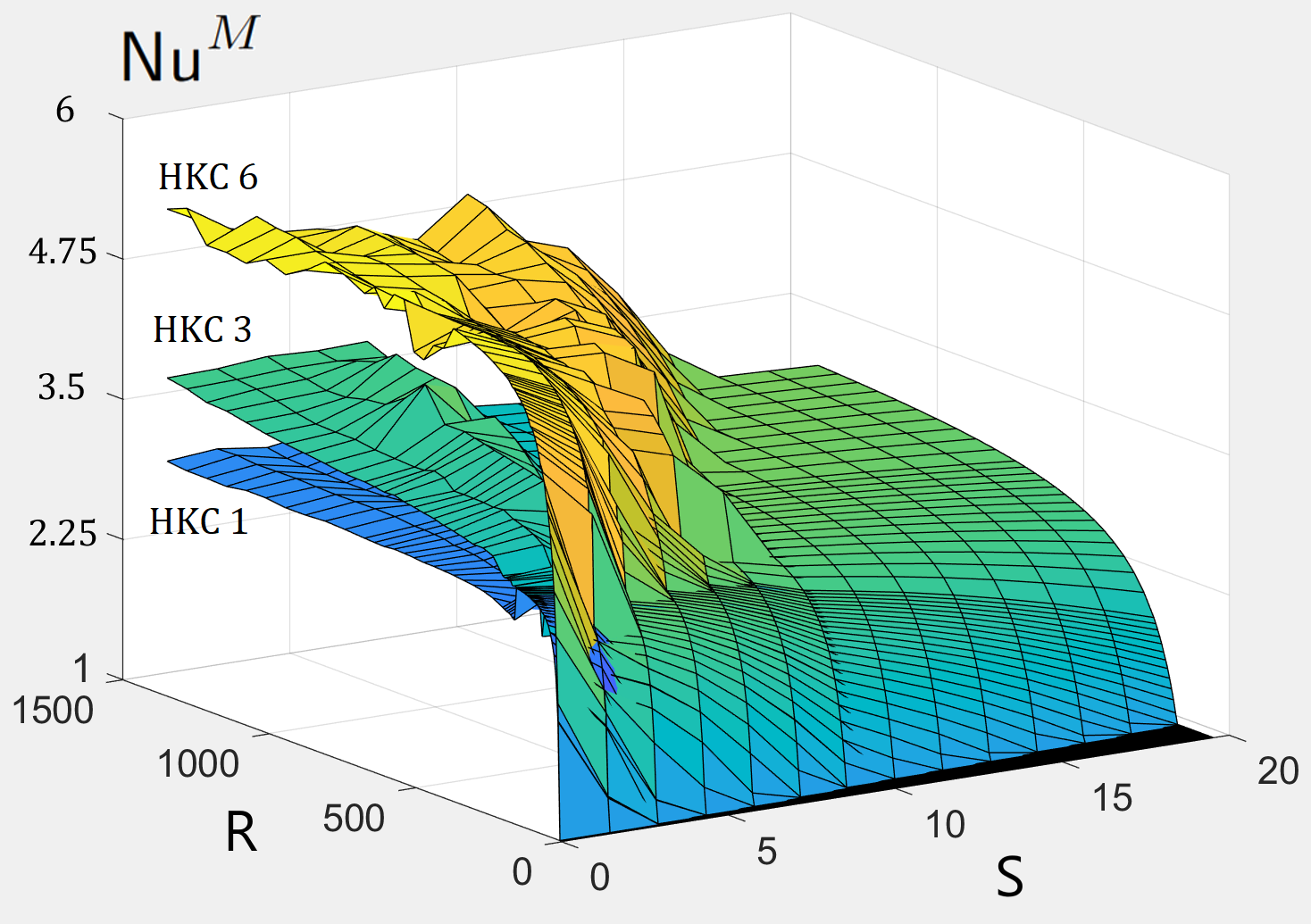}   \\
				(a) & \hspace{.05 cm} & (b) 
			\end{tabular}
			\caption{For $\Pra =10$, $\Sha = \frac{1}{\sqrt{2}}$, (a) The temperature field of the HKC 190 model with $\Ray = 7500, \Rot = 0$ and (b) Nusselt number approximations computed over $0 \leq \Ray \leq 1500$ and $0 \leq \Rot \leq 20$ using the HKC 1, HKC 3 and HKC 6 models. }
			\label{fig:FluidVisualization}
		\end{center}
	\end{figure}

	In this paper a hierarchy of truncated ODE models, referred to as the HKC hierarchy, is developed and it is proven analytically that each model in the hierarchy is energetically consistent.  One ascends the HKC hierarchy by adjoining more Fourier modes, obtaining higher spatial resolution by considering a higher dimensional ODE.  The HKC hierarchy is then implemented in MATLAB and used to compute approximations of $\mathsf{Nu}$ as a function of the Rayleigh number and rotation numbers.  The results of these numerical studies provide strong evidence of multi-stability of $\mathsf{Nu}$, suggest that the required spatial resolution increases as $\Ray^{1/2}$, and in general map out the dependence of $\mathsf{Nu}$ on $\Ray$ and $\Rot$.  Figure \ref{fig:FluidVisualization} displays the temperature field of a fluid flow computed from one of the models in this hierarchy, and different Nusselt number values computed using HKC models of increasing spatial resolution.
	
	It can be difficult to discern the sources and sizes of error only by comparing numerical methods to each other.  Most of this paper focuses on theoretical results to set the stage for an unambiguous intercomparison.  Beyond the energetic consistency mentioned above, these include well-posedness, regularity and the existence of a compact attractor for \eqref{BoussinesqCoriolis}, estimates on the attractor dimension, and a bifurcation analysis.  Indeed, while the well-posedness theory provides convergence results showing that the Nusselt number can in principle be computed, such theorems fall short of providing concrete information about the rate of convergence.  Thus bounds on the attractor dimension are obtained to address the required spatial resolution.
	
	Thus the structure of the paper can  be summarized as follows.  In section \ref{sec:HKCAnalysis} the HKC hierarchy is defined and energetic consistency is proven analytically.  In section \ref{sec:BCAnalysis} well-posedness, regularity and dynamical properties of the Boussinesq Coriolis model \eqref{BoussinesqCoriolis} are established analytically.  In section \ref{sec:AttractorAnalysis} the attractor dimension of \eqref{BoussinesqCoriolis} is studied, wherein a lower bound is obtained by studying the local bifurcations occurring at the origin and an upper bound is given by the analogue of Temam's result.  Section \ref{sec:NumAnalysis} presents the results of several numerical studies into heat transport that address various numerical issues and interesting dynamical behaviors.  In order to make this work more accessible, code for generating these models in \LaTeX, MATLAB and Mathematica was written, and can be found on GitHub at \url{https://github.com/rkwelter/HKC_CodeRepo}.  Finally, these results are interpreted in a larger context in section \ref{sec:Conclusions}.
	
	\subsection{Notation and basic definitions}
	
	In this paper, scalars are written in normal font, vectors $\textbf{v}$ in bold, matrices $\mathcal{M}$ and operators in calligraphic, parameters $\mathsf{P}$ in serif and sets $\mathscr{S}$ in script.  Generally subscripts will be used to denote the component of a vector or matrix, and superscripts will be used to indicate indexed quantities.  There will be a couple of exceptions to this rule, namely that some quantities will be indexed by the symbols $\textbf{u}, \theta$ as subscripts to indicate a different definition for the two different variables, and the subscript zero (e.g. $\textbf{u}_0$) will be used to indicate initial conditions.  This should not cause confusion, since these could not indicate a component of a vector or matrix.
	
	For a Banach space $B$, $\| \cdot \|_{B}$ will denote the norm, $B'$ the dual and $(\cdot, \cdot)_{B'}$ the dual pairing.  When $B$ is a Hilbert space $(\cdot, \cdot )_X$ denotes the inner product.  $C^{k}$, $W^{k,p}$, $H^k = W^{k,2}$ and $L^p$ will denote the space of $k$ times differentiable functions, the Sobolev spaces, the Hilbert Sobolev spaces and the Lebesgue spaces respectively.  For $k < 0$, $H^k$ is the well-known negative Sobolev space.  Unless otherwise specified these will denote functions defined on $\Omega$; when referring to functions defined on another domain $\Omega'$, the notation $C^k(\Omega')$ and so on will be used.  All functions defined on $\Omega$ are assumed periodic in $x_1$.  For $k \geq 0$ let $C_0^{k}$ be the corresponding subspace with zero trace, and for $k > 0$ let $W_0^{k,p}$, $H_0^k$ be the closure of $C^k_0$ in $W^{k,p}$, $H^k$.  Similarly let $\textbf{C}^{k}_{\sigma}$ be the space of incompressible vector fields $\textbf{u} = (u_1,u_2,u_3)^T$ with components in $C^{k}$ and satisfying the boundary conditions \eqref{BC_Velocity}, and let $\textbf{W}^{k,p}_{\sigma}$, $\textbf{H}^{k}_{\sigma}$, $\textbf{L}^{p}_{\sigma}$ denote the closure of $\textbf{C}^{\infty}_{\sigma}$ in $W^{k,p}$, $H^k$ and $L^p$, respectively.  Since one can translate to a moving frame if necessary, without loss of generality one can take these vector valued function spaces such that $u_1$, $u_2$ have zero mean.  We will use the following notation to denote the product spaces
	\[  \textbf{H}^0 := \textbf{L}_{\sigma}^2 \times L^2 \hspace{.5 cm} \text{ , } \hspace{.5 cm} \textbf{H}^k := \textbf{H}^k_{\sigma} \times H_0^k \hspace{.5 cm} \text{ , } \hspace{.5 cm} \textbf{C}^k := \textbf{C}_{\sigma}^k \times C^k_0 \text{ . } \]
	For a Banach space $B$, let $C^{k}( [0,T]; B)$ denote the space of mappings $f:[0,T]\mapsto B$ such that $t \mapsto \|f(t)\|_B$ belongs to $C^k([0,T])$, and so on for the spaces $W^{k,p}( [0,T]; B)$, $H^k( [0,T]; B)$ and $L^p( [0,T]; B)$.  For a semi-group of operators $\mathcal{S}(t):[0,\infty) \times B \mapsto B$, a subset $Y \subset B$ is said to be functionally invariant under $\mathcal{S}(t)$ if $\mathcal{S}(t) Y = Y$ for all $t > 0$, and is said to attract another set $Z \subset B$ if for all $z_0 \in Z$
	\[ \lim_{t \to \infty} \inf_{y\in Y} \|\mathcal{S}(t)z_0 - y\|_{B} = 0 . \]
	A set $\mathscr{A} \subset B$ is the global attractor for $\mathcal{S}(t)$ if it is the maximal functional invariant set which attracts all bounded subsets of $B$.  
	
	For a vector $\textbf{v} \in \mathbb{R}^d$ the notation $|\textbf{v}|_{p} = (\sum_i v_i^p)^{1/p}$ will be used and $|\textbf{v}| = |\textbf{v}|_{2}$.  $\mathcal{K} = \mathsf{diag}(\Sha,1)$ will denote the $2\times 2$ diagonal matrix with $\Sha,1$ on the diagonal and $\eta^{\textbf{m}}, V$ will denote the normalizing constants defined by
	\begin{equation} \label{Normalizers} \eta^{\textbf{m}} := \eta^{m_1}\eta^{m_3} , \hspace{1 cm} \eta^m := \left \{ \begin{array}{ll} 1 & \text{ if } m > 0 , \\ \frac{1}{\sqrt{2}} & \text{ if } m = 0 , \end{array} \right . \hspace{1 cm} V := \Big ( \frac{\pi^2}{2\Sha} \Big )^{1/2} . \end{equation}
	We will define $\delta$ to be an indicator function as follows:
	\[ \delta^{i,j} := \left \{ \begin{array}{cc} 1 & \text{ if } i=j , \\ 0 & \text{ if } i\neq j , \end{array} \right .  \hspace{1 cm} \delta^{\textbf{v},\tilde{\textbf{v}}} := \Pi_{i=1}^d  \delta^{v_i,\tilde{v}_i} . \]
	$\parameters = (\Ray,\Rot,\Pra, \Sha)$ will denote the vector of parameters for \eqref{BoussinesqCoriolis}, and $\parameters$ will be said to be admissible if $\Ray, \Rot \geq 0$, $\Pra, \Sha > 0$.
	
	\subsection{Reformulation as an evolution equation}
	
	It is desirable to write the system \eqref{BoussinesqCoriolis} as an evolution equation on a function space.  To this end, define $\theta, p_{dev}$ to be the temperature and pressure deviations from the pure conduction profiles:
	\[ \theta = \pi (T - 1) + x_3 \hspace{1 cm} \text{ , } \hspace{1 cm} p_{dev} = p - \pi \Ray (x_3 - \frac{1}{2\pi}x_3^2) . \]
	Inserting these into (\ref{BoussinesqCoriolis}) and applying the (Leray) projector, $\mathcal{P}$, onto divergence-free vector fields, one obtains the following evolution equation:
	\begin{equation} \label{EvolutionPDE} \begin{split}
			\partial_{t} \textbf{u} & = \Pra \Delta \textbf{u} + \mathcal{P} \big [ \Pra \Ray \theta \uVecThree - \Pra \Rot \hspace{.5 mm} \uVecThree \times \textbf{u} - \textbf{u} \cdot \nabla \textbf{u} \big ] , \\
			\partial_t \theta & =  \Delta \theta + u_3 - \textbf{u} \cdot \nabla \theta.
	\end{split}    \end{equation}
	The velocity field $\textbf{u}$ should have finite kinetic energy, hence a natural choice of function space for the velocity is $\textbf{L}^2_{\sigma}$, and the function $\theta$ is equal to zero at both boundaries, hence a natural choice of function space is $H^1_0$.
	
	\section{The HKC hierarchy of rotating convection models} 
	
	\label{sec:HKCAnalysis}
	
	\subsection{Generic reduced order rotating convection models}
	
	As described above, the first goal of this paper is to develop a hierarchy of ODE models by expanding the variables $\textbf{u},\theta$ in a Fourier series and projecting \eqref{EvolutionPDE} onto the corresponding Fourier coefficients.  In this section, we lay out the fundamentals of the projection process by which a generic Fourier truncation results in an ODE system.  We start by looking for a general expansion of the form 
	\begin{equation} \label{GeneralFourierExpansion} \begin{split} \textbf{u}(\textbf{x},t) = \sum_{\textbf{n}} u^{\textbf{n}}(t) \textbf{v}^{\textbf{n}}(\textbf{x}) \hspace{.5 cm} \text{ , } \hspace{.5 cm} \theta(\textbf{x},t) = \sum_{\textbf{n}} \theta^{\textbf{n}}(t) f^{\textbf{n}}(\textbf{x}) . \end{split} \end{equation}
	where $u^{\textbf{n}}(t),\theta^{\textbf{n}}(t)$ are Fourier coefficients corresponding to some sinusoidal functions $\textbf{v}^{\textbf{n}}(\textbf{x})$, $f^{\textbf{n}}(\textbf{x})$ yet to be defined.  One could naively choose the functions $\textbf{v}^{\textbf{n}}(\textbf{x})$, $f^{\textbf{n}}(\textbf{x})$ to be complex exponentials typical for Fourier expansions, but eventually the nonlinear analysis becomes much messier due to the complicated implicit relationships between the Fourier coefficients required to enforce the real-valued condition, boundary conditions and divergence free conditions.  Since we avoid complex exponentials, let us briefly motivate the choice of the definitions of $\textbf{v}^{\textbf{n}}(\textbf{x})$, $f^{\textbf{n}}(\textbf{x})$ by recalling some basic facts about Fourier expansions in one dimension.  In this case a real-valued, periodic function $f(x_1)$ defined on $[0,2\pi]$ can be expanded in terms of Fourier coefficients $\hat{f}^{m_1,p_1}$ via
	\[ f(x_1) = \frac{1}{(2\pi)^{1/2}} \hat{f}^{0,1} + \sum_{m_1 \in \mathbb{Z}_{>0}} \hat{f}^{m_1,1} \frac{\cos(m_1 x_1)}{(\pi)^{1/2}} + \hat{f}^{m_1,2} \frac{\sin(m_1 x_1)}{(\pi)^{1/2}} . \]
	Here we see that using sinusoids requires that the Fourier coefficients be indexed by two quantities, namely a wave number $m_1$ and a phase $p_1$, describing whether the coefficient corresponds to cosine or sine.  Furthermore, the number of phases depends on the wave number, e.g. for $m_1 = 0$ there is only one phase, corresponding to $\cos(0) = 1$, whereas $\sin(0) = 0$ is not included in the expansion.  Finally, for a Fourier expansion of a vector field one must additionally specify which component of the vector is being expanded, which is done here by using a "component" index $c$.
	
	With this in mind, the index $\textbf{n} = (m_1,m_3,p_1,c)$ will consist of a wave vector $\textbf{m}=(m_1,m_3)$, a phase index $p_1$ and a component index $c$.  Due to the horizontal alignment condition, $m_2$ is always zero and hence excluded.  We will slightly abuse the notation and write $\textbf{n} = (\textbf{m},p_1,c)$ at times.  Since $\theta$ is a scalar one can just use the usual Fourier expansion satisfying the Dirichlet boundary conditions.  Hence define the functions $f^{\textbf{n}}$ via 
	\begin{equation} \label{ThetaFourierBasis} \begin{split} & f^{(\textbf{m},1,c)} := \frac{\eta^{\textbf{m}}}{V}  \cos(\Sha m_1 x_1 ) \sin ( m_3 x_3 ) \hspace{.5 cm} \text{ , } \hspace{.5 cm} f^{(\textbf{m},2,c)} := \frac{\eta^{\textbf{m}}}{V} \sin(\Sha m_1 x_1 ) \sin ( m_3 x_3 ) \text{ . } \end{split} \end{equation}
	Since $\theta$ is a scalar only $c = 1$ is defined, but this gives a unified notation with the velocity field.  The expansion for $\textbf{u}$ is more complicated since $\textbf{u}$  is a divergence free vector field.  In this case, one can expand either the first or second component of the vector field, and the third is determined by the divergence free condition.  Hence define a family of vector fields 
	\begin{align}  \textbf{v}^{(\textbf{m},1,1)} := \frac{\eta^{\textbf{m}}}{|\mathcal{K}\textbf{m}| V} & \begin{pmatrix} m_3 \sin( \Sha m_1 x_1 ) \cos( m_3 x_3 ) \\ 0 \\ - \Sha m_1 \cos( \Sha m_1 x_1 ) \sin( m_3 x_3 ) \end{pmatrix} \text{ , } \notag \\ \label{VectorFieldDef_FreeSlip} \textbf{v}^{(\textbf{m},2,1)} := \frac{\eta^{\textbf{m}}}{|\mathcal{K}\textbf{m}| V} &  \begin{pmatrix} m_3 \cos( \Sha m_1 x_1 ) \cos( m_3 x_3 ) \\ 0 \\ \Sha m_1 \sin( \Sha m_1 x_1 ) \sin( m_3 x_3 ) \end{pmatrix} ,  \\ \textbf{v}^{(\textbf{m},1,2)} := \frac{\eta^{\textbf{m}}}{V} \uVecTwo \sin( \Sha m_1 x_1 ) \cos( & m_3 x_3 ) \hspace{.25 cm}  \text{ , } \hspace{.25 cm} \textbf{v}^{(\textbf{m},2,2)} :=  \frac{\eta^{\textbf{m}}}{V} \uVecTwo \cos( \Sha m_1 x_1 ) \cos( m_3 x_3 ) . \notag \end{align}
	One can quickly check that these vector fields are divergence free and satisfy the boundary conditions.  One can also check that the following orthonormality properties are satisfied:
	\begin{equation}
		\label{OrthonormalityProperties}
		\langle \textbf{v}^{\textbf{n}} \cdot \textbf{v}^{\tilde{\textbf{n}}} \rangle = \delta^{\textbf{n},\tilde{\textbf{n}}} \hspace{.5 cm} \text{ , } \hspace{.5 cm} \langle f^{\textbf{n}} f^{\tilde{\textbf{n}}} \rangle = \delta^{\textbf{n},\tilde{\textbf{n}}} .
	\end{equation}
	It can also be checked that the collection \eqref{VectorFieldDef_FreeSlip} is complete and hence nothing is missed by the expansion \eqref{GeneralFourierExpansion}.  This is slightly more involved and a distraction from the present goal, hence is carried out in \ref{app:FourierExpDeriv}.
	
	As above, some of the functions $\textbf{v}^{\textbf{n}}(\textbf{x}),f^{\textbf{n}}(\textbf{x})$ can equal zero due to the fact that $\sin(0) = 0$.  We therefore define the admissible wave, phase and component index sets by simply inspecting \eqref{ThetaFourierBasis},\eqref{VectorFieldDef_FreeSlip} for sine:
	\begin{align} \label{PhaseIndexSets}
		\mathscr{M}_{\theta} = \mathbb{Z}_{\geq 0} \times \mathbb{Z}_{> 0} \hspace{.3 cm} & \text{ , } \hspace{.3 cm} \mathscr{M}_{\textbf{u}} = \mathbb{Z}_{\geq 0}^2 \setminus \{ \textbf{0} \} \text{ , } \hspace{.3 cm} \\ \mathscr{P}^{\textbf{m}} = \left \{ \begin{array}{cl}
			\{ 1,2 \} & \text{ if } m_1 > 0, \\
			\{ 1 \} & \text{ otherwise, } 
		\end{array} \right . \hspace{.25 cm} & \mathscr{C}^{\textbf{m}}_{\textbf{u}} = \left \{ \begin{array}{cl}
			\{ 1,2 \}  & \text{ if } m_3 > 0, \\
			\{ 2 \}  & \text{ if } m_3 = 0 \text{ . } 
		\end{array} \right . \notag
	\end{align}
	Define also 
	\[ \mathscr{N}_{\theta} = \{ \textbf{n} : \textbf{m} \in \mathscr{M}_{\theta} \text{ , } p_1 \in \mathscr{P}^{\textbf{m}} \text{ , } c = 1 \} \hspace{.25 cm} \text{ , } \hspace{.25 cm} \mathscr{N}_{\textbf{u}} = \{ \textbf{n} : \textbf{m} \in \mathscr{M}_{\textbf{u}} \text{ , } p_1+1 \in \mathscr{P}^{\textbf{m}} \text{ , } c \in \mathscr{C}_{\textbf{u}}^{\textbf{m}} \} . \]
	Henceforth, whenever referring to $\textbf{v}^{\textbf{n}}(\textbf{x}),f^{\textbf{n}}(\textbf{x})$ or the corresponding coefficients it is always assumed that $\textbf{n} \in \mathscr{N}_{\textbf{u}}$ or $\textbf{n} \in \mathscr{N}_{\theta}$ respectively.  Furthermore, any calculation involving the phase index $p_1$ and component $c$ is taken to be mod $2$ (e.g. $u^{(\textbf{m},1,3)} = u^{(\textbf{m},1,1)}$).  Note $\textbf{m} = \textbf{0}$ is not included for either variable: for $\theta$ it is excluded due to the boundary conditions, and without loss of generality it is excluded for $\textbf{u}$ due to the Galilean invariance of \eqref{BoussinesqCoriolis}. 
	
	The above basis can be used to construct reduced order models via Galerkin truncation of the full PDE model.  In general, this can be done by considering expansions of the form
	\begin{equation} \label{FiniteFourierExpansion} \textbf{u}^M(\textbf{x},t) = \sum_{\textbf{n} \in \mathscr{N}_{\textbf{u}}^M} u^{\textbf{n},M}(t) \textbf{v}^{\textbf{n}}(\textbf{x}) \hspace{.5 cm} \text{ , } \hspace{.5 cm} \theta^M(\textbf{x},t) = \sum_{\textbf{n} \in \mathscr{N}_{\theta}^M} \theta^{\textbf{n},M}(t) f^{\textbf{n}}(\textbf{x}) . \end{equation}
	in which $\mathscr{N}^{M}_{\textbf{u}} \subseteq \mathscr{N}_{\textbf{u}}$, $\mathscr{N}^{M}_{\theta} \subseteq \mathscr{N}_{\theta}$ are some finite index sets.  Define the projection operators $\mathcal{P}_{\textbf{u}}^M,\mathcal{P}_{\theta}^M$ associated to these index sets via
	\[ \mathcal{P}_{\textbf{u}}^M \big [ \textbf{v} \big ] = \sum_{\textbf{n} \in \mathscr{N}_{\textbf{u}}^M} \langle \textbf{v} \cdot \textbf{v}^{\textbf{n}} \rangle \textbf{v}^{\textbf{n}}(\textbf{x}) \hspace{.5 cm} \text{ , } \hspace{.5 cm} \mathcal{P}_{\theta}^M \big [ g \big ] = \sum_{\textbf{n} \in \mathscr{N}_{\theta}^M} \langle g f^{\textbf{n}} \rangle f^{\textbf{n}}(\textbf{x}) . \] 
	One can consider the reduced problem obtained by projecting of the full problem \eqref{EvolutionPDE} onto this basis, as follows
	\begin{equation} \label{TruncatedEvolutionPDE} \begin{split}
			\partial_{t} \textbf{u}^M & = \Pra \Delta \textbf{u}^M + \Pra \Ray \mathcal{P}^M_{\textbf{u}} \big [ \theta^M \uVecThree \big ] - \Pra \Rot \mathcal{P}^M_{\textbf{u}} \big [ \uVecThree \times \textbf{u}^M \big ] - \mathcal{P}^M_{\textbf{u}} \big [ \textbf{u}^M \cdot \nabla \textbf{u}^M \big ] , \\
			\partial_t \theta^M & =  \Delta \theta^M + \mathcal{P}^M_{\theta} \big [ u_3^M \big ] - \mathcal{P}^M_{\theta} \big [ \textbf{u}^M \cdot \nabla \theta^M \big ].
	\end{split} \end{equation}
	This truncated PDE can also be viewed as a system of ODE's for the finite set of Fourier coefficients, where the equation for $u^{\textbf{n}}$ can be recovered by computing the inner product of the first equation with $\textbf{v}^{\textbf{n}}$, and the equation for $\theta^{\textbf{n}}$ can be recovered by computing the inner product of the second with $f^{\textbf{n}}$.  For each $\textbf{n} \in \mathscr{N}_{\textbf{u}}^{M}$ one thereby obtains an equation of the following form:
	\begin{equation} \label{GeneralBoussinesqODE_Vel}  \frac{d}{dt} u^{\textbf{n},M} = - \Pra |\mathcal{K}\textbf{m}|^2 u^{\textbf{n},M} + (-1)^{p_1} \Pra \Ray \frac{\Sha m_1}{|\mathcal{K}\textbf{m}|} \delta^{c,1} \theta^{\textbf{n},M} - (-1)^c \frac{\Pra \Rot m_3}{|\mathcal{K}\textbf{m}|} u^{(\textbf{m},p_1,c+1),M} - N^{\textbf{n},M}_{\textbf{u}} , \end{equation}
	and for each $\textbf{n} \in \mathscr{N}_{\theta}^{M}$ one obtains an equation of the following form:
	\begin{equation} \label{GeneralBoussinesqODE_Temp} \frac{d}{dt} \theta^{\textbf{n},M} = - |\mathcal{K}\textbf{m}|^2 \theta^{\textbf{n},M} + (-1)^{p_1} \frac{\Sha m_1}{|\mathcal{K}\textbf{m}|} u^{\textbf{n},M} - N^{\textbf{n},M}_{\theta} . \end{equation}
	The linear terms in the system \eqref{GeneralBoussinesqODE_Vel} - \eqref{GeneralBoussinesqODE_Temp} are obtained easily from the relations
	\begin{equation} \label{LinearBasisRelations} \langle v^{\textbf{n}}_3 f^{\tilde{\textbf{n}}} \rangle = \frac{(-1)^{p_1}\Sha m_1}{|\mathcal{K}\textbf{m}|} \delta^{\textbf{m},\tilde{\textbf{m}}} \delta^{p_1,\tilde{p}_1} \delta^{c,1} \text{ , } \hspace{.15 cm} \langle \textbf{v}^{\textbf{n}} \cdot \big ( \uVecThree \times \textbf{v}^{\tilde{\textbf{n}}}  \big ) \rangle = \frac{(-1)^{c}m_3}{|\mathcal{K}\textbf{m}|} \delta^{\textbf{m},\tilde{\textbf{m}}} \delta^{p_1,\tilde{p}_1} \delta^{c+1,\tilde{c}}, \end{equation}
	which can be established from the explicit formulas for $\textbf{v}^{\textbf{n}}$, $f^{\textbf{n}}$.  Abstractly, the nonlinear terms are given by 
	\begin{equation} \label{AbstractNonlinear} N_{\textbf{u}}^{\textbf{n},M} = \int_{\Omega} \textbf{v}^{\textbf{n}} \cdot \big [ \big ( \textbf{u}^M \cdot \nabla \big ) \textbf{u}^M \big ] d\textbf{x} \hspace{.5 cm} \text{ , } \hspace{.5 cm} N_{\theta}^{\textbf{n},M} = \int_{\Omega} f^{\textbf{n}} \cdot \big [ \textbf{u}^M \cdot \nabla \theta^M \big ] d\textbf{x} . \end{equation}
	However, to write \eqref{GeneralBoussinesqODE_Vel} - \eqref{GeneralBoussinesqODE_Temp} as an ODE system one must insert the expansions for $\textbf{u}^M$, $\theta^M$ to determine these sums explicitly in terms of the Fourier variables $u^{\textbf{n},M},\theta^{\textbf{n},M}$, namely one must find the time-independent constants $I_{\textbf{u}}^{\boldsymbol{\alpha}}, I_{\theta}^{\boldsymbol{\alpha}}$ such that
	\[ N_{\textbf{u}}^{\textbf{n},M} = \sum_{ \textbf{n}' \in \mathscr{N}_{\textbf{u}}^M } \sum_{\textbf{n}'' \in \mathscr{N}_{\textbf{u}}^M} I^{\boldsymbol{\alpha}}_{\textbf{u}} u^{\textbf{n}',M} u^{\textbf{n}'',M} \hspace{.5 cm} \text{ , } \hspace{.5 cm} N_{\theta}^{\textbf{n},M} = \sum_{ \textbf{n}'\in \mathscr{N}_{\textbf{u}}^M } \sum_{\textbf{n}''\in \mathscr{N}_{\theta}^M}   I^{\boldsymbol{\alpha}}_{\theta} u^{\textbf{n}',M} \theta^{\textbf{n}'',M} \text{ , } \]
	where the multi-index $\boldsymbol{\alpha} := (\textbf{n},\textbf{n}',\textbf{n}'')$ is used to indicate that $I_{\textbf{u}}^{\boldsymbol{\alpha}}, I_{\theta}^{\boldsymbol{\alpha}}$ depend on all of the indices involved in the triad interaction.  These constants are explicitly derived in \ref{app:NonlinearDerivation} by computing inner products of the basis elements inserted into \eqref{AbstractNonlinear}.  This derivation is somewhat involved, but much of the coming analysis only requires general properties of $I_{\textbf{u}}^{\boldsymbol{\alpha}}, I_{\theta}^{\boldsymbol{\alpha}}$ rather than their precise form, so the reader could skip this derivation on a first read.  The proof of Theorem \ref{thm:OriginBifurcations} as well as the numerical implementation do require the precise formulas for $I_{\textbf{u}}^{\boldsymbol{\alpha}}, I_{\theta}^{\boldsymbol{\alpha}}$ and so for a complete definition of the ODE model we give the general formulas here.  First, we introduce some notation to describe the triad interactions:
	\begin{align} \label{def:TriadNotation}  \boldsymbol{\phi} := (p_1,p_1',p_1'') \hspace{.5 cm} & \text{ , } \hspace{.5 cm} \boldsymbol{\mu}^{k} := (m_k,m_k',m_k'') \hspace{.5 cm} \text{ for } k =1,3, \\
		\boldsymbol{\xi}^1 := (1,1,1) \hspace{.5 cm} & \text{ , } \hspace{.5 cm} \boldsymbol{\xi}^2 := (1,2,2) \hspace{.5 cm}  \text{ , } \hspace{.5 cm} \boldsymbol{\xi}^3 := (2,1,2) \hspace{.5 cm} \text{ , } \hspace{.5 cm} \boldsymbol{\xi}^4 := (2,2,1) . \notag
	\end{align}
	The constants $I^{\boldsymbol{\alpha}}_{\textbf{u}},I^{\boldsymbol{\alpha}}_{\theta}$ are zero unless the wave number and phase triads $\boldsymbol{\mu}^{k},\boldsymbol{\phi}^{k}$ satisfy certain compatibility conditions, namely $\boldsymbol{\mu}^{k}$ must satisfy convolution type conditions, the phases $\boldsymbol{\phi}^{1}$ satisfy sinusoidal orthogonality conditions, and the component indices $c,c',c''$ satisfy pointwise orthogonality conditions:
	\begin{equation} \label{CompatibilityCond} m_k = |m_k' \pm m_k''| \hspace{.25 cm} \text{ for } \hspace{.25 cm} k = 1,3 \hspace{.5 cm} \text{ , } \hspace{.5 cm}  \boldsymbol{\phi} \in \{ \boldsymbol{\xi}^k \}_{k\leq 4} \hspace{.5 cm} \text{ , } \hspace{.5 cm} c' = 1  \text{ , }  c'' = c . \end{equation} 
	When these conditions are satisfied then one can write  \setlength\extrarowheight{2.5pt} 	 
	\begin{align} \label{def:NonlinCoeffs1} I_{\theta}^{\boldsymbol{\alpha}} = C^{\boldsymbol{\alpha}} \zeta^{\boldsymbol{\alpha},3} \hspace{.25 cm} & \text{ , } \hspace{.25 cm} I_{\textbf{u}}^{\boldsymbol{\alpha}} = C^{\boldsymbol{\alpha}} \left \{ \begin{array}{cl}
			\frac{ -m_3 m_3'' \zeta^{\boldsymbol{\alpha},1} + (-1)^{p_1 + p_1''} \Sha^2 m_1 m_1'' \zeta^{\boldsymbol{\alpha},3} }{|\mathcal{K} \textbf{m}||\mathcal{K} \textbf{m}''| } & \text{ if } c = 1, \\
			-\zeta^{\boldsymbol{\alpha},1} & \text{ if } c = 2 .
		\end{array} \right .  \end{align}
	where
	\[ C^{\boldsymbol{\alpha}} = \frac{ \Sha }{4\eta^{\textbf{m}}\eta^{\textbf{m}'}\eta^{\textbf{m}''} |\mathcal{K} \textbf{m}'| V} \text{ . } \]
	In order to specify the coefficients $\zeta^{\boldsymbol{\alpha},j}$, $j=1,3$, we introduce the phase maps
	\[ \rho^{1}(\boldsymbol{\phi}) := \boldsymbol{\phi} + (1,1,0) \text{ , } \hspace{.5 cm} \rho^{2}(\boldsymbol{\phi}) = \boldsymbol{\phi} + (1,0,1) \text{ , } \hspace{.5 cm} \rho^{3}(\boldsymbol{\phi}) = \boldsymbol{\phi} + (0,1,1) \text{ , } \hspace{.5 cm} \rho^{4}(\boldsymbol{\phi}) = \boldsymbol{\phi} \text{ . }  \]
	The coefficients $\zeta^{\boldsymbol{\alpha},j}$ are then given by
	\begin{equation} \label{def:Zeta} 
		\zeta^{\boldsymbol{\alpha},j} = (-1)^{p_1''} m_3' m_1'' s^{(\boldsymbol{\mu}^1,\rho^{j}(\boldsymbol{\phi}))} s^{(\boldsymbol{\mu}^3,\boldsymbol{\xi}^j)} + (-1)^{p_1'} m_1' m_3'' s^{(\boldsymbol{\mu}^1,\rho^{j+1}(\boldsymbol{\phi}))} s^{(\boldsymbol{\mu}^3,\boldsymbol{\xi}^{j+1})} \text{ , }
	\end{equation}
	in which the sign coefficients are defined as follows:
	\begin{equation} \label{def:SignCoefs} \begin{split}
			s^{(\boldsymbol{\mu},\boldsymbol{\xi}^1)} := 1 \hspace{.15 cm} \text{ , } \hspace{.15 cm} s^{(\boldsymbol{\mu},\boldsymbol{\xi}^2)} := S^{m,m',m''} & \hspace{.15 cm} \text{ , } \hspace{.15 cm} s^{(\boldsymbol{\mu},\boldsymbol{\xi}^3)} := S^{m',m'',m} \hspace{.15 cm} \text{ , } \hspace{.15 cm} s^{(\boldsymbol{\mu},\boldsymbol{\xi}^4)} := S^{m'',m,m'} \text{ , } \\
			S^{a_1,a_2,a_3} & := \left \{ \begin{array}{cc}
				-1 & \text{ if } a_1 = a_2+a_3, \\
				1 & \text{ otherwise. } \end{array} \right . 
	\end{split} \end{equation}
	
	\subsection{Mode selection criteria and the HKC hierarchy }
	
	\subsubsection{Energetic consistency}
	
	The formulas from the previous section enable one to write down a system of ODE's for any choice of a finite set of Fourier modes $\mathscr{N}_{\textbf{u}}^{M},\mathscr{N}_{\theta}^{M}$.  We turn now to the question of how to choose these sets to best retain the behavior of the full system.  Sufficiently regular solutions $\textbf{u}, \theta$ of \eqref{EvolutionPDE} should satisfy the following balance equations, where $\vorticity := \nabla \times \textbf{u}$: 
	\begin{align}
		\label{Balance_KinEnergy}
		\frac{d}{dt} \langle \frac{1}{2}|\textbf{u}|^2 \rangle & = - \Pra \langle |\nabla \textbf{u} |^2 \rangle + \Pra \Ray \langle u_3 \theta \rangle ,  \\
		\label{Balance_TempVar}
		\frac{d}{dt} \langle \frac{1}{2} \theta^2 \rangle & = - \langle |\nabla \theta |^2 \rangle + \langle u_3 \theta \rangle , \\
		\label{Balance_PotEnergy}
		\frac{d}{dt} \langle (1-\frac{x_3}{\pi}) \theta \rangle & = \langle (1-\frac{x_3}{\pi}) \partial_{x_3}^2 \theta \rangle - \frac{1}{\pi} \langle u_3 \theta \rangle , \\
		\label{Balance_Vort}
		\frac{d}{dt} \langle \vorticity \rangle & = \Pra \langle \partial_{x_3}^2 \vorticity \rangle + \langle ( \vorticity \cdot \nabla ) \textbf{u} \rangle + \Pra \Rot \langle  \partial_{x_3} \textbf{u} \rangle .
	\end{align}
	These balance relations are physically meaningful, since they express the time evolution of kinetic energy, temperature variance, potential energy and the total vorticity for the PDE \eqref{BoussinesqCoriolis}, respectively.  For any choice of index sets $\mathscr{N}_{\textbf{u}}^{M}$, $\mathscr{N}_{\theta}^{M}$ the same balance equations \eqref{Balance_KinEnergy}, \eqref{Balance_TempVar} hold also for the truncated fields $\textbf{u}^M,\theta^M$, since these balances are derived from inner products of the evolution equations for $\textbf{u}^M, \theta^M$ with themselves, hence projection operators $\mathcal{P}_{\textbf{u}}^M,\mathcal{P}_{\theta}^M$ in \eqref{TruncatedEvolutionPDE} have no effect in their derivation.  More concretely, the projectors $\mathcal{P}^{M}_{\textbf{u}},\mathcal{P}^{M}_{\theta}$ have a self-adjointness property, namely for any scalar functions $F,G \in L^2_0$ one has
	\begin{equation} \label{SelfAdjointProp} \langle \mathcal{P}^{M}_{\theta} [ F ] G \rangle = \Big \langle \sum_{\textbf{n} \in \mathscr{N}_{\theta}^M} \langle  F f^{\textbf{n}} \rangle f^{\textbf{n}} G \Big \rangle =  \Big \langle F \sum_{\textbf{n} \in \mathscr{N}_{\theta}^M} \langle f^{\textbf{n}} G \rangle f^{\textbf{n}} \Big \rangle = \langle F \mathcal{P}^{M}_{\theta} [ G ] \rangle , \end{equation}
	and a similar identity holds for $\mathcal{P}_{\textbf{u}}^M$.  Since also $\mathcal{P}^{M}_{\theta} [ \theta^M ] = \theta^M$ it follows that the projectors disappear when computing these "self" inner products, so for instance with the nonlinear terms one can then integrate by parts:
	\[ \langle \theta^M \cdot \mathcal{P}_{\theta}^M \big [ \textbf{u}^M \cdot \nabla \theta^M \big ] \rangle = \langle \theta^M \cdot \big [ \textbf{u}^M \cdot \nabla \theta^M \big ] \rangle = \langle \textbf{u}^M \cdot \nabla \frac{1}{2} (\theta^M)^2 \rangle = 0 . \]
	However, it has been recognized in several works that the balance equations \eqref{Balance_PotEnergy},\eqref{Balance_Vort} can fail to hold in a truncated ODE model unless the modes are selected in a particular way \cite{HermizGuzdarFinn_1995,ThiffeaultHorton_1996,GluhovskyTongAgee_2002}, and this failure can lead to unphysical dynamics such as unbounded trajectories \cite{HowardKrishnamurti_1986}.  We first consider the following mode-selection Criterion:
	\begin{Criterion}
		(Energy balance): For any pair $\textbf{n}' \in \mathscr{N}_{\textbf{u}}^{M}, \textbf{n}'' \in \mathscr{N}_{\theta}^{M}$ satisfying
		\begin{equation} \label{eq:EnergyCrit} m_1'=m_1'' > 0 \hspace{.25 cm} \text{ , } \hspace{.25 cm} p' = p'' \hspace{.25 cm} \text{ , } \hspace{.25 cm} c' = c'' = 1  \text{ , } \end{equation}
		and $m_3' = m_3''$, one must have $(0,2m_3',1,1) \in \mathscr{N}^{M}_{\theta}$. On the other hand, if any pair $\textbf{n}' \in \mathscr{N}_{\textbf{u}}^{M}, \textbf{n}'' \in \mathscr{N}_{\theta}^{M}$ satisfies \eqref{eq:EnergyCrit} and $m_3' \neq m_3''$, then one can have $(0,|m_3' - m_3''|,1,1) \in \mathscr{N}^{M}_{\theta}$ if and only if $(0,m_3' + m_3'',1,1) \in \mathscr{N}^{M}_{\theta}$.
		\label{Crit:EnergyCrit}
	\end{Criterion}
	\begin{figure}[H]
		\begin{center}
			\begin{tabular}{ccc}
				\includegraphics[height=46mm]{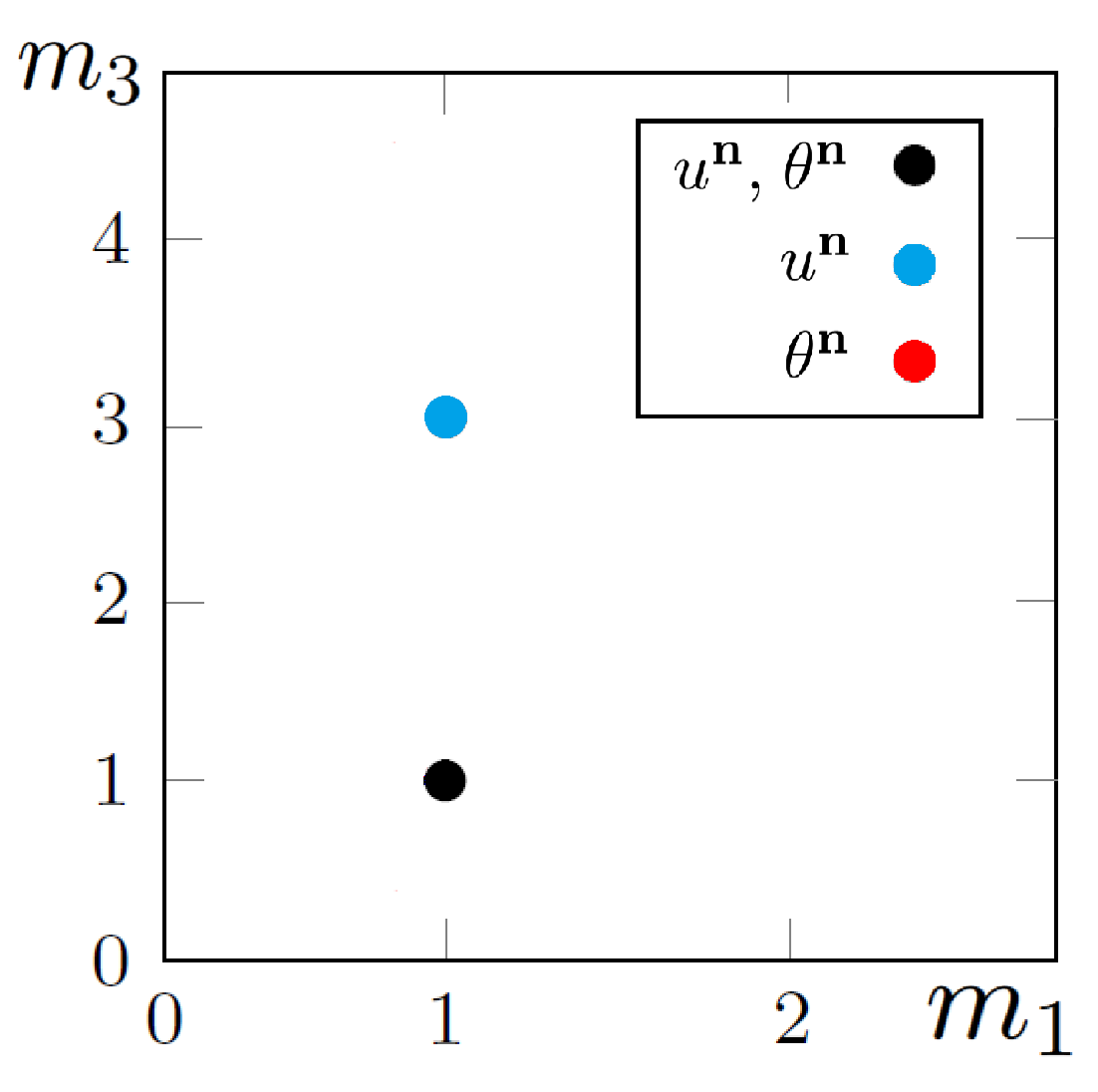}&
				\includegraphics[height=46mm]{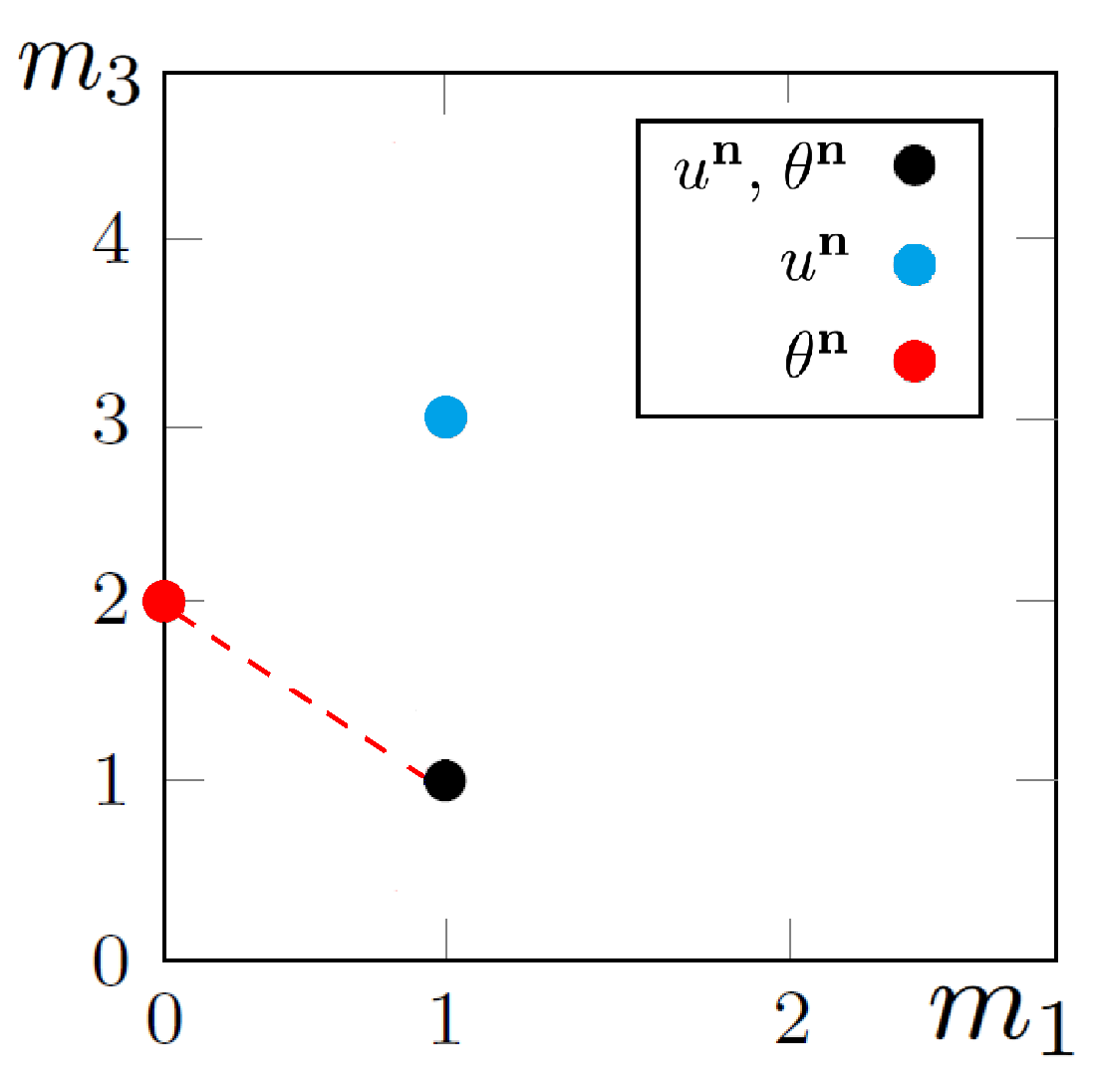}&
				\includegraphics[height=46mm]{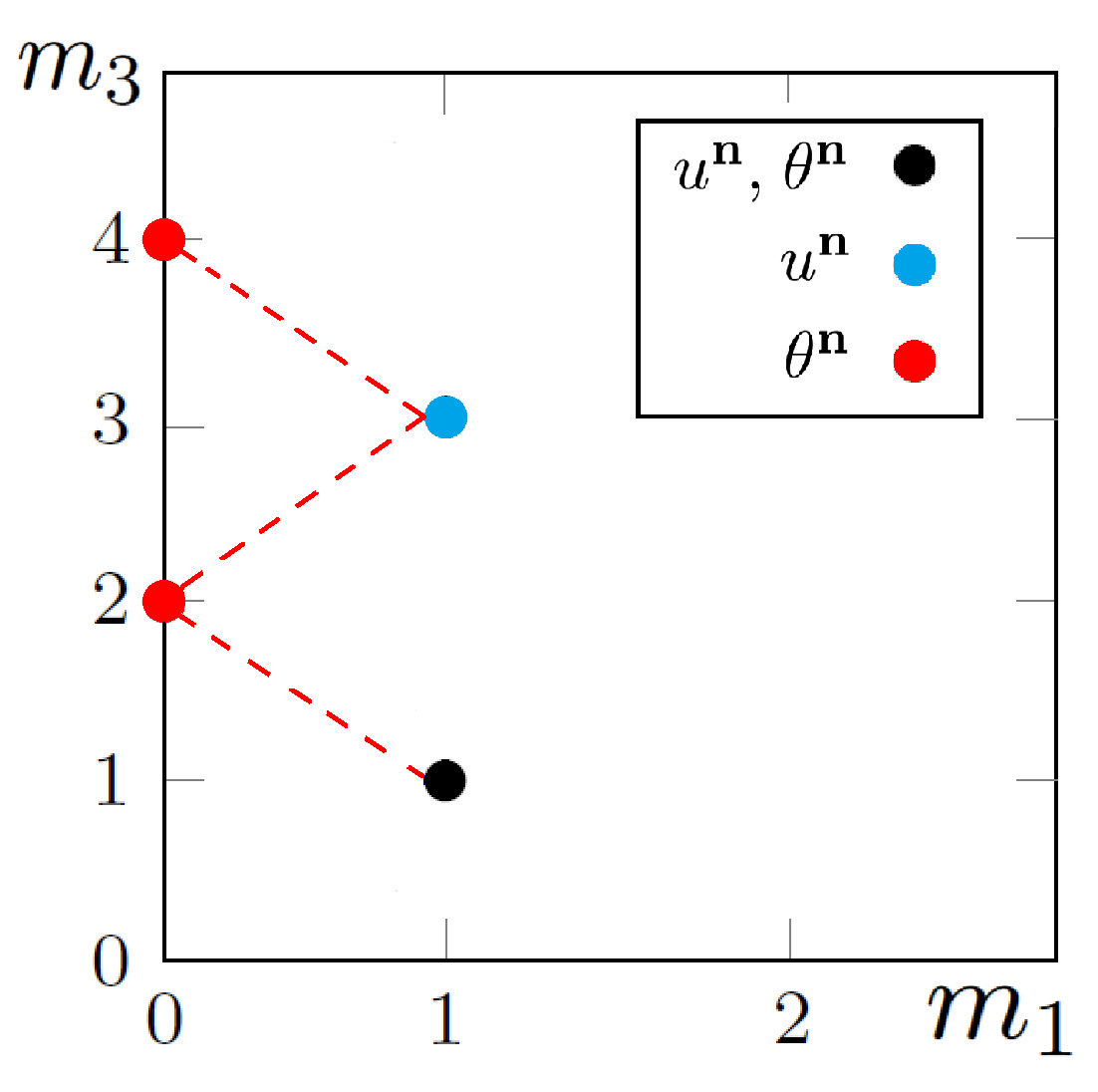} \\
				(a) & (b) & (c)
			\end{tabular}
			\caption{A depiction of the Fourier modes included in the minimal model which contains the modes in \eqref{Example:EnerCons} and satisfies Criterion \ref{Crit:EnergyCrit}. } 
			\label{fig:EnergyConsistCrit}
		\end{center}
	\end{figure}
	\revised{To illustrate this Criterion, consider the example depicted in Figure \ref{fig:EnergyConsistCrit}.  Suppose one wanted to study an energetically consistent truncated ODE which includes the following Fourier modes:
		\begin{equation} \label{Example:EnerCons} \{ (1,1,1,1),(1,3,1,1) \} \subset \mathscr{N}_{\textbf{u}}^M \hspace{.25 cm} \text{ , } \hspace{.25 cm} \{ (1,1,1,1) \}  \subset \mathscr{N}_{\theta}^M \text{ . }  \end{equation}
		Since both $\mathscr{N}_{\textbf{u}}^M$ and $ \mathscr{N}_{\theta}^M$ contain the mode $(1,1,1,1)$, Criterion \ref{Crit:EnergyCrit} says one must additionally include the mode $(0,2,1,1)$ in $\mathscr{N}_{\theta}^M$:
		\[ \{ (1,1,1,1),(1,3,1,1) \} \subset \mathscr{N}_{\textbf{u}}^M \hspace{.25 cm} \text{ , } \hspace{.25 cm} \{ (1,1,1,1), (0,2,1,1) \}  \subset \mathscr{N}_{\theta}^M \text{ . }  \]
		Considering now the pair $\textbf{n}' = (1,3,1,1) \in \mathscr{N}_{\textbf{u}}^M$ and $\textbf{n}'' = (1,1,1,1) \in  \mathscr{N}_{\theta}^M$, note that since $(0,2,1,1)$ in $\mathscr{N}_{\theta}^M$, Criterion \ref{Crit:EnergyCrit} says one must additionally include the mode $(0,4,1,1)$ in $\mathscr{N}_{\theta}^M$.  Thus the minimal model satisfying Criterion \ref{Crit:EnergyCrit} which contains the modes in \eqref{Example:EnerCons} is given by 
		\[ \mathscr{N}_{\textbf{u}}^M =  \{ (1,1,1,1),(1,3,1,1) \} \hspace{.25 cm} \text{ , } \hspace{.25 cm} \mathscr{N}_{\theta}^M = \{ (1,1,1,1), (0,2,1,1), (0,4,1,1) \}  \text{ . }  \]}
	
	In the following proposition, it is proven that Criterion \ref{Crit:EnergyCrit} is necessary and sufficient to prove \eqref{Balance_PotEnergy} also holds for the corresponding truncated model, and that \eqref{Balance_KinEnergy}-\eqref{Balance_PotEnergy} then imply the existence of a global attractor.
	
	\begin{Proposition}
		\label{prop:BalancePreservation}
		For any admissible parameters $\parameters$, index sets $\mathscr{N}_{\textbf{u}}^{M},\mathscr{N}_{\theta}^{M}$ and initial condition $\textbf{X}^M_0 = (\textbf{u}_{0}^{M}, \theta_{0}^M)$ a unique smooth solution $\textbf{X}^M(t) = (\textbf{u}^{M}(t), \theta^M(t))$ of the truncated ODE model \eqref{GeneralBoussinesqODE_Vel} - \eqref{GeneralBoussinesqODE_Temp} with $\textbf{X}(0) = \textbf{X}_0$ exists.  The following statements hold:
		\begin{enumerate}[label=\textbf{(\alph*)}]
			\item The potential energy equation \eqref{Balance_PotEnergy} holds with $\textbf{u},\theta$ replaced by $\textbf{u}^M,\theta^M$ if and only if Criterion \ref{Crit:EnergyCrit} (i) is satisfied.
			
			\item If \eqref{Balance_PotEnergy} holds, then there exists a forward invariant ball depending only on $\Pra, \Ray, \Sha$ and $|\mathscr{N}_{\theta}^{M}|$ which monotonically attracts all solutions lying outside. Thus all solutions are global and the semi-group $\mathcal{S}^M(t)$, defined by
			\begin{equation} \label{ODE_Semigroup} \mathcal{S}^M(t) \big [  \textbf{X}_0^M \big ] := \textbf{X}^M(t) , \end{equation}
			admits a compact global attractor $\mathscr{A}^M$.
		\end{enumerate}
	\end{Proposition}
	
	\begin{proof}
		The ODE models are locally Lipschitz, so the Cauchy-Lipschitz theorem is sufficient to guarantee existence of smooth, short-time solutions.  Global existence follows from the presence of an attracting ball, which can be proven using the balances.  Hence we consider the first balance \eqref{Balance_PotEnergy}.  Computing the inner product of the temperature equation in \eqref{TruncatedEvolutionPDE} with $(1-\frac{x_3}{\pi})$, one obtains
		\begin{equation} \label{EnerConsis_Eq2} \begin{split} \frac{d}{dt} \langle (1-\frac{x_3}{\pi}) \theta^M \rangle & = \langle (1-\frac{x_3}{\pi} ) \partial_{x_3}^2 \theta^M \rangle - \langle (1-\frac{x_3}{\pi} ) \mathcal{P}^{M}_{\theta} \big [ \textbf{u}^M \cdot \nabla \theta^M \big ] \rangle . \end{split} \end{equation}
		Let $l^{M}(x_3)$ denote the projection of the linear background state onto the modes included in the index set $\mathscr{N}_{\theta}^{M}$:
		\begin{equation} \label{ConductStateFourier} l^{M}(x_3) := \mathcal{P}^{M}_{\theta} \big [ 1-\frac{x_3}{\pi} \big ] = \sum_{\substack{ (0,m_3,0,1) \in \mathscr{N}_{\theta}^M \\ m_3 > 0 }} \frac{2}{\sqrt{\Sha} m_3}  \frac{\sin( m_3 x_3) }{\sqrt{2} V} . \end{equation}
		This function of course depends only on the vertical variable, and in fact the set of vertically stratified modes will occur so often in our analysis that we define 
		\begin{align} \label{VertStratNodes}
			\mathscr{M}^* = \{ \textbf{m} \in \mathbb{Z}^2_{\geq 0} & : m_1 = 0 \text{ , } m_3 > 0 \}  \text{ , } \\ \mathscr{N}^*_{\textbf{u}} = \{ \textbf{n} \in \mathscr{N}_{\textbf{u}} : \textbf{m} \in \mathscr{M}^* \} \hspace{.5 cm} & \text{ , } \hspace{.5 cm} \mathscr{N}^*_{\theta} = \{ \textbf{n} \in \mathscr{N}_{\theta} : \textbf{m} \in \mathscr{M}^* \} . \notag  
		\end{align}
		Focusing on the last term in \eqref{EnerConsis_Eq2}, one finds the following by using \eqref{SelfAdjointProp}, the mean-free property of $\theta$ and integrating by parts 
		\[ \langle (1-\frac{x_3}{\pi} ) \mathcal{P}^{M}_{\theta} \big [ \textbf{u}^M \cdot \nabla \theta^M \big ] \rangle = \langle l^{M}(x_3)  \textbf{u}^M \cdot \nabla \theta^M \rangle  = \langle \partial_{x_3} l^{M}(x_3)  u_3^M \theta^M  \rangle . \]
		Comparing with the expression in \eqref{Balance_PotEnergy}, one sees that the ODE model obeys a consistent potential energy balance iff
		\begin{equation} \label{EnerConsis_Condition} \langle \partial_{x_3} l^{M}(x_3)  u_3^M \theta^M  \rangle = \langle u_3^M \theta^M \rangle . \end{equation}
		Inserting the Fourier expansions into the left hand side of \eqref{EnerConsis_Condition} gives
		\begin{align} \langle \partial_{x_3} l^{M}(x_3)  u_3^M \theta^M  \rangle & = 2 \sum_{ \textbf{n} \in \mathscr{N}_{\theta}^M \cap \mathscr{N}_{\theta}^*} \langle  \cos(m_3 x_3) u_3^{M} \theta^{M} \rangle \notag \\ \label{EnerConsis_LHS} & = 2 \sum_{\textbf{n} \in \mathscr{N}_{\theta}^M \cap \mathscr{N}_{\theta}^*} \sum_{\textbf{n}' \in \mathscr{N}_{\textbf{u}}^M } \sum_{\textbf{n}'' \in \mathscr{N}_{\theta}^M} u^{\textbf{n}',M} \theta^{\textbf{n}'',M} \langle  \cos(m_3 x_3) v^{\textbf{n}'}_3 f^{\textbf{n}''} \rangle .  \end{align}
		The orthogonality properties of the sinusoidal functions imply that the only terms which survive in this sum satisfy $m_1' = m''_1$, $p_1' = p_1''$, $c' = c'' = 1$ and either $m_3 = m_3' + m''_3$,  $m_3 = m_3' - m''_3$ or  $m_3 = -m_3' + m''_3$.  On the other hand, inserting the Fourier expansions into the right hand side of \eqref{EnerConsis_Condition} one has
		\[ \langle u_3^M \theta^M \rangle = \sum_{\textbf{n}' \in \mathscr{N}_{\textbf{u}}^M } \sum_{\textbf{n}'' \in \mathscr{N}_{\theta}^M } u^{\textbf{n}',M} \theta^{\textbf{n}'',M} \langle v^{\textbf{n}'}_3 f^{\textbf{n}''} \rangle , \] 
		so using \eqref{LinearBasisRelations} one must have $\textbf{m}' = \textbf{m}''$, $p_1' = p_1''$, $c' = c'' = 1$:
		\begin{equation} \label{EnerConsis_RHS} \begin{split}  \langle u_3^M \theta^M \rangle & =  \sum_{\textbf{n}' \in \mathscr{N}_{\textbf{u}}^M \cap \mathscr{N}_{\theta}^M } u^{\textbf{n}',M} \theta^{\textbf{n}',M} (-1)^{p_1} \frac{\Sha m_1}{|\mathcal{K}\textbf{m}|} . \end{split} \end{equation}
		Assume now that Criterion \ref{Crit:EnergyCrit} (i) is satisfied.  Then for each term appearing in \eqref{EnerConsis_LHS} with $m_3' \neq m''_3$ the terms with $m_3 = m_3' + m''_3$ and $m_3 = |m_3' - m''_3|$ are both included in the sum, and one can easily check that
		\[ \langle \cos \big ( ( m_3' + m''_3 ) x_3 \big ) v_3^{\textbf{n}'} f^{\textbf{n}''} \rangle = - \langle \cos \big ( | m_3' - m''_3 | x_3 \big ) v_3^{\textbf{n}'} f^{\textbf{n}''} \rangle , \]
		since here only the integral in $x_3$ needs consideration.  Thus all such terms cancel out.  Furthermore, for each term appearing in \eqref{EnerConsis_LHS} with $m_3' = m''_3$ the index $m_3 = 2m_3'$ is included in the sum, and by inserting the definitions \eqref{ThetaFourierBasis}, \eqref{VectorFieldDef_FreeSlip} and computing the integrals one can similarly check that
		\[  -2 \langle  \cos( 2m_3' x_3) v^{\textbf{n}'}_3 f^{\textbf{n}'} \rangle = (-1)^{p_1+1} \frac{4}{\pi} \frac{\Sha m_1}{|\mathcal{K}\textbf{m}|} \int_0^{\pi} \cos( 2m_3' x_3) \sin^2( m_3 x_3) dx_3 = (-1)^{p_1} \frac{\Sha m_1}{|\mathcal{K}\textbf{m}|} . \]
		Thus \eqref{EnerConsis_LHS} and \eqref{EnerConsis_RHS} are equal and \eqref{EnerConsis_Condition} holds.  On the other hand if Criterion \ref{Crit:EnergyCrit} (i) is not satisfied, then it is clear that at least one term which appears on one side of \eqref{EnerConsis_LHS}, \eqref{EnerConsis_RHS} does not appear on the other side, either because it is not included in the sum, or because a cancellation fails to occur. 
		
		For part (b) note that by combining the balance equations \eqref{Balance_KinEnergy} - \eqref{Balance_PotEnergy} one obtains the following:
		\begin{equation} \label{ODE_AttractingBall} \frac{1}{2} \frac{d}{dt} \Big \langle \frac{|\textbf{u}^M|^2}{\Pra \Ray}  + \big ( \theta^M \big )^2 + 4\pi (1-\frac{x_3}{\pi} )\theta^M   \Big \rangle = - \Big \langle \frac{|\nabla \textbf{u}^M|^2}{\Ray} + |\nabla \theta^M |^2 - 2 \pi (1-\frac{x_3}{\pi} )\partial_{x_3}^2 \theta^M \Big \rangle .
		\end{equation}
		This alone implies the existence of an attracting ball $\mathscr{B}^M$ for the ODE system.  To see why, one can use \eqref{SelfAdjointProp} and by adding a constant term inside the time derivative to complete the square, the left hand side becomes
		\[ \frac{1}{2} \frac{d}{dt} \Big \langle \frac{|\textbf{u}^M|^2}{\Pra \Ray}  + \big ( \theta^M + 2\pi  l^{M}(x_3) \big )^2 \Big \rangle . \]
		For the right hand side one can use \eqref{SelfAdjointProp}, integrate by parts and add and subtract a term to obtain the following:
		\[ \Big \langle - |\nabla \theta^M |^2 + 2 \pi (1-\frac{x_3}{\pi} )\partial_{x_3}^2 \theta^M \Big \rangle =  \Big \langle ( \theta^M + \pi l^{M}) \Delta \theta^M + \theta^M \Delta \pi l^{M} + \tilde{T} - \tilde{T} \Big \rangle \text{ , } \]
		where $\tilde{T} = \pi l^{M}(x_3) \Delta \pi l^{M}(x_3) $.  Hence by integrating by parts, one obtains 
		\[ \Big \langle - |\nabla \theta^M |^2 + 2 \pi (1-\frac{x_3}{\pi} )\partial_{x_3}^2 \theta^M \Big \rangle =  \Big \langle - \big | \nabla \big ( \theta^M + \pi l^{M} \big ) \big  |^2 + \pi^2 \big | \nabla l^{M} \big |^2 \Big \rangle \text{ . } \]
		Letting $\rho := \frac{2\pi}{\sqrt{\Sha}} |\mathscr{N}_{\theta}^{M} \cap \mathscr{N}^{*}_{\theta}|^{1/2}$ , \eqref{ODE_AttractingBall} can thus be rewritten as 
		\begin{equation} \label{ODE_AttractingBall2} \frac{1}{2} \frac{d}{dt} \Big \langle \frac{|\textbf{u}^M|^2}{\Pra \Ray}  + \big ( \theta^M + 2\pi l^{M} \big )^2 \Big \rangle = - \Big \langle \frac{|\nabla \textbf{u}^M|^2}{\Ray} + \big | \nabla \big ( \theta^M + \pi l^{M} \big ) \big  |^2 \Big \rangle  + \rho^2 \text{ . } \end{equation}
		On the right hand side one sees the time derivative of a weighted $\textbf{H}^{0}$ norm of $(\textbf{u}^M,\theta^M + 2\pi l^{M})$, and on the left one has the a weighted $\textbf{H}^{1}$ norm of $(\textbf{u}^M,\theta^M + \pi l^{M}(x_3))$ subtracted from a positive constant.  Thus the $\textbf{H}^{0}$ norm decreases until the $\textbf{H}^{1}$ is smaller than the positive constant, and hence a ball containing the region where this $\textbf{H}^{1}$ norm is smaller will be attracting.  In fact, the weighted norms in \eqref{ODE_AttractingBall2} express ellipsoids in phase space with principle semi-axes $\Pra \Ray$, $1$ and $\Ray$, $1$ along the $\textbf{u}^{M}, \theta^M$ directions respectively.  To prove existence of an attractor one only needs an attracting ball, but we point out that using these ellipsoids gives a bound on $\theta^M$ which is independent of $\Ray$.
		
		Denoting such an attracting ball by $\mathscr{B}^M$, recall that the $\omega$-limit set of $\mathscr{B}^M$ with respect to $\mathcal{S}(t)$ is defined
		\[ \omega(\mathscr{B}^M) := \cap_{s > 0} \overline{\cup_{t > s} \mathcal{S}(t)\mathscr{B}^M} \text{ , } \]
		where the overline denotes the closure.  Since $\mathscr{B}^M$ is forward invariant, this is an intersection of compact sets, hence compact.  It is easy to check that it is forward invariant, attracts the bounded sets of the phase space, and that it is the maximal set with these properties.  Thus $\mathscr{A}^M := \omega(\mathscr{B}^M)$ is the global attractor.
		
		\vspace{-.75 cm}
		\[ \textcolor{white!100}{ . } \]
	\end{proof}
	
	Turning now to the vorticity balance \eqref{Balance_Vort}, consider the following mode selection Criteria:
	\begin{Criteria}
		\begin{enumerate}[label=(\roman*)]
			\item \revised{(Vorticity balance): Suppose $\textbf{n}', \textbf{n}'' \in \mathscr{N}_{\textbf{u}}^{M}$ satisfy
				\[ m_1'=m_1'' > 0 \hspace{.25 cm} \text{ , } \hspace{.25 cm} m_3'+m_3'' \text{ odd } \hspace{.25 cm} \text{ , } \hspace{.25 cm} p' = p'' + 1 \text{ . } \]
				If $c' = c'' = 1$, then one can have $(0,|m_3' - m_3''|,2,1) \in \mathscr{N}^{M}_{\textbf{u}}$ if and only if $(0,m_3' + m_3'',2,1) \in \mathscr{N}^{M}_{\textbf{u}}$.  On the other hand, if $(c',c'') = (1,2)$ or $(2,1)$ then one must have both $(0,m_3' + m_3'',2,2), (0,|m_3' - m_3''|,2,2) \in \mathscr{N}^{M}_{\textbf{u}}$.}
			
			\item (Rotating vorticity balance): If $\Rot \neq 0$ and $m_3 > 0$ is odd, $(0,m_3,p_1,c)  \in \mathscr{N}_{\textbf{u}}^{M} \Rightarrow (0,m_3,p_1,c+1)  \in \mathscr{N}_{\textbf{u}}^{M}$.
		\end{enumerate}
		\label{Crit:VortCrit} 
	\end{Criteria}
	In the following Proposition, it is proven that these Criteria are necessary and sufficient to ensure that \eqref{Balance_Vort} also holds for the corresponding truncated model.  Similar to Prop. \ref{prop:BalancePreservation}, the essence of the proof involves checking when the curl operator commutes with the projection operator $\mathcal{P}_{\textbf{u}}^{M}$.  However, the analysis is more technical due to the curl operator, and so the proof is given in 
	\ref{app:VortBal}.
	\begin{Proposition}
		\label{prop:VortBalancePreservation}
		For any admissible parameters $\parameters$, index sets $\mathscr{N}_{\textbf{u}}^{M},\mathscr{N}_{\theta}^{M}$ and initial condition $\textbf{X}^M_0 = (\textbf{u}_{0}^{M}, \theta_{0}^M)$ let $\textbf{X}^M(t) = (\textbf{u}^{M}(t), \theta^M(t))$ be the solution of the truncated ODE model \eqref{GeneralBoussinesqODE_Vel} - \eqref{GeneralBoussinesqODE_Temp} with $\textbf{X}(0) = \textbf{X}_0$.  The vorticity equation \eqref{Balance_Vort} holds for $\vorticity^M := \nabla \times \textbf{u}^M, \textbf{u}^M$ if and only if Criteria \ref{Crit:VortCrit} are satisfied.
	\end{Proposition}
	
	\subsubsection{Buoyancy criterion, invariant subspaces, phase locking and the HKC hierarchy}
	
	In this section we define the Howard-Krishnamurthy-Coriolis (HKC) hierarchy of ODE models, the main computational tool used herein.  In order to converge the PDE in the limit, the index sets in this hierarchy should be chosen to increase to eventually include all Fourier modes, ie one should have $\mathscr{N}_{\textbf{u}}^{M} \subseteq \mathscr{N}_{\textbf{u}}^{\tilde{M}}$ for $M < \tilde{M}$ and $\cup_{M} \mathscr{N}_{\textbf{u}}^{M} = \mathscr{N}_{\textbf{u}}$.  Furthermore at each step in the hierarchy Criteria \ref{Crit:EnergyCrit}, \ref{Crit:VortCrit} should be satisfied.  Two further considerations are involved in the definition of the HKC hierarchy: buoyancy instabilities and computational expense. 
	
	Based only on Criteria \ref{Crit:EnergyCrit}, \ref{Crit:VortCrit} one could include a velocity mode $\textbf{u}^{(\textbf{m},p_1,1)}$ in the model without including the corresponding temperature mode $\theta^{(\textbf{m},p_1,1)}$ or vice versa.  In the full PDE, these two modes are coupled via the buoyancy force, and in Section \ref{sec:OriginLocalBif} below we shall see that this coupling gives rise to a linear instability, which in turn produces a non-trivial sub-manifold of the attractor.  Since an accurate representation of the attractor is important for the Nusselt number, we adopt the following additional mode-selection Criterion so that the ODE models also contain such sub-manifolds:
	\begin{Criterion}
		(Buoyancy): For all $\textbf{m} \in \mathbb{Z}_{> 0}^2$, $(\textbf{m},p_1,1) \in \mathscr{N}_{\textbf{u}}^{M} \Leftrightarrow (\textbf{m},p_1,1) \in \mathscr{N}_{\theta}^{M}$.
		\label{Crit:BuoyCrit}
	\end{Criterion}	
	
	Invariant subspaces provide the opportunity to study solutions to \eqref{EvolutionPDE} which are lower dimensional, hence less computationally expensive.  One can deduce the existence of many invariant subspaces of the Boussinesq-Coriolis system from the compatibility conditions \eqref{CompatibilityCond}.  For example the fully 2d planar convection system can be seen to be invariant when the rotation is set to zero by checking the nonlinear couplings of the different phases $c = 1,2$.  While there are many possible interesting choices of invariant subspaces, this paper will focus on invariant subspaces coming from the phase condition in \eqref{CompatibilityCond} as done by Olson and Doering \cite{OlsonDoering2022}, so as to have results to compare against.  Note that the phase condition implies that every nonlinear term for the phase $p_1 = 2$ must always contain a term with either $p_1' = 2$ or $p_1'' = 2$, and hence if all such modes are initially zero then they will remain so for all time.  In fixing the phases of the sinusoids, one throws away the translation invariance of \eqref{BoussinesqCoriolis}, e.g. the convection rolls will occur at fixed locations, hence this is referred to as "phase locking".  Olson and Doering specifically use the following phase locking condition:
	\begin{equation} \label{OlsonDoering_PhaseCond} p_1 =  m_1 + m_3 + 1 \hspace{.25 cm} \text{ mod } \hspace{.1 cm} 2 . \end{equation}
	This is a somewhat unusual choice, but one can check via \eqref{CompatibilityCond} that it is indeed invariant.  This choice might have been made so that Olson and Doering could study models with shear flows while also setting things up such the Lorenz model appears as the first model in their hierarchy.  Note especially that the admissible phase conditions in \eqref{PhaseIndexSets} together with the phase lock condition \eqref{OlsonDoering_PhaseCond} imply that only certain Fourier variables are included depending on the wave vector $\textbf{m}$.  
	
\begin{figure}[H]
		\begin{center}
			\begin{tabular}{cc}
				\includegraphics[height=60mm]{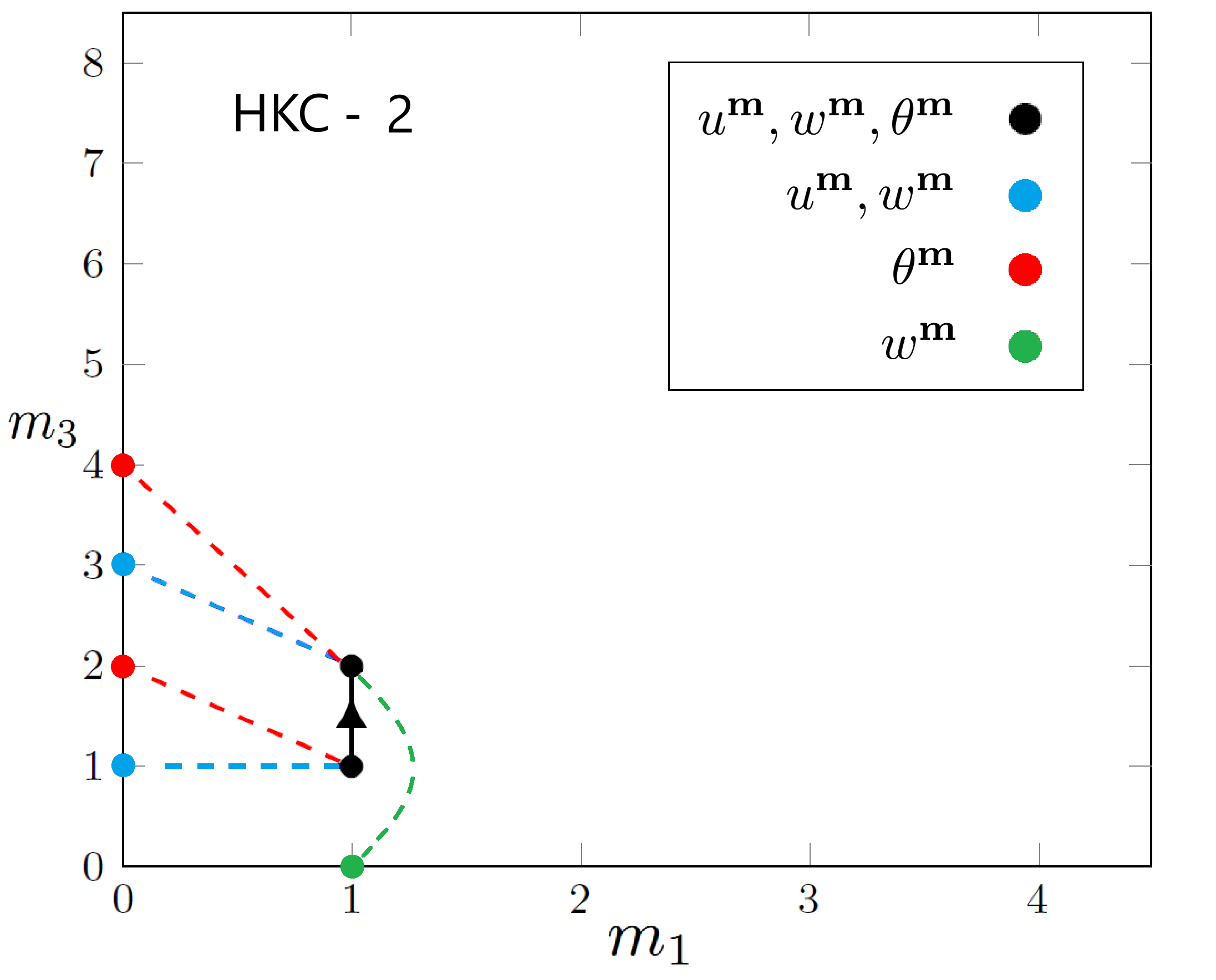}&\includegraphics[height=60mm]{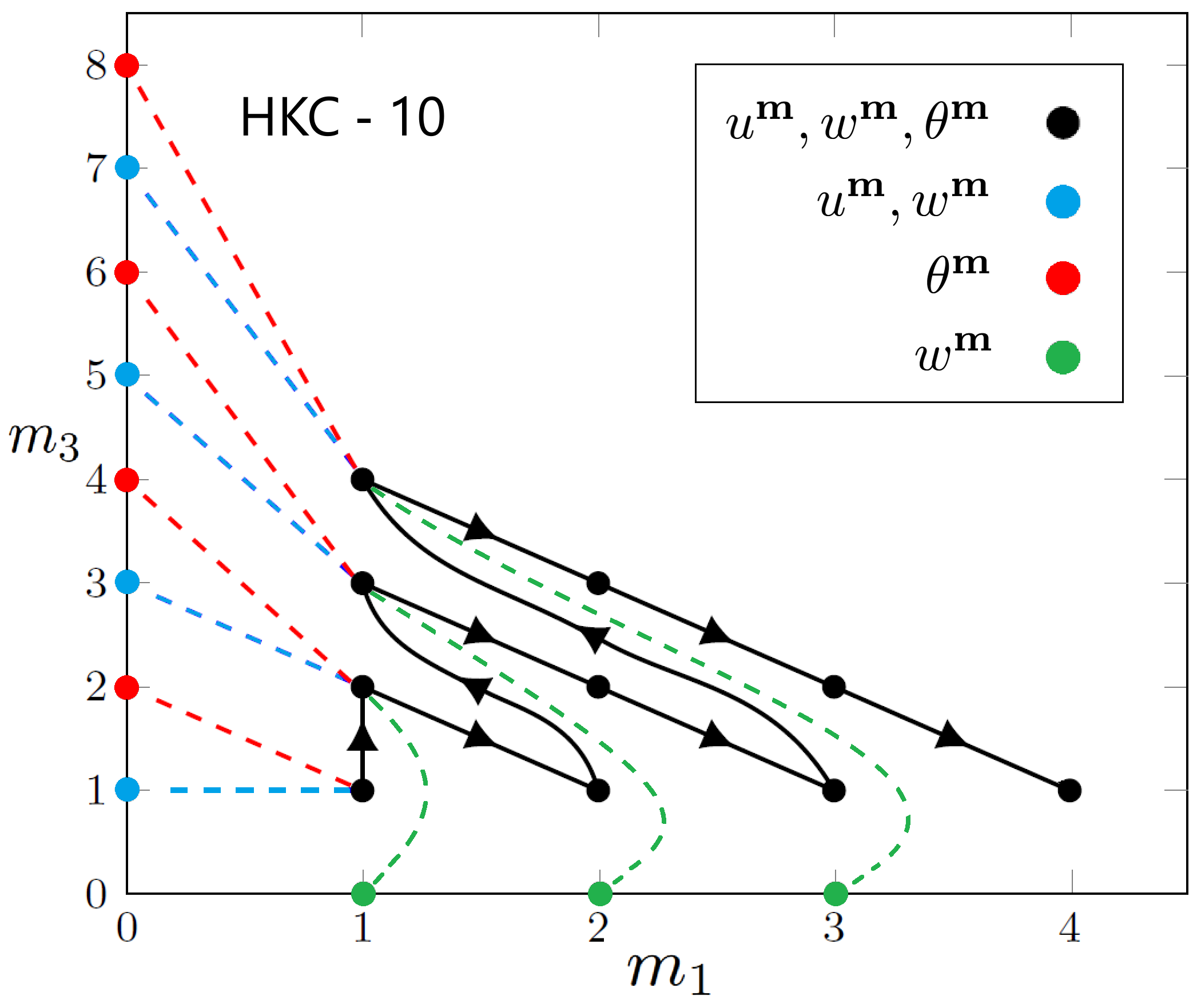} \\
				(a) & (b)
			\end{tabular}
			\caption{A depiction of the mode selection process, displaying which Fourier variables are included in the (a) HKC-2 and (b) HKC-10 models, respectively.  The black line indicates the progression of higher order models, the colors indicate which index sets $\mathscr{M}_{u}^{M},\mathscr{M}_{w}^{M},\mathscr{M}_{\theta}^{M}$ include the mode $\textbf{m}$, and the dashed lines indicating when the vertically stratified modes are adjoined to ensure that Criteria \ref{Crit:EnergyCrit}, \ref{Crit:VortCrit} are satisfied.} 
			\label{fig:HKC_ModeSelect}
		\end{center}
	\end{figure}

	We are therefore ready to define the HKC hierarchy.  For a given wave vector $\textbf{m} \in \mathbb{Z}_{\geq 0}^2$, the phases $p_1$ will always be chosen according to \eqref{OlsonDoering_PhaseCond} and all admissible component indices $c$ will always be included.  Since the phase $p_1$ is always determined by \eqref{OlsonDoering_PhaseCond} we will henceforth drop this from the notation and simply write $u^{\textbf{m}} = u^{(\textbf{m},p_1,1)}$, $w^{\textbf{m}} = u^{(\textbf{m},p_1,2)}$ and $\theta^{\textbf{m}} = \theta^{(\textbf{m},p_1,1)}$.  We also adopt a different notation for the index sets:
	\begin{align} \mathscr{M}_{u}^M := \{ \textbf{m} \in \mathbb{Z}_{\geq 0}^{2} : (\textbf{m},p_1,1) \in & \mathscr{N}_{\textbf{u}}^M \} \text{ , } \mathscr{M}_{w}^M = \{ \textbf{m} \in \mathbb{Z}_{\geq 0}^{2} : (\textbf{m},p_1,2) \in \mathscr{N}_{\textbf{u}}^M \} \text{ , } \notag \\ \mathscr{M}_{\theta}^M = \{ & \textbf{m} \in \mathbb{Z}_{\geq 0}^{2} : (\textbf{m},p_1,1) \in \mathscr{N}_{\theta}^M \} . \notag \end{align} 
	To specify the models we choose these wave vector sets via the following iterative algorithm.  Define an ordering on the wave vectors by saying $\textbf{m} > \tilde{\textbf{m}}$ iff 
	\begin{equation} \label{IndexOrdering} m_1 + m_3 > \tilde{m}_1 + \tilde{m}_3 \text{ , } \hspace{.5 cm} \text{ or } \hspace{.5 cm} m_1 + m_3 = \tilde{m}_1 + \tilde{m}_3 \text{ and } m_1 > \tilde{m}_1 . \notag \end{equation}
	Let $\textbf{m}^{i}$ denote the $i^{th}$ smallest wave vector in $\mathbb{Z}_{> 0}^2$ according to this order, e.g. $\textbf{m}^1 = (1,1), \textbf{m}^2=(2,1),$ and so on.  The first HKC model is defined via
	\[ \mathscr{M}_{u}^{1} := \mathscr{M}_{w}^{1} :=\{ (0,1), (1,1) \} \text{ , } \hspace{.5 cm} \mathscr{M}_{\theta}^{1} := \{ (1,1), (0,2) \} , \]
	For $M > 1$ and $m_1^M > 1$, the $M^{th}$ HKC model is defined by simply adjoining the next mode $\textbf{m}^M$ to the $M-1^{th}$ model, as displayed in Figure \ref{fig:HKC_ModeSelect}:
	\[ \mathscr{M}_u^{M} := 
	\mathscr{M}_u^{M-1} \cup \{ \textbf{m}^{M} \} \text{ , } \mathscr{M}_w^{M} := 
	\mathscr{M}_w^{M-1} \cup \{ \textbf{m}^{M} \} \text{ , } \mathscr{M}_{\theta}^{M} := \mathscr{M}_{\theta}^{M-1} \cup \{ \textbf{m}^{M} \} . \]
	For the case $m_1^M = 1$, the HKC-$M$ model is defined by adjoining $\textbf{m}^M$, but one must additionally include the vertically stratified temperature mode required for energy consistency as in Criterion \ref{Crit:EnergyCrit}.  Although not required by Criteria \ref{Crit:EnergyCrit}, \ref{Crit:VortCrit} the shear velocity modes $u^{(0,m_3)}$, $w^{(0,m_3)}$, $w^{(m_1,0)}$ are also included at this stage so as to have a more complete model.  In particular, the choice to adjoin the $w^{(m_1,0)}$ modes beginning from HKC-2 was made so that the lowest order model is the well-known Lorenz-Stenflo model, which is useful for testing.  Therefore in this case the HKC model is defined as follows:
	\begin{align} \mathscr{M}_u^{M} := \mathscr{M}_u^{M-1} \cup \{ \textbf{m}^{M}, (0,2m^{M}_3-1) \} \hspace{.5 cm} & \text{ , } \hspace{.5 cm} \mathscr{M}_{\theta}^{M} := \mathscr{M}_{\theta}^{M-1} \cup \{ \textbf{m}^{M}, (0,2m^{M}_3) \} \text{ , } \notag \\ \mathscr{M}_w^{M} := \mathscr{M}_w^{M-1} \cup \{ \textbf{m}^{M} , & (0,2m^{M}_3-1) , (m_3^M-1,0) \} \text{ . } \notag \end{align}
	
	Due to relative simplicity of the lowest level model of this hierarchy, the HKC-1 model is given here explicitly:
	\begin{equation}
		\label{def:HKC1}
		\begin{split} 
			\frac{d}{dt} \theta^{(0,2)} & = - 4 \theta^{(0,2)} +\frac{ \Sha}{\sqrt{2(\Sha^2+1)}V} u^{(1,1,1)} \theta^{(1,1)}  \\
			\frac{d}{dt} u^{(1,1)} & = - \Pra (\Sha^2 +1) u^{(1,1)} -  \frac{\Pra \Ray \Sha}{\sqrt{\Sha^2+1}} \theta^{(1,1)} + \frac{\Pra \Rot}{\sqrt{\Sha^2+1}} w^{(1,1)} \\
			\frac{d}{dt} w^{(1,1)} & = - \Pra (\Sha^2 +1) w^{(1,1)} - \frac{\Pra \Rot}{\sqrt{\Sha^2+1}} u^{(1,1)} \\
			\frac{d}{dt} \theta^{(1,1)} & = - (\Sha^2 +1) \theta^{(1,1)} -  \frac{\Sha}{\sqrt{\Sha^2+1}} u^{(1,1)} - \frac{\Sha}{\sqrt{2(\Sha^2+1)}V} u^{(1,1)} \theta^{(0,2)} .
		\end{split} 
	\end{equation}
	together with the equations
	\[ \frac{d}{dt} u^{(0,1)} = - \Pra u^{(0,1)} + \Pra \Rot w^{(0,1)} \hspace{.5 cm} \text{ , } \frac{d}{dt} w^{(0,1)} = - \Pra w^{(0,1)} - \Pra \Rot u^{(0,1)} . \] 
	which are decoupled from the rest of the system and hence $u^{(0,1)}, w^{(0,1)}$ decay exponentially.  As noted this model is essentially the Lorenz-Stenflo model \cite{LStenflo_1996}, aside from a linear change of variables and the trivial addition of $u^{(0,1)}, w^{(0,1)}$.
	
	
	In general the $M$th model in the HKC hierarchy is a system of ODE's of dimension
	\begin{equation} \label{HKC_Dimension} 3 M  + 4 \Big \lfloor \frac{\sqrt{8M +1 } -1}{2} \Big \rfloor - 1 , \end{equation}
	and one can define an analogous Nusselt number as follows:
	\[ \mathsf{Nu}^M = 1 + \frac{\Sha}{2 \pi^2} \overline{\langle u^{M}_3 \theta^M \rangle} . \]
	On the other hand, note that by taking the infinite time averages of the evolution equations in \eqref{EnergyEvolution}, one obtains the following equivalent expressions for the Nusselt number in terms of $\theta$:
	\[ \mathsf{Nu} = 1 + \frac{\Sha}{2\pi^2} \overline{\langle u_3 \theta \rangle} = 1 + \frac{\Sha}{2\pi^2} \overline{\langle | \nabla \theta |^2 \rangle} = 1 + \frac{\Sha}{2\pi} \overline{\langle (1-\frac{x_3}{\pi}) \partial_{x_3}^2 \theta \rangle} . \]
	These various expressions are useful in different contexts.  For instance while the first has a clear physical interpretation as the average vertical heat transport, it is clear from the second that the Nusselt number must be greater than or equal to 1, i.e. a convective flow transports more heat than the pure conductive state.  Since the same evolution equations hold for the HKC models, one has the analogous equivalent expressions.  In the computational work below, the third expression will be most useful, since it expresses the Nusselt number in terms of the least number of modes.  Furthermore the computational work involves only finite time approximations of the infinite time average.  Hence the finite time Nusselt number is defined via the third expression as follows:
	\begin{equation} \label{HKC_Nusselt} \mathsf{Nu}^M(t) := 1 - \sum_{\textbf{m} \in \mathscr{M}_{\theta}^M \cap \mathscr{M}^{*}} \frac{\sqrt{\Sha}m_3}{\pi} \frac{1}{t} \int_0^t \theta^{\textbf{m},M}(s)ds . \end{equation}
	
	\section{Existence of solutions and attractors} 
	
	\label{sec:BCAnalysis}
	
	\subsection{Well-posedness and attractor theory }
	
	\revised{In this section it is shown that the above PDE \eqref{EvolutionPDE} is well-posed and that solutions to the HKC hierarchy converge to solutions of \eqref{EvolutionPDE}.  Furthermore it shown that the PDE admits a compact global attractor, hence the Nusselt number \eqref{NusseltDef} is well-defined.}  Since the fluid is assumed to be horizontally aligned, the well-posedness theory for the model \eqref{EvolutionPDE} is similar to that in the two-dimensional, non-rotating case (for example \cite{temam_InfDimDynSys}).  However, the proof here is in a slightly different context (3d velocity, non-zero rotation), contains more details, and provides further regularity such that \eqref{Balance_KinEnergy} - \eqref{Balance_Vort} are satisfied in a classical sense.  First, some classical inequalities will be used in the existence proof, given in the following Lemma.  Since an explicit basis has been given, these inequalities can be established by referring to the Fourier basis, hence the proofs are omitted.
	\begin{Lemma}
		\label{lem:Ineqs}
		
		\textbf{(a) (Poincar\'e)} For $\textbf{u} \in \textbf{H}_{\sigma}^1$, $\theta \in H_0^1$ one has
		\begin{equation} \label{PoincareIneq} \|u_3 \|_{L^2} \leq \|\nabla u_3 \|_{L^2} \hspace{.15 cm} \text{ , } \hspace{.15 cm} \| \theta \|_{L^2} \leq \|\nabla \theta \|_{L^2} \hspace{.15 cm} \text{ , } \hspace{.15 cm} \|u_i \|_{L^2} \leq \frac{1}{\min(1,\Sha^2)} \|\nabla u_i \|_{L^2} \text{ for } i =1,2 .\end{equation}
		
		\noindent \textbf{(b) (Gagliardo-Nirenburg-Sobolev)} Let $p \in [2,\infty)$ and $\sigma = \frac{p-2}{p}$.  There exists a constant $C > 0$ such that for any $u \in H^1$ having mean-zero one has
		\begin{equation} \label{GagliardoNirenburgSobolevIneq} \| u\|_{L^p} \leq C \|u\|_{L^2}^{1-\sigma} \| \nabla u \|_{L^2}^{\sigma} . \end{equation}
		
		
		\noindent \textbf{(c) (Gronwall)} Suppose $v(t),a(t),\beta(t) \in C \big ( [0,\infty) ; \mathbb{R} \big )$ and $\beta(t)$ is non-negative.   If the left hand side of the following holds $\forall t \geq 0$, then the right hand side holds $\forall t \geq 0$:
		\begin{equation} \label{GronwallIneq} v(t) \leq a(t) + \int_0^t \beta(s) v(s) ds \hspace{.35 cm} \Rightarrow \hspace{.35 cm} v(t) \leq a(t) + \int_0^t a(s) \beta(s) \exp \big [ \int_s^t \beta(r) dr \big ] ds . \end{equation}
		
	\end{Lemma}
	
	The following theorem proves existence of smooth solutions for \eqref{EvolutionPDE} with Lipschitz dependence on initial conditions, which implies uniqueness.  A maximum principle is then proven, which together with the balances \eqref{Balance_KinEnergy} - \eqref{Balance_PotEnergy} imply the existence of an absorbing ball in $\textbf{H}^0$.  By bootstrapping to prove existence of an absorbing ball in $\textbf{H}^1$, we show that all solutions are in fact global and obtain a compact global attractor in $\textbf{H}^0$.
	
	\begin{Theorem}
		\label{thm:PDE_WellPosedness_Balances}
		
		\textbf{(a)} Given any admissible $\parameters$, $\textbf{X}_0 = (\textbf{u}_0, \theta_0) \in \textbf{H}^0$, $\tau > \delta > 0$ and $\kappa \in \mathbb{Z}_{\geq 0}$ there exists a unique $\textbf{X}(t) =(\textbf{u}(t),\theta(t)) \in C( [0,\tau];\textbf{H}^0) \cap L^2((0,\tau);\textbf{H}^1) \cap C^{\kappa}([\delta,\tau];\textbf{C}^{\kappa})$ such that $\textbf{X}(0) = \textbf{X}_0 $ which solves \eqref{EvolutionPDE} for $t \in (0,\tau)$, and satisfies \eqref{Balance_KinEnergy} - \eqref{Balance_Vort} in the classical sense.  Furthermore, one has Lipschitz dependence on initial conditions, i.e. letting $\textbf{X}(t), \textbf{X}^*(t)$ denote solutions with initial conditions $\textbf{X}_0,\textbf{X}_0^*$, one has
		\begin{equation} \label{Thm1_LipDepInitCond} \begin{split} \| \textbf{X}(t) - \textbf{X}^*(t) \|_{\textbf{H}^0}^2 & \leq \| \textbf{X}_0 - \textbf{X}^*_0 \|_{\textbf{H}^0}^2 \Big ( 1 + p(t) \exp \big [ p(t) \big ] \Big ) \hspace{.25 cm} \text{ for a.e. } t\in [0,\tau] \text{ , } \\  & \text{ where } \hspace{.25 cm} p(t) = (\Pra \Ray + 1)t + \frac{C \| \textbf{X}_0 \|_{\textbf{H}^0}^2 }{\mathsf{min}(\Pra^2,1)} e^{(\Pra \Ray +1 )t} \text{ . } \end{split} \end{equation}
		Thus the semi-group $ \mathcal{S}(t) \big [ \textbf{X}_0 \big ] :=  \textbf{X}(t)$ is well defined and continuous.

		\noindent \textbf{(b)} If $\| \theta_0 \|_{L^{\infty}} \leq \pi$ then $\| \theta (t) \|_{L^{\infty}} \leq \pi$ for a.e. $t > 0$.  On the other hand, if $\| \theta_0 \|_{L^{\infty}} \leq \pi$ is not assumed, then 
		\[ \theta(t) = \Theta(t) + \tilde{\theta}(t) , \]
		where $\| \Theta(t) \|_{L^{\infty}} \leq \pi$ for a.e. $t > 0$ and 
		\[ \| \tilde{\theta}(t) \|_{L^2} \leq C (1 + \|\theta_0 \|_{L^2} ) e^{-t} .  \]
		
		\noindent \textbf{(c)} There exists an absorbing ball in $\textbf{H}^1$ for the semi-group $\mathcal{S}(t)$.  Thus one can take $\tau = \infty$ and $\mathcal{S}(t)$ admits a compact global attractor $\mathscr{A}$.
		
	\end{Theorem}
	
	\begin{proof}
		
		For $\textbf{X}_0 \in \textbf{H}^0$ let $\textbf{X}^M_0 = (\textbf{u}_0^M, \theta_0^M)$ denote the projection onto the modes defined by the $M$th index set in the hierarchy.  Note that the Galerkin solutions satisfy
		\[ \frac{1}{2} \frac{d}{dt} \| \textbf{u}^M \|_{\textbf{L}^2}^2 + \Pra \|\nabla \textbf{u}^M \|_{\textbf{L}^2}^2 = \Pra \Ray \langle u_3^M \theta^M \rangle \hspace{.5 cm} \text{ , } \hspace{.5 cm}  \frac{1}{2} \frac{d}{dt} \| \theta^M \|_{L^2}^2 +  \|\nabla \theta^M \|_{L^2}^2  = \langle u_3^M \theta^M \rangle , \]
		hence integrating and applying Young's inequality one obtains 
		\begin{align} \frac{1}{2} \| \textbf{u}^M (t) \|_{\textbf{L}^2}^2 + \Pra \int_0^t \|\nabla \textbf{u}^M (s) \|_{\textbf{L}^2}^2ds & \leq \frac{1}{2} \| \textbf{u}^M_0 \|_{\textbf{L}^2}^2+ \frac{\Pra \Ray }{2} \int_0^t  \big [ \| u_3^M (s)\|_{L^2}^2 +  \| \theta^M (s)\|_{L^2}^2 \big ] ds , \notag \\ \label{EnergyInequality} \frac{1}{2} \| \theta^M(t) \|_{L^2}^2 + \int_0^t \|\nabla \theta^M(s) \|_{L^2}^2 ds & \leq \frac{1}{2} \| \theta^M_0 \|_{L^2}^2 + \frac{1}{2} \int_0^t  \big [ \| u_3^M (s)\|_{L^2}^2 +  \| \theta^M (s)\|_{L^2}^2 \big ] ds . \end{align}
		By adding these two inequalities and neglecting the integral terms on the left hand side, one obtains the following by using Gronwall's inequality \eqref{GronwallIneq}
		\begin{equation} \label{Thm1_Bound1} \| \textbf{X}^M (t) \|_{\textbf{H}^0}^2 \leq \| \textbf{X}^M_0 \|_{\textbf{H}^0}^2 \exp \big [ (\Pra \Ray + 1) t \big ] \leq \| \textbf{X}_0 \|_{\textbf{H}^0}^2 \exp \big [ (\Pra \Ray + 1) \tau \big ] . \end{equation}
		This bound is independent of $t$ and $M$, hence for any $\tau > 0$ all of the Galerkin solutions are uniformly bounded in $L^{\infty}( [0,\tau]; \textbf{H}^0 )$.  Inserting \eqref{Thm1_Bound1} into \eqref{EnergyInequality}, one obtains a bound on the $L^{2}( [0,\tau]; \textbf{H}^1 )$ norm which is independent of $t,M$:
		\begin{equation} \label{Thm1_Bound2} \int_0^t \| \textbf{X}^M (s) \|_{\textbf{H}^1}^2 ds \leq \frac{1}{2 \hspace{.5 mm} \mathsf{min}(\Pra,1)} \| \textbf{X}_0 \|_{\textbf{H}^0}^2 \Big ( 1 + (\Pra \Ray + 1) \tau \exp \big [ (\Pra \Ray + 1) \tau \big ] \Big ). \end{equation}
		Having proven $\textbf{X}^M$ is uniformly bounded in $L^{\infty}([0,\tau];\textbf{H}^0) \cap L^2((0,\tau);\textbf{H}^{1})$, one can inductively show that for any $k > 0$ the sequence $t^{\frac{k}{2}} \textbf{X}^M$ is uniformly bounded in $L^{\infty}((0,\tau);\textbf{H}^k) \cap L^{2}((0,\tau);\textbf{H}^{k+1})$, which implies that for any $0 < \delta < \tau$ the sequence $\textbf{X}^M$ is uniformly bounded in $L^{\infty}((\delta,\tau);\textbf{H}^k) \cap L^{2}((\delta,\tau);\textbf{H}^{k+1})$.  The method is the same as the derivation of \eqref{Thm1_Bound1}, \eqref{Thm1_Bound2}, but involves more complicated arithmetic, hence is carried out in \ref{app:SmoothSolns}.  From \eqref{TruncatedEvolutionPDE}, one sees that for any $\textbf{Y} = (\textbf{v},\varphi) \in \textbf{H}^1 $ one has
		\begin{equation} \label{UniformBoundTimeDerivative} \begin{split} \langle \frac{\partial\textbf{u}^M}{\partial t}  \cdot \textbf{v} \rangle & = - \Pra \langle \nabla \textbf{u}^M : \nabla \mathcal{P}_{\textbf{u}}^M [\textbf{v}] \rangle +  \Pra \Ray \langle \theta^M \uVecThree \cdot \mathcal{P}_{\textbf{u}}^M[\textbf{v}] \rangle \\ & \hspace{.5 cm} - \Pra \Rot \langle ( \uVecThree \times \textbf{u}^M ) \cdot \mathcal{P}_{\textbf{u}}^M[\textbf{v}]\rangle - \langle  \mathcal{P}_{\textbf{u}}^M [\textbf{v}] \cdot ( \textbf{u}^M \cdot \nabla \textbf{u}^M ) \rangle , \\ \langle \frac{\partial\theta^M}{\partial t} \varphi \rangle & = - \langle \nabla \theta^M \cdot \nabla \varphi \rangle + \langle u^M_3 \mathcal{P}_{\theta}^M[\varphi ] \rangle - \langle  \mathcal{P}_{\theta}^M [\varphi] ( \textbf{u}^M \cdot \nabla \theta^M ) \rangle . \end{split} \end{equation}
		Note that by integrating by parts, using Cauchy-Schwarz and \eqref{GagliardoNirenburgSobolevIneq} one obtains
		\begin{equation} \label{NonlinearTermBound} \begin{split} \big | \langle  \mathcal{P}_{\textbf{u}}^M [\textbf{v}] \cdot ( \textbf{u}^M \cdot \nabla \textbf{u}^M ) \rangle \big | & \leq  \| \textbf{u}^M \|_{L^4}^2 \| \nabla \textbf{v} \|_{L^2} \leq C \| \textbf{u}^M \|_{L^2} \| \nabla \textbf{u}^M \|_{L^2} \| \nabla \textbf{v} \|_{L^2}  , \\ \big | \langle  \mathcal{P}_{\theta}^M [\varphi] ( \textbf{u}^M \cdot \nabla \theta^M ) \rangle \big | & \leq  \| \textbf{u}^M \|_{L^4} \| \theta^M \|_{L^4}  \| \nabla \varphi \|_{L^2} \\ & \leq C \| \textbf{u}^M \|_{L^2}^{1/2} \| \nabla \textbf{u}^M \|_{L^2}^{1/2} \| \theta^M \|_{L^2}^{1/2} \| \nabla \theta^M \|_{L^2}^{1/2} \| \nabla \varphi \|_{L^2} . \end{split} \end{equation} 
		One can then square both sides of \eqref{UniformBoundTimeDerivative}, integrate over $[0,\tau]$, and bound all of the terms on the right hand side using \eqref{Thm1_Bound1},\eqref{Thm1_Bound2}, \eqref{NonlinearTermBound}.  One thus obtains a uniform bound for $\partial_t \textbf{X}^M $ in the space $L^{2}( [0,\tau]; \textbf{H}^{-1} )$, hence by the Aubin-Lions lemma \cite{Simon_1986} the collection of Galerkin solutions is relatively compact in $L^{2}( [0,\tau]; \textbf{H}^0 )$.  Thus a sub-sequence converges strongly to some $\textbf{X}(t)$ in $L^{2}( [0,\tau]; \textbf{H}^0 )$.  In the same way, for any $k>1$ one can obtain a uniform bound for $\partial_t \textbf{X}^M $ in the space $L^{2}( [\delta,\tau]; \textbf{H}^{-k} )$, and together with the uniform bound in $L^{2}( [\delta,\tau]; \textbf{H}^{k} )$ the Aubin-Lions Lemma says that any sub-sequence which converges strongly in $L^{2}( [0,\tau]; \textbf{H}^0 )$ must have a further sub-sequence which converges strongly in $L^{2}( [\delta,\tau]; \textbf{H}^{k-1} )$. 
		
		In fact, every sub-sequential limit must belong to the more regular subspace $L^{\infty}( [0,\tau];\textbf{H}^0) \cap L^{2}( [0,\tau];\textbf{H}^{1})$.  Let $\textbf{X}^{M_i}$ be a sub-sequence which converges strongly to some $\textbf{X}(t)$ in $L^{2}( [0,\tau]; \textbf{H}^0 )$.  Since this sub-sequence is again uniformly bounded in $L^{\infty}( [0,\tau]; \textbf{H}^0)$, it is relatively compact in the weak-* topology due to Banach-Alaoglu.  Hence there must be some point $\textbf{X}^*(t) \in L^{\infty}( [0,\tau];\textbf{H}^0)$ and a further sub-sequence such that for all $\textbf{Y}(t) \in L^{1}( [0,\tau];\textbf{H}^0)$ one has
		\begin{equation} \label{Thm1_SubseqConv1} \lim_{j \to \infty} \int_0^{\tau} \big \langle \big ( \textbf{X}^{M_{i_j}}(t) - \textbf{X}^{*}(t) \big ) \cdot \textbf{Y}(t) \big \rangle dt = 0 . \end{equation}
		Since one has strong convergence this gives that for all $\textbf{Y}(t) \in L^{2}( [0,\tau];\textbf{H}^0)$ one has
		\begin{equation} \label{Thm1_SubseqConv2} \int_0^{\tau} \big \langle \big  ( \textbf{X}(t) - \textbf{X}^{*}(t) \big ) \cdot \textbf{Y}(t) \big \rangle dt  , \end{equation}
		and hence $\textbf{X}(t) = \textbf{X}^*(t) $ belongs to the subspace $L^{\infty}( [0,\tau];\textbf{H}^0)$.  Similarly, since this sub-sequence is uniformly bounded in $L^{2}( [0,\tau];\textbf{H}^1)$, it is relatively compact in its weak-* topology.  However, since this is a Hilbert space this is the same as its weak topology, hence there must be some (apriori possibly distinct) point $\textbf{X}^{*}(t) \in L^{2}( [0,\tau];\textbf{H}^1)$, and a further sub-sequence such that for all $\tilde{\textbf{Y}}(t) \in L^{2}( [0,\tau];\textbf{H}^1)$ one has the following, where $:$ denotes the Frobenius product:
		\[ \lim_{j \to \infty} \int_0^{\tau} \big \langle \nabla \big ( \textbf{X}^{M_{i_j}}(t) - \textbf{X}^{*}(t) \big ) : \nabla \tilde{\textbf{Y}}(t) \big \rangle dt = 0 . \]
		Since $\Delta:H^1 \mapsto H^{-1}$ is an isometry this implies \eqref{Thm1_SubseqConv1} holds for all $\textbf{Y}(t) = \Delta \tilde{\textbf{Y}}(t) \in L^{2}( [0,\tau];\textbf{H}^{-1})$, thus \eqref{Thm1_SubseqConv2} holds for such $\textbf{Y}(t)$ and hence $\textbf{X}(t) = \textbf{X}^*(t) \in L^{2}( [0,\tau];\textbf{H}^1)$.  In the same way, one can show that the uniform bounds of $\textbf{X}^M$ in $L^{\infty}((\delta,\tau),\textbf{H}^{k}) \cap L^{2}((\delta,\tau),\textbf{H}^{k+1}) $ imply that the sub-sequential limits belong to these spaces as well.
		
		Testing \eqref{TruncatedEvolutionPDE} with $\textbf{Y} \in C^{\infty}([0,\tau]; \textbf{C}^{\infty} )$ one can then use the respective strong, weak and weak-* convergence to show that \eqref{EvolutionPDE} holds for any sub-sequential limit $\textbf{X}(t)$, where $\partial_t \textbf{u}$, $\partial_t \theta$ are the distributional derivatives in time.  Since \eqref{EvolutionPDE}, and all of the terms on the right hand side have been shown to have certain regularity, it then becomes clear that the distributional derivatives $\partial_t \textbf{X}(t)$ belong to $L^2([0,\tau]; \textbf{H}^{-1})$, and using the Lions-Magenes lemma one can conclude that $\textbf{X}(t) \in C( [0,\tau];\textbf{H}^0)$, and that the following holds in a distributional sense
		\[ \frac{1}{2} \frac{d}{dt} \| \textbf{u}\|_{\textbf{L}^2}^2 = \langle \partial_t \textbf{u}(t) , \textbf{u}(t) \rangle \hspace{.5 cm} \text{ , } \hspace{.5 cm} \frac{1}{2}  \frac{d}{dt} \| \theta\|_{L^2}^2 = \langle \partial_t \theta(t) , \theta(t) \rangle \text{ . } \]
		In fact, since these sub-sequential limits belong to $L^{\infty}((\delta,\tau),\textbf{H}^k)$ for any $k$ implies that they belong to $C^{\iota}((\delta,\tau);\textbf{C}^{\iota})$ due to the Sobolev embedding and the fact that the time derivatives can be expressed in terms of the spatial derivatives by virtue of the differential equation \eqref{EvolutionPDE}.  Since these sub-sequential limits solve \eqref{EvolutionPDE} in a strong sense for $t \in [\delta, \tau]$, it follows that \eqref{Balance_KinEnergy} - \eqref{Balance_Vort} hold as well.
		
		Having proven existence and regularity, one can show the Lipschitz dependence on initial conditions \eqref{Thm1_LipDepInitCond}, from which uniqueness clearly follows.  This is proven in \ref{app:LipDepOnIC_MaxPrinc}, along with the maximum principle in Theorem \ref{thm:PDE_WellPosedness_Balances}, part (b).
		
		Next, it is proven that the evolution defined by \eqref{EvolutionPDE} admits an attracting, forward invariant ball in $\textbf{H}^0$.  Using (b) one can prove that for any $\epsilon > 0$ there exists $t^*(\epsilon) > 0$ such that $\| \theta(t) \|_{L^{2}} \leq \sqrt{ \pi |\Omega| } + \epsilon $ for all $t > t^{*}(\epsilon)$.  One can therefore prove the velocity eventually enters a ball in $L^2$, where it remains ever after as follows.  Testing the momentum equation with $\textbf{u}$, one obtains
		\begin{align} \frac{1}{2 \Pra} \frac{d}{dt} \|\textbf{u}\|_{\textbf{L}^2}^2 & = - \langle |\nabla \textbf{u}|^2 \rangle + \Ray \langle u_3 \theta \rangle \leq -  \min(1,\Sha^2) \langle u_1^2 + u_2^2 \rangle -  \langle u_3^2 \rangle  + \Ray \|\theta \|_{L^{2}} \|u_3\|_{L^2} \notag \\ & \leq - \min(1,\Sha^2) \langle u_1^2 + u_2^2 \rangle - \langle u_3^2 \rangle + \Ray ( \sqrt{ \pi |\Omega| } +\epsilon) \|u_3 \|_{L^2} \notag \\ \label{PDE_L2AbsorbingSet} & \leq - \min(1,\Sha^2) \langle u_1^2 + u_2^2 \rangle - \big ( \|u_3\|_{L^2} - \frac{1}{2}\Ray (\sqrt{ \pi |\Omega| }+\epsilon) \big )^2 + \frac{\Ray^2 (\sqrt{ \pi |\Omega| }+\epsilon)^2}{4} . \end{align}
		Thus the radius of $(\|u_1\|_{L^2},\|u_2\|_{L^2},\|u_3\|_{L^2})$ decreases unless the norms of the solution lie in the ellipsoid centered at $(\|u_1\|_{L^2},\|u_2\|_{L^2},\|u_3\|_{L^2}) = (0,0,\frac{1}{2}\Ray (\sqrt{ \pi |\Omega| } + \epsilon))$.  Therefore any larger ball centered at the origin containing this ellipsoid in its interior is an attracting, forward invariant set, for instance the ball $\mathscr{B}(0, \Ray (\sqrt{ \pi |\Omega| } + 2\epsilon ) )$.
		
		Finally, one can bootstrap the result to prove the evolution defined by \eqref{EvolutionPDE} admits an attracting, forward invariant ball in $\textbf{H}^1_{\sigma}\times H_0^1$.  First, by integrating the equations for the $L^2$ norms from time $t$ to $t+1$, one obtains
		\begin{equation} \notag 
			\begin{split} \frac{1}{2\Pra} \big ( \|\textbf{u}(t+1)\|_{\textbf{L}^2}^2 - \|\textbf{u}(t)\|_{\textbf{L}^2}^2 \big ) & = -\int_{t}^{t+1} \|\nabla \textbf{u}(s) \|_{\textbf{L}^2}^2 ds + \Ray \int_t^{t+1} \langle u_3(s) \theta(s) \rangle ds , \\ \frac{1}{2} \big ( \|\theta(t+1)\|_{L^2}^2 - \|\theta(t)\|_{L^2}^2 \big ) & = -\int_{t}^{t+1} \|\nabla \theta(s) \|_{L^2}^2 ds + \int_t^{t+1} \langle u_3(s) \theta(s) \rangle ds  ,
			\end{split} 
		\end{equation}
		hence for $t > t^*(\epsilon)$ one can use Cauchy-Schwarz, Young's inequality and the bounds on the $L^2$ norms of $\textbf{u}, \theta$ to obtain
		\begin{equation} \label{Attractor_IntegralBound} \begin{split}
				\int_{t}^{t+1} \|\nabla \textbf{u}(s) \|_{\textbf{L}^2}^2 ds & \leq \frac{\Ray}{2} ( \Ray +1+\frac{2}{\Pra} )(\sqrt{ \pi |\Omega| } + 2\epsilon ) \text{ , } \\ \int_{t}^{t+1} \|\nabla \theta(s) \|_{L^2}^2 ds & \leq \frac{1}{2} (\Ray + 3)( \sqrt{ \pi |\Omega| } + 2\epsilon ) .	\end{split} \end{equation}
		Next, by testing the momentum equation with $-\Delta \textbf{u}$, then using Cauchy-Schwarz, Poincar\'e and Young's inequalities together with taking $\textbf{v} = \Delta \textbf{u}$ in \eqref{NonlinearTermBound}, one obtains 
		\begin{equation} \notag
			\begin{split} \frac{1}{2 \Pra} \frac{d}{dt} \| \nabla \textbf{u}\|_{\textbf{L}^2}^2 & = - \|\Delta \textbf{u} \|_{L^2}^2 + \frac{1}{\Pra} \langle \big ( (\textbf{u}\cdot \nabla ) \textbf{u} \big ) \cdot \Delta \textbf{u} \rangle + \Ray \langle \nabla u_3 \cdot \nabla \theta \rangle \\ & \leq - \| \Delta \textbf{u} \|_{L^2}^2 + \frac{1}{\Pra} \| (\textbf{u}\cdot \nabla ) \textbf{u} \|_{L^2} \| \Delta \textbf{u} \|_{L^2} + \Ray \| \nabla u_3 \|_{L^2} \| \nabla \theta \|_{L^2} \\ & \leq - \frac{1}{2} \| \Delta \textbf{u} \|_{L^2}^2 + \frac{1}{\Pra} \| \textbf{u} \|_{L^2}^{1/2} \| \nabla \textbf{u} \|_{L^2} \| \Delta \textbf{u} \|_{L^2}^{3/2} + \frac{\Ray^2}{2 \min (1,\Sha^4)} \| \nabla \theta \|_{L^2}^2 \\ & \leq - \frac{1}{4} \| \Delta \textbf{u} \|_{L^2}^2 + \frac{27}{4\Pra^4} \| \textbf{u} \|_{L^2}^{2} \| \nabla \textbf{u} \|_{L^2}^4 + \frac{\Ray^2}{2 \min (1,\Sha^4)} \| \nabla \theta \|_{L^2}^2 .
			\end{split} 
		\end{equation}
		Hence one has
		\[ \frac{d}{dt} \| \nabla \textbf{u}\|_{\textbf{L}^2}^2 \leq \Big ( \frac{27}{2\Pra^3} \| \textbf{u} \|_{L^2}^{2} \| \nabla \textbf{u} \|_{L^2}^2 \Big ) \| \nabla \textbf{u} \|_{L^2}^2 + \frac{\Pra \Ray^2}{\min (1,\Sha^4)} \| \nabla \theta \|_{L^2}^2 , \]
		and together with \eqref{Attractor_IntegralBound} it follows from the Gronwall inequality, that for $t > t^*(\epsilon) + 1$ one has
		\[ \| \nabla \textbf{u} (t) \|_{\textbf{L}^2}^2 \leq \big ( \frac{\Pra \Ray^2(\Ray +3)}{2 \min(1,\Sha^4)}  + \frac{\Ray}{2} (\Ray +1+\frac{2}{\Pra} ) \big ) (\sqrt{ \pi |\Omega| } + 2\epsilon ) \exp \Big [ \frac{27 \Ray^3}{4 \Pra^3} ( \pi |\Omega|  + 2\epsilon )^3 (\Ray + 1 + \frac{2}{\Pra} ) \Big ] . \]
		The same argument can be used to obtain a bound on $\|\nabla \theta\|_{L^2}^2$, thus proving the result.
		
		Finally, since the attracting ball $\mathscr{B} \subset \textbf{H}^1$ is compactly embedded in $\textbf{H}^0$, the proof $\mathscr{A} = \omega(\mathscr{B})$ has the properties of the global attractor is the same as in the finite dimensional case above.  
	\end{proof}

\section{The dimension of the attractor and an analysis of bifurcations }

\label{sec:AttractorAnalysis}

As the above proof of Theorem \ref{thm:PDE_WellPosedness_Balances} has shown, one has the convergence result that for any $\tau > 0$ and any $\epsilon > 0$ there exists a $M^*(\tau,\epsilon) > 0$ such that for all $M > M^*(\tau,\epsilon)$ one has
\begin{equation} \label{StrongL2Convergence} \|\textbf{X} - \textbf{X}^M \|_{L^2([0,\tau];\textbf{H}^0)} < \epsilon . \end{equation}
Hence with a solution to a sufficiently large HKC model one can approximate the Nusselt number via 
\[ \big | \frac{1}{\tau} \int_0^{\tau} \int_{\Omega} \big [ u_3 \theta - u_3^M \theta^M \big ] dx dt \big | \leq \frac{1}{\tau^{1/2}}\|\textbf{X} \|_{C([0,\tau];\textbf{H}^0)} \|\textbf{X} - \textbf{X}^M \|_{L^2([0,\tau];\textbf{H}^0)} < \frac{C \epsilon}{\tau^{1/2}} , \]
where the constant $C$ depends only on the initial condition and the radii of the absorbing balls above.  However, \revised{note that the above proof makes use of the Aubin-Lions lemma.  Unlike contraction mapping techniques, this abstract result does not provide any information about the rate of convergence}.  Hence even for finite time averages the number of necessary modes $M^*(\tau,\epsilon)$ is unknown.  Furthermore, it is unclear whether the above result about finite time approximation is sufficient for an infinite time average such as the Nusselt number, since in the worst case one could have $M^*(\tau,\epsilon) \to \infty$ as $\tau \to \infty$.


This section therefore aims at studying the properties of the global attractor for \eqref{EvolutionPDE} in order to provide a stronger theoretical basis for the numerical heat transport analysis in section \ref{sec:NumAnalysis}.  This global attractor likely has a highly complicated, fractal structure which can be difficult to study, so first we study the local bifurcations which give rise to the unstable manifold at the origin.  This is a smooth sub-manifold of the attractor about which analytical statements can be made.  In this way one obtains a nice lower bound on the dimension of the attractor which can help inform how large $M^*(\tau,\epsilon)$ might be.  Next, we adapt Temam's result on the Hausdorff dimension of the attractor to the present case, giving an upper bound on the complexity of the attractor as a function of the Rayleigh number.

\subsection{Local bifurcations at the origin and a lower bound on the attractor dimension}

\label{sec:OriginLocalBif}

The analysis of the unstable manifold at the origin begins with a detailed analysis of the linearization of \eqref{EvolutionPDE} about the origin, denoted $\mathcal{L}^{\textbf{0}}$.  Note only terms with the same wave vector $\textbf{m}$ are coupled via $\mathcal{L}^{\textbf{0}}$, hence $\mathcal{L}^{\textbf{0}}$ acts on $\textbf{L}^2_0\times L^2$ by multiplying each $\hat{\textbf{X}}^{\textbf{m}}$ by a corresponding matrix $\hat{\mathcal{L}}^{\textbf{0},\textbf{m}}$.  In the case $\textbf{m} \in \mathbb{Z}_{>0}^2$ these matrices are given explicitly as follows:
\begin{equation} \label{ExplicitLinearization}  \hat{\mathcal{L}}^{\textbf{0},\textbf{m}} =  \begin{pmatrix} -\Pra | \mathcal{K}\textbf{m}|^2 & \Pra \Rot \frac{m_3}{|\mathcal{K}\textbf{m}|} & (-1)^{|\textbf{m}|_1 +1} \Pra \Ray \frac{\Sha m_1}{|\mathcal{K}\textbf{m}|} \\ - \Pra \Rot \frac{m_3}{|\mathcal{K}\textbf{m}|} & -\Pra |\mathcal{K}\textbf{m}|^2 & 0 \\ (-1)^{|\textbf{m}|_1 +1} \frac{\Sha m_1}{|\mathcal{K}\textbf{m}|} & 0 & -|\mathcal{K}\textbf{m}|^2 \end{pmatrix} \text{ , } \end{equation}
whereas when either $m_1$ or $m_3$ is zero it is an appropriate submatrix (recall \eqref{PhaseIndexSets}).  Due to the block diagonal action, the spectrum of $\mathcal{L}^{\textbf{0}}$ is just the set of eigenvalues of $\hat{\mathcal{L}}^{\textbf{0},\textbf{m}}$, ranging across $\textbf{m} \in \mathbb{Z}_{\geq 0}^2$.  The interesting case for these eigenvalues is $\textbf{m} \in \mathbb{Z}_{>0}^2$, since in all other cases the eigenvalues have strictly negative real part for all parameter values.  In this case the eigenvalues $\lambda^{\textbf{m},j}$ must belong to at least one of the following discs in the complex plane, due to Gershgorin's theorem:
\begin{equation} \label{GershgorinDiscs}  |\lambda^{\textbf{m},j} + \Pra |\mathcal{K}\textbf{m}|^2| \leq \Pra (\Ray + \Rot) \hspace{.25 cm} \text{ , } \hspace{.25 cm} |\lambda^{\textbf{m},j} + \Pra |\mathcal{K}\textbf{m}|^2| \leq \Pra \Rot \hspace{.25 cm} \text{ , } \hspace{.25 cm} |\lambda^{\textbf{m},j} + |\mathcal{K}\textbf{m}|^2| \leq 1 \text{ . } \end{equation}
For fixed $\parameters$, the radii are fixed while the centers of these discs tend to $-\infty$ as $-|\mathcal{K}\textbf{m}|^2$ when $|\textbf{m}|$ becomes large.  

\begin{figure}[H]
	\begin{center}
		\begin{tabular}{cc}
			\includegraphics[height=45mm]{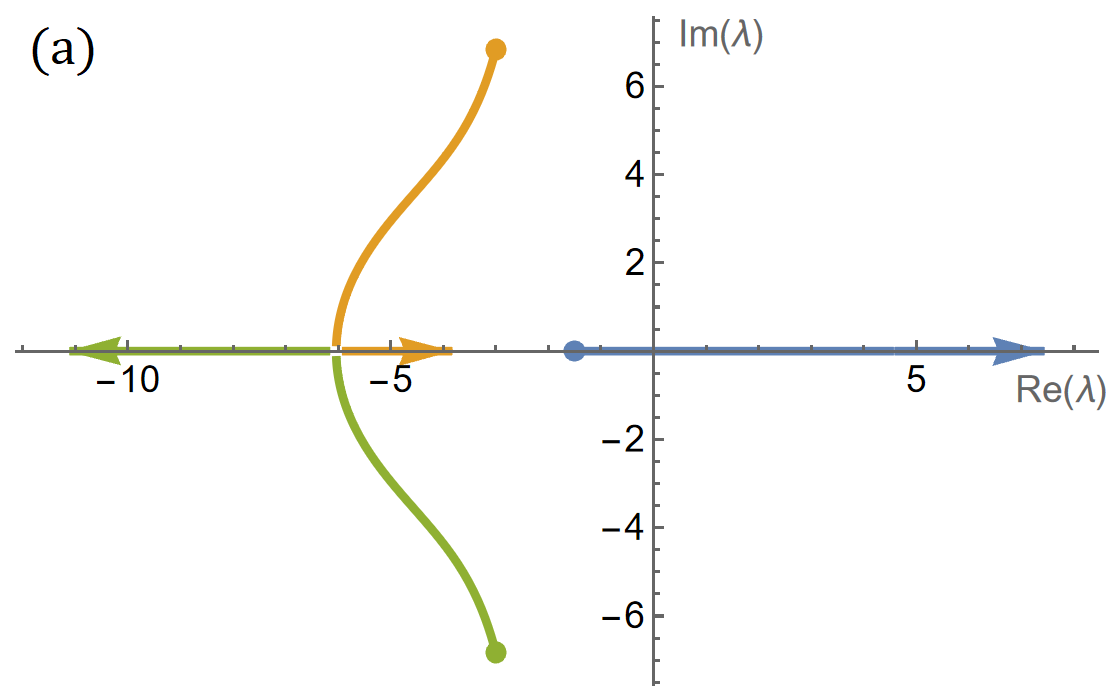}
			&\includegraphics[height=45mm]{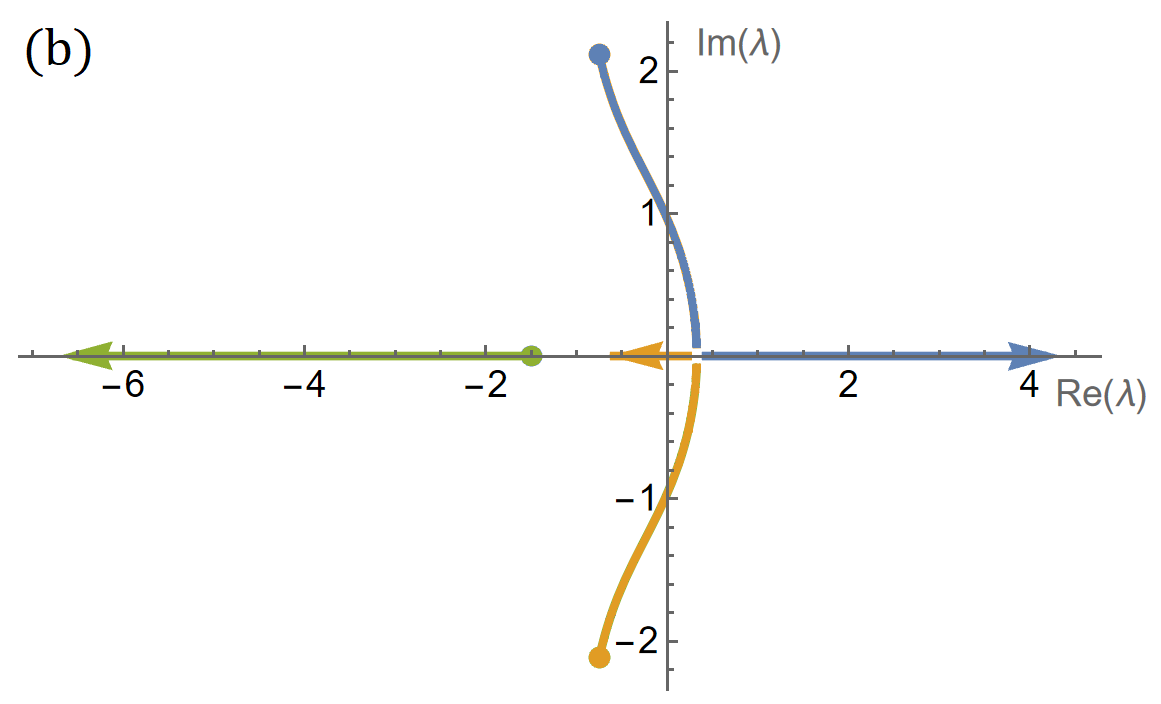}
		\end{tabular}
		\caption{Eigenvalues of the matrix $\mathcal{L}^{\textbf{0},(1,1)}$ as a function of $\Ray$, for $\Sha = 1/\sqrt{2}$.  (a) For $\Pra = 2$, $\Rot = \frac{9}{2\sqrt{2}}+1$, $0 \leq \Ray \leq 200$, a type 1 crossing occurs at $\Ray^{(1,1),1}$. (b) For $\Pra = 1/2$, $\Rot = \frac{9}{2\sqrt{2}}+2$, $0 \leq \Ray \leq 200$, a type 2 crossing occurs at $\Ray^{(1,1),2}$, and a type 3 crossing occurs at $\Ray^{(1,1),1}$. } 
	\label{fig:Eigenvalues}
\end{center}
\end{figure}

The following Lemma studies how these eigenvalues behave as $\parameters$ are varied.   The statements in this Lemma prove that as the Rayleigh number increases the eigenvalues of $\hat{\mathcal{L}}^{\textbf{0},\textbf{m}}$ must either appear as depicted in Figure \ref{fig:Eigenvalues} panel (a), or as in Figure \ref{fig:Eigenvalues} panel (b) in the special case of low Prandtl and high rotation.  The proof is given in \ref{app:OriginEigs}.

\begin{Lemma}
	\label{lem:OriginEigenvalues}  For any admissible parameters $\parameters$ and $m_1,m_3 \geq 1$, let the critical Rayleigh numbers $\Ray^{\textbf{m},1}, \Ray^{\textbf{m},2}$ be defined by
	\begin{equation} \label{CriticalRayleighNumber} \Ray^{\textbf{m},1} = \frac{|\mathcal{K}\textbf{m}|^6 + \Rot^2 m_3^2}{\Sha^2m_1^2} \text{ , } \hspace{.5 cm} \Ray^{\textbf{m},2} = \frac{2(\Pra +1)|\mathcal{K}\textbf{m}|^6 + \frac{2\Pra^2}{\Pra+1} \Rot^2 m_3^2 }{\Sha^2m_1^2} \text{ , } \end{equation}
	and let $\Ray^{\textbf{m},c} = \min (\Ray^{\textbf{m},1} , \Ray^{\textbf{m},2})$.  Furthermore let the rotation threshold $\Rot^{\textbf{m}}$ be defined by 
	\begin{equation} \label{RotThresh} \Rot^{\textbf{m}} = \left \{ \begin{array}{cl}
			\sqrt{\frac{1+\Pra}{1-\Pra}} \frac{|\mathcal{K}\textbf{m}|^3}{m_3} & \text{ if } 0 < \Pra < 1 , \\ 
			\frac{|\mathcal{K}\textbf{m}|^3}{2 m_3 \sqrt{\Pra (\Pra - 1)}} & \text{ if } \Pra > 1 , 
		\end{array}  \right . \end{equation}
	The following statements hold:
	\begin{enumerate}[label=(\roman*)]
		\item For all $ 0 \leq \Ray < \Ray^{\textbf{m},c}$, the matrix $\hat{\mathcal{L}}^{\textbf{0},\textbf{m}}$ has eigenvalues with strictly negative real parts, and for all $\Ray > \Ray^{\textbf{m},c}$ at least one has positive real part.
		\item For all $\Ray > \Ray^{\textbf{m},1}$, the matrix $\hat{\mathcal{L}}^{\textbf{0},\textbf{m}}$ has two (possibly complex valued) eigenvalues with strictly negative real parts, and a third real valued eigenvalue with positive real part.
		\item There exists a unique $\Ray^* := \Ray^*(\Pra,\Rot,\Sha,\textbf{m}) \geq 0$ such that $\mathcal{L}^{\textbf{0},\textbf{m}}$ has two complex conjugate eigenvalues for $0 \leq \Ray < \Ray^*$, and has three distinct real eigenvalues for $\Ray > \Ray^*$.  Furthermore $\Ray^* = 0$ if and only if $\Rot = 0$.
		\item If $0 < \Pra < 1$, then the rotation $\Rot$ determines whether a (real) eigenvalue or a conjugate pair of eigenvalues cross the imaginary axis:
		\begin{align}
			|\Rot| < \Rot^{\textbf{m}} \hspace{1 cm} & \Rightarrow \hspace{1 cm} \Ray^* < \Ray^{\textbf{m},1} < \Ray^{\textbf{m},2} \text{ , } \notag \\ |\Rot| = \Rot^{\textbf{m}}  \hspace{1 cm} & \Rightarrow \hspace{1 cm} \Ray^* = \Ray^{\textbf{m},1} = \Ray^{\textbf{m},2} \text{ , } \notag \\ |\Rot| > \Rot^{\textbf{m}}  \hspace{1 cm} & \Rightarrow \hspace{1 cm}  \Ray^{\textbf{m},2} < \Ray^* < \Ray^{\textbf{m},1} . \notag 
		\end{align}
		\item If $\Pra \geq 1$ then $\Ray^{\textbf{m},1} < \Ray^{\textbf{m},2}$.  When $\Pra = 1$, $\Ray^* < \Ray^{\textbf{m},1}$ whereas for $\Pra > 1$ one has
		\[ |\Rot| > \Rot^{\textbf{m}} \hspace{.25 cm} \Leftrightarrow \hspace{.25  cm} \Ray^* > \Ray^{\textbf{m},1} \hspace{.5 cm} \text{ , } \hspace{.5 cm} |\Rot| = \Rot^{\textbf{m}} \hspace{.25 cm} \Leftrightarrow \hspace{.25  cm} \Ray^* = \Ray^{\textbf{m},1} \text{ . } \]
	\end{enumerate}
\end{Lemma}

In particular, Lemma \ref{lem:OriginEigenvalues} characterizes how the eigenvalues can cross the imaginary axis as the Rayleigh number is increased, which will determine what kind of bifurcations can occur.  In particular, there are four types of crossing scenarios, which will simply be referred to as type 1,2,3,4.  The crossing depicted in Figure \ref{fig:Eigenvalues} panel (a) will be referred to as type 1, which occurs when $\Ray^{\textbf{m},1} < \Ray^{\textbf{m},2}$.  Lemma \ref{lem:OriginEigenvalues} (iv),(v) show that this is the more common crossing scenario.  The crossing of the conjugate eigenvalues in panel (b) will be referred to as type 2, which occurs in the more special case when $\Ray^{\textbf{m},2} < \Ray^{\textbf{m},1}$.  As shown in panel (b), subsequent to a type 2 crossing the complex eigenvalues meet and one of the real eigenvalues returns across the imaginary axis.  This will be referred to as type 3.  Finally, the crossing that occurs when $\Ray^{\textbf{m},1} = \Ray^{\textbf{m},2}$ will be referred to as type 4.

With the spectral picture in hand, the next ingredient required to classify the bifurcations at the origin is further regularity for the semi-group $\mathcal{S}(t)$ obtained in Theorem \ref{thm:PDE_WellPosedness_Balances}.  Although Theorem \ref{thm:OriginBifurcations} only uses the Frechet regularity of $\mathcal{S}(t)$ at the origin, the following Proposition provides global Frechet regularity used in Theorem \ref{thm:AttractorUppderBound}.  The proof is given in \ref{app:FurtherRegularity}.
\begin{Proposition}
\label{prop:FurtherRegularity}
For admissible $\parameters$ and $\textbf{X}_0 \in \textbf{H}^0$, let $\textbf{X}(t)$ be the solution of \eqref{EvolutionPDE} given in Theorem \ref{thm:PDE_WellPosedness_Balances}, and let $\mathcal{L}^{\textbf{X}_0}(t)$ be the linearization of \eqref{EvolutionPDE} along $\textbf{X}(t)$, given explicitly by 
\begin{equation} \label{LinearizedOperator} \hspace{.5 cm} \mathcal{L}^{ \textbf{X}_0}(t) \begin{pmatrix}
		\textbf{g} \\ \psi 
	\end{pmatrix} := \begin{pmatrix} \Pra \Delta \textbf{g} + \mathcal{P} \big [ \Pra \Ray \psi \uVecThree - \Pra \Rot \hspace{.5 mm} \uVecThree \times \textbf{g} - \textbf{u} \cdot \nabla \textbf{g} - \textbf{g} \cdot \nabla \textbf{u} \big ] , \\
		\Delta \psi + g_3 - \textbf{u} \cdot \nabla \psi - \textbf{g} \cdot \nabla \theta. \end{pmatrix} , \end{equation}
where $\mathcal{P}$ is the Leray projector.  For any  $\tau > \delta > 0$, $\kappa \in \mathbb{Z}_{\geq 0}, \textbf{Y}_0 \in \textbf{H}^0$ there exists a unique $\textbf{Y}(t) =(\textbf{g}(t),\psi(t)) \in C( [0,\tau];\textbf{H}^0) \cap L^2((0,\tau);\textbf{H}^1) \cap C^{\kappa}([\delta,\tau];\textbf{C}^{\kappa})$ with $\textbf{Y}(0) = \textbf{Y}_0$ which solves the following equation:
\begin{equation} \label{LinearizedEquation}  \frac{d}{dt} \textbf{Y}(t) = \mathcal{L}^{\textbf{X}_0}(t) \textbf{Y}(t) \text{ . } \end{equation}
The semi-group $\mathcal{S}(t)$ is uniformly Frechet differentiable at $\textbf{X}_0$, with Frechet derivative $\mathcal{F}^{\textbf{X}_0}(t)$ given by
\[ \mathcal{F}^{\textbf{X}_0}(t) \textbf{Y}_0 := \textbf{Y}(t) \text{ . }  \]
\end{Proposition}

In order to understand how the unstable manifold unfolds as the Rayleigh number is increased, we need a holistic view of the critical Rayleigh numbers $\Ray^{\textbf{m},c}$ across all wave numbers $\textbf{m}$.  We must therefore consider the level curves of $\Ray^{\textbf{m},1},\Ray^{\textbf{m},2}$ in the $m_1, m_3$ plane, as well as the boundary $(\Pra+1)|\mathcal{K}\textbf{m}|^6 = (1-\Pra)\Rot^2m_3^2 $ along which these are equal.  Specifically, let the level curves $m_3(m_1,\parameters), \tilde{m}_3(m_1,\parameters)$ solve
\begin{equation} \label{CriticalLevelSet} \begin{split} \Ray \Sha^2m_1^2 & = ( \Sha^2 m_1^2 + m_3^2 )^3 + \Rot^2 m_3^2 \hspace{3.65 cm} \text{ for } \hspace{.25 cm} 0 < \Sha m_1 <  \Ray^{1/4} , \\  \Ray \Sha^2m_1^2 & =  2(\Pra +1) ( \Sha^2 m_1^2 + \tilde{m}_3^2 )^3 + \frac{2\Pra^2}{\Pra +1 } \Rot^2 \tilde{m}_3^2 \hspace{1.4 cm} \text{ for } \hspace{.25 cm} 0 < \Sha m_1 <  \frac{\Ray^{1/4}}{(2\Pra +2)^{1/4}}. \end{split} \end{equation}
At any given parameter value $\parameters$, each of the variables $u^{\textbf{m}},\theta^{\textbf{m}}$ with $m_3$ lying below these level curves give rise to a linear (buoyancy) instability, which we will see corresponds one dimension of the unstable manifold.  Figure \ref{fig:CriticalLevelSets} below depicts the progression of instabilities as the Rayleigh number increases, corresponding to an increasingly high dimensional unstable manifold.  For $\Pra \geq 1$, the progression must be similar to that shown in (a), where an intersection a level curve of $\Ray^{\textbf{m},1}$ with an integer vertex $\textbf{m}$ corresponds to an eigenvalue of $\mathcal{L}^{\textbf{0},\textbf{m}}$ crossing the imaginary axis.  Here all of the crossings are of type 1.  For $0 < \Pra < 1$, the progression must be similar to that shown in (b), where the red curve depicts the boundary along which $\Ray^{\textbf{m},1}, \Ray^{\textbf{m},2}$ are equal.  Inside the red boundary, intersections of the solid curves with the integer vertices correspond to crossings of type 2, whereas intersections with the dashed curve correspond to a crossing of type 3.  If any integer vertices lie on the red boundary curve, then these correspond to a crossing of type 4 at the appropriate Rayleigh number.  All other crossings are of type 1.

\begin{figure}[H]
\begin{center}
	\begin{tabular}{cc}
		\includegraphics[height=48mm]{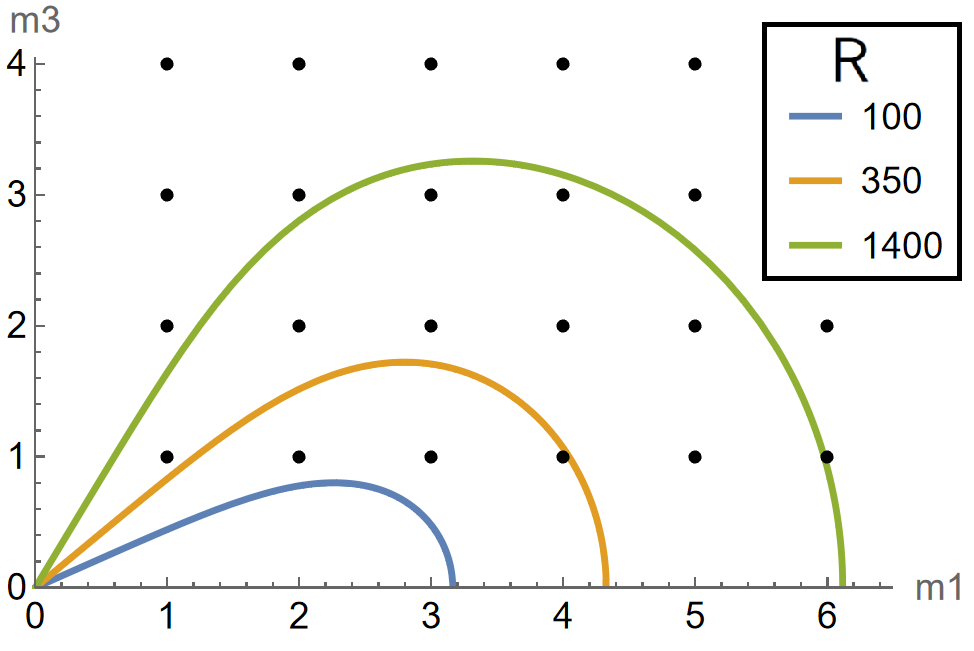}
		&\includegraphics[height=48mm]{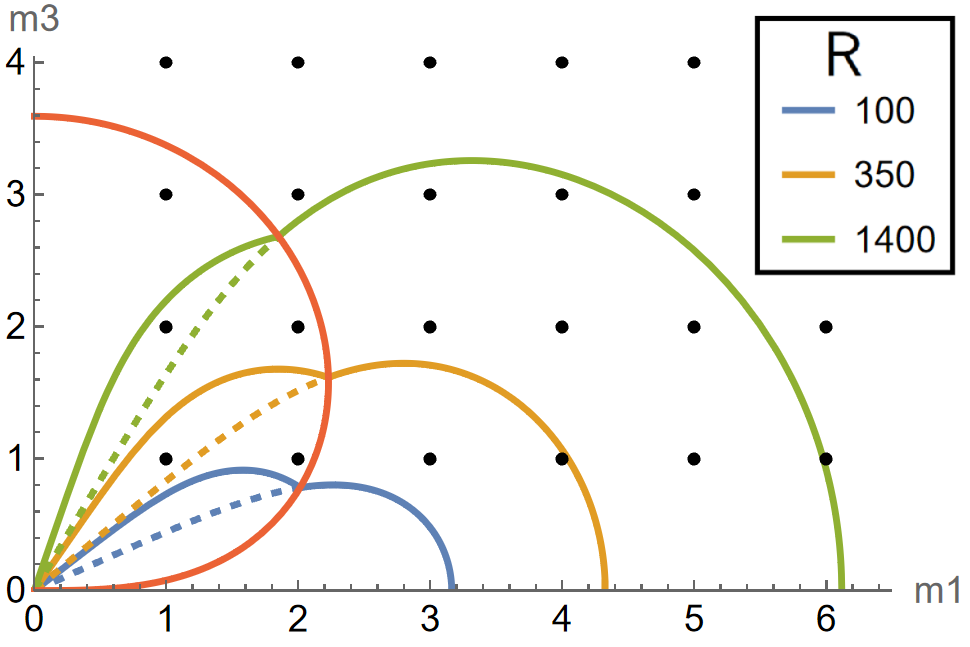} \\
		(a) & (b)
	\end{tabular}
	\caption{(a) For $\Pra = 2$, $\Rot = 10\sqrt{5}$, $\Sha = 1$, the level curves of $\Ray^{\textbf{m},c}$ as $\Ray$ increases. (b) For $\Pra = 1/2$, $\Rot = 100\sqrt{5}$, $\Sha = 1$, the boundary curve $\Ray^{\textbf{m},1} = \Ray^{\textbf{m},2}$ (red), as well as level curves of $\Ray^{\textbf{m},c}$.  For wave vectors $\textbf{m}$ inside the red boundary, the solid curves indicate the initial (Hopf) type 2 crossing and dashed curves indicate the subsequent type 3 crossing. }
	\label{fig:CriticalLevelSets}
\end{center}
\end{figure}

It is possible that multiple critical Rayleigh numbers can be equal, i.e. $\mathcal{R}^{\tilde{\textbf{m}},c} = \mathcal{R}^{\textbf{m},c}$ with $\tilde{\textbf{m}} \neq \textbf{m}$, in which case more complicated bifurcations could occur.  On the other hand, in \ref{app:LevelCurves} it is proven that $m_3,\tilde{m}_3$ are strictly concave functions of $m_1$, and the boundary curve is a strictly concave function of $m_3$.  This puts a nice constraint on the complexity of possible scenarios for simultaneous critical Rayleigh numbers.  For example, when $\Pra \geq 1$, one can have a maximum of two critical Rayleigh numbers with the same value of $m_3$.

Before we study the bifurcations associated to these crossings, we need to introduce some notation.  First, we introduce $\mathscr{M}^{c}_{\parameters} \subset \mathbb{Z}^2_{> 0}$ be the set of all wave vectors $\textbf{m}$ which are critical at parameter value $\parameters$, ie $\mathscr{M}^{c}(\parameters) := \cup_{\ell \leq 4} \mathscr{M}^{c,\ell}(\parameters)$, where $\mathscr{M}^{c,\ell}(\parameters) \subset \mathbb{Z}_{>0}^2$ are defined to be the set of all wave vectors $\textbf{m} \in \mathbb{Z}_{>0}^2$ such that the eigenvalues of $\mathcal{L}^{\textbf{0},\textbf{m}}$ have a crossing of type $\ell$ at a given Rayleigh number $\Ray$.  Explicitly, 
\[ \mathscr{M}^{c,1} = \{ \textbf{m} \in \mathbb{Z}^2_{> 0}: \Ray^{\textbf{m},1} = \Ray \} \hspace{.5 cm} \text{ , } \hspace{.5 cm} \mathscr{M}^{c,2}(\parameters) = \mathscr{M}^{c,3}(\parameters) = \mathscr{M}^{c,4}(\parameters) = \varnothing , \]
for $\Pra \geq 1$, whereas for $0 < \Pra < 1$, these are defined as
\begin{equation} \notag
\begin{split} \mathscr{M}^{c,1}(\parameters) = \{ \textbf{m} : \Ray^{\textbf{m},1} = \Ray \text{ , } |\Rot| < \Rot^{\textbf{m}} \} \hspace{.1 cm} & \text{ , } \hspace{.1 cm} \mathscr{M}^{c,2}(\parameters) = \{ \textbf{m}: \Ray^{\textbf{m},2} = \Ray \text{ , } |\Rot| > \Rot^{\textbf{m}} \} , \\ \mathscr{M}^{c,3}(\parameters) = \{ \textbf{m} : \Ray^{\textbf{m},1} = \Ray \text{ , } |\Rot| > \Rot^{\textbf{m}} \} \hspace{.1 cm} & \text{ , } \hspace{.1 cm} \mathscr{M}^{c,4}(\parameters) = \{ \textbf{m} : \Ray^{\textbf{m},1} = \Ray \text{ , } |\Rot| = \Rot^{\textbf{m}} \} . \end{split} 
\end{equation}
Similarly, we let $\mathscr{M}^{s},\mathscr{M}^{u} \subset \mathbb{Z}^2_{\geq 0}$ denote the sets of all wave vectors which are stable and unstable at parameter $\parameters$, respectively.  Explicitly, these sets are given by
\[ \mathscr{M}^{s}(\parameters) = \{ \textbf{m} : \Ray^{\textbf{m},c} > \Ray \} \hspace{.5 cm} \text{ , } \hspace{.5 cm} \mathscr{M}^{u}(\parameters) = \{ \textbf{m} : \Ray^{\textbf{m},c} < \Ray \} . \]
Let $\textbf{X}^s,\textbf{X}^c,\textbf{X}^u$ denote the projections of $\textbf{X}$ onto the corresponding index sets $\mathscr{M}^s,\mathscr{M}^c,\mathscr{M}^u$.

In the following theorem, it is proven that the unstable manifold of the origin $\mathscr{W}^u(\textbf{0})$ exists, and its dimension $d( \mathscr{W}^{u}(\textbf{0}))$ is estimated in terms of the parameters $\parameters$.  In particular, by holding $\Rot$ constant $d( \mathscr{W}^{u}(\textbf{0}))$ increases as $\Ray^{1/2}$ as $\Ray \to \infty$, and by holding $\Ray$ constant $d( \mathscr{W}^{u}(\textbf{0}))$ collapses as $\Rot^{-1}$ as $\Rot \to \infty$.  The bifurcations giving rise to the unstable manifold are then classified in the special case when only one mode is critical and all eigenvalues are distinct, hence $|\Rot^c| \neq \Rot^{\textbf{m}}$.  Note that this theorem is not meant to be an exhaustive summary of all possible bifurcations, but only a proof regarding typical bifurcations, and the case of simultaneous critical modes and degenerate eigenvalues should be quite rare in parameter space.  The proof is given in \ref{app:ProofOfOriginBifur}.

\begin{Theorem}
\label{thm:OriginBifurcations}

\textbf{(a)} For admissible $\Pra, \Rot, \Sha$, let $\Ray^{*} = \min_{\textbf{m} \in \mathbb{Z}^2_{>0}} \Ray^{\textbf{m},c}$.  For $\Ray > \Ray^*$ a non-trivial unstable manifold of the origin $\mathscr{W}^u(\textbf{0})$ exists, and one has
\[  d( \mathscr{W}^{u}(\textbf{0})) \sim \frac{ \Ray^{1/2} }{\Sha ( 1 + \frac{\Rot^2}{\Ray})^{1/2}} \text{  . }  \]
More precisely, there exists constants $C^{\textbf{0},u}_1, C^{\textbf{0},u}_2 > 0$ such that one has the following bounds, and furthermore one has the following limit
\begin{equation} \label{UnstableManifoldDimension} \big \lfloor \frac{ C^{\textbf{0},u}_1 \Ray^{1/2} }{\Sha ( 1 + \frac{\Rot^2}{\Ray})^{1/2}} \big \rfloor \leq d( \mathscr{W}^{u}(\textbf{0})) \leq \big \lfloor \frac{ C^{\textbf{0},u}_2 \Ray^{1/2} }{\Sha ( 1 + \frac{\Rot^2}{\Ray})^{1/2}} \big \rfloor  \hspace{.5 cm} \text{ , } \hspace{.5 cm}  \lim_{\Ray \to \infty} \frac{\Sha}{\Ray^{1/2}} d( \mathscr{W}^{u}(\textbf{0})) = \frac{1}{2} \text{ . } \end{equation}
Since $\mathscr{W}^{u}(\textbf{0}) \subset \mathscr{A}$, one immediately obtains the following lower bound:
\begin{equation} \label{DimensionLowerBound} d_{Haus}( \mathscr{A}) \geq d( \mathscr{W}^{u}(\textbf{0})) \text{ . } \end{equation}

\noindent \textbf{(b)} For admissible $\parameters^c$, suppose $\mathscr{M}^{c}(\parameters^c)$ is non-empty.  Then the PDE \eqref{EvolutionPDE} undergoes a local bifurcation at $(\textbf{u},\theta,\parameters) = (\textbf{0},0,\parameters^c)$.  If $|\mathscr{M}^{c}| = 1$ and $|\Rot^c| \neq \Rot^{\textbf{m}}$ these bifurcations are classified as follows:  
\begin{enumerate}[label=(\roman*)]
	\item If $\textbf{m} \in \mathscr{M}^{c,1}$ and $\Pra^c \geq \frac{m_3}{\Sha m_1}$ then a supercritical pitchfork bifurcation occurs as $\Ray$ is increased through $\Ray^c$ with the other parameters fixed.  If $\Pra^c < \frac{m_3}{\Sha m_1}$ then let
	\begin{equation} \label{CriticalityRotationThreshold} C^{\textbf{m}} =  \frac{1}{m_3} |\mathcal{K}\textbf{m}|^3 \big ( (\frac{m_3}{\Sha m_1 \Pra})^2-1 \big )^{-1/2} \text{ . }  \end{equation}
	If $|\Rot^c| <  C^{\textbf{m}}$ then a supercritical pitchfork bifurcation occurs, whereas if $|\Rot^c| >  C^{\textbf{m}}$ a subcritical pitchfork bifurcation occurs.  If $\textbf{m} \in \mathscr{M}^{c,3}$, then the same bifurcations occur as $\Ray$ is decreased through $\Ray^c$.  
	\item \revised{If $\textbf{m} \in \mathscr{M}^{c,2}$, $\Pra \geq \frac{1}{2}$ and $m_3 \geq \sqrt{2}\Sha m_1$ then a supercritical Hopf bifurcation occurs.}
\end{enumerate}
\end{Theorem}

\begin{Remark}
Note that the rotation number $\Rot$ could just as well be treated as the control parameter, rather than the Rayleigh number, and the same bifurcations occur.  For example, as one decreases the rotation $\Rot$ through $\Rot^c$, the critical Rayleigh number $\Ray^c$ defined in \eqref{CriticalRayleighNumber} decreases through the fixed Rayleigh number $\Ray$.
\end{Remark}

\subsection{Upper bound on the attractor dimension for the Boussinesq Coriolis model}

Here we apply the techniques of Temam \cite{temam_InfDimDynSys} to obtain an upper bound on the dimension of the attractor.  Since the attractor presumably has a fractal structure, it is more precise to use alternative notions of dimension.  Therefore recall that the $d$-dimensional Hausdorff measure of $\mathscr{A}$ is the number
\[ \mu_{H}(\mathscr{A},d) := \sup_{\epsilon > 0} \inf_{ C \in \mathcal{C}_{\epsilon}} \sum \nolimits_{i=1}^{|C|} r_i^d , \]
where $\mathcal{C}_{\epsilon}$ is the collection of all coverings $C$ of $\mathscr{A}$ by balls $B_i$ of radius $0 < r_i \leq \epsilon$.  The Hausdorff dimension of $\mathscr{A}$, $d_{Haus}(\mathscr{A})$, is defined as the unique number such that for all $d \in [0,\infty)$ one has $\mu_{H}(\mathscr{A},d) = 0$ if $d > d_{Haus}(\mathscr{A})$ and $\mu_{H}(\mathscr{A},d) = \infty$ if $d < d_{Haus}(\mathscr{A})$.  On the other hand the box-counting dimension of $\mathscr{A}$ is defined
\[ d_{box}(\mathscr{A}) := \limsup_{\epsilon \to 0} \frac{\log n_{\mathscr{A}}(\epsilon)}{\log 1/\epsilon} , \]
where $n_{\mathscr{A}}(\epsilon)$ is the number of balls of radius $\leq \epsilon$ which is necessary to cover $\mathscr{A}$.  Note that Temam refers to the box-counting dimension as the fractal dimension, whereas others use the term fractal dimension more generally to refer to any index which characterizes fractal sets by quantifying their complexity.

Temam et al's result is for non-rotating, planar flows, hence $\Rot = 0$, and $u_2(\textbf{x},t) = 0$ for all $t \geq 0$.  For $\Rot = 0$, the condition $u_2(\textbf{x},t) = 0$ is invariant for the flow defined by \eqref{BoussinesqCoriolis}, their result can be stated in our notation as follows:
\begin{Theorem}[Foias, Manley, Temam \cite{FoiasManleyTemam_1987}  Theorem 5.1]
\label{thm:AttractorUppderBound}
For $\Rot = 0$, $u_{2,0}(\textbf{x}) = 0$, let $\mathscr{A}$ be the global attractor of $\mathcal{S}(t)$ defined in Theorem \ref{thm:PDE_WellPosedness_Balances} (a).  There exists constants $C_{Haus},C_{box}$ depending only on $\Pra,\Sha$ such that
\begin{equation} \label{DimensionUpperBound} d_{Haus}(\mathscr{A}) \leq C_{Haus} ( 1 + \Ray ) \hspace{.5 cm} \text{ , } \hspace{.5 cm} d_{box}(\mathscr{A}) \leq C_{box} ( 1 + \Ray ) . \end{equation}
\end{Theorem}
In fact, one can check that this result carries over to the general case $\Rot \geq 0$, $u_{2,0}(\textbf{x}) \neq 0$, with different constants $\tilde{C}_{Haus}, \tilde{C}_{box}$.  Hence \eqref{DimensionUpperBound} is used in the analysis in section \ref{sec:NumAnalysis} with $\tilde{C}_{Haus}$.  The proof for the general case is so similar to that in  \cite{FoiasManleyTemam_1987} that it is not repeated here in full detail.  However, by carefully analyzing the proof, one can extract the constant $\tilde{C}_{Haus}$, hence a sketch of the proof is provided in \ref{App:AttractDimUpperBound} to show how these arise.  The result is as follows:
\begin{equation} \label{DimensionUpperBound_Constants} \tilde{C}_{Haus} = \frac{320\pi^3}{\Pra (1+\Pra)\min(1,\Sha^2)} . \end{equation}

\begin{Remark}
Note that while \eqref{DimensionUpperBound} is independent of the rotation number, this bound is not necessarily sharp.  Indeed, the analysis of the local bifurcations at the origin and the numerical work in section \ref{sec:NumAnalysis} strongly suggest that the rotation plays a significant role.  However, a rotation dependent bound ended up being too difficult to prove for the present work, since the rotation is energy neutral and hence doesn't enter into \eqref{Attractor_IntegralBound}, hence one must use new methods.  This same problem arises in background field arguments, described in the conclusion section.
\end{Remark}

\section{Numerical analysis}

\label{sec:NumAnalysis}

This section summarizes some numerical investigations into heat transport in the Boussinesq Coriolis model via the HKC models.  First, an overview of some heat transport phenomenology is given, leading to the goals for the numerical studies.  A summary of the codes written to achieve these goals is then given.  The results of several numerical studies regarding convergence and dynamics are then presented.  Note that all of the computations were done using MATLAB on a standard laptop or desktop, and as mentioned above the codes are available on GitHub.

\subsection{Context, goals and code inventory}  

\begin{figure}[H]
\begin{center}
	\begin{tabular}{cc}
		\includegraphics[height=52mm]{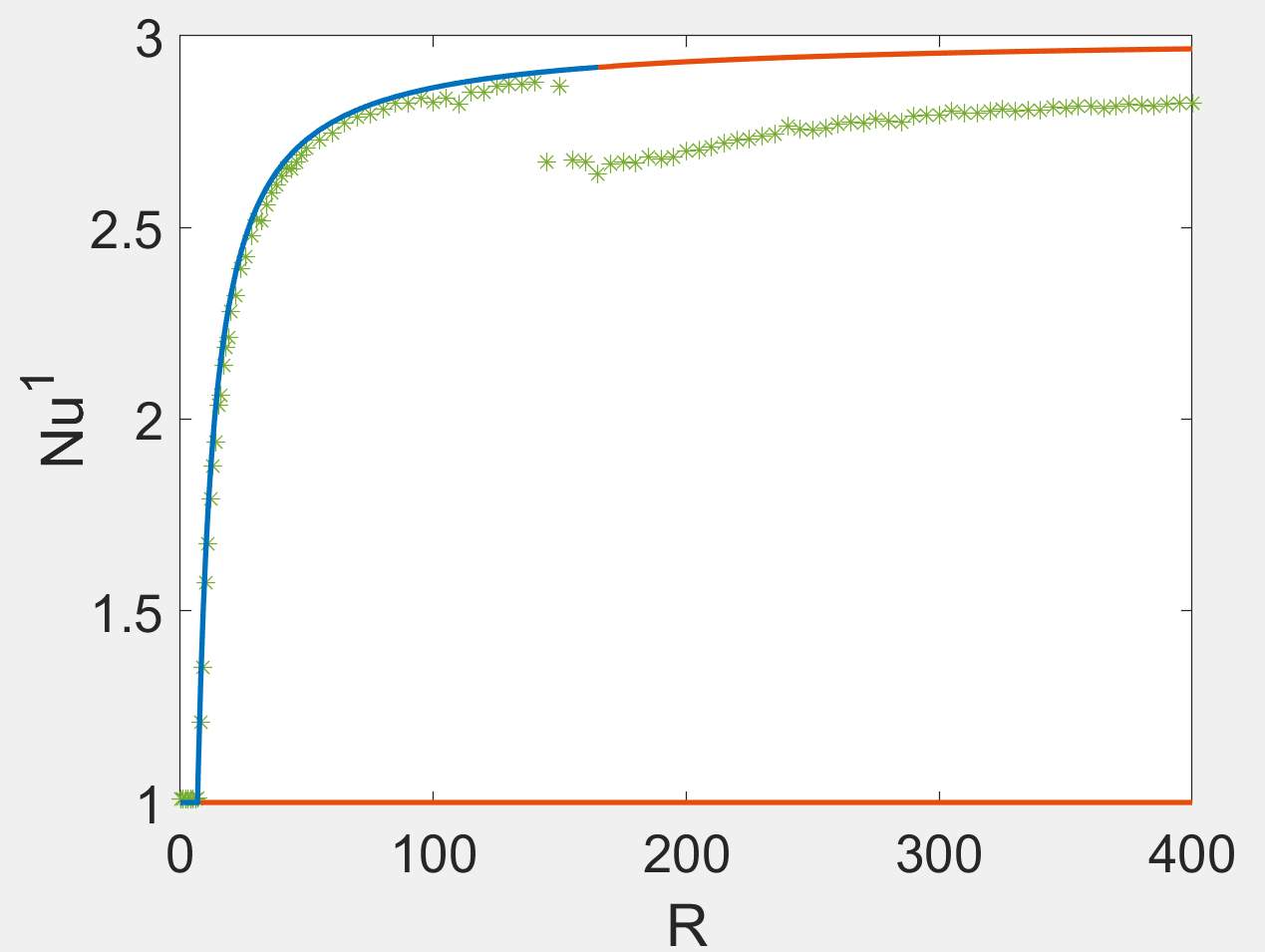}
		& \includegraphics[height=52mm]{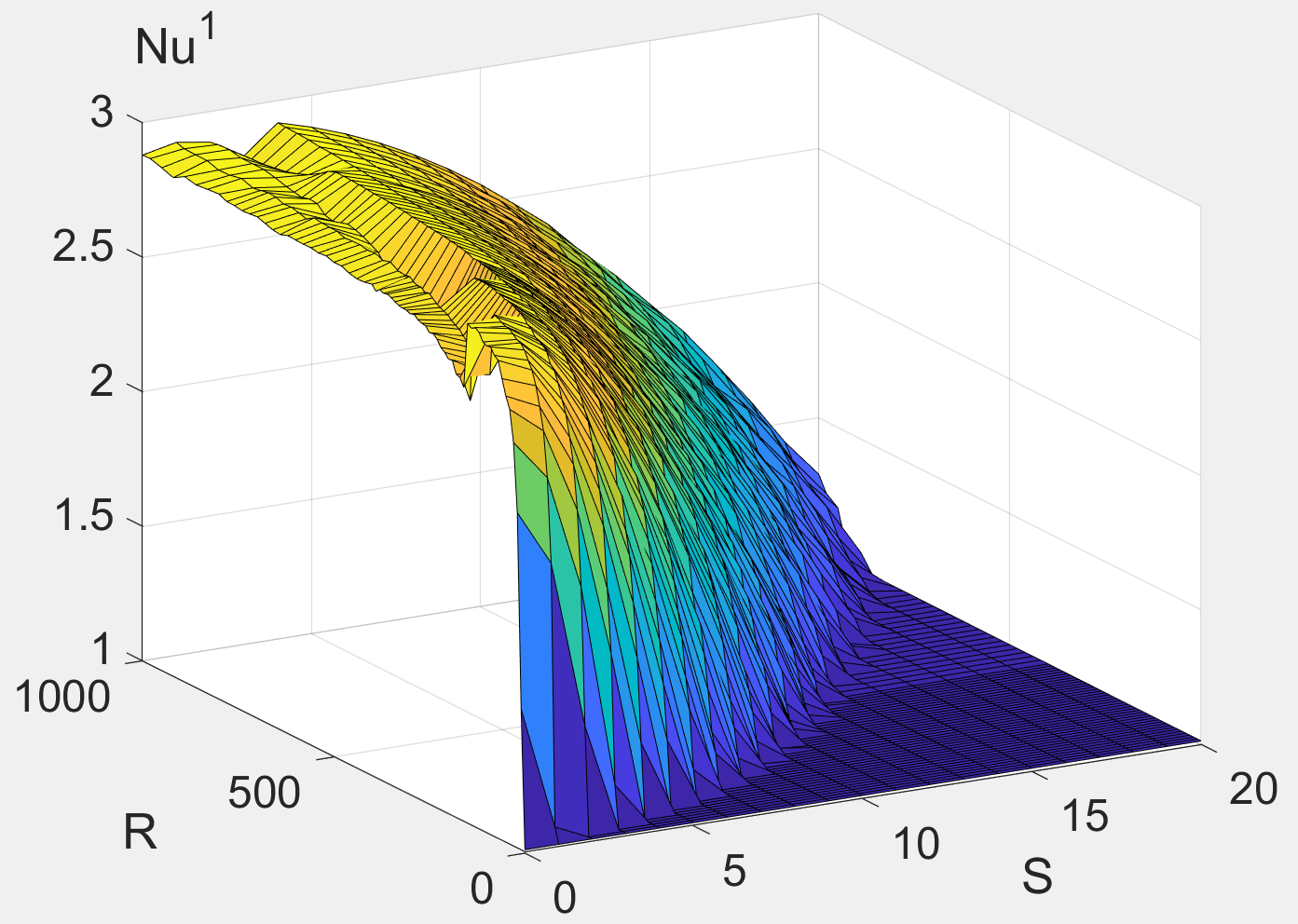} \\
		(a) & (b)
	\end{tabular}
	
	\caption{For $\Pra =10$, $\Sha = \frac{1}{\sqrt{2}}$, Nusselt number approximations for the HKC-1 model (a) for $\Rot = 0$ with random initial conditions (green) and the fixed point values, with blue indicating local stability and red indicating instability and (b) for random initial conditions $0 \leq \Rot \leq 20$, $0 \leq \Ray \leq 1000$. }
	\label{fig:BasicPhenomena}
\end{center}
\end{figure}

In order to provide context for the numerical results below, we summarize here some heat transport phenomena which are fairly well understood, and give a more detailed account in \ref{app:HeatTransportSurvey}.  Already one can gain many qualitative insights by carefully considering the HKC-1 model.  Figure \ref{fig:BasicPhenomena} depicts the heat transport for the HKC-1 model, in which one can see an intimate correspondence between the dynamics of this model and its heat transport.  Here we chose the parameters $(\Pra,\Sha )$ equal their classical Lorenz values $(10,1/\sqrt{2})$, since with this choice there is a huge body of literature to compare to.

Panel (a) depicts a bifurcation diagram of the Nusselt number as $\Ray$ increases with $\Rot$ fixed at $0$, where blue and red curves indicate locally stable and locally unstable equilibria, respectively.  Solutions with random initial conditions are depicted in green.  For very small Rayleigh numbers, the origin is the global attractor and hence all solutions have $\mathsf{Nu}^1 = 1$.  At $\Ray = 6.75$, a pair of nontrivial fixed points emerge from the origin in a pitchfork bifurcation, representing steady convection rolls.  They are locally stable for $ 6.75 \leq \Ray \lesssim 166.97$, and solutions with random I.C.'s track their Nusselt number.  At $\Ray \approx 166.97$ these fixed points lose stability in a subcritical Hopf bifurcation, whereupon solutions tend to the famous chaotic attractor.  One then observes that trajectories exhibit less heat transport than those of the unstable fixed points.  In this case, the fluid seems to be spending more energy moving around chaotically and hence transports less energy.

Panel (b) depicts the Nusselt number computed for random initial conditions as a function of both $\Ray$ and $\Rot$.  The effect of rotation seems to always be to stabilize simpler dynamical behavior.  This stabilization can increase $\mathsf{Nu}^1$ only as it moves the dynamics out of the chaotic region, but typically the rotation decreases the heat transport, and it never increases $\mathsf{Nu}^1$ above the maximal rate realized among solutions with $\Rot = 0$.  As the rotation becomes very large $\mathsf{Nu}^1$ must decrease down to $1$, since eventually the stabilizing effect becomes so strong that vertical motion is impossible and the origin becomes the global attractor.  

Several of the above phenomena seem to carry over the PDE.  For example, it is proven that vertical motion is suppressed to zero in the limit of infinite rotation (see for example \cite{Chemin2006Mathematical}), so vertical convective heat transport will be zero.  Since the rotation drops out of the energy balance equation, the rotation cannot increase the heat transport above the bound for the $\Rot =0$ case.  There is numerical evidence suggesting that the maximal heat transport for the PDE is also realized by time-independent solutions \cite{wen_goluskin_doering_2022}, although this remains unproven.

The PDE may have qualitative similarities to the HKC-1 model, but one must use an ODE of appropriately high spatial resolution to have any chance of accurately representing the PDE quantatively.  For instance, it is widely believed that the PDE exhibits unbounded heat transport as $\Ray \to \infty$, whereas it can be proven analytically that the heat transport in the HKC-1 model tends toward a constant value \cite{OvsyRadeWelterLu_2023_LorenzLargeRayleigh}.  The required spatial resolution to accurately capture the heat transport of the PDE at given Rayleigh and rotation numbers is unknown and presents one of the main questions studied herein.  In the 6 dimensional HKC-1 model the short time steps and long integration times required by the chaotic dynamics are expensive but not prohibitive.  If the dimension required for accurate heat transport were to increase as in \eqref{DimensionUpperBound}, then the computational expense would very quickly become prohibitive.

Therefore the goals of the numerical studies below aim at addressing issues which fall roughly into two categories:
\begin{enumerate}
\item \textbf{Convergence issues:}  How long must one integrate a trajectory to obtain an approximation of the infinite time average?  How large are the errors associated with step-size, especially for solutions which appear chaotic (ie persistent errors)?  For a given $\Ray,\Rot$, how large must $M$ be such that the heat transport values from the HKC-$M$ model appear to converge?  How does this compare to the bounds in section \ref{sec:AttractorAnalysis}?  How large is the error associated with using an insufficiently large model?

\item \textbf{Dynamical issues:} What are the dynamically stable values of the heat transport (ie heat transport obtained by all initial conditions in a set of positive measure)?  At a given Rayleigh number, are there multiple stable heat transport values?  What heat transport do stationary solutions exhibit, and how much does this differ from time dependent solutions?  How can chaotic solutions be identified?
\end{enumerate}

The codes begin with the script \verb|ModelConstructor.m| which takes as input a positive integer $M$, and writes a file called \verb|HKCM.m| containing the right hand side of \eqref{GeneralBoussinesqODE_Vel} - \eqref{GeneralBoussinesqODE_Temp}.  Solutions of the HKC-$M$ model can then be simulated with the function \verb|FluidSolver.m|, which takes as input the desired parameters $\parameters$ and model number $M$, uses \verb|ode45| to call \verb|HKCM.m| and saves the trajectory along with the parameters.  \verb|FluidSolver.m| was written so that it can be given an initial condition (either random or near-uniform) and start at $t = 0$, or it can extend an already computed trajectory for a specified time.  Several codes are then available to visualize the trajectory, either in phase space, or via the corresponding velocity or temperature fields.  Several codes are also available to compute the cumulative average heat transport corresponding to these trajectories using the expression \eqref{HKC_Nusselt}.  The function \verb|HeatTransport_Trajectory.m| takes a single trajectory and computes its heat transport, whereas the script \verb|HeatTransport_Iterator.m| is initialized with an array of Rayleigh and rotation numbers, computes trajectories for each parameter value by calling \verb|FluidSolver.m| and then computes the corresponding heat transport.

\subsection{Convergence issues}

The first convergence issue was that of convergence in time, namely the step size for the numerical integration must be chosen sufficiently small to approximate the ODE solution well, and the total integration time must be chosen sufficiently long to obtain a good approximation for the infinite time average.  Regarding this first issue, most of the time some small step size such as $\Delta t = 10^{-3}$ was chosen, and then \verb|ode45|'s built in error control was used to determine whether smaller time steps were needed by setting the relative and absolute error.  

Regarding the long integration time, the function \verb|HeatTransport_Iterator.m| was designed to first compute a trajectory for some long integration time such as $10^{5}$ time increments.  The function \verb|HeatTransport_Iterator.m| then enters into a loop which tests the fluctuations in the heat transport value, and either continues the trajectory if the fluctuations are larger than a specified threshold or saves the data if the fluctuations are acceptably small.  Specifically, the trajectory is extended by an additional $10^{3}$ time units until the standard deviation in the second half of the trajectory was less than 2\% of the cumulative heat transport at the end of the trajectory.  Most commonly the solutions exhibited quick convergence toward an apparent infinite-time average, although in rare cases a trajectory can spend a long time near solutions which are weakly unstable (see for example Figure \ref{fig:HeatTransportTimeConvergence}), hence long integration times may be required.

\begin{figure}[H]
\begin{center}
	\begin{tabular}{cc}
		\includegraphics[height=66mm]{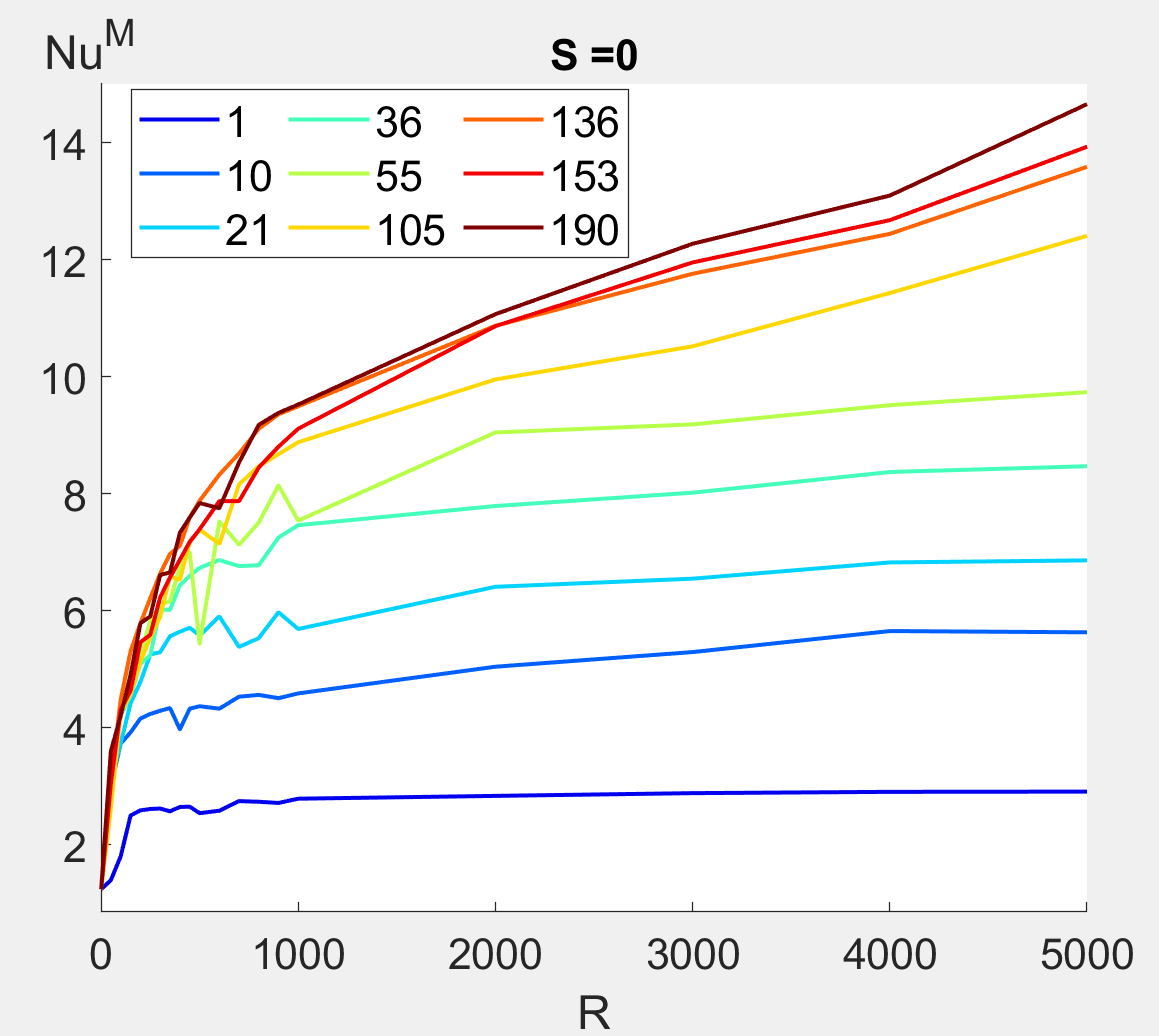}
		&\includegraphics[height=66mm]{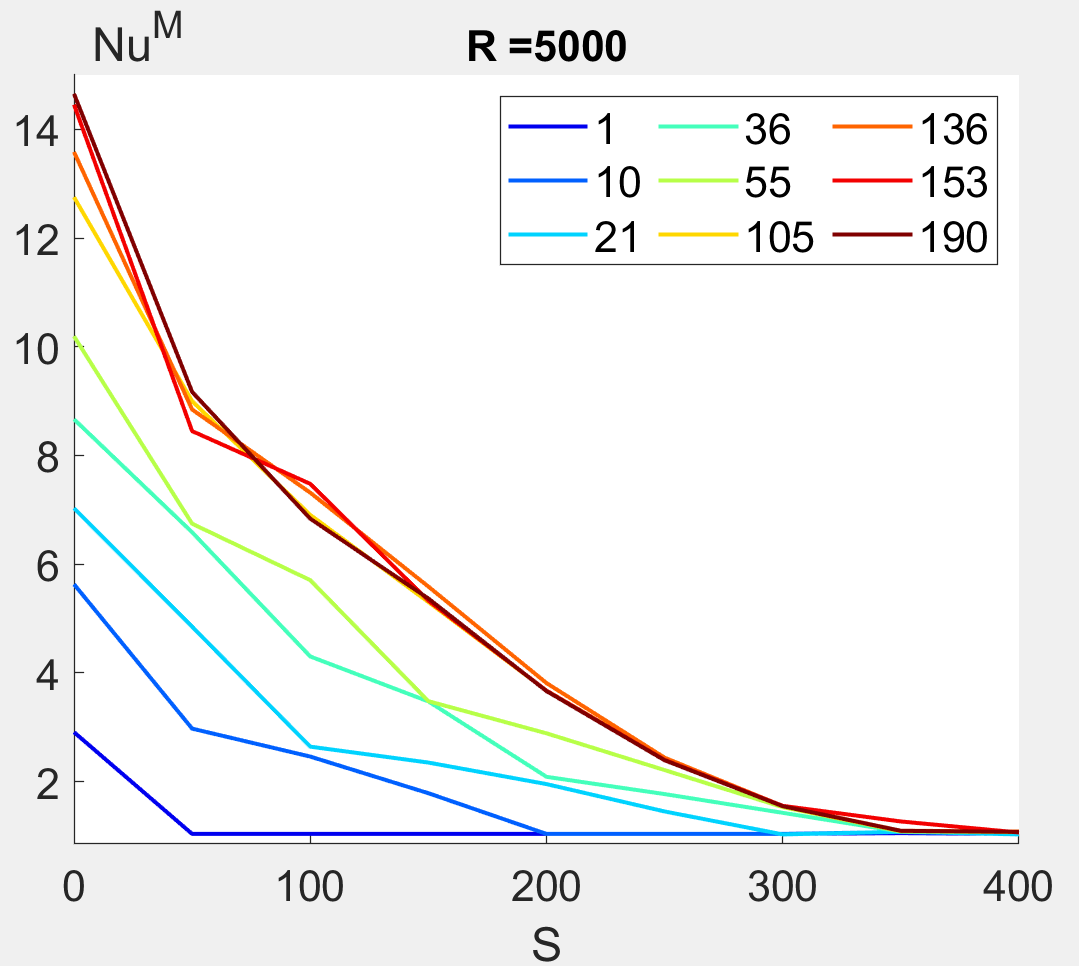}  \\
		(a) & (b)
	\end{tabular}
	\caption{Heat transport model comparison (a) for all $0 \leq \Ray \leq 5000, \Rot = 0$ and (b) for $\Ray = 5000$, $0 \leq \Rot \leq 400$.}
	\label{fig:HeatTransportModelComparison}
\end{center}
\end{figure}

In order to study the convergence of the ODE models toward the PDE, trajectories were generated using \verb|HeatTransport_Iterator.m| for each the following HKC-$M$ models, beginning from a random perturbation of the uniform initial state $(\textbf{u}_0(\textbf{x}),T_0(\textbf{x})) = (0,1/2)$:
\begin{equation} \label{NumericalStudy_ModelNumbers} M = 1 \text{ , } 3 \text{ , } 6 \text{ , } 10 \text{ , } 21 \text{ , } 36 \text{ , } 45 \text{ , } 55 \text{ , } 66 \text{ , } 78 \text{ , } 91 \text{ , } 105 \text{ , } 120 \text{ , } 136 \text{ , } 153 \text{ , } 171 \text{ , } 190 \text{ . } \end{equation}
These corresponding to the $1^{st}$ through $19^{th}$ completed shell.  For each model, a trajectory was generated  with $\Pra = 10,$ $\Sha = 1/\sqrt{2}$, for each pair of $\Ray,\Rot$ in the following range:
\begin{equation} \label{HKC_ModelComparison_ParameterRange}  \Ray = [ 1 , 50:50:500 , 600:100:1000 , 2000:1000:5000  ]  \hspace{.5 cm} \text{ , } \hspace{.5 cm}  \Rot =  0:50:400  . \end{equation}
where $x_0 : \Delta x: x_f$ is the MATLAB notation denoting the set of all integers beginning from $x_0$, incrementing by $\Delta x$ and ending with $x_f$.  For these trajectories, the relative and absolute error in \verb|ode45| were set to $10^{-10}$, the time increment was set to $10^{-4}$, the trajectories were initially integrated for $10^{4}$ time increments, or $10^{3}$ more time increments until the fluctuations in heat transport met the 2\% error threshold.  The results are shown in Figure \ref{fig:HeatTransportModelComparison}, although for clarity only a smaller selection of models is displayed and the dependence on the parameters $\Ray, \Rot $ is shown separately.

Panel (a) in Figure \ref{fig:HeatTransportModelComparison} displays the dependence of the heat transport on $\Ray$ for $\Rot = 0$.  Most prominently, this figure displays that in each HKC model the heat transport increases rather rapidly for low Rayleigh numbers, then the rate of increase levels off for higher Rayleigh numbers.  This increase is fairly monotonic, although with some fluctuations.  The higher dimensional HKC models almost always exhibit more heat transport than the lower dimensional models.  We point out one interesting exception where at the value $\Ray = 500$ the HKC-55 model appears to exhibit less heat transport than many of the others.  Rather than a peculiarity of this model, this might be due to the random initial conditions, and we will return to this point later.  All HKC models tend to agree fairly well for low Rayleigh numbers, but as one moves toward higher Rayleigh numbers the lower dimensional models one by one exhibit significantly less heat transport than the higher dimensional models.  As described above, one sees that the HKC-1 model levels off to an apparently constant value, as one would expect from analytical results, and already by $\Ray = 5000$ exhibits around $5$ times less transport than the HKC-190 model.  It might be that the failure mechanism for the HKC-1 model described above also appears in each truncated model, where the heat transport eventually levels off to a constant value as the Rayleigh number begins to vastly dominate over the viscosity of the smallest scale of the truncated model.  Even the HKC-105 model seems to lie significantly below the other models beginning around $\Ray = 2000$.  On the other hand, the HKC-136, HKC-153 and HKC-190 all seem to agree quite well in Figure \ref{fig:HeatTransportModelComparison} (a), so on this range it seems that these models are accurately representing the PDE.  

Panel (b) in Figure \ref{fig:HeatTransportModelComparison} shows the dependence of the heat transport on the rotation number for $\Ray = 5000$.  The heat transport decreases with increasing rotation for all models, although more slowly for the higher dimensional HKC models.  This decrease appears to be monotonic, although one sees that the slopes are not strictly increasing.  One sees a more interesting relationship across the HKC models, where for $\Rot \geq 100$ the HKC-190 model seems to exhibit slightly less heat transport than the HKC-136 and HKC-153 models.  This again could be due to the random initial conditions, but it is somewhat curious that it is so consistent.  Finally, in this case the HKC-105, HKC-136, HKC-153 and HKC-190 all seem to agree quite well over the whole range, and one sees that all models tend to agree for higher rotation numbers.  Hence for higher rotation, the PDE seems well represented by lower dimensional models. 

In order to quantify the apparent convergence across models of different spatial resolutions, the heat transport for each model in \eqref{NumericalStudy_ModelNumbers} was compared to the heat transport from the HKC-190 model at each Rayleigh and rotation number in \eqref{HKC_ModelComparison_ParameterRange}.  For definiteness, the value $10\%$ was chosen arbitrarily, and if the heat transport from a given model was within $10\%$ of the heat transport value from the HKC-190 model, then that the model was considered to represent the PDE "well".  The dimension of the smallest model found in this way will be referred to as the empirical Nusselt dimension, $d_{Nusselt}$.  In Figure \ref{fig:HeatTransportModelComparison2} below, one can see the plot of $d_{Nusselt}$, indicating the increasing computational cost associated with resolving the flow as the Rayleigh number increases.  

\begin{figure}[H]
\begin{center}
	\includegraphics[height=60mm]{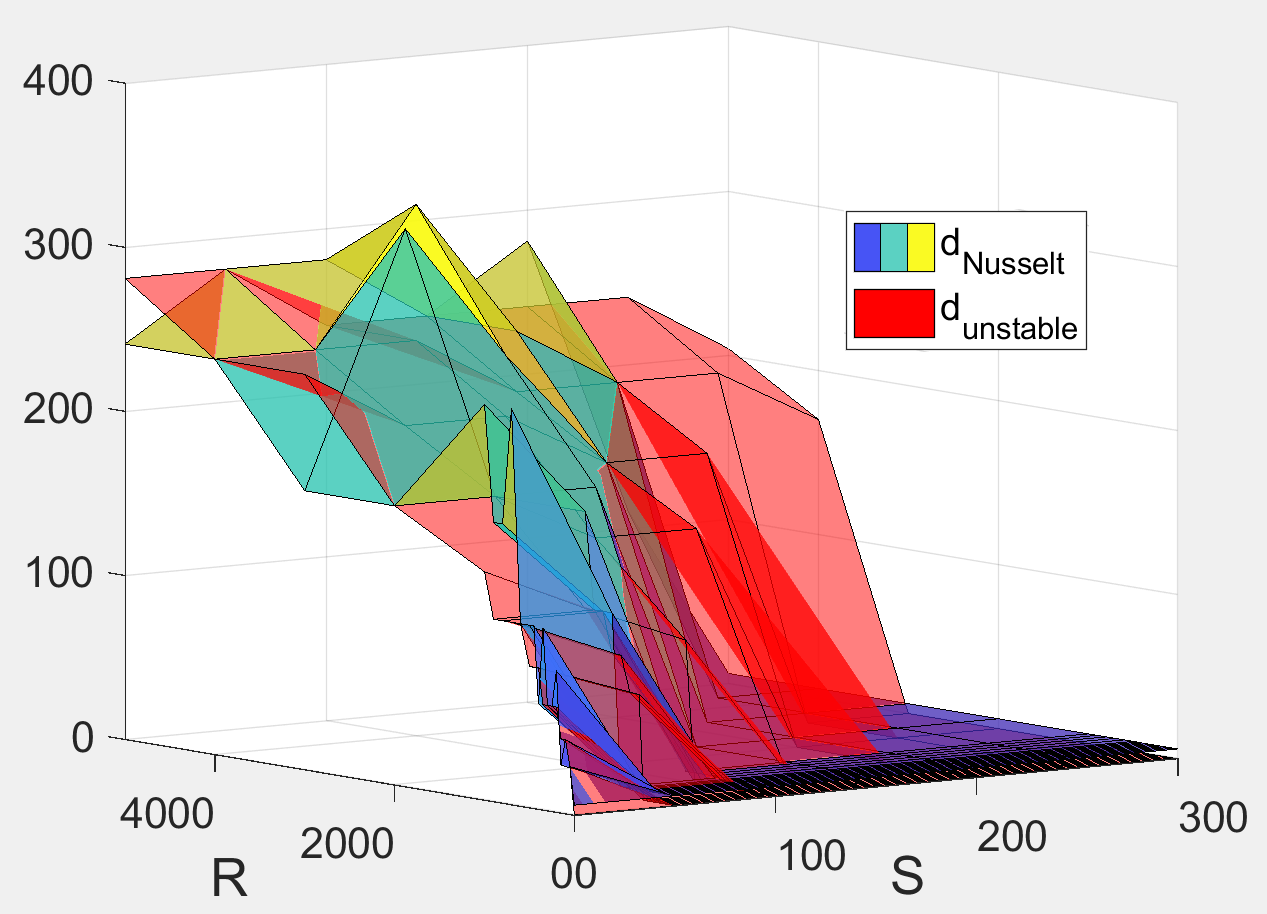}  
	\caption{The computational cost measured by the empirical Nusselt dimension $d_{Nusselt}$ vs. the dimension of the minimal HKC model containing all unstable modes $d_{unstable}$. }
	\label{fig:HeatTransportModelComparison2}
\end{center}
\end{figure}

The lower bound on the Hausdorff dimension found from \eqref{DimensionLowerBound} was found to be quite small compared to the plot of $d_{Nusselt}$.  This is not at all surprising, and in fact this is not a meaningful comparison because the attractor must of course be embedded in a larger ambient space.  In order to provide a more meaningful comparison, the quantity $d_{unstable}$ was plotted in red in Figure \ref{fig:HeatTransportModelComparison2} panel, defined as the minimal HKC model which contains all of the unstable modes at the origin.  The upper bound \ref{DimensionUpperBound} was not depicted in \ref{fig:HeatTransportModelComparison2}, because this is far larger, on the order of $10^5$ for most of the range depicted in Figure \ref{fig:HeatTransportModelComparison2}. 

The two graphs plotted in Figure \ref{fig:HeatTransportModelComparison2} seem to agree fairly well, although for larger rotation numbers it is noticeable that $d_{unstable}$ jumps up immediately to high dimension when the Rayleigh number crosses a critical threshold, whereas $d_{Nusselt}$ ramps up over a larger range of Rayleigh numbers.  This is not surprising because the structures developing after a critical threshold are initially small and do not immediately produce a large change in transport.  It is also noticeable that $d_{Nusselt}$ fluctuates in a non-monotonic way.  This is surprising, and possibly just a result of the crude nature of this comparison.  This study should be interpreted as a rough idea about how much resolution is needed, since it is limited by the fact that only the restricted number of HKC models in \eqref{NumericalStudy_ModelNumbers} was used, rather than every HKC model from 1 to 190.  Furthermore this study does not necessarily compare across different HKC models for analogous initial conditions.  As we shall see in the coming section this dependence of the heat transport is nontrivial due to the presence of multiple stable values of transport, hence the comparison across models may be skewed by the fact that trajectories for each model were initialized using random perturbations of the uniform state. 



\subsection{Dynamical issues}

While stationary solutions are important for understanding the dynamics of any system, they are particularly important in Rayleigh-B\'enard convection since it is thought that they might maximize heat transport \cite{wen_goluskin_doering_2022,OlsonDoering2022}.  Hence this subsection presents results on heat transport and stability properties of all equilibria bifurcating from the origin tracked using MATCONT, a MATLAB software project for numerical continuation and bifurcation analysis of parameterized dynamical systems \cite{Dhooge_2008_MatCont}.  

After some experimentation, the HKC21 model was selected for analysis using MATCONT, since this model was found to be large enough to represent the PDE over a decent range of Rayleigh numbers, but also small enough where a host of numerical continuations could be performed.  Fixing $\Rot,\Pra,\Sha$ at $(0,10,1/\sqrt{2})$ and starting at the origin, the Rayleigh number was increased from $0$ to $2000$, and the equilibria bifurcating from the origin were tracked.  The second generation equilibria arising at each branch point were then tracked through $\Ray = 2000$, as well as the third generation.  The heat transport along each branch was then computed, and the result is depicted in Figure \ref{fig:HKC21_BifDiag_Ray} below.  This figure displays the stability information as well, where the heat transport corresponding to locally stable equilibria is depicting in blue vs the heat transport for unstable equilibria is depicted in red.  The branches of equilibria also exhibited many Hopf bifurcations, hence these were continued as well.  Some of these branches of periodic solutions were locally stable for a range of Rayleigh numbers, hence locally stable periodic solutions are depicted in purple vs locally unstable periodic solutions in orange.  These branches of periodic solutions typically lost stability via torus bifurcations, and hence the stable heat transport values were not tracked further due to the extreme complexity.  Furthermore, a time dependent solution beginning from a random perturbation of the uniform state was computed for each of the following Rayleigh numbers, and depicted in green:
\[  \Ray =  0:5:2000 \text{ . } \]
A very complicated bifurcation structure is evident in Figure \ref{fig:HKC21_BifDiag_Ray}, and in fact Figure \ref{fig:HKC21_BifDiag_Ray} only displays the branches which were most directly related to the stable transport values.  Dozens of branches had to be left out to avoid overcrowding the plot.  Most interestingly, there appear branches which emerge from the origin as unstable, but later restabilize for a range of Rayleigh numbers by ejecting further unstable equilibria, leading to the presence of multiple stable values of heat transport at a given Rayleigh number.  From Figure \ref{fig:HeatTransportModelComparison} panel (a) one sees that the HKC-21 model might only be accurate for the PDE up to around $\Ray = 500$, but the multi-stability appears well before this, and the dynamics of the HKC-21 model are interesting nevertheless.  In the region of multi-stability, note that the green points do not fall exactly on the blue equilibria branches, and the points seem roughly evenly distributed between the two branches.  This might indicate that these two branches have comparably sized basins of attraction which are fairly close together, and the values of the time dependent trajectories are tainted by finite integration time errors.

\begin{figure}[H]
	\begin{center}
		\includegraphics[height=90mm]{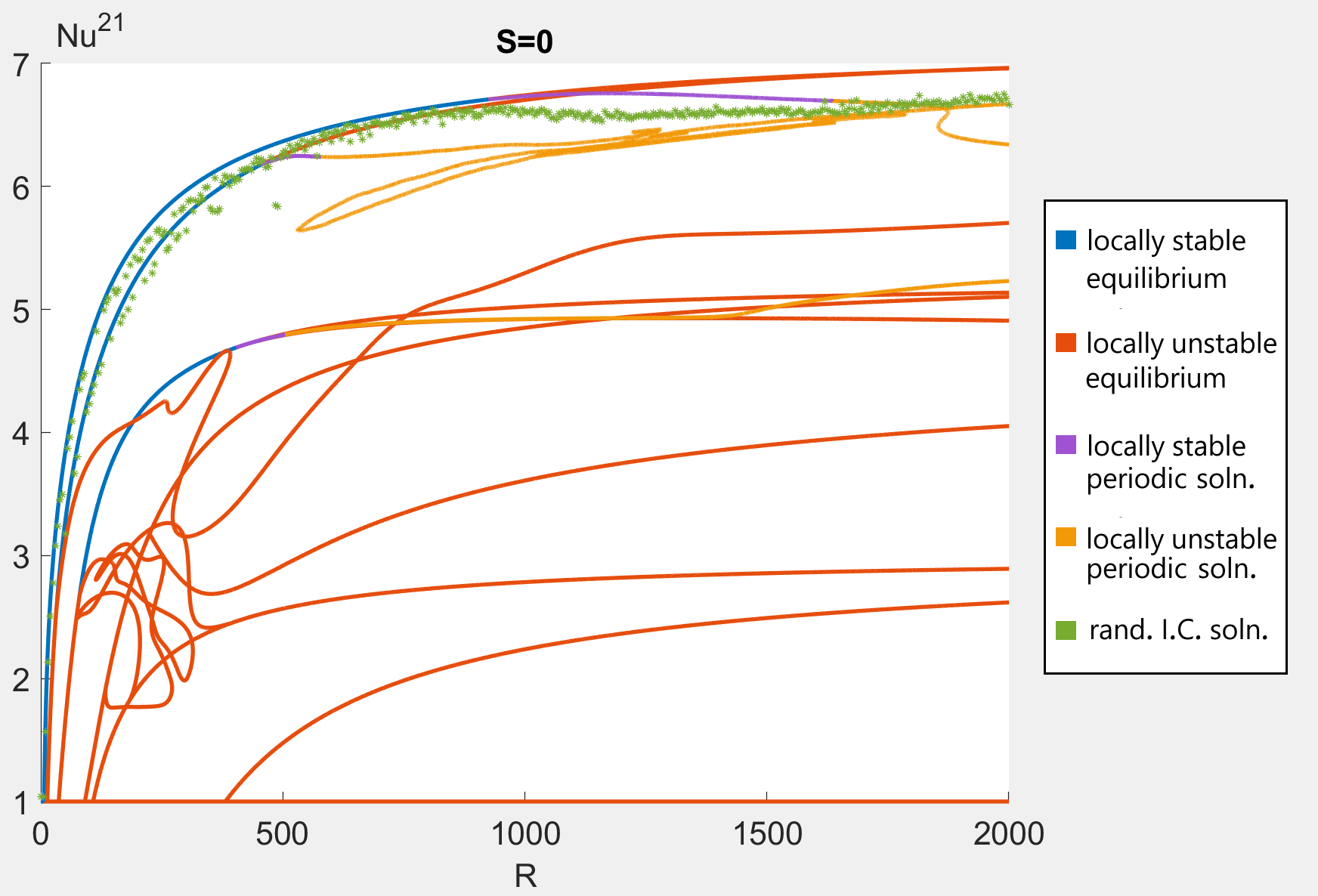}
		\caption{Heat transport bifurcation diagram for HKC-21 with the Rayleigh number as the control parameter at $\Rot =0$.}
		\label{fig:HKC21_BifDiag_Ray}
	\end{center}
\end{figure}

\begin{figure}[H]
	\begin{center}
		\includegraphics[height=80mm]{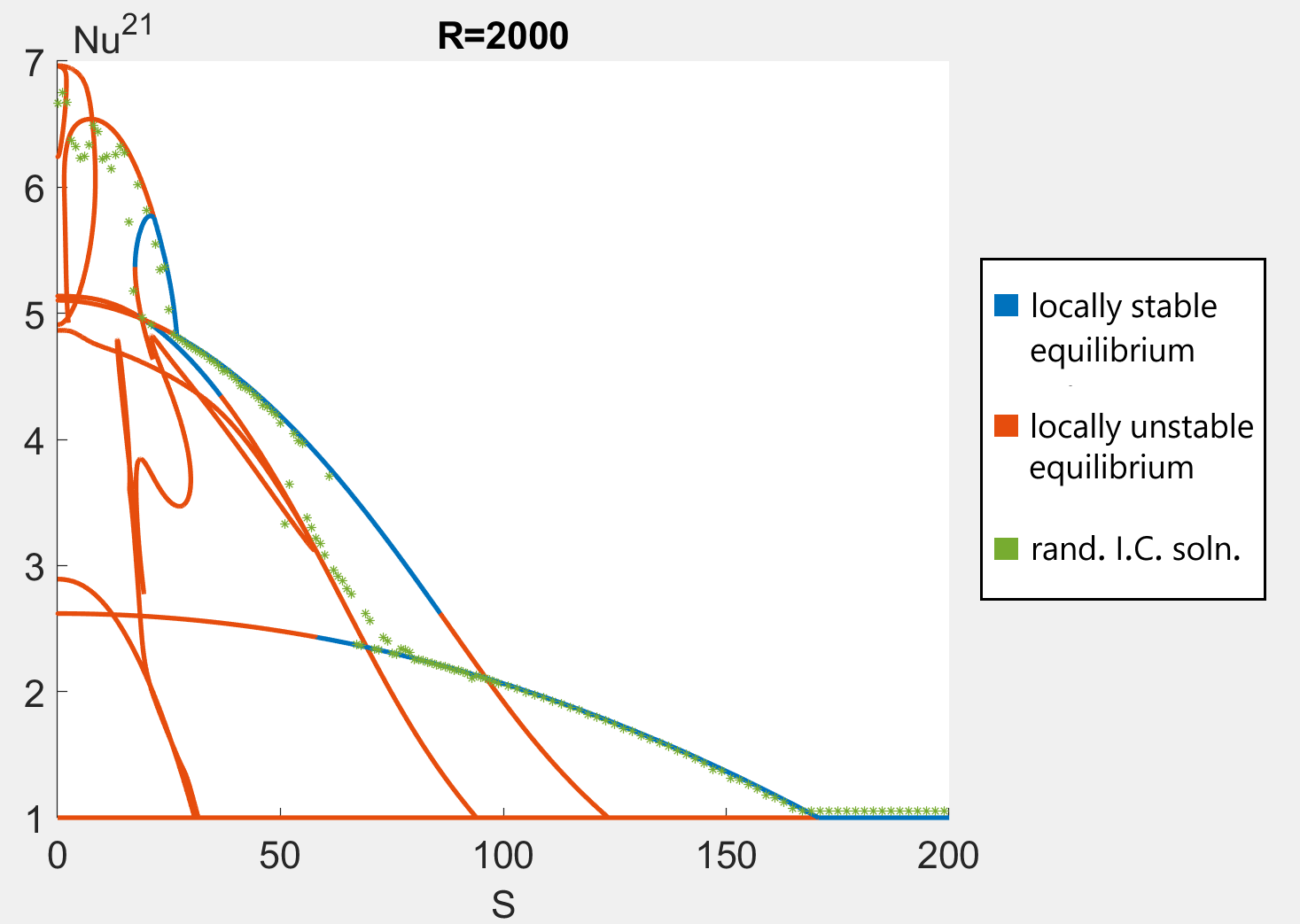} 
		\caption{Heat transport bifurcation diagram for HKC-21 with the rotation number as the control parameter at $\Ray =2000$.  }
		\label{fig:HKC21_BifDiag_Rot}
	\end{center}
\end{figure}

Figure \ref{fig:HKC21_BifDiag_Rot} depicts the analogous study using the rotation number as the control parameter at $\Ray = 2000$.  In this case, a time dependent solution beginning from a random perturbation of the uniform state was computed for each of the following rotation numbers:
\[  \Rot =  0:1:200 \text{ . } \]
These were again depicted in green.  At this high Rayleigh number the branches of periodic solutions were found to be too expensive to compute, hence they are not depicted.  However, the equilibria apparently regain stability fairly quickly and seem to explain the observed values of heat transport after $\Rot \sim 25$.  As one might expect, the bifurcation diagram was found to be symmetric in $\Rot$, hence for $\Rot < 0$ the bifurcation diagram is the mirror image of that shown in panel (b).  However, it is interesting because several of the branches form closed loops with their mirror image which cannot be accessed by tracing branches from the origin.  This is for instance the case for the branch with maximal heat transport at $\Rot = 0$.  In particular the stability might not be explained by tracing the stability to the origin.  For instance, in the region $60 \leq \Rot \leq 80$, note that the green points fall fairly far away from the two stable branches.  While this could again just be a finite integration time error, the trajectories seemed to show good convergence in time here.  It could be that there is a stable branch here which could not be found by tracking the bifurcations either beginning from the origin or beginning from one of the branches found with $\Rot = 0$.

\begin{figure}[H]
	\begin{center}
		\includegraphics[height=70mm]{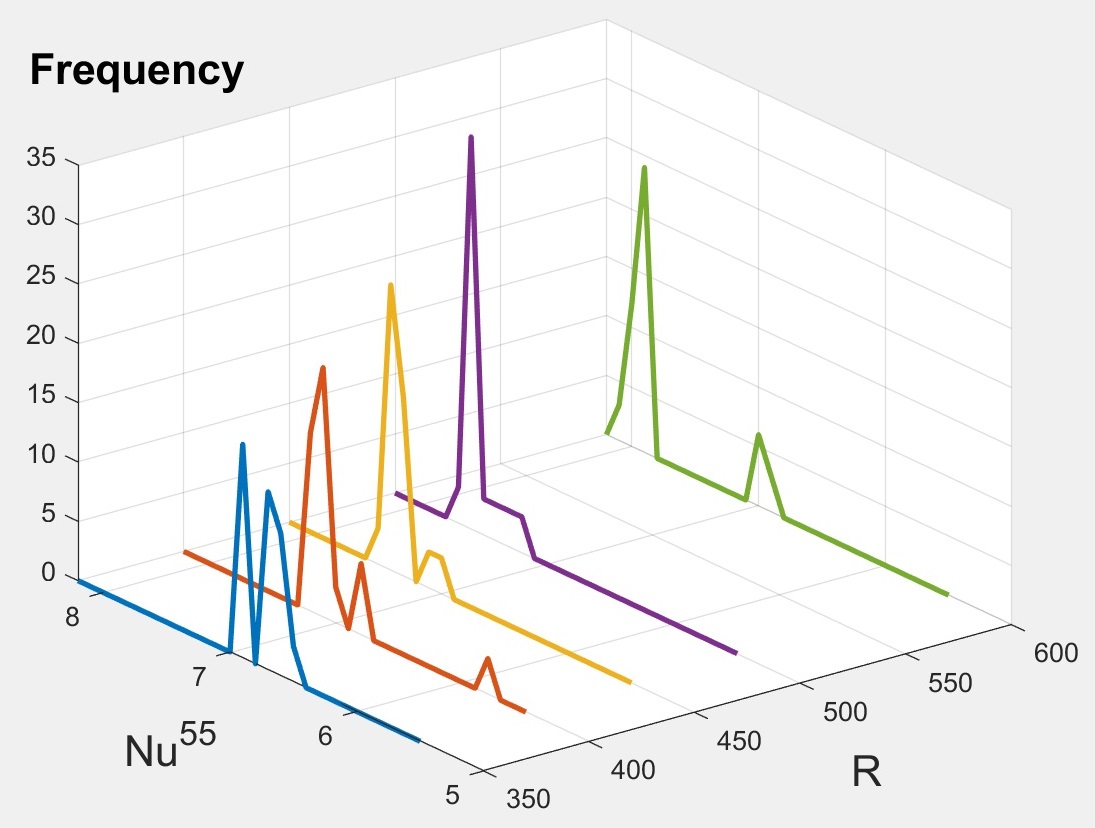}
		\caption{A frequency plot counting the number of initial conditions exhibiting various Nusselt number values for $350 \leq \Ray \leq 600$, $\Rot = 0$.}
		\label{fig:HKC55_IC_Statistics}
	\end{center}
\end{figure}

In order to study this apparent multi-stability further, the HKC-55 model was selected.  From Figure \ref{fig:HeatTransportModelComparison}, one sees that this model could be sufficiently high dimensional to represent the PDE up to roughly $\Ray = 1000$, but also sufficiently low dimensional where many simulations could be performed for different initial conditions.  To gather statistical information regarding the multi-stability, trajectories were simulated at 48 initial conditions obtained from randomly perturbing from the uniform state.  This was done for each of the following $\Ray$ and $\Rot$ values:
\begin{equation} \label{HKC55_MultiStabStats_ParameterRange}  \Ray = [0:50:500 , 600:100:1000 , 1250:250:3000  ]  \hspace{.5 cm} \text{ , } \hspace{.5 cm}  \Rot =  0:50:300  . \end{equation}
The results are plotted in Figure \ref{fig:HKC55_IC_Dependence} below.  The trajectories were then binned binned into groups according to their Nusselt number, where a trajectory was put into group $i$ if its Nusselt number was closest to the $i^{th}$ element of the following vector:
\begin{equation} \label{HKC55_NusseltBin}  \mathsf{Nu}^{55} = 0:0.5:10 . \end{equation}
The number of initial conditions falling into each bin is depicted in Figure \ref{fig:HKC55_IC_Statistics} at $\Rot = 0$ for a few Rayleigh number values.

\begin{figure}[H]
	\begin{center}
		\begin{tabular}{cc}
			\includegraphics[height=65mm]{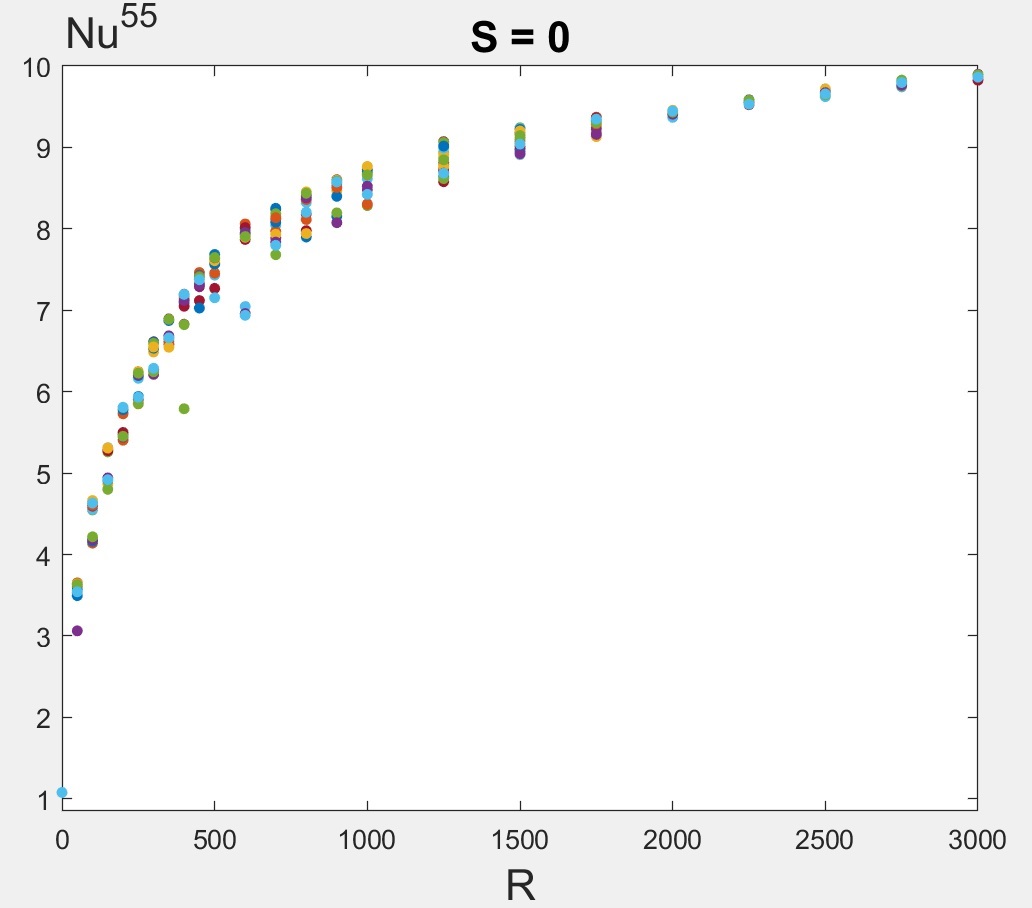}
			&
			\includegraphics[height=65mm]{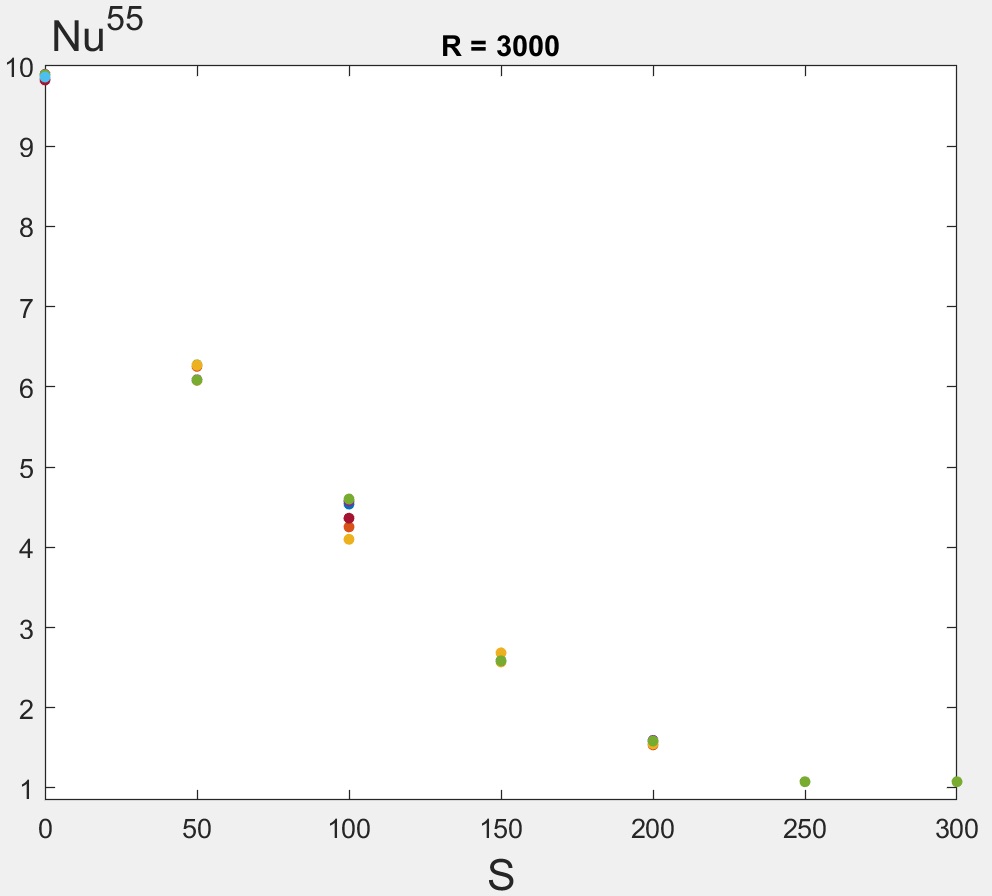} \\
			(a) & (b)
		\end{tabular}
		\caption{(a) The dependence of the heat transport on the Rayleigh number for 48 random initial conditions in the HKC-55 model (b) A frequency plot counting the number of initial conditions exhibiting various Nusselt number values.}
		\label{fig:HKC55_IC_Dependence}
	\end{center}
\end{figure}

At $\Ray = 350$ the distribution of Nusselt values is distinctly bimodal with two peaks of roughly even frequency, whereas for other values there seems to be a more dominant heat transport value.  Nevertheless at all values of the Rayleigh in Figure \ref{fig:HKC55_IC_Statistics} there is a non-trivial spread of Nusselt values, reflecting a multi-modal distribution.

\section{Conclusions}

\label{sec:Conclusions}

This paper has provided several results to aide both theoretical and applied researchers.   Explicitly computing the formulas for the Galerkin truncations enabled the HKC hierarchy to be implemented in MATLAB, which can be of use for applied researchers.  Criteria \ref{Crit:EnergyCrit}, \ref{Crit:VortCrit} were shown to not only guarantee the energetic consistency desired by climate researchers, but compact dynamical behavior as well.  The well-posedness theory and upper bound on the attractor dimension carry over from the two-dimensional case, and the rotation in the present case apparently does not ruin the smoothness properties or increase the radius or dimension of the attractor.  By studying the local bifurcations at the origin, one mechanism by which the attractor increases in complexity was identified and a rate of increase was given.  Furthermore, several heuristic steps toward an analysis of heat transfer were given in the numerical results, including a basic outline of the dependence of the heat transport on the Rayleigh and rotation numbers, convergence results as one looks at higher dimensional HKC models, and results on bifurcations and multi-stability.  In particular, multiple stable values of heat transport presents an interesting challenge for representing convection as a sub-grid process, since one would in fact need to know what is going on beneath the subgrid scale to get the correct value of transport.  The studies conducted herein could inform a stochastic representation of sub-grid convection, where the stochastic measure is chosen to be consistent with the frequency statistics in Figure \ref{fig:HKC55_IC_Statistics}.

These results invite further analysis.  On the applied side, one natural next step is to extend the hierarchy to allow fully 3-dimensional flows, rather than horizontally aligned flows.  In this paper a careful study of the convergence across models was conducted for Rayleigh and rotation numbers only up to around $5000$ and $400$ respectively, but with more computational resources one could certainly use the HKC models herein at higher parameter values.  One could push this even further by exploiting the convolution type conditions in \eqref{CompatibilityCond}, for example one could study the special class of solutions to \eqref{BoussinesqCoriolis} such that $\textbf{m} \in 100 \mathbb{Z}^2_{\geq 0}$ for all $u^{\textbf{n}},\theta^{\textbf{n}}$.  Even using relatively few modes one could study solutions with very fine scale structures using this class.  A sequel paper is already underway to address these issues.  The results could then be used to compare heat transport in established climate models such as ICON-A \cite{Giorgetta2018ICON} or to specific heat transport parameterizations.  For the model studied herein such direct comparisons would only be reasonable in a restricted physical regime.  One could add moisture into the model and it seems very likely that the energetic consistency Criteria would remain the same.  A further challenge would be to consider more physically realistic boundary conditions, for instance with no-slip conditions or allowing the upper boundary to behave as a free surface. 

On the theoretical side, the results herein open up a new approach toward obtaining universal upper bounds on the Nusselt number.  The standard approach for such bounds is the background field method \cite{ConstantinDoering_1996}, which for instance provides the sharpest known bound in the non-rotating, free slip case \cite{WhiteheadDoering_2011}, namely that for fixed $\Pra,\Sha$ there exists a constant $C_{\Pra,\Sha}$ such that 
\begin{equation} \label{HeatTransportBound_BackgroundFieldMethod} \sup_{(\textbf{u}_0,T_0)} \mathsf{Nu} \leq 1 + C_{\Pra,\Sha} \Ray^{5/12} \hspace{1 cm } \text{ for all} \hspace{.25 cm } \Ray \geq 0 . \end{equation}
To the best knowledge of the author the only analytical result in the rotating case was obtained by Constantin et al. \cite{Constantin_1999,Constantin_2001} in the special case where $\Pra = \infty$, very far from relevant for atmospheric flows.  The background field method appears to fail for the rotating case when $\Pra < \infty$, since many of the global analytical estimates use quadratic expressions such as the kinetic energy $\| \textbf{u} \|_{L^2}^2$ formed by computing an inner product of a term with itself, and due to the skew-symmetry of the Coriolis term the dependence on rotation disappears.  In this case one might have to apply more advanced techniques such as the more general auxiliary function method \cite{Chernyshenko2022Relationship}, which when combined with convex optimization techniques has been successful in obtaining bounds for several ODE systems \cite{FantuzziGoluskin,Goluskin2018,Goluskin2021}.  The results of this paper allow one to apply such techniques to arbitrary Fourier truncations of the BC model.  

While the bounds on the attractor dimension obtain here give a rough idea of the required spatial resolution, the reduction of the BC model to the ODE system is still not fully rigorous.  For more ideal equations, such as reaction-diffusion or Navier-Stokes-alpha models, one can make a fully rigorous reduction to an ODE system by showing the existence of an inertial manifold, which loosely speaking is a graph $\Gamma: \mathbb{R}^M \mapsto \textbf{L}^2$ such that the solution $\textbf{u}^M(t)$ to a sufficiently high dimensional ODE (say, of dimension $M$) satisfies
\[ \| \textbf{u}(t) - ( \textbf{u}^M(t) + \Gamma \textbf{u}^M(t) ) \|_{\textbf{L}^2} \leq C e^{- t} . \]
With an inertial manifold in hand one could recover infinite time averages from the PDE (such as the Nusselt number) by solutions of the ODE system.  Progress continues in the general theory of inertial manifolds \cite{Zelik_2014}, but challenges remain to prove existence for Navier-Stokes type models.  However, studying inertial manifolds in the models developed herein may provide insights regarding how this could be done for the full BC system.

\appendix

\section{Non-dimensionalization and values of the physical parameters}

\label{app:PhysParam}

Using variables with physical units, the three dimensional Boussinesq equations with a Coriolis force are given as:
\begin{equation}
	\label{PhysicalBoussinesqCoriolis}
	\begin{split}
		\partial_{\tilde{t}} \textbf{v} + \textbf{v} \cdot \nabla \textbf{v} + \frac{1}{\rho} \nabla \tilde{p} & = \nu \Delta \textbf{v} + g \alpha (\tilde{T} -\tilde{T}_r) \uVecThree - 2 s \hspace{.5 mm} \uVecThree\times \textbf{v} \text{ , }  \\
		\nabla \cdot \textbf{v} & = 0 \text{ , } \\
		\partial_t \tilde{T} + \textbf{v}\cdot \nabla \tilde{T} & =  \kappa \Delta \tilde{T} .
	\end{split}
\end{equation}
Here $\textbf{v} = (v_1,v_2,v_3)$, $\tilde{T}$, $\tilde{p}$ are the velocity, temperature and pressure of a fluid depending on $\textbf{y} = (y_1,y_2,y_3) \in [0,L_1]\times [0,L_2] \times [0,H]$ and $\tilde{t} \geq 0$.  The parameters $\rho, \nu, g, \alpha, T_r, s, \kappa$ are the density, kinematic viscosity, gravitational constant, thermal expansion coefficient, local angular speed of the rotating frame, thermal diffusivity and reference temperature.  The temperature is held constant at the top and bottom, hence $\tilde{T}(y_1,y_2,0) = \tilde{T}_b$, $\tilde{T}(y_1,y_2,H) = \tilde{T}_t$, and it is assumed that $\tilde{T}_b > \tilde{T}_t$.  Without loss of generality one can choose the reference temperature $\tilde{T}_r = \tilde{T}_t$ and rescale as follows:
\[ \textbf{x} = \frac{\pi \textbf{y}}{H} \hspace{.5 cm} \text{ , } \hspace{.5 cm} t = \frac{H^2}{\pi^2 \kappa} \tilde{t} \hspace{.5 cm} \text{ , } \hspace{.5 cm} \textbf{u} = \frac{H \textbf{v}}{\pi \kappa} \hspace{.5 cm} \text{ , } \hspace{.5 cm} T = \frac{\tilde{T}-\tilde{T}_t}{\tilde{T}_b-\tilde{T}_t} \hspace{.5 cm} \text{ , } \hspace{.5 cm} p = \frac{H^2}{\pi^2 \kappa \nu \rho} \tilde{p} \hspace{.5 cm} \text{ , } \]
and obtain the non-dimensionalized system \eqref{BoussinesqCoriolis}.  The dimensionless parameters in \eqref{ParameterVals} are given explicitly by
\begin{equation} \label{ParameterDefs} \Pra = \frac{\nu}{\kappa} \hspace{1 cm} \text{ , } \hspace{1 cm}  \Ray = \frac{g \alpha H^3 (T_t-T_b)}{ \pi^4 \kappa \nu}  \hspace{1 cm} \text{ , } \hspace{1 cm} \Rot = \frac{2 s H^3}{\pi^3 \nu \kappa} . \end{equation}	
The gravitational constant near the surface of the Earth is given by $g \approx 9.8 \frac{m}{s^2}$.  Taking the fluid to represent air, some values for the physical parameters are as follows:
\[ \rho \approx 1.2 \frac{kg}{m^3} \hspace{.5 cm} \text{ , } \hspace{.5 cm} \nu \approx 4.3 \cdot 10^{-5} \frac{m^2}{s} \hspace{.5 cm} \text{ , } \hspace{.5 cm} \alpha \approx 1.5 \cdot 10^{-3} \frac{1}{K} \hspace{.5 cm} \text{ , } \hspace{.5 cm} \kappa \approx 1.9 \cdot 10^{-5} \frac{m^2}{s} . \]
There is some choice involved in the other parameters, depending on what is to be modelled.  For instance, the local angular speed is zero at the equator or $s \approx 1.6 \cdot 10^{-5} \frac{1}{s}$ near the poles.  If the box $\Omega$ is taken to represent the troposphere then some approximate values are given by 
\[ T_b \approx 300 K \hspace{.5 cm} \text{ , } \hspace{.5 cm} T_t \approx 222 \hspace{.5 cm} \text{ , } \hspace{.5 cm} H \approx 10^4 m . \]
Inserting these values into \eqref{ParameterDefs} one obtains \eqref{ParameterVals}.









\section{Model construction}

\subsection{Completeness of the general Fourier expansion for divergence-free vector fields}

\label{app:FourierExpDeriv}

First, artificially extend the domain by defining $\textbf{u}$ for $x_3 \in [-\pi,0)$ by an odd extension for $u_3$ and an even extension for $u_1,u_2$ in order to satisfy the boundary conditions.  The sinusoidal functions are a complete basis for $L^2([0,2\pi])$, hence one can expand each component of the vector field.  As mentioned above one can assume each component has zero mean, and furthermore due to the even/odd properties above $u_1$, $u_2$ must be given only in terms of $\cos( m_3 x_3)$ and $u_3$ must be given only in terms $\sin( m_3 x_3)$.  Defining the following index sets
\[ \mathscr{I}_1 = \{  m_1 , m_3 > 0 \} \hspace{.25 cm} \text{ , } \hspace{.25 cm} \mathscr{I}_2 = \{  m_1  > 0 , m_3 = 0\} \hspace{.25 cm} \text{ , } \hspace{.25 cm} \mathscr{I}_3 = \{  m_3  > 0 , m_1 = 0\} , \]
one therefore has the following general expansion 
\begin{align}
	\textbf{u}(\textbf{x},t) = & \sum_{\textbf{m} \in \mathscr{I}_1} \big [ \big ( \hat{u}_{1}^{\textbf{m},1} \uVecOne + \hat{u}_{2}^{\textbf{m},1} \uVecTwo  \big ) \cos( m_3 x_3) +  \hat{u}_{3}^{\textbf{m},1} \uVecThree \sin( m_3 x_3) \big ] \cos ( \Sha m_1 x_1 ) \notag \\ & \hspace{2 cm} + \big [ \big ( \hat{u}_{1}^{\textbf{m},2} \uVecOne + \hat{u}_{2}^{\textbf{m},2} \uVecTwo  \big ) \cos( m_3 x_3) +  \hat{u}_{3}^{\textbf{m},2} \uVecThree \sin( m_3 x_3) \big ] \sin ( \Sha m_1 x_1 ) \notag \\ & + \sum_{\textbf{m} \in \mathscr{I}_2} \big ( \hat{u}_{1}^{\textbf{m},1} \uVecOne + \hat{u}_{2}^{\textbf{m},1} \uVecTwo  \big ) \cos ( \Sha m_1 x_1 ) + \big ( \hat{u}_{1}^{\textbf{m},2} \uVecOne + \hat{u}_{2}^{\textbf{m},2} \uVecTwo  \big ) \sin ( \Sha m_1 x_1 ) \notag \\ & + \sum_{\textbf{m} \in \mathscr{I}_3} \big ( \hat{u}_{1}^{\textbf{m},1} \uVecOne + \hat{u}_{2}^{\textbf{m},1} \uVecTwo  \big ) \cos( m_3 x_3) +  \hat{u}_{3}^{\textbf{m},1} \uVecThree \sin( m_3 x_3)  , \notag
\end{align}
where for $\textbf{m} \in \mathscr{I}_1$ the Fourier coefficients are given by
\begin{align} \hat{u}^{\textbf{m},1}_{\ell}(t) & = \frac{1}{V^2} \int_{0}^{\frac{2\pi}{\Sha}}\int_{-\pi}^{\pi} u_\ell(\textbf{x},t) \cos(\Sha m_1 x_1) \cos( m_3 x_3) d\textbf{x} \hspace{.1 cm} \text{ , } \notag \\ \hat{u}^{\textbf{m},2}_{\ell}(t) & = \frac{1}{V^2} \int_{0}^{\frac{2\pi}{\Sha}}\int_{-\pi}^{\pi} u_\ell(\textbf{x},t) \sin(\Sha m_1 x_1) \cos( m_3 x_3) d\textbf{x} \text{ , } \notag \end{align}
for $\ell = 1,2$ and 
\begin{align} \hat{u}^{\textbf{m},1}_{3}(t) & = \frac{1}{V^2} \int_{0}^{\frac{2\pi}{\Sha}}\int_{-\pi}^{\pi} u_3(\textbf{x},t) \cos(\Sha m_1 x_1) \sin( m_3 x_3) d\textbf{x} \hspace{.1 cm} \text{ , } \notag \\ \hat{u}^{\textbf{m},2}_{3}(t) & = \frac{1}{V^2} \int_{0}^{\frac{2\pi}{\Sha}}\int_{-\pi}^{\pi} u_3(\textbf{x},t) \sin(\Sha m_1 x_1) \sin( m_3 x_3) d\textbf{x} \text{ . } \end{align}
For $\textbf{m} \in \mathscr{I}_2 ,  \mathscr{I}_3$ the Fourier coefficients are defined similarly.

For $\textbf{m} \in \mathscr{I}_1$, one can use the definitions of the $\hat{u}^{\textbf{m},j}_{\ell}$, integrate by parts and apply the incompressibility condition to establish that for each $j$ one must have the relation
\begin{equation} \label{FourierBasis_DivFreeCond} (-1)^{j} \Sha m_1 \hat{u}^{\textbf{m},j+1}_{1} + m_3 \hat{u}^{\textbf{m},j}_{3} = 0 . \end{equation}
This relation then reduces the number of coefficients to be considered for each $\textbf{m}$.  In this case each component of the vector field has two Fourier coefficients for a total of six, namely $\hat{u}^{\textbf{m},j}_{\ell}$ for $j \leq 2$, $\ell \leq 3$.  However, the condition \ref{FourierBasis_DivFreeCond} implies there can be at most one independent Fourier coefficients for each $j$, hence one ends up with four independent coefficients.  Therefore there is some transformation carrying the four corresponding vector fields to the four vector fields in \eqref{VectorFieldDef_FreeSlip}.  For $\textbf{m} \in \mathscr{I}_2, \mathscr{I}_3$, the analogue of \eqref{FourierBasis_DivFreeCond} shows that one of the Fourier coefficients must be zero, hence in these cases there is some transformation carrying the two remaining linearly independent vector fields to the two vector fields in \eqref{VectorFieldDef_FreeSlip}.

\subsection{Derivation of the explicit form of the nonlinear terms}

\label{app:NonlinearDerivation}

We drop the superscript $M$ for this derivation.  Expanding the expressions \eqref{AbstractNonlinear} using the Fourier expansions \eqref{GeneralFourierExpansion}, one obtains the following:
\begin{equation} \notag
	\begin{split}  \int_{\Omega} \textbf{v}^{\textbf{n}} \cdot \big [ \big ( \textbf{u} \cdot \nabla \big )  \textbf{u} \big ]  d\textbf{x} & = \int_{\Omega} \textbf{v}^{\textbf{n}} \cdot \big [ \big ( \sum_{\textbf{n}'} u^{\textbf{n}'} \textbf{v}^{\textbf{n}'}  \cdot \nabla \big )  \sum_{\textbf{n}''} u^{\textbf{n}''} \textbf{v}^{\textbf{n}''}  \big ] d\textbf{x}  \\ & = \sum_{\textbf{n}'} \sum_{\textbf{n}''}  u^{\textbf{n}'} u^{\textbf{n}''} \int_{\Omega} \textbf{v}^{\textbf{n}} \cdot \big [ (\textbf{v}^{\textbf{n}'} \cdot \nabla ) \textbf{v}^{\textbf{n}''} \big ] d\textbf{x} , \\ \int_{\Omega} f^{\textbf{n}} \cdot \big [ \textbf{u} \cdot \nabla \theta  \big ] d\textbf{x} & = \int_{\Omega} f^{\textbf{n}} \cdot \big [ \sum_{\textbf{n}'} u^{\textbf{n}'} \textbf{v}^{\textbf{n}'} \cdot \nabla \sum_{\textbf{n}''} \theta^{\textbf{n}''} f^{\textbf{n}''} \big ] d\textbf{x} \\ & = \sum_{\textbf{n}'} \sum_{\textbf{n}''}  u^{\textbf{n}'} \theta^{\textbf{n}''} \int_{\Omega} f^{\textbf{n}} \cdot \big [ \textbf{v}^{\textbf{n}'} \cdot \nabla f^{\textbf{n}''} \big ] d\textbf{x} . \end{split}
\end{equation}
Having extracted the time dependent terms, denote the (time-independent) integrals via
\[ I_{\textbf{u}}^{\boldsymbol{\alpha}} = \int_{\Omega} \textbf{v}^{\textbf{n}} \cdot \big [ (\textbf{v}^{\textbf{n}'} \cdot \nabla ) \textbf{v}^{\textbf{n}''} \big ] d\textbf{x} \hspace{.5 cm} \text{ , } \hspace{.5 cm} I_{\theta}^{\boldsymbol{\alpha}} = \int_{\Omega} f^{\textbf{n}} \cdot \big [ \textbf{v}^{\textbf{n}'} \cdot \nabla f^{\textbf{n}''} \big ] d\textbf{x} , \]
where we have used $\boldsymbol{\alpha} := (\textbf{n},\textbf{n}',\textbf{n}'')$.  Next, we expand the dot products in terms of the vector components via
\[ I_{\textbf{u}}^{\boldsymbol{\alpha}} =  \sum_{i = 1,3} \sum_{j \leq 3 } J_{i,j}^{\boldsymbol{\alpha}} \hspace{.5 cm} \text{ , } \hspace{.5 cm} I_{\theta}^{\boldsymbol{\alpha}} =  \sum_{i = 1,3} \tilde{J}_{i}^{\boldsymbol{\alpha}} \text{ , } \]
in which
\begin{equation} \label{NonlinDeriv_ExpDotProd} J_{i,j}^{\boldsymbol{\alpha}} = \int_{\Omega} v_{j}^{\textbf{n}}  \cdot \big [ v^{\textbf{n}'}_{i} \partial_{x_{i}} v^{\textbf{n}''}_{j} \big ] d\textbf{x} \hspace{.5 cm} \text{ , } \hspace{.5 cm} \tilde{J}_{i}^{\boldsymbol{\alpha}} = \int_{\Omega} f^{\textbf{n}} \cdot \big [ v^{\textbf{n}'}_{i} \partial_{x_{i}} f^{\textbf{n}''} \big ] d\textbf{x}  \text{ . } \end{equation}
Note the value $i = 2$ is excluded since the derivative in $x_2$ is zero.  Since the component of $v_{i}^{\textbf{n}'}$ must match the derivative $\partial_{x_{i}}$ in \eqref{NonlinDeriv_ExpDotProd}, and only the second component of $v_{i}^{(\textbf{m}',p_1',2)}$ is non-zero, one must have $c' = 1$. Similarly one must have $c = c''$, since otherwise $\textbf{v}^{\textbf{n}}, \textbf{v}^{\textbf{n}''}$ have no matching non-zero components, and hence the expression for $I_{\textbf{u}}^{\boldsymbol{\alpha}}$ simplifies:
\[  I_{\textbf{u}}^{\boldsymbol{\alpha}} = \left \{ \begin{array}{ll}
	J_{1,1}^{\boldsymbol{\alpha}} + J_{3,1}^{\boldsymbol{\alpha}} + J_{1,3}^{\boldsymbol{\alpha}} + J_{3,3}^{\boldsymbol{\alpha}} & \text{ if } c = 1 , \\
	J_{1,2}^{\boldsymbol{\alpha}} + J_{3,2}^{\boldsymbol{\alpha}}  & \text{ if } c = 2 .
\end{array} \right . \]

In order to explicitly resolve the integrals \eqref{NonlinDeriv_ExpDotProd} in a concise way, let $\mathsf{k}_3 = 1$, and let the sinusoidal functions $h^{m,p_1}$ and phase matrices $\mathcal{G}^{\textbf{n}}$ be defined
\begin{align} 
	h^{m,1}(y) := \cos ( m y ) & \hspace{.5 cm} \text{ , } \hspace{.5 cm}  h^{m,2}(y) := \sin ( m y ) , \notag \\
	\mathcal{G}^{(\textbf{m},p_1,1)} := \mathsf{diag} \big [  \frac{m_3}{|\mathcal{K}\textbf{m}|} , 0 , (-1)^{p_1} \frac{\Sha m_1}{|\mathcal{K}\textbf{m}|} \big ] & \hspace{.5 cm} \text{ , } \hspace{.5 cm} \mathcal{G}^{(\textbf{m},p_1,2)} := \mathsf{diag} \big [ 0,1,0 \big ] . \notag
\end{align}
Using this notation one can write all of the cases in \eqref{ThetaFourierBasis},\eqref{VectorFieldDef_FreeSlip} in the following general expressions
\begin{equation} \label{NonlinDeriv_Shorthand} 
	\textbf{v}^{\textbf{n}} = \frac{\eta^{\textbf{m}}}{V} \mathcal{G}^{\textbf{n}} \begin{pmatrix}
		h^{m_1,p_1+1}(\Sha x_1) \cos(m_3 x_3 ) \\ h^{m_1,p_1+1}(\Sha x_1) \cos(m_3 x_3 ) \\ h^{m_1,p_1}(\Sha x_1) \sin(m_3 x_3 )
	\end{pmatrix} \hspace{.15 cm} \text{ , } \hspace{.15 cm} f^{\textbf{n}} = \frac{\eta^{\textbf{m}}}{V} h^{m_1,p_1}(\Sha x_1) \sin (m_3 x_3 ) ,   
\end{equation}
where we recall that the phase and component indices $p_1,c$ are taken mod $2$.  Define also the sign function $\sigma_{i,j}^{p_1''}$ to give the signs produced by taking the derivative $\partial_{x_i}$ of the $j^{th}$ component in \eqref{NonlinDeriv_ExpDotProd}:
\[  \sigma_{i,j}^{p_1''} :=  \sigma_{i}^{p_1''} \sigma_{j} \hspace{.5 cm} \text{ , } \hspace{.5 cm} \sigma_{i}^{p_1''} := \left \{ \begin{array}{cl}
	(-1)^{p_1''+1} & \text{ if } i = 1 , \\
	-1 & \text{ if } i = 3 ,
\end{array} \right . \hspace{.5 cm} \text{ , } \hspace{.5 cm} \sigma_{j} := \left \{ \begin{array}{cl}
	1 & \text{ if }  j = 1,2, \\
	-1 & \text{ if }  j = 3 .
\end{array} \right .  \]
Finally, for $\boldsymbol{\mu} = (m,m',m'')$, $\boldsymbol{\phi} = (p_1,p_1',p_1'')$, let $E^{(\boldsymbol{\mu},\boldsymbol{\phi})}$ be the evaluation of the sinusoidal integral:
\[ E^{(\boldsymbol{\mu},\boldsymbol{\phi})} := \int_0^{2\pi} h^{m,p_1}(y) h^{m',p_1'}(y) h^{m'',p_1''}(y) dy . \]
We can now pull out all parameter dependence and explicit dependence on $\textbf{m},\textbf{m}',\textbf{m}''$ from the integrals by inserting \eqref{NonlinDeriv_Shorthand} into \eqref{NonlinDeriv_ExpDotProd}:
\begin{equation} \label{NonlinDeriv_WaveDep}  J_{i,j}^{\boldsymbol{\alpha}} =  \frac{\eta^{\textbf{m}}\eta^{\textbf{m}'}\eta^{\textbf{m}''}}{V^3} \mathcal{G}^{\textbf{n}}_{j,j} \mathcal{G}^{\textbf{n}'}_{i,i} \mathcal{G}^{\textbf{n}''}_{j,j} \sigma_{i,j}^{p_1''} \mathsf{k}_{i} m_{i}'' \hat{I}_{i,j}^{\boldsymbol{\alpha}} \hspace{.5 cm} \text{ , } \hspace{.5 cm} \tilde{J}_{i}^{\boldsymbol{\alpha}} =  \frac{\eta^{\textbf{m}}\eta^{\textbf{m}'}\eta^{\textbf{m}''}}{V^3} \mathcal{G}^{\textbf{n}'}_{i,i} \sigma_{i,3}^{p_1''} \mathsf{k}_{i} m_{i}'' \hat{I}_{i,3}^{\boldsymbol{\alpha}} , \end{equation}
in which the integrals $\hat{I}_{i,j}^{\boldsymbol{\alpha}}$ are given as follows (after using the substitutions $y = \Sha x_1, y = 2 x_3$):  	
\begin{align} 
	\hat{I}_{1,1}^{\boldsymbol{\alpha}} = \hat{I}_{1,2}^{\boldsymbol{\alpha}} = \frac{1}{2 \Sha} E^{(\boldsymbol{\mu}^1,p_1+1,p_1'+1,p_1'')}  E^{(\frac{\boldsymbol{\mu}^3}{2},1,1,1)} \hspace{.1 cm} & \text{ , } \hspace{.1 cm} \hat{I}_{3,1}^{\boldsymbol{\alpha}} = \hat{I}_{3,2}^{\boldsymbol{\alpha}} = \frac{1}{2 \Sha} E^{(\boldsymbol{\mu}^1,p_1+1,p_1',p_1''+1)}  E^{(\frac{\boldsymbol{\mu}^3}{2},1,2,2)} \text{ , } \notag \\  \label{NonlinDeriv_SinIntegrals}
	\hat{I}_{1,3}^{\boldsymbol{\alpha}} = \frac{1}{2 \Sha} E^{(\boldsymbol{\mu}^1,p_1,p_1'+1,p_1''+1)}  E^{(\frac{\boldsymbol{\mu}^3}{2},2,1,2)} \hspace{.1 cm} & \text{ , } \hspace{.1 cm} \hat{I}_{3,3}^{\boldsymbol{\alpha}} = \frac{1}{2 \Sha} E^{(\boldsymbol{\mu}^1,p_1,p_1',p_1'')}  E^{(\frac{\boldsymbol{\mu}^3}{2},2,2,1)} \text{ , } 
\end{align}
where $\boldsymbol{\mu}^k$ are as in \eqref{def:TriadNotation}  Note if either $\boldsymbol{\mu}^1 = \textbf{0}$ or $\boldsymbol{\mu}^3 = \textbf{0}$ then the integrals don't need to be evaluated, since $J_{i,j}^{\boldsymbol{\alpha}} = \tilde{J}_{i}^{\boldsymbol{\alpha}} = 0$ already follows from the fact that $m_k'= m_k'' = 0$, $c' = 1$ in \eqref{NonlinDeriv_WaveDep}.  Hence we can assume at least one wave number $m_k,m_k',m_k''$ is non-zero for both $k=1,3$.  

It remains to evaluate the eight integrals $E^{(\boldsymbol{\mu},\boldsymbol{\phi})}$ for $\boldsymbol{\phi} \in \{1,2\}^3$.  First, for any $\boldsymbol{\mu} \in \mathbb{Z}^3_{\geq 0}$ one has
\[ E^{(\boldsymbol{\mu},1,1,2)} = E^{(\boldsymbol{\mu},1,2,1)} = E^{(\boldsymbol{\mu},2,1,1)} = E^{(\boldsymbol{\mu},2,2,2)} = 0 \text{ . } \]
since these are integrals of odd functions over the domain $[-\pi,\pi]$ (after using a substitution).  This then gives the phase condition \eqref{CompatibilityCond} on $p_1,p_1',p_1''$.  So one need only consider the remaining four integrals, and as in \eqref{NonlinDeriv_SinIntegrals} $\boldsymbol{\mu}$ can have either integer or half integer components.  In order to evaluate these, one can repeatedly apply the angle addition formulas:

\vspace{-.2 cm}
\small

\begin{equation} \label{NonlinDeriv_TrigAddition} \begin{split}  \cos (m' y) \cos (m'' y) & = \frac{\cos ( (m' + m'')y) + \cos ((m' - m'')y )}{2} \text{ , } \\ \cos (m' y) \sin (m'' y) & = \frac{\sin ((m'+m'')y) - \sin ( (m' -m'')y )}{2} \text{ , } \\ \sin (m'y) \sin (m''y) & = \frac{\cos ((m'-m'')y ) -\cos ((m'+m'')y)}{2} \text{ , } \\ \sin (m'y) \cos (m''y) & = \frac{\sin ((m'+m'')y) + \sin ((m'-m'')y )}{2}. \end{split} \end{equation}

\normalsize

\noindent To see how the evaluation procedure works, consider the specific example 
\[ E^{(\boldsymbol{\mu},1,2,2)} = \int_0^{2\pi} \cos(my) \sin (m'y) \sin( m''y) dy . \]
In this case, $E^{(\boldsymbol{\mu},1,2,2)} = 0$ if $m' = 0$ or $m'' = 0$ since $\sin(0) = 0$, so $E^{(\boldsymbol{\mu},1,2,2)}$ should be proportional to the amplitude factor $A^{(\boldsymbol{\mu},1,2,2)}:=\chi^{m',m''}$, where
\[ \chi^{m,m'} := \left \{ \begin{array}{cc} 0 & \text{ if } m = 0 \text{ or } m' = 0 , \\ 1 & \text{ otherwise. } \end{array} \right . \]
Expanding via \eqref{NonlinDeriv_TrigAddition}, one obtains 
\begin{align}
	E^{(\boldsymbol{\mu},1,2,2)} = \frac{1}{4} \int_0^{2\pi} \big [ & \cos \big ( (m + m'-m'') y \big ) + \cos \big ( (m - m' + m'') y \big ) \notag \\ & - \cos \big ( (m - m'-m'') y \big ) - \cos \big ( (m+ m'+m'') y \big ) \big ] dy . \notag 
\end{align} 
If no convolution condition is satisfied (i.e. $m \neq |m'\pm m''|$) then in the case where $m,m',m''$ are integers one obtains $w^{-1} \sin(wx)$ for each of the integrands (for some $w \neq 0$), which then evaluates to zero due to periodicity.  If $m,m',m''$ are only half integers, then a convolution-type condition must still be satisfied, since the function $w^{-1} \sin(wx)$ is zero at both boundaries if $w \neq 0$ is a half integer.  On the other hand, if $m',m''> 0$ and a convolution-type condition is satisfied, then the integral is non-zero.  If $m = m' + m''$ then it is negative, otherwise it is positive.  The absolute value is equal to $\frac{\pi}{2}$ if $m,m',m''$ are all non-zero, or $\pi$ if $m = 0$.  By repeating this procedure for each of the remaining integrals $E^{(\boldsymbol{\mu},\boldsymbol{\phi})}$, one arrives at the following general statements:
\begin{enumerate}
	\item $E^{(\boldsymbol{\mu},\boldsymbol{\phi})} = 0$ if the triad $\boldsymbol{\mu}$ doesn't satisfy a convolution-type condition of the form $m = | m' \pm m''|$.  
	\item If $\boldsymbol{\mu} \neq \textbf{0}$ satisfies a convolution condition and $\boldsymbol{\phi} \in \{\boldsymbol{\xi}^i \}_{i \leq 4}$ then 
	\begin{equation} \label{NonlinDeriv_EvalResult} E^{(\boldsymbol{\mu},\boldsymbol{\phi})} = \frac{\pi}{2(\eta^{m}\eta^{m'}\eta^{m''})^2} s^{(\boldsymbol{\mu},\boldsymbol{\phi})}  A^{(\boldsymbol{\mu},\boldsymbol{\phi})} , \end{equation}
	for the sign coefficients $s^{(\boldsymbol{\mu},\boldsymbol{\phi})}$ given in \eqref{def:SignCoefs}, and amplitude factors $A^{(\boldsymbol{\mu},\boldsymbol{\phi})} $ given as follows:
	\[ A^{(\boldsymbol{\mu},\boldsymbol{\xi}^1)} = 1 \hspace{.25 cm} \text{ , } \hspace{.25 cm} A^{(\boldsymbol{\mu},\boldsymbol{\xi}^2)} = \chi^{m',m''} \hspace{.25 cm} \text{ , } \hspace{.25 cm} A^{(\boldsymbol{\mu},\boldsymbol{\xi}^3)} = \chi^{m,m''}\hspace{.25 cm} \text{ , } \hspace{.25 cm} A^{(\boldsymbol{\mu},\boldsymbol{\xi}^4)} = \chi^{m,m'} \text{ . } \]
\end{enumerate}
Note however that the amplitude factors $A^{(\boldsymbol{\mu},\boldsymbol{\phi})}$ enter into the analysis to handle the case $\sin ( 0 ) = 0$.  However, in our choice of admissible phases \eqref{PhaseIndexSets} we have specifically excluded these cases, and by this reason one should always have $\sin(mx), \sin(m'x) \neq 0$.  The sinusoid $\sin(m''x)$ for $m'' = 0$ might not be avoided by \eqref{PhaseIndexSets}, but since this sinusoid is obtained by taking a derivative one always has the factor $m'' =0 $ in \eqref{NonlinDeriv_WaveDep} which enforces $J^{\boldsymbol{\alpha}}_{i,j} = 0$ or $\tilde{J}_{i}^{\boldsymbol{\alpha}} = 0$.  Hence for all $\boldsymbol{\mu},\boldsymbol{\phi}$ chosen via \eqref{PhaseIndexSets} one can take $A^{(\boldsymbol{\mu},\boldsymbol{\phi})} = 1$, and \eqref{def:NonlinCoeffs1} - \eqref{def:Zeta} then follow from collecting \eqref{NonlinDeriv_ExpDotProd} - \eqref{NonlinDeriv_EvalResult} and simplifying.

\section{Proof of vorticity balance}

\label{app:VortBal}

\begin{proof}
	Note that the vorticity $\vorticity^M = \nabla \times \textbf{u}^M$ must satisfy
	\begin{equation} \label{TruncatedVorticity} \partial_{t} \vorticity^M = \Pra \Delta \vorticity^M +  \nabla \times \mathcal{P}^{M}_{\textbf{u}} \big [ \Pra \Ray \theta^M \uVecThree - \Pra \Rot ( \uVecThree \times \textbf{u}^M ) - \textbf{u}^M \cdot \nabla \textbf{u}^M  \big ] , \end{equation}
	Computing the volume integral of \eqref{TruncatedVorticity}, one can eliminate the horizontal derivatives in the viscous term using periodicity, and only $\partial_{x_3}^2 \vorticity$ cannot be eliminated due to the boundary conditions, as for \eqref{Balance_Vort}.  There are then three terms to consider, the buoyancy term, the nonlinear term and the Coriolis term, as follows:
	\[ T_1 =  \langle \nabla \times \mathcal{P}^{M}_{\textbf{u}} \big [ \theta^M \uVecThree \big ] \rangle \text{ , } \hspace{.25 cm} T_2 =  -\langle \nabla \times \mathcal{P}^{M}_{\textbf{u}} \big [ \textbf{u}^M \cdot \nabla \textbf{u}^M \big ] \rangle \text{ , } \hspace{.25 cm} T_3 =  -\langle \nabla \times \mathcal{P}^{M}_{\textbf{u}} \big [ \uVecThree \times \textbf{u}^M \big ] \rangle \text{ . } \]
	In order to transform these terms into those in \eqref{Balance_Vort}, we consider their expression in Fourier space.  First, consider the effect of the curl and projection operators.  For a vector valued function $\textbf{F}$ one can use the definition of the projection operator $\mathcal{P}_{\textbf{u}}^M$, apply the curl operator and compute the spatial integral, and then eliminate all of the terms involving derivatives in $x_1$ due to the periodicity:
	\begin{align} \langle \nabla \times \mathcal{P}^{M}_{\textbf{u}} \big [ \textbf{F} \big ] \rangle & = \sum_{ \textbf{n} \in \mathscr{N}_{\textbf{u}}^M} \big \langle \textbf{F} \cdot \textbf{v}^{\textbf{n}} \big \rangle \big \langle \nabla \times \textbf{v}^{\textbf{n}} \big \rangle \notag \\ \label{VortConsis_CurlProj} & = \sum_{ \textbf{n} \in \mathscr{N}_{\textbf{u}}^M} \big \langle \textbf{F} \cdot \textbf{v}^{\textbf{n}} \big \rangle \big ( \langle  \partial_{x_3} v_1^{\textbf{n}} \rangle \uVecTwo - \langle \partial_{x_3} v_2^{\textbf{n}} \rangle \uVecOne \big ) . \end{align} 
	Furthermore note that only the terms for which $m_1 = 0$ are non-zero, since $\partial_{x_3} v_1^{\textbf{n}}, \partial_{x_3} v_2^{\textbf{n}}$ otherwise have zero mean.  Since $m_1 = 0$, it follows from \eqref{PhaseIndexSets} that $p = 2$ and from the definition \eqref{VectorFieldDef_FreeSlip} one has 
	\[ \textbf{v}^{(0,m_3,2,1)} = \frac{\uVecOne}{\sqrt{2}V} \cos ( m_3 x_3 ) \hspace{.5 cm} \text{ , } \hspace{.5 cm} \textbf{v}^{(0,m_3,2,2)} = \frac{\uVecTwo}{\sqrt{2}V} \cos ( m_3 x_3 ) \text{ . } \]
	One can therefore evaluate the integrals in \eqref{VortConsis_CurlProj}, and one finds that the only non-zero terms in the sum must have $m_3$ odd.  Thus the curl of the projection gives the following sum: 
	\begin{align} \label{VortConsis_Sum} \langle \nabla \times \mathcal{P}^{M}_{\textbf{u}} \big [ \textbf{F} \big ] \rangle = \sum_{\substack{\textbf{n} \in \mathscr{N}_{\textbf{u}}^M \cap \mathscr{N}_{\textbf{u}}^{*} \\  m_3 \text{ odd }}} \frac{4}{ \sqrt{\Sha}} \big \langle \textbf{F} \cdot \textbf{v}^{\textbf{n}}  \big \rangle \big ( \delta^{c,2} \uVecOne - \delta^{c,1} \uVecTwo \big ) .\end{align}
	Thus each of the terms $T_1,T_2,T_3$ are sums over such restricted indices $\textbf{n}$.  The buoyancy term $T_1$ can therefore be eliminated since from \eqref{LinearBasisRelations} one has
	\[ \langle  \theta^M \uVecThree \cdot \textbf{v}^{\textbf{n}} \rangle = \sum_{\tilde{\textbf{n}} \in \mathscr{N}_{\theta}^M} \frac{(-1)^{p_1}\Sha m_1}{|\mathcal{K}\textbf{m}|} \delta^{\textbf{m},\tilde{\textbf{m}}} \delta^{p_1,\tilde{p}_1} \delta^{c,1} \theta^{\tilde{\textbf{n}},M} = 0 \text{ , }  \] 
	since $m_1 = 0$.  By expanding the nonlinear term, integrating by parts and using $m_1 = 0$, one finds 
	\begin{equation} \label{VortConsis_Nonlin} \langle \big [ ( \textbf{u}^M \cdot \nabla ) \textbf{u}^M \big ] \cdot \textbf{v}^{\textbf{n}} \rangle =  - \sum_{\textbf{n}' \in \mathscr{N}_{\textbf{u}}^M } \sum_{\textbf{n}'' \in \mathscr{N}_{\textbf{u}}^M } u^{\textbf{n}',M} u^{\textbf{n}'',M} \big \langle v^{\textbf{n}'}_3 \textbf{v}^{\textbf{n}''} \cdot \partial_{x_3} \textbf{v}^{\textbf{n}} \big  \rangle .  \end{equation}
	Due to the orthogonality of the sinusoids, one sees that the only terms which remain in this sum must satisfy $m_1' = m_1''$ and $m_3 = |m_3' \pm m_3''|$.  Clearly $v_3^{\textbf{n}'}$ must be non-zero, hence $m_1', m_3' > 0$ and $c' = 1$.  Furthermore, since only one component of $\textbf{v}^{\textbf{n}}$ is non-zero, the dot product in \eqref{VortConsis_Nonlin} involves only one component of $\textbf{v}^{\textbf{n}''}$, so the phase $p'$ of $v_3^{\textbf{n}}$ must match the phase of this component, hence one finds $c'' = c$ and $p' + p'' = 1$ mod 2.  One can then insert the definitions of $v_3^{\textbf{n}'}, \textbf{v}^{\textbf{n}''}$ and since the phases match the horizontal integral is easily evaluated.  The vertical integral has fixed phase and can be evaluated as in \ref{app:NonlinearDerivation} via \eqref{NonlinDeriv_EvalResult}.  For $c'' = 1$ one obtains the following:
	\[ \big \langle v^{\textbf{n}'}_3 \textbf{v}^{\textbf{n}''} \cdot \partial_{x_3} \textbf{v}^{\textbf{n}} \big  \rangle = \frac{(-1)^{p'+1}\Sha m_1' m_3''m_3}{2\sqrt{2} |\mathcal{K}\textbf{m}'|\mathcal{K}\textbf{m}''|V} \frac{\delta^{m_3'}\delta^{c',1}\delta^{c'',c}}{\eta^{m_3'}\eta^{m_3''}} s^{(\boldsymbol{\mu}^3,\boldsymbol{\xi}^4)} \text{ , }  \]
	whereas for $c'' = 2$ one obtains
	\[ \big \langle v^{\textbf{n}'}_3 \textbf{v}^{\textbf{n}''} \cdot \partial_{x_3} \textbf{v}^{\textbf{n}} \big  \rangle = \frac{(-1)^{p'+1}\Sha m_1' m_3}{2\sqrt{2} |\mathcal{K}\textbf{m}'|V} \frac{\delta^{m_3'}\delta^{c',1}\delta^{c'',c}}{\eta^{m_3'}\eta^{m_3''}} s^{(\boldsymbol{\mu}^3,\boldsymbol{\xi}^4)}  \text{ . }  \]
	In order to avoid repeated terms in the expression for $T_2$, one can impose $m_3' > m_3''$ by adding the corresponding term with $\textbf{n}',\textbf{n}''$ reversed, and hence by collecting the above results, one obtains the following expression:
	\begin{equation} \label{VortBal_LHS} T_2 = \sum_{\substack{\textbf{n}',\textbf{n}''  \in \mathscr{N}_{\textbf{u}}^M \\  m_3' + m_3'' \text{ odd } \\ m_3' > m_3'' } }  (-1)^{p'} \frac{2}{\pi} \Sha m_1' \delta^{m_1',m_1''}\delta^{p',p''+1} \textbf{a}^{\textbf{n}',\textbf{n}''} u^{\textbf{n}',M} u^{\textbf{n}'',M} \text{ , }  \end{equation}
	in which
	\begin{align} \label{VortBal_LHS2} 
		\textbf{a}^{\textbf{n}',\textbf{n}''} = \sum_{\substack{\textbf{n}\in \mathscr{N}_{\textbf{u}}^M \cap \mathscr{N}_{\textbf{u}}^{*} \\  m_3 = m_3' \pm m_3'' }} m_3 \begin{pmatrix} \delta^{c,2} \big ( \frac{\delta^{c',2}\delta^{c'',1}\delta^{m_3''}}{|\mathcal{K}\textbf{m}''|} s^{(\boldsymbol{\mu}^3,\boldsymbol{\xi}^3)} -  \frac{\delta^{c',1}\delta^{c'',2}}{|\mathcal{K}\textbf{m}'|\eta^{m_3''}} \big ) \\ \frac{\delta^{c,1}\delta^{c',1}\delta^{c'',1}\delta^{m_3''}}{|\mathcal{K}\textbf{m}'||\mathcal{K}\textbf{m}''|}\big ( m_3'' - m_3's^{(\boldsymbol{\mu}^3,\boldsymbol{\xi}^3)} \big ) \\ 0
		\end{pmatrix} \text{ . }
	\end{align}
	On the other hand, for the vorticity balance to hold one must have
	\begin{equation} \label{VortBal_T2Equal} T_2 = \langle ( \vorticity^M \cdot \nabla ) \textbf{u}^M \rangle \text{ , } \end{equation}
	hence we consider the expression of the right hand side in Fourier space.  Expanding the right hand side one obtains:
	\begin{align}
		\langle ( \vorticity^M \cdot \nabla ) \textbf{u}^M \rangle = \sum_{\textbf{n}' \in \mathscr{N}_{\textbf{u}}^M } \sum_{\textbf{n}'' \in \mathscr{N}_{\textbf{u}}^M } u^{\textbf{n}',M} u^{\textbf{n}'',M} \big \langle \big ( (\nabla \times \textbf{v}^{\textbf{n}'} )\cdot \nabla \big ) \textbf{v}^{\textbf{n}''}  \big \rangle \text{ . } \notag
	\end{align}
	For the third component $v_3^{\textbf{n}''}$, one can simply integrate by parts to cancel the various terms from the curl, whereas for the first two components one cannot integrate by parts in the vertical direction without obtaining a boundary term.  Thus one obtains the following:
	\begin{align}
		\big \langle \big ( (\nabla \times \textbf{v}^{\textbf{n}'} )\cdot \nabla \big ) \textbf{v}^{\textbf{n}''} \big  \rangle = \sum_{j = 1,2} \big \langle ( -\partial_{x_3} v_2^{\textbf{n}'} \partial_{x_1} + \partial_{x_1} v_2^{\textbf{n}'} \partial_{x_3}  ) v_j^{\textbf{n}''} \big \rangle  \uVecJ \notag
	\end{align}
	One can then insert the definitions of $\textbf{v}^{\textbf{n}'},\textbf{v}^{\textbf{n}''}$.  Since $v_2^{\textbf{n}'}$ must be non-zero, one has $c' = 2$, the orthogonality of sinusoids implies $m_1' = m_1''$ and since each term involves a derivative in $x_1$ one has $m_1' > 0$.  The phases of the components must match, hence one finds $p' + p'' =1$, and then the horizontal integral is easily evaluated.  By evaluating the vertical integral as well one find $m_3' + m_3''$ must be odd.  In the case $c'' = 1$ one obtains the following:
	\[ \big \langle \big ( (\nabla \times \textbf{v}^{\textbf{n}'} )\cdot \nabla \big ) \textbf{v}^{\textbf{n}''} \big  \rangle = 
	\delta^{c',2} \delta^{m_1',m_1''}\delta^{p',p''+1} \frac{(-1)^{p'}4 \Sha m_1'm_3''}{\pi |\mathcal{K}\textbf{m}''|} \eta^{m_3'}\eta^{m_3''} \delta^{c'',1} \uVecOne  .  \]
	whereas in the case $c'' = 2$ one obtains the following:
	\[ \big \langle \big ( (\nabla \times \textbf{v}^{\textbf{n}'} )\cdot \nabla \big ) \textbf{v}^{\textbf{n}''} \big  \rangle = 
	\delta^{c',2} \delta^{m_1',m_1''}\delta^{p',p''+1} \frac{(-1)^{p'}4 \Sha m_1'm_3''}{\pi |\mathcal{K}\textbf{m}''|} \eta^{m_3'}\eta^{m_3''} (\delta^{c'',1} \uVecOne + \delta^{c'',2} \uVecTwo ) .  \]
	In order to avoid repeated terms one can again impose $m_3' > m_3''$ by adding the corresponding term with $\textbf{n}',\textbf{n}''$ reversed, and hence by collecting the above results, one obtains the following expression:
	\begin{equation} \label{VortBal_RHS} \langle ( \vorticity^M \cdot \nabla ) \textbf{u}^M \rangle = \sum_{\substack{\textbf{n}',\textbf{n}''  \in \mathscr{N}_{\textbf{u}}^M \\ m_3' > m_3''} }  (-1)^{p'} \frac{4}{\pi} \Sha m_1' \delta^{m_1',m_1''}\delta^{p',p''+1}  \textbf{b}^{\textbf{n}',\textbf{n}''} u^{\textbf{n}',M} u^{\textbf{n}'',M} \text{ , }  \end{equation}
	in which
	\begin{align} \label{VortBal_RHS2} 
		\textbf{b}^{\textbf{n}',\textbf{n}''} =  \eta^{m_3''} \big ( \frac{m_3'' \delta^{c'',1}\delta^{c',2}}{|\mathcal{K}\textbf{m}''|} - \frac{m_3'\delta^{c',1}\delta^{c'',2}}{|\mathcal{K}\textbf{m}'|} \big ) \uVecOne \text{ . }
	\end{align}
	Hence in order for the vorticity balance to hold one must have
	\[ \textbf{a}^{\textbf{n}',\textbf{n}''} = 2 \textbf{b}^{\textbf{n}',\textbf{n}''} \text{ , } \]
	for all $\textbf{n}',\textbf{n}'' \in \mathscr{N}_{\textbf{u}}^{M}$ such that $m_1' = m_1'' > 0$, $p' = p''+1$, $m_3' + m_3''$ odd.  There are three cases to consider:
	\begin{enumerate}
		\item If $c' = c'' = 2$, it is immediate from the definitions \eqref{VortBal_LHS2}, \eqref{VortBal_RHS2} that both $\textbf{a}^{\textbf{n}',\textbf{n}''}$ and $\textbf{b}^{\textbf{n}',\textbf{n}''}$ are zero.  
		
		\item If $c' = c'' = 1$, note that the first component of both $\textbf{a}^{\textbf{n}',\textbf{n}''}$ and $\textbf{b}^{\textbf{n}',\textbf{n}''}$ are zero, so only the second component needs to be considered.  The second component of $\textbf{b}^{\textbf{n}',\textbf{n}''}$ is explicitly zero.  The second component of $\textbf{a}^{\textbf{n}',\textbf{n}''}$ is zero if neither $(0,m_3'+m_3'',2,1)$ and $(0,|m_3'-m_3''|,2,1)$ are included in $\mathscr{N}_{\textbf{u}}^{M}$.  On the other hand, if both $(0,m_3'+m_3'',2,1)$ and $(0,|m_3'-m_3''|,2,1)$ are included in $\mathscr{N}_{\textbf{u}}^{M}$, then one has
		\[ \sum_{\substack{\textbf{n}\in \mathscr{N}_{\textbf{u}}^M \cap \mathscr{N}_{\textbf{u}}^{*} \\  m_3 = m_3' \pm m_3'' }} m_3 \big ( m_3'' - m_3's^{(\boldsymbol{\mu}^3,\boldsymbol{\xi}^3)} \big ) = ( m_3' + m_3'' )( m_3'' - m_3' ) +  ( m_3' - m_3'')( m_3'' + m_3' ) = 0 \text{ . } \]
		On the other hand, if only one of $(0,m_3'+m_3'',2,1)$ or $(0,|m_3'-m_3''|,2,1)$ is included in $\mathscr{N}_{\textbf{u}}^{M}$, then the above cancellation does not occur.
		
		\item If either $(c',c'')=(1,2)$ or $(c',c'')=(2,1)$ then the second component of both $\textbf{a}^{\textbf{n}',\textbf{n}''}$ and $\textbf{b}^{\textbf{n}',\textbf{n}''}$ are zero, hence one only need consider the first.  Here it is clear that both $(0,m_3'+m_3'',2,2)$ and $(0,|m_3'-m_3''|,2,2)$ must be included in $\mathscr{N}_{\textbf{u}}^{M}$, since one has
		\begin{align}
			\sum_{\substack{\textbf{n}\in \mathscr{N}_{\textbf{u}}^M \cap \mathscr{N}_{\textbf{u}}^{*} \\  m_3 = m_3' \pm m_3'' }} m_3 \delta^{m_3''} s^{(\boldsymbol{\mu}^3,\boldsymbol{\xi}^3)} & = ( m_3' + m_3'' )-  ( m_3' - m_3'') = 2m_3'' \text{ , } \notag \\ \sum_{\substack{\textbf{n}\in \mathscr{N}_{\textbf{u}}^M \cap \mathscr{N}_{\textbf{u}}^{*} \\  m_3 = m_3' \pm m_3'' }} \frac{m_3}{ \eta^{m_3''}} & = \left \{ \begin{array}{cc} ( m_3' + m_3'' )+  ( m_3' - m_3'') & m_3'' > 0 \\ \sqrt{2} m_3' & m_3'' = 0 \end{array} \right . = 2\eta^{m_3''} m_3' \text{ . } \notag 
		\end{align}
	\end{enumerate}
	Thus Criterion \ref{Crit:VortCrit} (i) is exactly what is required for \eqref{VortBal_T2Equal} to hold.
	
	Finally, we consider $T_3$.  For the balance to hold one requires 
	\begin{equation} \label{VortBal_T3Equal} -\Pra \Rot \langle \nabla \times \mathbb{P}_{\textbf{u}}^{M} \big [ \uVecThree \times \textbf{u}^M \big ] \rangle = \Pra \Rot \langle \partial_{x_3} \textbf{u}^M \rangle \text{ , } \end{equation}
	which is trivial for $\Rot = 0$, hence we consider $\Rot \neq 0$.  Considering the left hand side of \eqref{VortBal_T3Equal}, one recalls from \eqref{LinearBasisRelations} that 
	\[ \big \langle (\uVecThree \times \textbf{v}^{\tilde{\textbf{n}}} )\cdot \textbf{v}^{\textbf{n}}  \big \rangle = \frac{(-1)^{c}m_3}{|\mathcal{K}\textbf{m}|} \delta^{\textbf{m},\tilde{\textbf{m}}} \delta^{p_1,\tilde{p}_1} \delta^{c+1,\tilde{c}} , \]
	hence the only non-zero terms in the sum \eqref{VortConsis_Sum} must satisfy $\textbf{m} = (0,m_3) = \tilde{\textbf{m}}$, $m_3$ odd, and $p = \tilde{p}$.  Therefore $T_3$ has the form
	\[ T_3 = - \sum_{\substack{\textbf{n} \in \mathscr{N}^M_{\textbf{u}}  \cap \mathscr{N}^*_{\textbf{u}} \\ m_3 \text{ odd }}} \sum_{\substack{\tilde{\textbf{n}} \in \mathscr{N}^M_{\textbf{u}} \\ \tilde{\textbf{n}} = (0,m_3,p,c+1)}} \frac{4}{\sqrt{\Sha}} u^{\tilde{\textbf{n}}} \big ( \delta^{c,2} \uVecOne + \delta^{c,1} \uVecTwo \big ) . \]
	On the other hand, using the periodicity in $x_1$ and evaluating the vertical integral, one finds
	\begin{align}
		\langle \partial_{x_3} \textbf{u}^M \rangle = \sum_{\textbf{n} \in \mathscr{N}^M_{\textbf{u}}} u^n \langle \partial_{x_3} \textbf{v}^{\textbf{n}} \rangle \notag = - \sum_{\substack{\textbf{n} \in \mathscr{N}^M_{\textbf{u}}  \cap \mathscr{N}^*_{\textbf{u}} \\ m_3 \text{ odd }}} \frac{4}{\sqrt{\Sha}} u^{\textbf{n}} \big ( \delta^{c,1} \uVecOne + \delta^{c,2} \uVecTwo \big ) . 
	\end{align}
	Comparing these two expressions one sees that, for \eqref{VortBal_T3Equal} to hold, the mode $u^{(0,m_3,p,c)}$ can be included in $\mathscr{N}^M_{\textbf{u}}$ if and only if the mode $u^{(0,m_3,p,c)}$ is included as well, as in Criterion \ref{Crit:VortCrit} (ii).
\end{proof}

\section{Well posedness and regularity proofs}

\subsection{Additional proofs in Theorem \ref{thm:PDE_WellPosedness_Balances}}

\subsubsection{Further regularity}

\label{app:SmoothSolns}

\begin{proof}
	We have already obtained uniform bounds for $\{ \textbf{X}^{M} \}_{M \geq 0}$ in $L^{\infty}((0,\tau);\textbf{H}^{0})$ and $L^2((0,\tau);\textbf{H}^{1})$, so one can assume that for some $k \geq 0$ one has that for $0 \leq \tilde{k} \leq k$ the sequence $\{ t^{\frac{\tilde{k}}{2}} \textbf{X}^{M} \}_{M \geq 0}$ is uniformly bounded in $L^{\infty}((0,\tau);\textbf{H}^{\tilde{k}}) \cap L^2((0,\tau);\textbf{H}^{\tilde{k}+1})$.  We show that this implies the sequence $\{ t^{\frac{k+1}{2}} \textbf{X}^{M} \}_{M \geq 0}$ is uniformly bounded in $L^{\infty}((0,\tau);\textbf{H}^{k+1}) \cap L^2((0,\tau);\textbf{H}^{k+2})$.  Let $\boldsymbol{\alpha} \in \mathbb{Z}^2_{\geq 0}$ be such that $|\boldsymbol{\alpha}| = k+1$.  As above the strategy is to show that for some non-negative $A(t) \in L^{\infty}((0,\tau))$, $B(t) \in L^{1}((0,\tau))$ and $C > 0$ one has the bound
	\begin{equation} \label{AppSmooth_EnergyIneq} t^{k+1} \|\partial_x^{\boldsymbol{\alpha}} \textbf{X}^M (t) \|_{\textbf{H}^0}^2 +  C \int_0^t s^{k+1} \| \partial_x^{\boldsymbol{\alpha}} \textbf{X}^M (s) \|_{\textbf{H}^1}^2 ds \leq  A(t) + \int_0^t B(s) s^{k+1} \|\partial_x^{\boldsymbol{\alpha}} \textbf{X}^M(s) \|_{\textbf{H}^0}^2 ds \text{ , }  \end{equation}
	and hence Gronwall's inequality gives the $L^{\infty}((0,\tau);\textbf{H}^{k+1})$ bound.  Inserting the $L^{\infty}((0,\tau);\textbf{H}^{k+1})$ bound into \eqref{AppSmooth_EnergyIneq} gives the $L^{2}((0,\tau);\textbf{H}^{k+2})$ bound.  By differentiating \eqref{TruncatedEvolutionPDE} and computing an inner product with $t^{k+1} \partial_x^{\boldsymbol{\alpha}} \textbf{X}^M$ one obtains
	\begin{align} \frac{1}{2} t^{k+1} \frac{d}{dt} \big [ \|\partial_x^{\boldsymbol{\alpha}} \textbf{X}^M \|_{\textbf{H}^0}^2 \big ] = t^{k+1} & \big \langle \partial_x^{\boldsymbol{\alpha}} \textbf{X}^M \cdot \mathcal{L}^{\textbf{0}} \partial_x^{\boldsymbol{\alpha}} \textbf{X}^M \big \rangle \notag \\ & - \sum_{\boldsymbol{\alpha}_1 + \boldsymbol{\alpha}_2 = \boldsymbol{\alpha}} \binom{\boldsymbol{\alpha}}{\boldsymbol{\alpha}_1} t^{k+1} \Big \langle \partial_x^{\boldsymbol{\alpha}} \textbf{X}^{M} \cdot \Big ( \big ( \partial_x^{\boldsymbol{\alpha}_1} \textbf{u}^M \cdot \nabla \big ) \partial_x^{\boldsymbol{\alpha}_2} \textbf{X}^{M} \Big ) \Big \rangle \text{ , } \notag \end{align}
	where $\mathcal{L}^{\textbf{0}}$ is the linearization of \eqref{EvolutionPDE} about the origin.  Using Young's inequality as in \eqref{EnergyInequality} one has 
	\begin{equation} \label{AppSmooth_LinearOperatorBound} \big \langle \partial_x^{\boldsymbol{\alpha}} \textbf{X}^M \cdot \mathcal{L}^{\textbf{0}} \partial_x^{\boldsymbol{\alpha}} \textbf{X}^M \big \rangle \leq \frac{1}{2} (\Pra \Ray + 1 ) \| \partial_x^{\boldsymbol{\alpha}} \textbf{X}^M \|_{\textbf{H}^0}^2 - \mathsf{min}(\Pra ,1) \| \partial_x^{\boldsymbol{\alpha}} \textbf{X}^M \|_{\textbf{H}^1}^2  \end{equation} 
	and by also adding the term $\frac{k+1}{2} t^{k} \|\partial_x^{\boldsymbol{\alpha}} \textbf{X}^M \|_{\textbf{H}^0}^2$ to both sides and integrating in time, one obtains
	\begin{align} t^{k+1} & \|\partial_x^{\boldsymbol{\alpha}} \textbf{X}^M (t) \|_{\textbf{H}^0}^2  +  2 \hspace{1 mm} \mathsf{min}(\Pra ,1) \int_0^t s^{k+1} \| \partial_x^{\boldsymbol{\alpha}} \textbf{X}^M (s) \|_{\textbf{H}^1}^2 ds \notag \\ \label{AppSmooth_EnergyIneq2} & \leq  \int_0^t \big ( (k+1)s^{k} + (\Pra \Ray +1)s^{k+1} \big )  \|\partial_x^{\boldsymbol{\alpha}} \textbf{X}^M(s) \|_{\textbf{H}^0}^2 ds + \sum_{\boldsymbol{\alpha}_1 + \boldsymbol{\alpha}_2 = \boldsymbol{\alpha}} 2 \binom{\boldsymbol{\alpha}}{\boldsymbol{\alpha}_1} \textbf{N}^{\boldsymbol{\alpha}_1 , \boldsymbol{\alpha}_2 } \text{ , } \end{align}
	where
	\[ \textbf{N}^{\boldsymbol{\alpha}_1 , \boldsymbol{\alpha}_2 } = \int_0^t s^{k+1} \Big | \Big \langle \partial_x^{\boldsymbol{\alpha}} \textbf{X}^{M}(s) \cdot \Big ( \big ( \partial_x^{\boldsymbol{\alpha}_1} \textbf{u}^M(s) \cdot \nabla \big ) \partial_x^{\boldsymbol{\alpha}_2} \textbf{X}^{M}(s) \Big ) \Big \rangle \Big | ds  \text{ . } \]
	Due to the inductive assumption, one can include the first terms on the right hand side of \eqref{AppSmooth_EnergyIneq2} in the functions $A(t)$ and $B(t)$ respectively, so one need only consider the terms in the sum.  Note that when $\boldsymbol{\alpha}_1 = \textbf{0}$ one has
	\[ \Big \langle \partial_x^{\boldsymbol{\alpha}} \textbf{X}^{M}(s) \cdot \Big ( \big ( \textbf{u}^M(s) \cdot \nabla \big ) \partial_x^{\boldsymbol{\alpha}} \textbf{X}^{M}(s) \Big ) \Big \rangle  = \frac{1}{2} \Big \langle \textbf{u}^M(s) \cdot \nabla  |\partial_x^{\boldsymbol{\alpha}} \textbf{X}^{M}(s)|^2 \Big \rangle = 0 \text{ , } \]
	hence one need only bound the terms for which $|\boldsymbol{\alpha}_1| > 0$.  For such $\boldsymbol{\alpha}_1$, one can integrate by parts, apply Cauchy-Schwarz twice and then apply \eqref{GagliardoNirenburgSobolevIneq} to obtain the bound
	\begin{align} & \textbf{N}^{\boldsymbol{\alpha}_1 , \boldsymbol{\alpha}_2 } = \int_0^t s^{k+1} \Big | \Big \langle \partial_x^{\boldsymbol{\alpha}_2} \textbf{X}^{M}(s) \cdot \Big ( \big ( \partial_x^{\boldsymbol{\alpha}_1} \textbf{u}^M(s) \cdot \nabla \big ) \partial_x^{\boldsymbol{\alpha}} \textbf{X}^{M}(s) \Big ) \Big \rangle \Big | ds \notag \\ & \leq \int_0^t s^{k+1} \| \partial_x^{\boldsymbol{\alpha}_2} \textbf{X}^{M}(s) \|_{\textbf{L}^2\times L^4} \| \partial_x^{\boldsymbol{\alpha}_1} \textbf{X}^M(s) \|_{\textbf{L}^4 \times L^4} \| \partial_x^{\boldsymbol{\alpha}} \textbf{X}^{M}(s) \|_{\textbf{H}^1} ds \notag  \\ & \leq \int_0^t s^{k+1} \| \partial_x^{\boldsymbol{\alpha}_2} \textbf{X}^{M}(s) \|_{\textbf{H}^0}^{\frac{1}{2}} \| \partial_x^{\boldsymbol{\alpha}_2} \textbf{X}^{M}(s) \|_{\textbf{H}^1}^{\frac{1}{2}} \| \partial_x^{\boldsymbol{\alpha}_1} \textbf{X}^{M}(s) \|_{\textbf{H}^0}^{\frac{1}{2}} \| \partial_x^{\boldsymbol{\alpha}_1} \textbf{X}^{M}(s) \|_{\textbf{H}^1}^{\frac{1}{2}} \|\partial_x^{\boldsymbol{\alpha}} \textbf{X}^{M}(s) \|_{\textbf{H}^1} ds \text{ . } \notag  \end{align}
	For a sufficiently small $\varepsilon > 0$ (yet to be chosen), one then obtains the following from Young's inequality
	\begin{align} \textbf{N}^{\boldsymbol{\alpha}_1 , \boldsymbol{\alpha}_2 } \leq \frac{1}{\varepsilon} & \int_0^t s^{k+1} \| \partial_x^{\boldsymbol{\alpha}_2} \textbf{X}^{M}(s) \|_{\textbf{H}^0} \| \partial_x^{\boldsymbol{\alpha}_2} \textbf{X}^{M}(s) \|_{\textbf{H}^1} \| \partial_x^{\boldsymbol{\alpha}_1} \textbf{X}^{M}(s) \|_{\textbf{H}^0} \| \partial_x^{\boldsymbol{\alpha}_1} \textbf{X}^{M}(s) \|_{\textbf{H}^1} ds \notag \\ \label{AppSmooth_Bound1} & + \varepsilon \int_0^t s^{k+1} \|\partial_x^{\boldsymbol{\alpha}} \textbf{X}^{M}(s) \|_{\textbf{H}^1}^2 ds \text{ . }  \end{align}
	The first term on the right hand side of \eqref{AppSmooth_Bound1} is then dealt with in one of two possible ways, depending on $\boldsymbol{\alpha}_1,\boldsymbol{\alpha}_2$:
	\begin{enumerate}
		\item \underline{\textbf{Case 1:} $0<|\boldsymbol{\alpha}_1|_1 \leq k$ and $0 < |\boldsymbol{\alpha}_2|_1 \leq k$ (if $k > 0$)}
		
		In this case one can directly apply the $L^{\infty}((0,\tau);\textbf{H}^{\tilde{k}})$ and $L^2((0,\tau);\textbf{H}^{\tilde{k}+1})$ bounds from the induction hypothesis to the $\textbf{H}^0$ and $\textbf{H}^1$ terms respectively, hence here the first term on the right hand side of \eqref{AppSmooth_Bound1} is included in $A(t)$.
		
		\item \underline{\textbf{Case 2:} $\boldsymbol{\alpha}_1 = \boldsymbol{\alpha}$}
		
		Here one must instead apply Young's inequality again to obtain the bound
		\begin{align}  N^{\boldsymbol{\alpha}_1 , \boldsymbol{\alpha}_2 } \leq \frac{1}{\varepsilon^3} & \int_0^t s^{k+1} \| \textbf{X}^{M}(s) \|_{\textbf{H}^0}^2 \| \textbf{X}^{M}(s) \|_{\textbf{H}^1}^2 \| \partial_x^{\boldsymbol{\alpha}} \textbf{X}^{M}(s) \|_{\textbf{H}^0}^2 ds \notag \\ \label{AppSmooth_Bound2} & + 2 \varepsilon \int_0^t s^{k+1} \|\partial_x^{\boldsymbol{\alpha}} \textbf{X}^{M}(s) \|_{\textbf{H}^1}^2 ds \text{ , } \end{align}
		Due to the $L^{\infty}((0,\tau);\textbf{H}^0)$ and $L^{2}((0,\tau);\textbf{H}^1)$ bounds one has a uniformly bound on  
		\[ \int_0^t \| \textbf{X}^{M}(s) \|_{\textbf{H}^0}^2 \| \textbf{X}^{M}(s) \|_{\textbf{H}^1}^2 ds \text{ , } \]
		hence the first term on the right hand side of \eqref{AppSmooth_Bound2} can be included in $B(t)$.
		
	\end{enumerate}
	Finally, by choosing $\varepsilon$ such that
	\[ \sum_{\boldsymbol{\alpha}_1 + \boldsymbol{\alpha}_2 = \boldsymbol{\alpha}} 2 \varepsilon \binom{\boldsymbol{\alpha}}{\boldsymbol{\alpha}_1} = 2^{|\boldsymbol{\alpha}|_1+1} \varepsilon < \frac{\mathsf{min} (\Pra,1)}{2} \text{ , } \]
	the sum over the second terms in the right hand sides of \eqref{AppSmooth_Bound1}, \eqref{AppSmooth_Bound2} involving $\|\partial_x^{\boldsymbol{\alpha}} \textbf{X}^{M}(s) \|_{\textbf{H}^1}$ can be absorbed in the corresponding term on the left hand side of $\eqref{AppSmooth_EnergyIneq2}$, and hence $\eqref{AppSmooth_EnergyIneq}$ holds with $C = \mathsf{min}(\Pra,1)/2$.
	
\end{proof}

\subsubsection{Lipschitz dependence on initial conditions and maximum principle}

\label{app:LipDepOnIC_MaxPrinc}

\begin{proof}
	First we prove the Lipschitz dependence on initial conditions in \eqref{Thm1_LipDepInitCond}.  Suppose that $\textbf{X}(t)$, $\textbf{X}^{*}(t)$ are two different sub-sequential limits with initial conditions $\textbf{X}_0,\textbf{X}_0^*$ (not necessarily distinct).  The difference $\tilde{\textbf{X}} = \textbf{X} - \textbf{X}^{*}$ must satisfy
	\[ \frac{1}{2} \frac{d}{dt} \| \tilde{\textbf{X}}  \|_{\textbf{H}^0}^2 + \Pra \|  \tilde{\textbf{u}}  \|_{\textbf{H}^1}^2 + \| \tilde{\theta} \|_{\textbf{H}^1}^2 =  \text{ } \big ( \Pra \Ray + 1 \big ) \big \langle \tilde{u}_3 \tilde{\theta} \big \rangle - \big \langle \tilde{\textbf{u}} \cdot \big ( \textbf{u} \cdot \nabla \textbf{u} - \textbf{u}^{*} \cdot \nabla \textbf{u}^{*} \big ) \big \rangle - \big \langle \tilde{\theta} \big ( \textbf{u} \cdot \nabla \theta - \textbf{u}^{*}  \cdot \nabla \theta^{*}  \big ) \big \rangle .  \]
	Integrating over $[0,t]$, one can use
	\begin{equation} \label{ThmExist_NonlinearDiffTrick} \textbf{u} \cdot \nabla \textbf{u} - \textbf{u}^{*} \cdot \nabla \textbf{u}^{*} = \tilde{\textbf{u}} \cdot \nabla \textbf{u} + \textbf{u}^{*} \cdot \nabla \tilde{\textbf{u}} \hspace{.25 cm} \text{ , } \hspace{.25 cm} \textbf{u} \cdot \nabla \theta - \textbf{u}^{*} \cdot \nabla \theta^{*} = \tilde{\textbf{u}} \cdot \nabla \theta + \textbf{u}^{*} \cdot \nabla \tilde{\theta} , \end{equation}
	and together with the Cauchy-Schwartz inequality and a bound similar to \eqref{NonlinearTermBound} one obtains
	\begin{align} \label{Thm1_DifferenceIneq} \| & \tilde{\textbf{X}}(t) \|_{\textbf{H}^0}^2 + 2 \hspace{.5 mm} \mathsf{min}(\Pra,1) \int_0^t \| \tilde{\textbf{X}}(s) \|_{\textbf{H}^1}^2 ds \\ & \leq \| \tilde{\textbf{X}}_0 \|_{\textbf{H}^0}^2 + \int_0^t \Big ( (\Pra \Ray +1) \| \tilde{\textbf{X}}(s) \|_{\textbf{H}^0}^2 + C \| \textbf{X}(s) \|_{\textbf{H}^1} \| \tilde{\textbf{X}}(s) \|_{\textbf{H}^1} \| \tilde{\textbf{X}}(s) \|_{\textbf{H}^0} \Big )  ds \text{ , } \notag \end{align}
	for a.e. $t \in [0,\tau]$.  One can then use Young's inequality on this last term as follows
	\begin{align}  C \int_0^t & \| \textbf{X}(s) \|_{\textbf{H}^1} \| \tilde{\textbf{X}}(s) \|_{\textbf{H}^1} \| \tilde{\textbf{X}}(s) \|_{\textbf{H}^0} ds \notag \\ & \leq \frac{C^2}{8 \hspace{.5 mm} \mathsf{\min}(\Pra,1)} \int_0^t \| \textbf{X}(s) \|_{\textbf{H}^1}^2 \| \tilde{\textbf{X}}(s) \|_{\textbf{H}^0}^2 ds + 2 \hspace{.5 mm} \mathsf{\min}(\Pra,1) \int_0^t  \| \tilde{\textbf{X}}(s) \|_{\textbf{H}^1}^2 ds \text{ . } \notag \end{align}
	This last term can be cancelled with the identical term on the left hand side of \eqref{Thm1_DifferenceIneq}, and then \eqref{Thm1_DifferenceIneq} together with Gronwall's inequality and \eqref{Thm1_Bound2} gives \eqref{Thm1_LipDepInitCond}.  It follows that if $\tilde{\textbf{X}}_0 = \textbf{0}$ then $\tilde{\textbf{X}}(t) = \textbf{0}$, hence the limit is unique.
	
	The maximum principle in Theorem \ref{thm:PDE_WellPosedness_Balances}, part (b) is the analogue of Temam's Lemma 3.2 in Chapter 3 of \cite{temam_InfDimDynSys}, and the proof therein carries through identically for the present problem.  To briefly recall, the proof is given in terms of the temperature $T$ in order to eliminate the additional $u_3$ term, since part (a) provides existence for $T$ with the same regularity as $\theta$, although with different boundary conditions.  Set $\tilde{T}_+ = \max(T - 1,0)$, $\tilde{T}_- = \max(-T,0)$.  By testing the temperature in \eqref{BoussinesqCoriolis} with $\tilde{T}_+$, one obtains
	\[ \frac{1}{2} \frac{d}{dt} \| \tilde{T}_+ \|_{L^2}^2 + \| \tilde{T}_+ \|_{H^1}^2 = 0  , \]
	and by the Poincar\'e inequality and Gronwall's inequality one has
	\[ \| \tilde{T}_+(t) \|_{L^2} \leq \| \tilde{T}_+(0) \|_{L^2} e^{-t} , \]
	for a constant $C$.  The same reasoning applies to $\tilde{T}_-$, and taking $\tilde{T} = \tilde{T}_+ - \tilde{T}_-$ is sufficient to prove (b).
	
	\vspace{-.75 cm}
	\[ \textcolor{white!100}{ . } \]
\end{proof}

\subsection{Further regularity of the semi-group}

\label{app:FurtherRegularity}

\begin{proof}[Proof of Prop. \ref{prop:FurtherRegularity}]
	For $\textbf{Y}_0 \in \textbf{H}^0$, The existence of the solution $\textbf{Y}(t)$ of the linearized equation \eqref{LinearizedEquation} can obtained by essentially the same Galerkin argument as for the existence of $\textbf{X}(t)$ in Theorem \ref{thm:PDE_WellPosedness_Balances}, treating now $\textbf{X}(t)$ as a forcing term.  In so doing, one obtains constants $C^{\tau,\textbf{X}_0}_1,C^{\tau,\textbf{X}_0}_2$ and the following bounds analogous to \eqref{Thm1_Bound1}, \eqref{Thm1_Bound2}:
	\begin{equation} \label{AppFrechet_Bound1} \sup_{0\leq t \leq \tau} \| \textbf{Y} (t) \|_{\textbf{H}^0}^2 \leq C^{\tau,\textbf{X}_0}_1 \| \textbf{Y}_0 \|_{\textbf{H}^0}^2 \hspace{.5 cm} \text{ , } \hspace{.5 cm} \int_0^t \| \textbf{Y} (s) \|_{\textbf{H}^1}^2 ds \leq C^{\tau,\textbf{X}_0}_2 \| \textbf{Y}_0 \|_{\textbf{H}^0}^2 . \end{equation}
	This is left to the reader to verify.  To verify the uniform Frechet differentiability for $\textbf{X}_0 \in \textbf{H}^0$, let $\textbf{X}_0' = \textbf{X}_0 + \textbf{Y}_0$ and let $\textbf{X}(t),\textbf{X}'(t)$ be the solutions of \eqref{EvolutionPDE} corresponding to $\textbf{X}_0,\textbf{X}_0'$.  Define
	\[ \textbf{Z}(t) = \begin{pmatrix}
		\textbf{u}^{\textbf{Z}} \\ \theta^{\textbf{Z}}
	\end{pmatrix} = \textbf{X}'(t) - \textbf{X}(t) - \textbf{Y}(t) \text{ . } \]
	This solves the following equation
	\begin{align} \frac{d}{dt} \textbf{Z}(t) & = \mathcal{L}^{\textbf{0}} \textbf{Z}(t) - (\textbf{u}' \cdot \nabla ) \textbf{X}'(t) + (\textbf{u} \cdot \nabla ) \textbf{X}(t) + (\textbf{u} \cdot \nabla ) \textbf{Y}(t) + (\textbf{g} \cdot \nabla ) \textbf{X}(t)  \notag \\
		& = \mathcal{L}^{\textbf{0}} \textbf{Z}(t) - (\textbf{u}^{\textbf{Z}} \cdot \nabla ) \textbf{X}'(t) - \big ( (\textbf{u} + \textbf{g} ) \cdot \nabla ) \textbf{Z}(t) - (\textbf{g} \cdot \nabla ) \textbf{Y}(t) \text{ . } \notag \end{align}
	Computing inner products and using the analogue of \eqref{AppSmooth_LinearOperatorBound} one arrives at
	\begin{align} \frac{1}{2} \frac{d}{dt} \|\textbf{Z}(t) \|_{\textbf{H}^0}^2 + \mathsf{min}(\Pra,1) \big \| \textbf{Z}(t) \big \|_{\textbf{H}^1}^2  \leq  & \frac{1}{2} (\Pra \Ray + 1 ) \| \textbf{Z}(t) \|_{\textbf{H}^0}^2 \notag \\ & + \big | \big \langle \textbf{Z} \cdot \big ( (\textbf{u}^{\textbf{Z}} \cdot \nabla ) \textbf{X}'(t) \big ) \big \rangle \big | + \big | \big \langle \textbf{Z} \cdot \big ( (\textbf{g} \cdot \nabla ) \textbf{Y}(t) \big ) \big \rangle \big | \text{ . } \notag \end{align}
	Integrating in time, using the fact $\textbf{Z}(0) = 0$, integrating by parts and applying Cauchy-Schwarz gives
	\begin{align} \|\textbf{Z}(t) \|_{\textbf{H}^0}^2 + 2 \hspace{1 mm} \mathsf{min}(\Pra,1) \int_0^t \big \| \textbf{Z}(s) \big \|_{\textbf{H}^1}^2 ds \leq  \int_0^t \Big [ & (\Pra \Ray + 1 ) \| \textbf{Z}(s) \|_{\textbf{H}^0}^2 + 2 \| \textbf{Z}(s) \|_{\textbf{L}^4\times L^4}^2 \| \textbf{X}'(s) \|_{\textbf{H}^1} \notag \\ & + 2\| \textbf{Z}(s) \|_{\textbf{H}^1} \| \textbf{Y}(s) \|_{\textbf{L}^4 \times L^4}^2 \Big ] ds \text{ . } \notag  \end{align}
	Using \eqref{GagliardoNirenburgSobolevIneq} and Young's inequality these last two terms can be bounded as follows 
	\begin{align}
		\int_0^t \| \textbf{Z}(s) \|_{\textbf{L}^4\times L^4}^2 \| \textbf{X}'(s) \|_{\textbf{H}^1} ds & \leq \frac{\mathsf{min}(\Pra , 1)}{4} \int_0^t \| \textbf{Z}(s) \|_{\textbf{H}^1}^2 ds  \notag \\ & \hspace{1 cm}  + \frac{1}{\min(\Pra,1)} \int_0^t \| \textbf{X}'(s) \|_{\textbf{H}^1}^2 \| \textbf{Z}(s) \|_{\textbf{H}^0}^2 ds \text{ , } \notag \\ 
		\int_0^t \| \textbf{Z}(s) \|_{\textbf{H}^1} \| \textbf{Y}(s) \|_{\textbf{L}^4 \times L^4}^2 ds & \leq  \frac{\mathsf{min}(\Pra , 1)}{4} \int_0^t \| \textbf{Z}(s) \|_{\textbf{H}^1}^2 ds \notag \\ & \hspace{1 cm}  + \frac{1}{\min(\Pra,1)} \int_0^t \| \textbf{Y}(s) \|_{\textbf{H}^1}^2 \| \textbf{Y}(s) \|_{\textbf{H}^0}^2 ds \text{ . } \notag 
	\end{align}
	Subtracting the integral of the $\textbf{H}^1$ norm from both sides and using \eqref{AppFrechet_Bound1}, one obtains the bound
	\begin{align} \|\textbf{Z}(t) \|_{\textbf{H}^0}^2 + \mathsf{min}(\Pra,1) \int_0^t \big \| \textbf{Z}(s) \big \|_{\textbf{H}^1}^2 ds \leq & \frac{2C^{\tau,\textbf{X}_0}_1 C^{\tau,\textbf{X}_0}_2 \| \textbf{Y}_0 \|_{\textbf{H}^0}^2 }{\min(\Pra,1)} \| \textbf{Y}_0 \|_{\textbf{H}^0}^4 \notag \\ & +  \int_0^t \Big [ (\Pra \Ray + 1 ) + \frac{2\| \textbf{X}'(s) \|_{\textbf{H}^1}^2}{\min(\Pra,1)} \Big ] \| \textbf{Z}(s) \|_{\textbf{H}^0}^2 \text{ . } \notag \end{align} 
	Finally, using Gronwall's lemma one obtains
	\[ \sup_{0 \leq t \leq \tau} \frac{\|\textbf{Z}(t) \|_{\textbf{H}^0}^2}{\| \textbf{Y}_0 \|_{\textbf{H}^0}^2} \leq C_3^{\tau,\textbf{X}_0} \| \textbf{Y}_0 \|_{\textbf{H}^0}^2 \text{ , } \]
	for some constant $C_3^{\tau,\textbf{X}_0}$.
\end{proof}

\section{Proofs regarding the local bifurcations at the origin}

\subsection{Proof of Lemma \ref{lem:OriginEigenvalues}}

\label{app:OriginEigs}

\begin{proof}
	
	For $\textbf{m} \in \mathbb{Z}^2_{\geq 1}$, the eigenvalues of $\hat{\mathcal{L}}^{\textbf{0},\textbf{m}}$, denoted $\lambda^{\textbf{m},j}$, $j=1,2,3$, must solve the following characteristic equation:
	\begin{align} 
		(\lambda^{\textbf{m}})^3 & + c_2 (\lambda^{\textbf{m}})^2 + c_1 (\lambda^{\textbf{m}})^2 + c_0  = 0 \text{ , } \notag \\
		\label{OriginCharEq} c_2 & = (2\Pra +1)|\mathcal{K}\textbf{m}|^2 \text{ , }  \\
		c_1 & = \Pra \Big ( (\Pra + 2) |\mathcal{K}\textbf{m}|^4 + \Pra \Rot^2 \frac{m_3^2}{|\mathcal{K}\textbf{m}|^2} - \Ray \frac{\Sha^2m_1^2}{|\mathcal{K}\textbf{m}|^2} \Big ) \text{ , } \notag \\ 
		c_0 & = \Pra^2 \big ( |\mathcal{K}\textbf{m}|^6 + \Rot^2 m_3^2 - \Ray \Sha^2m_1^2 \big ) \text{ . } \notag
	\end{align} 
	For convenience, define the $\lambda^{\textbf{m},j}$ to be descending by lexicographic order, so that $\lambda^{\textbf{m},1}$ has largest real part, or with largest real part and positive imaginary part if this is ambiguous.  Note that the sum $\sum_j \lambda^{\textbf{m},j}$ must equal $-(2\Pra+1)|\mathcal{K}\textbf{m}|^2$, hence for all values of $\parameters$ at least one of the eigenvalues must have negative real part.  According to the Routh-Hurwitz criteria, the characteristic equation has roots with negative real parts iff the coefficients $c_k$ in \eqref{OriginCharEq} are positive and $c_2 c_1 - c_0 > 0$.  These quantities are clearly positive at $\Ray = 0$, so one can check that the coefficient $c_0$ in \eqref{OriginCharEq} changes sign at $\mathcal{R}^{\textbf{m},1}$, the quantity $c_2 c_1 - c_0$ changes sign at at $\mathcal{R}^{\textbf{m},2}$, and the coefficient $c_1$ changes sign at $\mathcal{R}^{\textbf{m},3}$ defined by
	\[ \mathcal{R}^{\textbf{m},3} =  \frac{(\Pra + 2) |\mathcal{K}\textbf{m}|^4 + \Pra \Rot^2 m_3^2}{\Sha^2 m_1^2} \text{ . } \]
	For $\Pra \geq 1$ one can verify $\Ray^{\textbf{m},1} < \mathsf{min}( \Ray^{\textbf{m},2}, \Ray^{\textbf{m},3})$.  On the other hand, for $0 < \Pra < 1$ one has $\Ray^{\textbf{m},1} < \mathsf{min}( \Ray^{\textbf{m},2}, \Ray^{\textbf{m},3})$ only for $|\Rot| < \Rot^{\textbf{m}}$.  For $|\Rot| = \Rot^{\textbf{m}}$ one has $\Ray^{\textbf{m},1} = \Ray^{\textbf{m},2} = \Ray^{\textbf{m},3}$, and $|\Rot| > \Rot^{\textbf{m}}$ one has $\Ray^{\textbf{m},2} < \mathsf{min}( \Ray^{\textbf{m},1}, \Ray^{\textbf{m},3})$.  Hence claim (i) follows.  Since the last coefficient is equal to the product $\Pi_{j} \lambda^{\textbf{m},j}$ and at least one of the roots must have negative real part, claim (ii) follows.  When $\Ray^{\textbf{m},2} < \Ray^{\textbf{m},1}$ the last coefficient in the characteristic equation is strictly positive on an interval around $\Ray = \Ray^{\textbf{m},c}$, hence claim (iv) above follows if claim (iii) is proven.
	
	To prove claim (iii), note that the discriminant $D(\lambda)$ of the characteristic equation is given by
	\begin{equation} \notag
		\begin{split} |\mathcal{K}\textbf{m}|^6 D(\lambda ) = & 4 \tilde{\Ray}^3 + \big ( (\Pra -1)^2 |\mathcal{K}\textbf{m}|^6 - 12 m_3^2 \Pra^2 \Rot^2 \big ) \tilde{\Ray}^2 \\ & - 4m_3^2 \Pra^2 \Rot^2 \big ( 5(\Pra -1)^2 |\mathcal{K}\textbf{m}|^6 - 3m_3^2 \Pra^2 \Rot^2 \big ) \tilde{\Ray} \\ & - 4m_3^2 \Pra^2 \Rot^2 \big ( (\Pra -1)^2 |\mathcal{K}\textbf{m}|^6 + m_3^2\Pra^2\Rot^2 \big )^2  \text{ , } 
		\end{split} 
	\end{equation}
	in which $\tilde{\Ray} = \Sha^2m_1^2 \Pra \Ray$.  For $\Rot = 0$ this is always positive, whereas for $\Rot \neq 0$ this is negative for $\Ray =0$ and positive for $\Ray$ sufficiently large, hence claim (iii) follows if this discriminant has only one positive root in $\Ray$.  To this end let $\rho_1,\rho_2,\rho_3$ denote the roots of the discriminant, and note that in the special case $\Pra = 1$, the roots of the characteristic equation are 
	\begin{align} \label{PrandtlOneEVals} \lambda^{\textbf{m},1} = - |\mathcal{K}\textbf{m}|^2 & + \sqrt{\frac{\Sha^2m_1^2 \Ray - m_3^2 \Rot^2}{|\mathcal{K}\textbf{m}|^2}} \hspace{.25 cm}  \text{ , } \hspace{.25 cm} \lambda^{\textbf{m},2} = - |\mathcal{K}\textbf{m}|^2 \text{ , } \\ \lambda^{\textbf{m},3} = & - |\mathcal{K}\textbf{m}|^2 - \sqrt{\frac{\Sha^2m_1^2 \Ray - m_3^2 \Rot^2}{|\mathcal{K}\textbf{m}|^2}} \text{ , } \notag \end{align}
	so here it is obvious that the discriminant has only one positive root in $\Ray$, and in fact one has $\rho_1 = \rho_2 = \rho_3 = \frac{m_3^2\Rot^2}{\Sha^2m_1^2}$.  Note also that $\rho_1,\rho_2,\rho_3$ depend continuously on $\parameters$.  One can therefore check how the roots $\rho_i$ behave for $0 < |\Pra - 1| \ll 1$ by checking the discriminant of the discriminant $D( D (\lambda ))$, given by
	\[ D( D (\lambda )) = \frac{16}{|\mathcal{K}\textbf{m}|^{12}} m_3^2\Pra^2 \Rot^2 (\Pra-1)^4 \big ( (\Pra-1)^2 |\mathcal{K}\textbf{m}|^6 - 27 m_3^2\Pra^2 \Rot^2 \big )^3 \text{ . } \]
	This is of course zero for $\Pra =1$, but for $0 < |\Pra - 1| < \frac{27 m_3^2\Pra^2\Rot^2}{|\mathcal{K}\textbf{m}|^6}$ it is negative, i.e. there is one real root $\rho_1$, which must be positive by continuity, and two complex roots $\rho_2,\rho_3$.  Note that the product of the roots must be
	\[ \rho_1 \rho_2 \rho_3 = 4m_3^2 \Rot^2 \big ( (\Pra -1)^2 |\mathcal{K}\textbf{m}|^6 + m_3^2\Pra^2\Rot^2 \big )^2 \text{ . } \]
	Since this is strictly positive for $\Rot \neq 0$, the roots $\rho_j$ can only cross the imaginary axis as a pair of complex conjugates, hence $\rho_1$ must remain positive on $0 < |\Pra - 1| < \frac{27 m_3^2\Pra^2\Rot^2}{|\mathcal{K}\textbf{m}|^6}$.  The complex roots $\rho_2,\rho_3$ join on the real axis at $|\Pra - 1| = \frac{27 m_3^2\Pra^2\Rot^2}{|\mathcal{K}\textbf{m}|^6}$, and apriori one could then see additional positive roots.  However, note that the sum of the roots must equal
	\[ \rho_1 + \rho_2 + \rho_3 = \frac{1}{|\mathcal{K}\textbf{m}|^6} \big ( 12 m_3^2 \Pra^2 \Rot^2 - (\Pra -1)^2 |\mathcal{K}\textbf{m}|^6 \big ) \text{ . } \]
	This is positive only for $|\Pra - 1| < \frac{12 m_3^2\Pra^2\Rot^2}{|\mathcal{K}\textbf{m}|^6}$, so the pair of complex conjugate roots $\rho_2,\rho_3$ must cross the imaginary axis and have negative real part before the conjugate roots join on the real axis.  Finally, for $|\Pra - 1| > \frac{27 m_3^2\Pra^2\Rot^2}{|\mathcal{K}\textbf{m}|^6}$ the roots $\rho_2,\rho_3$ are trapped on the negative real axis, and it follows that $\Ray^* := \rho_1$ is the unique positive root of the discriminant.
	
	Finally, by inserting $\Ray = \Ray^{\textbf{m},1}$ into the characteristic equation and solving the quadratic equation for the non-trivial roots, one obtains the following expressions, from which (v) plainly follows:
	\begin{align} \label{NonCriticalEigenvalues} \lambda^{\textbf{m},2} & = - (\Pra + \frac{1}{2})|\mathcal{K}\textbf{m}|^2 + \sqrt{\frac{1}{4} |\mathcal{K}\textbf{m}|^4 + \Pra (1-\Pra)\frac{\Rot^2m_3^2}{|\mathcal{K}\textbf{m}|^2}} \text{ , } \\ \lambda^{\textbf{m},3} & = - (\Pra + \frac{1}{2})|\mathcal{K}\textbf{m}|^2 - \sqrt{\frac{1}{4} |\mathcal{K}\textbf{m}|^4 + \Pra (1-\Pra)\frac{\Rot^2m_3^2}{|\mathcal{K}\textbf{m}|^2}} . \notag \end{align}
	
	
\end{proof}

\subsection{Analysis of the level curves of the critical Rayleigh numbers}

\label{app:LevelCurves}

Both equations in \eqref{CriticalLevelSet} can be brought to the form
\[ x^2 = (x^2 + y^2)^3 + r y^2 \hspace{.25 cm} \text{ for } \hspace{.25 cm} 0 < x < 1 \text{ , } \]
via the transformations
\[  x = \frac{\Sha m_1}{\Ray^{1/4}} \text{ , } y = \frac{m_3}{\Ray^{1/4}} \text{ , } r = \frac{\Rot^2}{\Ray} \hspace{.25 cm} \text{ , } \hspace{.25 cm} x = (2\Pra +2)^{1/4} \frac{\Sha m_1}{\Ray^{1/4}} \text{ , } y = (2\Pra +2)^{1/4} \frac{\tilde{m}_3}{\Ray^{1/4}} \text{ , } r = \frac{2\Pra^2 \Rot^2}{\Ray(\Pra + 1)} , \]
respectively.  The above is cubic in $y^2$ with discriminant $- (1+r) x^2$.  Since this is strictly negative, there is precisely one real root.  For the case $r = 0$ it is trivial to prove $y(x)$ is concave, since there is a nice explicit solution.  For $r > 0$, let $\xi = x^2$, $\eta = y^2$ and implicitly differentiate:
\begin{equation} \label{AppLevelCurves_ImplicitDiff} \xi = (\xi+\eta)^3 + r\eta \hspace{.5 cm} \Rightarrow \hspace{.5 cm} \frac{d \eta}{d \xi} = \frac{1-3(\xi+\eta)^2}{3(\xi+\eta)^2+r} \hspace{.25 cm} \text{ , } \hspace{.25 cm} \frac{d^2 \eta}{d \xi^2} = - \frac{6(\xi+\eta)(1+\frac{d\eta}{d\xi})^2}{3(\xi+\eta)^2+r} . \end{equation}
Using the chain rule, one finds
\[  \frac{dy}{dx} = \frac{x}{y} \frac{d\eta}{d\xi} \hspace{.5 cm} \text{ , } \hspace{.5 cm}  \frac{d^2y}{dx^2} = \frac{1}{y^3} \big [ 2 \xi \eta \frac{d^2 \eta}{d \xi^2} -  \eta \frac{d \eta}{d \xi} - \xi \big ( \frac{d \eta}{d \xi} \big )^2 \big ] , \]
which must be negative to prove concavity.  Since $0 < x < 1$ and $\xi = x^2$, $\eta = y^2$ one must have $0 < \xi < 1$ and $\eta > 0$.  Furthermore since $\xi,\eta$ solve \eqref{AppLevelCurves_ImplicitDiff}, one has $\eta < \mathsf{min}(\xi^{1/3}-\xi,\frac{1}{r}\xi) < 1 -\xi$.  Note that $\frac{d \eta}{d\xi}$ is strictly positive for $\xi + \eta < \frac{1}{\sqrt{3}} =: U^{0,1}$, and zero for $\xi + \eta = U^{0,1}$, hence $\frac{d^2y}{dx^2}$ is easily seen to be strictly negative on this region.

We must therefore check that $\frac{d^2y}{dx^2}$ is negative on the region $U^{0,1} < \xi + \eta < 1$.  For this we expand the above and obtain
\begin{align} ( 3 (\xi + \eta)^2 + r) \big ( 2 \xi \eta \frac{d^2 \eta}{d \xi^2} - \eta \frac{d \eta}{d \xi} - \xi \big ( \frac{d \eta}{d \xi} \big )^2 \big ) =  & 3\eta (\xi+\eta)^2 - \eta -12 \xi \eta (\xi+ \eta) \Big ( \frac{1+r}{3(\xi+\eta)^2+r} \Big )^2 \notag \\ & - \xi \Big ( \frac{1-3(\xi+\eta)^2}{3(\xi+\eta)^2+r} \Big )^2 . \notag \end{align}
On the right hand side there is one positive term and three negative terms, hence we split the positive term between the negative terms using some constants $f_1^{i,j},f_2^{i,j} >0$ yet to be chosen:
\begin{align} ( 3 (\xi + \eta)^2 + r) \big ( 2 \xi \eta & \frac{d^2 \eta}{d \xi^2} - \eta \frac{d \eta}{d \xi} - \xi \big ( \frac{d \eta}{d \xi} \big )^2 \big ) = \eta T_1 + \eta (\xi + \eta) T_2 + T_3 ,  \notag \\  T_1 & =  (3-f_1^{i,j}-f_2^{i,j}) (\xi + \eta)^2 -1 \text{ , } \notag \\ T_2 & = f_1^{i,j}(\xi+\eta) - 12 \xi ( \frac{1+r}{3(\xi+\eta)^2+r} )^2 \text{ , } \notag \\ T_3 & = f_2^{i,j} \eta (\xi+\eta)^2 - \xi \frac{ (3(\xi+\eta)^2-1)^2 }{ 3(\xi+\eta)^2+r} . \notag \end{align}
One easily finds where these terms are non-positive:
\begin{align} \label{App_Ineq1} & T_1 \leq 0  \hspace{.5 cm} \Leftrightarrow \hspace{.5 cm} \xi + \eta \leq \big ( 3-f_1^{i,j}-f_2^{i,j} \big )^{-1/2} =: U^{i,j} \text{ , } \\ \label{App_Ineq2} & T_2 \leq 0  \hspace{.5 cm} \Leftrightarrow \hspace{.5 cm} f_1^{i,j} (\xi + \eta) \big ( 3(\xi+\eta)^2+r \big )^2 \leq 12 \xi (1+r)^2 , \\ \label{App_Ineq3} & T_3 \leq 0  \hspace{.5 cm} \Leftrightarrow \hspace{.5 cm} f_2^{i,j} \eta (\xi + \eta)^2 \big ( 3(\xi+\eta)^2+r \big ) \leq \xi \big ( 3(\xi+\eta)^2 -1 \big )^2 .  \end{align}
In the region $U_1^{i-1,j} < \xi + \eta \leq U^{i,j}$ one can use $\eta < \xi^{1/3} - \xi$ to prove $\xi > (U^{i-1,j})^3$.  One can therefore ensure the inequality for $T_2$ on this region by choosing $f_1^{i,j},f_2^{i,j}$ such that:
\begin{equation} \label{App_Ineq4} 0 \leq 12(U^{i-1,1})^3 (1+r)^2 -f_1^{i,j}U^{i,j} \big ( 3(U^{i,j})^2+r \big )^2 =: Q^{i,j}(r) . \end{equation}
On the other hand, the right hand side of the inequality in \eqref{App_Ineq3} is zero at the boundary $\xi + \eta = U_{1}^{0}$, hence first we must take $f_2^{1,1} = 0$ and then $T_3 < 0$ is trivially satisfied.  Choosing $f_1^{1,1} = 1.16$ for instance, the resulting $Q^{1,1}(r)$ is a quadratic polynomial with strictly positive coefficients, hence \eqref{App_Ineq4} is satisfied and $\frac{d^2y}{dx^2} < 0$ for $\xi + \eta \leq U^{1,1} \approx 0.74$ for all $r > 0$.   One can then iteratively expand the region on which $\frac{d^2y}{dx^2}$ is known to be negative by choosing $f_2^{i,1} = 0$ and choosing $f_1^{i,1}$ such that $Q^{i,1}(r)$ has strictly positive coefficients.  For example, one could choose
\[ f_1^{2,1} = 1.49 \hspace{.25 cm} \text{ , } \hspace{.25 cm} f_1^{3,1} = 1.61 \hspace{.25 cm} \text{ , } \hspace{.25 cm} f_1^{4,1} = 1.66 \hspace{.25 cm} \text{ , } \hspace{.25 cm} f_1^{5,1} = 1.68 \hspace{.5 cm} \Rightarrow \hspace{.5 cm} U^{5,1} \approx 0.87 \text{ . } \]
Further iterations only yield smaller improvements, so we must choose $f_2^{i,j} > 0$.  However, note the left hand side of the last inequality in \eqref{App_Ineq3} increases unboundedly with $r$, hence first we choose $f_1^{5,2} = 2, f_2^{5,2} = 0$.  This gives $T_3 < 0$, $T_1 \leq 0$ on the whole domain $\xi + \eta \leq 1$, but then $Q^{5,2}$ has some negative coefficients:
\[ Q^{5,2}(r) \approx 5.91 r^2 +3.83 r -10.09 \hspace{1 cm} \Rightarrow \hspace{1 cm} Q^{5,2}(r) > 0 \text{ for } r \geq 1.03 \text{ . } \]
Therefore $T_2 \leq 0$ on $\xi + \eta \leq 1$ for $r \geq r^5 := 1.03$.

Thus it remains to prove $\frac{d^2y}{dx^2} < 0$ on the region $U^{5,1} < \xi + \eta < 1$ for $r < r^{5}$.  To do so we can iteratively choose $f_2^{i,1}, f_2^{i,1} > 0$ to find an extended region $U^{i-1,1} < \xi + \eta \leq U^{i,1}$ where $\frac{d^2y}{dx^2}$ is negative for all $0 < r < r^{i-1}$, although it is now nontrivial to enforce $T_3 \leq 0$.  To do so one can use the upper bound $\eta \leq \xi^{1/3} - \xi$ to prove $\eta \leq \eta^{i-1,j} :=  U^{i-1,1} - \big (  U^{i-1,1} \big )^3$ on this new domain, and hence one has $T_3 \leq 0$ if $f_1^{i,1},f_2^{i,1}$ are chosen such that
\[ 0 \leq  \big ( U^{i-1,1} \big )^3 \big ( 3 \big ( U^{i-1,1} \big )^2 -1  \big )^2 - f^{i,1}_2 \eta^{i-1,1} \big ( U^{i,1}  \big )^2 \big ( 3 \big ( U^{i,1} \big )^2 + r^{i-1} \big ) =: V^{i,1} , \]
Choosing $f_1^{6,1} = 1.24$, $f_2^{6,1} = 0.61$ gives $V^{6,1} > 0$ and $Q^{6,1}(r)$ strictly positive, and finally choosing $f_1^{7,1} = 0$, $f_2^{7,1} = 2$ proves negativity on the remaining part of the region $\xi + \eta \leq 1$.

On the other hand, the boundary curve $m_1(m_3,\parameters)$ at which $\Ray^{\textbf{m},1} = \Ray^{\textbf{m},2}$ is easily shown to be concave.  For $\Pra \geq 1$, there is no boundary curve, whereas for $0 < \Pra < 1$ $m_1(m_3,\parameters)$ must solve:
\begin{equation} \label{CriticalBoundary} \Rot^2 (1-\Pra ) m_3^2 = (\Pra + 1) ( \Sha^2 m_1^2 + \tilde{m}_3^2 )^3 \hspace{.5 cm} \text{ for } \hspace{.5 cm} 0 < m_3 < \left ( \frac{\Rot^2(1-\Pra)}{\Pra +1} \right )^{1/4}. \end{equation}
Using $y = ( \frac{(1-\Pra)}{1+\Pra}\Rot^2 )^{-1/4} m_3$, $x = ( \frac{(1-\Pra)}{1+\Pra}\Rot^2 )^{-1/4}\Sha m_1$ one finds a solution which is clearly strictly concave on this interval:
\[ x = \sqrt{y^{2/3}-y^2} \hspace{1 cm} \text{ for } \hspace{1 cm} 0 < y < 1 \text{ . }  \]

\subsection{Proof of Theorem \ref{thm:OriginBifurcations}}

\label{app:ProofOfOriginBifur}

\begin{proof}[Proof of Theorem \ref{thm:OriginBifurcations}]
	
	There are plenty of tools available for proving the existence of the unstable manifold, but for concreteness we invoke specifically the invariant manifold theorem of Chen, Hale and Tan \cite{Chen1997Invariant}.  For a semi-group $\mathcal{S}(t)$ of operators acting on a Banach space $B$, this theorem requires the following hypotheses:
	\begin{enumerate}[label=(\roman*)]
		\item \textbf{(Regularity)} $\mathcal{S}(t): [0,\infty) \times B \mapsto B$ is continuous in $t$, Lipschitz in $B$ and for each $\tau > 0$ bounded as follows:
		\[ \sup_{0 \leq t \leq \tau} \sup_{\substack{u,v \in B \\ u \neq v}} \frac{\| \mathcal{S}(t) \big [ u \big ]  - \mathcal{S}(t) \big [ v \big ] \|_{B} }{\|u-v \|_B} < \infty . \]
		For each $t \in [0, \tau]$ one can also decompose $\mathcal{S}(t) = \mathcal{L}(t) + \mathcal{R}(t)$, where $\mathcal{L}(t)$ is a bounded linear operator and $\mathcal{R}(t)$ is globally Lipschitz and Frechet differentiable at $\textbf{0}$, and satisfies $\mathcal{R}(t)[\textbf{0}] = \textbf{0}$, $D\mathcal{R}(t)[\textbf{0}] = \textbf{0}$.
		\item \textbf{(Spectral gap)} There are subspaces $B_1,B_2$ with $B = B_1 \oplus B_2$ which are invariant with respect to $\mathcal{L}(t)$ and are associated with continuous projections $\mathcal{P}_i: B \mapsto B_i$.  Also $\mathcal{L}(t)$ commutes with $\mathcal{P}_i$.  If $\mathcal{L}_i(t) = \mathcal{L}(t)|_{B_i}$, then $\mathcal{L}_1(t)$ has bounded inverse and there exist constants $C_i$ and $\alpha_i$ such that $\alpha_1 > \alpha_2 \geq 0$ and for $t > 0$ one has
		\begin{equation} \label{Thm3_GrowthDecayEst} \| \mathcal{L}_1(-t) \mathcal{P}_1 \|_{C(B;B)} = \| \big ( \mathcal{L}_1(t) \big )^{-1}  \mathcal{P}_1 \|_{C(B;B)} \leq C_1 \alpha_1^{-t} \hspace{.15 cm} \text{ , } \hspace{.15 cm} \| \mathcal{L}_2(t) \mathcal{P}_2 \|_{C(B;B)} \leq C_1 \alpha_2^{t} \text{ . } \end{equation}
		\item \textbf{(Nonlinearity small compared to spectral gap)} The constants $C_i,\alpha_i$ are related to the nonlinearity via
		\begin{equation} \label{Thm3_SpecGapCond}  \frac{(\sqrt{C_1}+\sqrt{C_2})^2}{\alpha_1 - \alpha_2} \sup_{\substack{u,v \in B \\ u \neq v}} \frac{\| \mathcal{R}(t) \big [ u \big ] - \mathcal{R}(t) \big [ v \big ] \|_{B} }{\|u-v \|_B} < 1 \text{ . } \end{equation}
	\end{enumerate}
	In fact, (i) slightly overstates the required regularity, but this is fine for the present case.  When these are satisfied, the theorem gives the existence of a globally Lipschitz function $\Psi: B_1 \to B_2$ with $\Psi(\textbf{0}) = \textbf{0}$ whose graph is invariant under the semi-group $\mathcal{S}(t)$ and satisfies certain growth/decay estimates.  One can verify that the hypotheses hold for the present case as follows, and hence the unstable manifold exists. 
	\begin{enumerate}[label=(\roman*)]
		\item The required regularity for $\mathcal{S}(t)$ was obtained in Theorem \ref{thm:PDE_WellPosedness_Balances}.  Defining
		\[ \mathcal{L}(t) := \exp \big [ \mathcal{L}^{\textbf{0}} t \big ] \hspace{.5 cm} \text{ , } \hspace{.5 cm } \mathcal{R}(t) := \mathcal{S}(t) - \mathcal{L}(t) \text{ , } \]
		it follows from Lemma \ref{lem:OriginEigenvalues} that $\mathcal{L}(t)$ is a bounded linear operator, and the Lipschitz property and $\mathcal{R}(\textbf{0}) = \textbf{0}$ follow.  Of course, the Lipschitz constant in \eqref{Thm1_LipDepInitCond} is only bounded on bounded subsets of $\textbf{H}^0$, but since one only needs to prove existence in a neighborhood of the origin, one can introduce a cutoff function such that $\mathcal{R}(t)$ is zero outside of some ball around the origin, and hence obtain a global Lipschitz constant.  One can also check that $\mathcal{S}(t)$ is Frechet differentiable at the origin and that $D\mathcal{S}(t)[\textbf{0}] = \mathcal{L}(t)$.
		\item The required subspaces are given by $B_1 = \mathsf{span}$ $(\textbf{X}^u)$, $B_2 = \mathsf{span}$ $(\textbf{X}^s) \oplus \mathsf{span}$ $(\textbf{X}^c)$, which are clearly invariant under $\mathcal{L}(t)$ since they are formed by the spectral projections of $\mathcal{L}^{\textbf{0}}$, and similarly the commutability property is obvious.  Since $\mathcal{L}_1(t)$ acts on only a finite number of Fourier modes, it is invertible.  Letting $\beta_1 = \max_{\textbf{m} \in \mathscr{M}^{u}} \mathsf{Re}(\lambda^{\textbf{m},j}), \beta_2 = \max_{\textbf{m} \in \mathscr{M}^{s} \cup \mathscr{M}^c} \mathsf{Re}(\lambda^{\textbf{m},j})$ be the largest real parts among the stable/center and unstable eigenvalues, respectively, the growth/decay estimates \eqref{Thm3_GrowthDecayEst} hold with $\alpha_1 = e^{\beta_1}$, $\alpha_2 =  e^{\beta_2}$.  
		\item One can check using Duhamel's formula that with $\mathcal{R}(t)$ defined as above one has the Lipschitz bound
		\[ \| \mathcal{R}(t) \big [ \textbf{X}_{0} \big ] - \mathcal{R}(t) \big [ \textbf{X}_{0}^{*} \big ] \|_{\textbf{H}^0} \leq C(t) \| \textbf{X}_{0} \|_{\textbf{H}^0} \| \textbf{X}_{0} - \textbf{X}_{0}^* \|_{\textbf{H}^0} \text{ , } \] 
		where  $C(t)$ is independent of $\textbf{X}_{0}, \textbf{X}_{0}^*$.  As above, this Lipschitz constant is unbounded on the whole domain, but by considering short times and using a cutoff function on a sufficiently small neighborhood of the origin the spectral gap condition \eqref{Thm3_SpecGapCond} is satisfied. 
	\end{enumerate}
	
	The dimension of the unstable manifold is found simply by counting the number of positive eigenvalues of each $\mathcal{L}^{\textbf{0},\textbf{m}}$ at a given set of parameters.  In the case $\Pra \geq 1$, each integer vertex $\textbf{m}$ under the level curve $\Ray^{\textbf{m},1}$ corresponds to exactly one positive eigenvalue, hence the number of positive eigenvalues is comparable to the area under this curve.  This area can be computed by using the substitution $\Sha m_1 = \Ray^{1/4}x$, $m_3(m_1) = (\Ray + \Rot^2)^{1/4}y(x)$, $r = (1+\frac{\Rot^2}{\Ray})^{-3/4}$, and then using Cardano's formula to solve \eqref{CriticalLevelSet} for $y$:
	\begin{align} & \int_0^{\Ray^{1/4}/k_1} m_3(m_1) dm_1 = \frac{\Ray^{1/2}(1+\frac{\Rot^2}{\Ray})^{1/4}}{k_1} \int_0^{1} y(x,r)dx = \frac{\Ray^{1/2}(1+\frac{\Rot^2}{\Ray})^{1/4}}{k_1} F(r) \text{ , } \notag \\ & F(r) := \int_0^{1} \Big [ \Big ( \sqrt{p} + \frac{r^{\frac{2}{3}}x^2}{2} \Big )^{\frac{1}{3}} - \Big ( \sqrt{p} - \frac{r^{\frac{2}{3}}x^2}{2} \Big )^{\frac{1}{3}} - r^{\frac{2}{3}} x^2 \Big ]^{\frac{1}{2}} dx \hspace{.2 cm} \text{ , } \hspace{.2 cm} p := \frac{(1-r^{\frac{4}{3}})^3}{27} + \frac{r^{\frac{4}{3}}x^4}{4} \text{ . } \notag \end{align}
	This integrand is smooth for $r \in (0,1]$, continuous for $r \in [0,1]$ and equal to zero at $r = 0$.  In order to determine the asymptotic behavior as $r \to 0$, note the integrand is not smooth but one can expand using the identity
	\[ \sqrt{1-z} = 1 - \frac{z}{2} -\frac{z^2}{8} - \frac{z^3}{16} + \frac{\frac{z^4}{256}(20+4z+z^2)}{1 - \frac{z}{2} -\frac{z^2}{8} - \frac{z^3}{16} + \sqrt{1-z}} \text{ , } \]
	and analogous identities for $(1\pm z)^{1/3}$, $(1\pm z)^{2/3}$.  Note that the last term in this expansion can be uniformly bounded in a neighborhood of $z = 0$ by a term of order $z^4$.  Using these expansions one can obtain the following identity in a neighborhood of $r = 0$:
	\[ F(r) = r \int_0^1 \sqrt{x^2-x^6} dx  + \mathscr{O}(r^{5/3}) =  \frac{\pi^2}{8} r  + \mathscr{O}(r^{5/3}) \text{ . } \]
	Thus $F(r)/r$ is bounded above and below by some constants on $r \in [0,1]$, and the bounds in \eqref{UnstableManifoldDimension} follow.  On the other hand, when $0 < \Pra < 1$ the integer vertices $\textbf{m}$ which undergo type 2 crossings temporarily increase the dimension by 2, before the subsequent type 3 crossing decreases the dimension by 1.  As in Figure \ref{fig:CriticalLevelSets} the set of such integer vertices is bounded by a convex curve and is therefore finite, so in this case the bound \eqref{UnstableManifoldDimension} still holds with different constants.  One can check $F(1) = 1/2$, and since the number of integer vertices contained beneath the level curve will be asymptotically equal to the area under the curve, one obtains the limit in \eqref{UnstableManifoldDimension}.

	The result of part (b) follows by proving the existence of a parameter dependent invariant manifold in a neighborhood of the origin and a small neighborhood of $\parameters^c$, determining its Taylor coefficients out to second order and analyzing the reduced system.  The theorem of Chen, Hale and Tan doesn't address higher orders of smoothness required for order two Taylor coefficients, hence we use here Theorem 3.3 of Chapter 2 of Haragus and Iooss's book \cite{haragus2011local}, which gives existence of a smooth, parameter dependent center manifold under some conditions regarding the equation itself rather than the semigroup.  For Banach spaces $\mathscr{Z} \hookrightarrow \mathscr{Y} \hookrightarrow \mathscr{X}$ with continuous embeddings and an equation of the form
	\[ \frac{du}{dt} = \textbf{L}u + \textbf{R}(u,\mu) \text{ , } \]
	the Theorem of Haragus and Iooss requires the following hypotheses:
	\begin{enumerate}[label=(\roman*)]
		\item $\textbf{L}$ is a bounded linear operator from $\mathscr{Z}$ to $\mathscr{X}$, and $\textbf{R}$ is a $C^k$ map from some neighborhood of $(u,\mu) = (0,0)$ in $\mathscr{Z} \times \mathbb{R}^m$ into $\mathscr{Y}$.  Furthermore $\textbf{R}(0,0) = 0$, $D_u\textbf{R}(0,0) = 0$.
		\item Decomposing the spectrum $\boldsymbol{\sigma}$ of $\textbf{L}$ as $\boldsymbol{\sigma} = \boldsymbol{\sigma}_+ \cup \boldsymbol{\sigma}_0 \cup \boldsymbol{\sigma}_-$, where
		\[ \boldsymbol{\sigma}_{\pm} = \{ \lambda \in \boldsymbol{\sigma}; \mathsf{Re}(\lambda) \gtrless 0 \} \hspace{.25 cm} \text{ , } \hspace{.25 cm}  \boldsymbol{\sigma}_0 = \{ \lambda \in \boldsymbol{\sigma}; \mathsf{Re}(\lambda) = 0 \} \text{ , } \]
		the set $\boldsymbol{\sigma}_0$ must consist of a finite number of eigenvalues with finite multiplicity and there must exist a positive constant $\gamma > 0$ such that
		\[ \inf_{\lambda \in \boldsymbol{\sigma}_+ } \mathsf{Re}(\lambda) > \gamma \hspace{.25 cm} \text{ , } \hspace{.25 cm} \sup_{\lambda \in \boldsymbol{\sigma}_- } \mathsf{Re}(\lambda) < \gamma \text{ . } \]
		\item There exist positive constants $\omega_0 > 0, c > 0,$ and $\alpha \in [0,1)$ such that for all $\omega \in \mathbb{R}$, with $|\omega| > \omega_0$ we have that $i \omega$ belongs to the resolvent set of $\textbf{L}$, and 
		\[ \| (i\omega I - \textbf{L})^{-1} \|_{\mathscr{L}(\mathscr{Z})} \leq \frac{c}{|\omega|} \hspace{.5 cm} \text{ , } \hspace{.5 cm} \| (i\omega I - \textbf{L})^{-1} \|_{\mathscr{L}(\mathscr{Y},\mathscr{Z})} \leq \frac{c}{|\omega|^{1-\alpha}} \text{ . } \]
	\end{enumerate}
	In particular, hypothesis (i) requires that the nonlinear term must be $C^2$.  This is easily verified if one chooses $\mathscr{Z} = \textbf{C}^2$, $\mathscr{Y} = \textbf{C}^1$, $\mathscr{X} = \textbf{C}^0$.  Since the semi-group has been shown to have smoothing properties, one can restrict to these subspaces, especially since the goal is only to classify the bifurcations, rather than prove growth/decay estimates for solutions beginning from general initial data.  The other hypotheses have already been verified, for instance (iii) follows from the application of Gershgorin's theorem \eqref{GershgorinDiscs}.  Hence a $C^2$ center manifold exists.
	
	In order to analyze the bifurcation structure, recall the projections $\textbf{X}^s,\textbf{X}^c,\textbf{X}^u$ introduced in before the statement of Theorem \ref{thm:OriginBifurcations}.  Furthermore let $\textbf{X}^{\textbf{m}}$ denote the vector of all of the Fourier variables with a wave vector $\textbf{m} \in \mathbb{Z}^2_{\geq 0}$, for one has the following explicit expression:
	\begin{equation} \label{def:FourierProjectedVars} \textbf{X}^{\textbf{m}} = \left \{ \begin{array}{cl}
			(u^{\textbf{m}},w^{\textbf{m}},\theta^{\textbf{m}}) & m_1 , m_3 > 0 , \\
			(u^{\textbf{m}},w^{\textbf{m}}) & m_1 = 0 \text{ , } m_3 > 0 \text{ , } m_3 \text{ odd, } \\
			\theta^{\textbf{m}} & m_1 = 0 \text{ , } m_3 > 0 \text{ , } m_3 \text{ even, }  \\
			w^{\textbf{m}} & m_1 > 0 \text{ , } m_3 = 0 .
		\end{array} \right . \end{equation}
	The PDE \eqref{EvolutionPDE} can be rewritten as an equation for each $\textbf{X}^{\textbf{m}}$, analogous to \eqref{GeneralBoussinesqODE_Vel} - \eqref{GeneralBoussinesqODE_Temp}  but now with $\textbf{m}$ ranging over all $\textbf{m} \in \mathbb{Z}^2_{\geq 0 } \setminus \{ \textbf{0} \}$.  In particular we write \eqref{EvolutionPDE} as 
	\begin{equation} \label{EvolutionPDE_FourierVersion} \frac{d}{dt} \textbf{X}^{\textbf{m}} = \mathcal{L}^{\textbf{0},\textbf{m}} \textbf{X}^{\textbf{m}} - \textbf{N}^{\textbf{m},c}(\textbf{X}^{c},\textbf{X}^s,\textbf{X}^u) - \tilde{\textbf{N}}^{\textbf{m},c}(\textbf{X}^s,\textbf{X}^u) \text{ , } \end{equation}
	where $\textbf{N}^{\textbf{m},c}$ denotes the sum of all nonlinear terms which involve at least one variable from $\textbf{X}^{c}$ and $\tilde{\textbf{N}}^{\textbf{m},c}$ is all other nonlinear terms involving only variables from $\textbf{X}^{s},\textbf{X}^{u}$.  Since we only consider the case $|\mathscr{M}^c| = 1$, for $\textbf{m} \in \mathscr{M}^c$, the only nonlinear terms involving variables from $\textbf{X}^c$ are as follows:
	\[ \textbf{N}^{\textbf{m},c} = \begin{pmatrix}
		I_{u}^{(\textbf{m},\textbf{m},2\textbf{m})} u^{\textbf{m}} u^{2\textbf{m}} + I_{u}^{(\textbf{m},2\textbf{m},\textbf{m})} u^{2\textbf{m}} u^{\textbf{m}} \\ 
		I_{w}^{(\textbf{m},\textbf{m},2\textbf{m})} u^{\textbf{m}} w^{2\textbf{m}} + I_{w}^{(\textbf{m},2\textbf{m},\textbf{m})} u^{2\textbf{m}} w^{\textbf{m}} + I_{w}^{(\textbf{m},\textbf{m},(2m_1,0))} u^{\textbf{m}} w^{(2m_1,0)}  \\  I_{\theta}^{(\textbf{m},\textbf{m},2\textbf{m})} u^{\textbf{m}} \theta^{2\textbf{m}} + I_{\theta}^{(\textbf{m},2\textbf{m},\textbf{m})} u^{2\textbf{m}} \theta^{\textbf{m}} + I_{\theta}^{(\textbf{m},\textbf{m},(0,2m_3))} u^{\textbf{m}} \theta^{(0,2m_3)} 
	\end{pmatrix} , \]
	where we have written $I_{u}^{\boldsymbol{\mu}},I_{w}^{\boldsymbol{\mu}},I_{\theta}^{\boldsymbol{\mu}}$ for the coefficients $I_{\textbf{u}}^{\boldsymbol{\alpha}},I_{\theta}^{\boldsymbol{\alpha}}$ defined in \eqref{def:NonlinCoeffs1}-\eqref{def:Zeta}.  This is unambiguous, since the component information $c,c',c''$ is determined by the variables $u,w,\theta$, and the phase information $\boldsymbol{\phi}$ is determined by the wave vector triad $\boldsymbol{\mu}$, since we are considering phase locked solutions as in \eqref{OlsonDoering_PhaseCond}.  Many of these terms drop out, for instance one can check by inserting $\textbf{m}'=\textbf{m},\textbf{m}'' = 2\textbf{m}$ into \eqref{def:Zeta} that $I^{(\textbf{m},\textbf{m},2\textbf{m})}_u = 0$.  By evaluating $I_{u}^{\boldsymbol{\mu}},I_{w}^{\boldsymbol{\mu}},I_{\theta}^{\boldsymbol{\mu}}$, one finds 
	\[ \textbf{N}^{\textbf{m},c} = I^{\textbf{m}} \begin{pmatrix}
		0 \\ 
		- u^{\textbf{m}} w^{(2m_1,0)}  \\ u^{\textbf{m}} \theta^{(0,2m_3)} 
	\end{pmatrix} \hspace{.5 cm} \text{ for } \hspace{.5 cm} I^{\textbf{m}} := \frac{(-1)^{|\textbf{m}|_1} \Sha m_1 m_3}{\sqrt{2}|\mathcal{K}\textbf{m}|V} \text{ . } \]
	By integrating by parts, one can check the following anti-symmetry property
	\[ I^{(\textbf{m},\textbf{m}',\textbf{m}'')}_{w} = \int_{\Omega} \textbf{v}^{\textbf{n}} \cdot \big [ (\textbf{v}^{\textbf{n}'} \cdot \nabla ) \textbf{v}^{\textbf{n}''} \big ] d\textbf{x} = - \int_{\Omega} \textbf{v}^{\textbf{n}''} \cdot \big [ (\textbf{v}^{\textbf{n}'} \cdot \nabla ) \textbf{v}^{\textbf{n}} \big ] d\textbf{x} = - I^{(\textbf{m}'',\textbf{m}',\textbf{m})}_{w} \text{ , } \]
	and one finds the same property for $I_{\theta}^{\boldsymbol{\mu}}$, hence one has the explicit equations
	\begin{align} \label{BifThm_NonlinInteractEq1} 
		\frac{d}{dt} w^{(2m_1,0)} & = - 4 \Pra \Sha^2 m_1^2 w^{(2m_1,0)} - I^{\textbf{m}} u^{\textbf{m}} w^{\textbf{m}} - \tilde{N}^{(2m_1,0),c}(\textbf{X}^s,\textbf{X}^u) \text{ , } \\ \label{BifThm_NonlinInteractEq2} 
		\frac{d}{dt} \theta^{(0,2m_3)} & = -4m_3^2 \theta^{(0,2m_3)} + I^{\textbf{m}} u^{\textbf{m}} \theta^{\textbf{m}} - \tilde{N}^{(0,2m_3),c}(\textbf{X}^s,\textbf{X}^u) \text{ . } \end{align}
	Next, as proven in Lemma \ref{lem:OriginEigenvalues} some of the eigenvalues of $\mathcal{L}^{0,\textbf{m}}$ have strictly negative real part, hence the corresponding eigenprojections of $\textbf{X}^{\textbf{m}}$ belong to the stable subspace at $\parameters^c$.  On the other hand, the eigenprojections corresponding to the eigenvalue(s) which pass through the imaginary axis will be called "critical variables".  
	
	In order to separate the critical variables from the stable subspace, we need to consider the eigenvectors of $\mathcal{L}^{0,\textbf{m}}$.  For admissible $\parameters$ and $\textbf{m} \in \mathbb{Z}^2_{> 0}$, let $\lambda^{\textbf{m},j}$ denote the eigenvalues of $\mathcal{L}^{\textbf{0},\textbf{m}}$ in lexicographic order, as described after \eqref{OriginCharEq}.   By inserting $- \Pra |\mathcal{K}\textbf{m}|^2$ into \eqref{OriginCharEq} one obtains $\Sha^2 m_1^2 \Pra \Ray (1-\Pra)$, hence for $\Pra \neq 1$, $\Rot \neq 0$ none of the eigenvalues $\lambda^{\textbf{m},j}$ are equal to $- \Pra |\mathcal{K}\textbf{m}|^2$.   For such $\parameters$ one can therefore define 
	\begin{equation} \label{def:EigVecs} 
		\tilde{\textbf{v}}^{\textbf{m},j} := \begin{pmatrix}
			|\mathcal{K}\textbf{m}|^2 + \lambda^{\textbf{m},j} \\
			-\Pra \Rot \frac{m_3 (|\mathcal{K}\textbf{m}|^2 + \lambda^{\textbf{m},j})}{|\mathcal{K}\textbf{m}|( \Pra |\mathcal{K}\textbf{m}|^2 + \lambda^{\textbf{m},j})} \\
			(-1)^{|\textbf{m}|_1+1} \frac{\Sha m_1}{|\mathcal{K}\textbf{m}|}
		\end{pmatrix} \text{ . } \end{equation}
	We claim that for arbitrary admissible $\parameters$ (also for $\Pra = 1$, $\Rot = 0$), one can then define the following matrix $\mathcal{V}^{\textbf{m}}$
	\begin{equation} \label{def:EigVecs2} {\normalfont \mathcal{V}^{\textbf{m}}_{i,j} := \lim_{\substack{\tilde{\parameters} \to \parameters \\ \tilde{\Pra} \neq 1 \text{ , } \tilde{\Rot} \neq 0}} \frac{\tilde{v}^{\textbf{m},j}_i}{|\tilde{\textbf{v}}^{\textbf{m},j}|} \text{ . } }     \end{equation}
	For $\Pra  =1$ and $\Rot = 0$, one has $\lambda^{\textbf{m},2} = -\Pra |\mathcal{K}\textbf{m}|^2$, hence \eqref{def:EigVecs} is of indeterminate form.  The eigenvalues $\lambda^{\textbf{m},j}$ are always continuous functions of $\parameters$, and using the characteristic equation \eqref{OriginCharEq} one obtains the identity
	\begin{equation} \label{EigVecIdentity} \frac{\Pra \Rot m_3 ( |\mathcal{K}\textbf{m}|^2 + \lambda^{\textbf{m},j})}{\Pra |\mathcal{K}\textbf{m}|^2 + \lambda^{\textbf{m},j}} = \frac{(\Pra \Ray \Sha^2 m_1^2 - |\mathcal{K}\textbf{m}|^2(\Pra |\mathcal{K}\textbf{m}|^2 + \lambda^{\textbf{m},j})( |\mathcal{K}\textbf{m}|^2 + \lambda^{\textbf{m},j}))}{\Pra \Rot m_3 } \text{ . } \end{equation}
	Using this identity and the normalization one sees that the limit in \eqref{def:EigVecs2} is also well defined, and the matrix $\mathcal{V}^{\textbf{m}}$ defined in this way is then automatically has continuous dependence on $\parameters$.  The vectors $\tilde{\textbf{v}}^{\textbf{m},j}$ are easily verified to be eigenvectors and hence the columns of $\mathcal{V}^{\textbf{m}}$ inherit this property.  Since we assume $|\Rot^c| \neq \Rot^{\textbf{m}}$ the eigenvalues are distinct in a neighborhood of $\parameters^c$, hence the eigenvectors are linearly independent.

	Since the matrix $\mathcal{V}^{\textbf{m}}$ has full rank one can define $\mathcal{U}^{\textbf{m}} = (\mathcal{V}^{\textbf{m}})^{-1}$ and separate the stable variables from the critical variables by projecting $\textbf{X}^{\textbf{m}}$ onto the eigenspaces:
	\[ \textbf{Z}^{\textbf{m}} := \mathcal{U}^{\textbf{m}} \textbf{X}^{\textbf{m}} \hspace{.5 cm} \text{ , } \hspace{.5 cm} \mathcal{D}^{\textbf{m}} := \mathcal{U}^{\textbf{m}} \mathcal{L}^{0,\textbf{m}} \mathcal{V}^{\textbf{m}} = \mathsf{diag}(\lambda^{\textbf{m},1},\lambda^{\textbf{m},2},\lambda^{\textbf{m},3}) \text{ . } \]
	Then \eqref{EvolutionPDE_FourierVersion} becomes:
	\begin{equation} \label{BifThm_GenCritSys} \frac{d}{dt} \textbf{Z}^{\textbf{m}} = \mathcal{D}^{\textbf{0},\textbf{m}} \textbf{Z}^{\textbf{m}} - I^{\textbf{m}} \mathcal{U}^{\textbf{m}} \sum_{j \leq 3} \begin{pmatrix}
			0 \\ 
			- \mathcal{V}_{1,j}^{\textbf{m}} Z^{\textbf{m}}_j w^{(2m_1,0)}  \\ \mathcal{V}_{1,j}^{\textbf{m}} Z^{\textbf{m}}_j \theta^{(0,2m_3)} 
		\end{pmatrix} - \mathcal{U}^{\textbf{m}}\tilde{\textbf{N}}^{\textbf{m},c}(\textbf{X}^s,\textbf{X}^u) \text{ . } \end{equation}
	Furthermore we will denote the critical variables by $\textbf{Z}^{\textbf{m},c}$, and note that $\textbf{Z}^{\textbf{m},c} = Z_1^{\textbf{m}}, Z_2^{\textbf{m}},$  or $(Z_1^{\textbf{m}},Z_2^{\textbf{m}})$, depending on the type of eigenvalue crossing.  Now we apply the center manifold theorem, which says that all non-critical variables can be written as a graph over the critical variables, and since $\Psi(\textbf{0}) = \textbf{0}$, $\nabla_{\textbf{Z}^c} \Psi ( \textbf{0}) = \textbf{0}$ all non-critical variables are at least quadratic in the critical variables.  In order to classify the type of bifurcation that occurs, \eqref{BifThm_GenCritSys} must be rewritten solely in terms of the critical variables out to cubic order.  The terms belonging to $\tilde{\textbf{N}}^{\textbf{m},c},\tilde{N}^{(2m_3,0),c} ,\tilde{N}^{(0,2m_3),c}$ consist of products of non-critical variables, hence are quartic order and can therefore be neglected.  So one needs only determine quadratic the Taylor coefficients for $w^{(2m_1,0)},\theta^{(0,2m_3)}$.  While in general the center manifold may not be unique, the Taylor coefficients of $\Psi$ at $\textbf{0}$ are unique, hence one can determine them using the invariance condition.  More explicitly, in the case where $\textbf{Z}^{\textbf{m},c} = Z_1^{\textbf{m}}$, one writes $w^{(2m_1,0)}$, $\theta^{(0,2m_3)}$ in terms of Taylor coefficients $w^{(2m_1,0)}_{1,1}$, $\theta^{(0,2m_3)}_{1,1}$ as follows:
	\[ w^{(2m_1,0)} = w_{1,1}^{(2m_1,0)} (Z_1^{\textbf{m}})^2 + \mathscr{O}(|Z_1^{\textbf{m}}|^3) \hspace{.5 cm} \text{ , } \hspace{.5 cm} \theta^{(0,2m_3)} = \theta_{1,1}^{(0,2m_3)} (Z_1^{\textbf{m}})^2 + \mathscr{O}(|Z_1^{\textbf{m}}|^3) \text{ , } \]
	One can then solve for $w^{(2m_1,0)}_{1,1}$, $\theta^{(0,2m_3)}_{1,1}$  by matching the left and right hand sides of \eqref{BifThm_NonlinInteractEq1}, \eqref{BifThm_NonlinInteractEq2}, and the result can be inserted into \eqref{BifThm_GenCritSys} to determine the type of bifurcation that occurs.
	
	From this point, the cases above must be considered more or less individually.  In the case $\textbf{m} \in \mathscr{M}^{c,1}$ one has $\textbf{Z}^{\textbf{m},c} = Z_1^{\textbf{m}}$, and hence $Z_2^{\textbf{m}},Z_3^{\textbf{m}}$ are non-critical and can be ignored. One obtains the following Taylor coefficients:
	\[ w_{1,1}^{(2m_1,0)} = \frac{- I^{\textbf{m}}\mathcal{V}_{1,1}^{\textbf{m}}\mathcal{V}_{2,1}^{\textbf{m}}}{2(\lambda^{\textbf{m},1} + 2\Pra \Sha^2 m_1^2)} \hspace{.5 cm} \text{ , } \hspace{.5 cm} \theta^{(0,2m_3)} = \frac{I^{\textbf{m}}\mathcal{V}_{1,1}^{\textbf{m}}\mathcal{V}_{3,1}^{\textbf{m}} }{2(\lambda^{\textbf{m},1} + 2m_3^2)} \text{ , }  \]
	where we note that since $\lambda^{\textbf{m},1} = 0$ at $\parameters^c$ the denominators are strictly positive in a neighborhood of $\parameters^c$.  Inserting these into \eqref{BifThm_GenCritSys} one obtains
	\begin{equation} \label{BifThm_Type1} \frac{d}{dt} Z_1^{\textbf{m}} = \lambda^{\textbf{m},1} Z_1^{\textbf{m}} - \frac{1}{2}(I^{\textbf{m}})^2 (\mathcal{V}_{1,1}^{\textbf{m}})^2 \Big ( \frac{\mathcal{U}_{1,2}^{\textbf{m}}\mathcal{V}_{2,1}^{\textbf{m}}}{\lambda^{\textbf{m},1}+2 \Pra \Sha^2m_1^2} + \frac{\mathcal{U}_{1,3}^{\textbf{m}}\mathcal{V}_{3,1}^{\textbf{m}}}{\lambda^{\textbf{m},1}+2 m_3^2} \Big )  \big ( Z_1^{\textbf{m}} \big )^3 + \mathscr{O}(|Z_1^{\textbf{m}}|^4) . \end{equation}
	By using the explicit expressions \eqref{def:EigVecs} one obtains the following identity at $\parameters = \parameters^c$:
	\begin{equation} \label{EigVecIdentity_CaseM1} \frac{\mathcal{U}_{1,2}^{\textbf{m}}\mathcal{V}_{2,1}^{\textbf{m}}}{\lambda^{\textbf{m},1}+2 \Pra \Sha^2m_1^2} + \frac{\mathcal{U}_{1,3}^{\textbf{m}}\mathcal{V}_{3,1}^{\textbf{m}}}{\lambda^{\textbf{m},1}+2 m_3^2} = \frac{ \Sha^2m_1^2 \Pra^2 (|\mathcal{K}\textbf{m}|^6 + m_3^2 \Rot^2 ) - m_3^4 \Rot^2 }{2\Sha^2m_1^2m_3^2\Pra ( (1+\Pra) |\mathcal{K}\textbf{m}|^6 + (\Pra - 1) \Rot^2 m_3^2  ) }  \text{ . }  \end{equation}
	Note this expression also holds when $\Pra^c =1$ and/or $\Rot^c = 0$ due to continuity.  Hence if $\Pra^c \geq \frac{m_3}{\Sha m_1}$ or $\Pra^c < \frac{m_3}{\Sha m_1}$ and $|\Rot^c| < C^{\textbf{m}}$ the coefficient of $(Z_1^{\textbf{m}})^3$ is strictly positive in a neighborhood of $\parameters^c$.  Therefore if one increases $\Ray$ through $\Ray^c$ one has $\lambda^{\textbf{m},1}$ passing from the negative real axis to the positive real axis, hence a supercritical pitchfork bifurcation occurs.  Alternatively if $\Pra^c < \frac{m_3}{\Sha m_1}$ and $|\Rot^c| > C^{\textbf{m}}$ the coefficient of $(Z_1^{\textbf{m}})^3$ is strictly negative in a neighborhood of $\parameters^c$, hence a subcritical pitchfork bifurcation occurs.
	
	The case $\textbf{m} \in \mathscr{M}^{c,3}$ is similar, although here  one has $\textbf{Z}^{\textbf{m},c} = Z_2^{\textbf{m}}$, and hence $Z_1^{\textbf{m}},Z_3^{\textbf{m}}$ are non-critical and can be ignored.  In this case one obtains the following Taylor expansion to quadratic order, where again the denominators are strictly positive in a neighborhood of $\parameters^c$:
	\begin{align} w_{1,1}^{(2m_1,0)} & = \frac{- I^{\textbf{m}}\mathcal{V}_{1,2}^{\textbf{m}}\mathcal{V}_{2,2}^{\textbf{m}}}{2(\lambda^{\textbf{m},2} + 2\Pra \Sha^2 m_1^2)} \big ( Z_2^{\textbf{m}} \big )^2 + \mathscr{O} (|Z_2^{\textbf{m}}|^3) \text{ , } \notag \\ \theta^{(0,2m_3)} & = \frac{I^{\textbf{m}}\mathcal{V}_{1,2}^{\textbf{m}}\mathcal{V}_{3,2}^{\textbf{m}} }{2(\lambda^{\textbf{m},2} + 2m_3^2)} \big ( Z_2^{\textbf{m}} \big )^2 + \mathscr{O} (|Z_2^{\textbf{m}}|^3) \text{ , } \notag \end{align}
	Inserting these into \eqref{BifThm_GenCritSys} one obtains
	\begin{equation} \label{BifThm_Type3} \frac{d}{dt} Z_2^{\textbf{m}} = \lambda^{\textbf{m},2} Z_2^{\textbf{m}} - \frac{1}{2}(I^{\textbf{m}})^2 (\mathcal{V}_{1,2}^{\textbf{m}})^2 \Big ( \frac{\mathcal{U}_{2,2}^{\textbf{m}}\mathcal{V}_{2,2}^{\textbf{m}}}{\lambda^{\textbf{m},2}+2 \Pra \Sha^2m_1^2} + \frac{\mathcal{U}_{2,3}^{\textbf{m}}\mathcal{V}_{3,2}^{\textbf{m}}}{\lambda^{\textbf{m},2}+2 m_3^2} \Big )  \big ( Z_2^{\textbf{m}} \big )^3 + \mathscr{O}(|Z_2^{\textbf{m}}|^4). \end{equation}
	By using the explicit expressions \eqref{def:EigVecs} one obtains the following identity at $\parameters = \parameters^c$:
	\begin{equation} \label{EigVecIdentity_CaseM3} \frac{\mathcal{U}_{2,2}^{\textbf{m}}\mathcal{V}_{2,2}^{\textbf{m}}}{\lambda^{\textbf{m},2}+2 \Pra \Sha^2m_1^2} + \frac{\mathcal{U}_{2,3}^{\textbf{m}}\mathcal{V}_{3,2}^{\textbf{m}}}{\lambda^{\textbf{m},2}+2 m_3^2}  = \frac{ \Sha^2m_1^2 \Pra^2 (|\mathcal{K}\textbf{m}|^6 + m_3^2 \Rot^2 ) - m_3^4 \Rot^2 }{2\Sha^2m_1^2m_3^2\Pra ( (1+\Pra) |\mathcal{K}\textbf{m}|^6 + (\Pra - 1) \Rot^2 m_3^2  ) }  \text{ . }  \end{equation}
	The right hand side is the same as in \eqref{EigVecIdentity_CaseM1}, however in this case if one decreases $\Ray$ through $\Ray^c$ one has $\lambda^{\textbf{m},2}$ passing from the negative real axis to the positive real axis, hence the same bifurcations occur by decreasing $\Ray$.
	
	Next, consider the case $\textbf{m} \in \mathscr{M}^{c,2}$.  In this case one has $\textbf{Z}^{\textbf{m},c} = (Z_1^{\textbf{m}},Z_1^{\textbf{m}})$, and $Z_3^{\textbf{m}}$ is non-critical and can be ignored.  Note that $\lambda^{\textbf{m},1},\lambda^{\textbf{m},2}$ are complex conjugates, hence their eigenvectors are also complex conjugates, and since $(u^{\textbf{m}},w^{\textbf{m}},\theta^{\textbf{m}})$ are real $Z_1^{\textbf{m}},Z_2^{\textbf{m}}$ are complex conjugates as well.  Therefore define $\rho,\phi$ such that $Z_1^{\textbf{m}} = \rho e^{i\phi}$, $Z_2^{\textbf{m}} = \rho e^{-i\phi}$.  One obtains the following Taylor expansions to quadratic order:
	\begin{align}
		w^{(2m_1,0)} & = \frac{-I^{\textbf{m}}\mathcal{V}_{1,1}^{\textbf{m}}\mathcal{V}_{2,1}^{\textbf{m}}}{2(\lambda^{\textbf{m},1} + 2\Pra \Sha^2 m_1^2)} \big ( Z_1^{\textbf{m}} \big )^2 - \frac{I^{\textbf{m}} ( \mathcal{V}_{1,1}^{\textbf{m}}\mathcal{V}_{2,2}^{\textbf{m}}+ \mathcal{V}_{1,2}^{\textbf{m}}\mathcal{V}_{2,1}^{\textbf{m}})}{\lambda^{\textbf{m},1} + \lambda^{\textbf{m},2} + 4\Pra \Sha^2 m_1^2} Z_1^{\textbf{m}} Z_2^{\textbf{m}}  \notag \\ & \hspace{1 cm} - \frac{I^{\textbf{m}}\mathcal{V}_{1,2}^{\textbf{m}}\mathcal{V}_{2,2}^{\textbf{m}}}{2(\lambda^{\textbf{m},2} + 2\Pra \Sha^2 m_1^2)} \big ( Z_2^{\textbf{m}} \big )^2 \text{ , }  \notag \\
		\theta^{(0,2m_3)} & = \frac{I^{\textbf{m}}\mathcal{V}_{1,1}^{\textbf{m}}\mathcal{V}_{3,1}^{\textbf{m}}}{2(\lambda^{\textbf{m},1} + 2m_3^2)} \big ( Z_1^{\textbf{m}} \big )^2 + \frac{I^{\textbf{m}} ( \mathcal{V}_{1,1}^{\textbf{m}}\mathcal{V}_{3,2}^{\textbf{m}}+ \mathcal{V}_{1,2}^{\textbf{m}}\mathcal{V}_{3,1}^{\textbf{m}})}{\lambda^{\textbf{m},1} + \lambda^{\textbf{m},2} + 4m_3^2} Z_1^{\textbf{m}} Z_2^{\textbf{m}} + \frac{I^{\textbf{m}}\mathcal{V}_{1,2}^{\textbf{m}}\mathcal{V}_{3,2}^{\textbf{m}}}{2(\lambda^{\textbf{m},2} + 2m_3^2)} \big ( Z_2^{\textbf{m}} \big )^2 \text{ . }  \notag
	\end{align}
	One can then insert these into \eqref{BifThm_GenCritSys}.  Using the chain rule and Taylor expanding the resulting sinusoids, it follows that $\rho$ and $\phi$ must solve
	\begin{equation} \label{Bifurcation_SingleType2} \frac{d}{dt} \rho = \mathsf{Re}(\lambda^{\textbf{m},1}) \rho  - (I^{\textbf{m}})^2 C(\phi) \rho^3 + \mathscr{O}( \rho^4 ) \hspace{.5 cm} \text{ , } \hspace{.5 cm} \frac{d}{dt} \phi = \mathsf{Im}(\lambda^{\textbf{m},1})  + \mathscr{O}( \rho^2 ) \text{ , } \end{equation}
	where
	\begin{align} 
		C(\phi) = & |\mathcal{V}_{1,1}^{\textbf{m}}| |\mathcal{U}_{2,2}^{\textbf{m}}| \Big ( \cos \big ( \mathsf{Arg}[ \mathcal{V}_{1,1}^{\textbf{m}} ] - \mathsf{Arg}[ \mathcal{U}_{2,2}^{\textbf{m}} ] \big ) + \cos \big ( 2\phi + \mathsf{Arg}[ \mathcal{V}_{1,1}^{\textbf{m}} ] + \mathsf{Arg}[ \mathcal{U}_{2,2}^{\textbf{m}} ]  \big ) \Big ) \notag \\ & \hspace{1 cm} \times \Big ( \frac{\mathsf{Re}(\mathcal{V}_{1,1}^{\textbf{m}}\mathcal{V}_{2,2}^{\textbf{m}})}{\mathsf{Re}(\lambda^{\textbf{m},1}) + 2 \Pra \Sha^2 m_1^2} + \left | \frac{\mathcal{V}_{1,1}^{\textbf{m}}\mathcal{V}_{2,2}^{\textbf{m}}}{\lambda^{\textbf{m},1} + 2 \Pra \Sha^2 m_1^2} \right | \cos \big ( 2\phi + \mathsf{Arg}[ \frac{\mathcal{V}_{1,1}^{\textbf{m}}\mathcal{V}_{2,2}^{\textbf{m}}}{\lambda^{\textbf{m},1} + 2 \Pra \Sha^2 m_1^2} ] \big ) \Big ) \notag \\ & + |\mathcal{V}_{1,1}^{\textbf{m}}| |\mathcal{U}_{2,3}^{\textbf{m}}| \Big ( \cos \big ( \mathsf{Arg}[ \mathcal{V}_{1,1}^{\textbf{m}} ] - \mathsf{Arg}[ \mathcal{U}_{2,3}^{\textbf{m}} ] \big ) + \cos \big ( 2\phi + \mathsf{Arg}[ \mathcal{V}_{1,1}^{\textbf{m}} ] + \mathsf{Arg}[ \mathcal{U}_{2,3}^{\textbf{m}} ]  \big ) \Big ) \notag \\ & \hspace{1.5 cm} \times \Big ( \frac{\mathsf{Re}(\mathcal{V}_{1,1}^{\textbf{m}}\mathcal{V}_{3,2}^{\textbf{m}})}{\mathsf{Re}(\lambda^{\textbf{m},1}) + 2 m_3^2} + \left | \frac{\mathcal{V}_{1,1}^{\textbf{m}}\mathcal{V}_{3,2}^{\textbf{m}}}{\lambda^{\textbf{m},1} + 2 m_3^2} \right | \cos \big ( 2\phi + \mathsf{Arg}[ \frac{\mathcal{V}_{1,1}^{\textbf{m}}\mathcal{V}_{3,2}^{\textbf{m}}}{\lambda^{\textbf{m},1} + 2 m_3^2} ] \big ) \Big ) \text{ . } \notag
	\end{align}
	For any $t > 0$ one can integrate the first equation in \eqref{Bifurcation_SingleType2} to obtain the following:
	\begin{equation} \label{Bifurcation_PicardIterate} \rho(t) = \rho_0 + \int_0^{t} \Big [ \mathsf{Re}(\lambda^{\textbf{m},1}) \rho(s) - (I^{\textbf{m}})^2 C(\phi) \rho^3(s) + \mathscr{O}( \rho^4 ) \Big ] ds \text{ . } \end{equation}
	From the second equation in \eqref{Bifurcation_SingleType2}, one sees that in a neighborhood $\mathscr{B}_1(0)$ of $\rho = 0$ the time derivative of the phase is bounded below, for instance by $\frac{1}{2}\mathsf{Im}(\lambda^{\textbf{m},1})$.  For some (possibly smaller) neighborhood $\mathscr{B}_2(0)$, solutions with $\rho_0 \in \mathscr{B}_2(0)$ will remain inside $\mathscr{B}_1(0)$ for at least time $t = \frac{4\pi}{\mathsf{Im}(\lambda^{\textbf{m},1})}$, and hence these solutions will make a full orbit from $\phi = 0$ to $2\pi$.  Hence these solutions have a well defined first return time $t^*(\rho_0)$, and one can therefore define the Poincar\'e map $\rho_0 \mapsto \rho(t^*)$ by evaluating \eqref{Bifurcation_PicardIterate} at $t^*$.  By iterating the formula \eqref{Bifurcation_PicardIterate}, using a first order Taylor expansion in $\rho_0$, and changing variables from $t$ to $\phi$ one obtains the following formula for the Poincar\'e map:
	\[ \rho(t^*) = \rho_0 + \mathsf{Re}(\lambda^{\textbf{m},1}) (1 + \mathsf{Re}(\lambda^{\textbf{m},1}) B(\rho_0))  \rho_0 - (I^{\textbf{m}})^2 \rho_0^3 \int_0^{2\pi} C(\phi) d\phi + \mathscr{O}( \rho_0^4 ) \text{ . } \]
	where $B(\rho_0)$ is a term from the Taylor expansion involving $\frac{d\rho}{d\rho_0}$, hence can be bounded.  We look for a nontrivial solution where $\rho(t^*) - \rho_0 = 0$, corresponding to a periodic orbit.  Since the term $B(\rho_0)\mathsf{Re}(\lambda^{\textbf{m},1})$ can be made arbitrarily small in a neighborhood of $\parameters = \parameters^c$, we look for solutions of
	\[  \mathsf{Re}(\lambda^{\textbf{m},1}) - (I^{\textbf{m}})^2 \rho_0^2 \int_0^{2\pi} C(\phi) d\phi + \mathscr{O}( \rho_0^3 ) = 0 \text{ . } \]
	In this case, one has $\lambda^{\textbf{m},3} = -(2\Pra +1 ) |\mathcal{K}\textbf{m}|^2$ at $\parameters = \parameters^c$, hence one can insert this into \eqref{OriginCharEq} to obtain the purely imaginary eigenvalues:
	\[ \lambda^{\textbf{m},1} = i \frac{\Pra}{|\mathcal{K}\textbf{m}|} \sqrt{\frac{1-\Pra}{1+\Pra} m_3^2 \Rot^2 -|\mathcal{K}\textbf{m}|^6} \hspace{.5 cm} \text{ , } \hspace{.5 cm} \lambda^{\textbf{m},2} = - i \frac{\Pra}{|\mathcal{K}\textbf{m}|} \sqrt{\frac{1-\Pra}{1+\Pra} m_3^2 \Rot^2 -|\mathcal{K}\textbf{m}|^6} \text{ . } \] 
	Note that $\textbf{m} \in \mathscr{M}^{c,2}$ only occurs for $\Pra < 1$, $\Rot >0$, hence one has $|\mathcal{V}_{1,1}^{\textbf{m}}| > 0$, and by using the explicit expressions \eqref{def:EigVecs} one obtains the following identity at $\parameters = \parameters^c$:
	\[ \int_0^{2\pi} C(\phi) d \phi = \pi f^{*} ( T_1 + T_2 + T_3 + \frac{ \Pra ( m_3^2 - \Sha^2 m_1^2 ) }{4 m_3^4 + |\lambda^{\textbf{m},1}|^2}T_4) \text{ , } \]
	where
	\[ f^{*} = \frac{ \Pra \Rot^2 m_3^2 |\mathcal{K}\textbf{m}|^5 }{ |\tilde{\textbf{v}}^{\textbf{m},1}|^2\big ( |\mathcal{K}\textbf{m}|^6 (1+\Pra)^2 (1+3 \Pra ) + (1-\Pra) \Pra^2 m_3^2 \Rot^2 \big ) (\Pra^2 |\mathcal{K}\textbf{m}|^4 + |\lambda^{\textbf{m},1}|^2)} \text{ , } \]
	which is clearly strictly positive for $\Pra < 1$, and 
	\begin{align} T_1 & = \frac{(1+\Pra)(\Pra |\mathcal{K}\textbf{m}|^4 + |\lambda^{\textbf{m},1}|^2)}{2\Pra \Sha^2 m_1^2} - \Pra \Big ( \frac{1}{m_3^2} + \frac{1}{ |\mathcal{K}\textbf{m}|^2} \Big ) \big ( (2\Pra^2 -1 )  |\mathcal{K}\textbf{m}|^4 + |\lambda^{\textbf{m},1}|^2 \big ) \text{ , } \notag \\  T_2 & =  \frac{1+\Pra}{2} \Big [ \frac{ 2\Pra \Sha^2 m_1^2 ( \Pra |\mathcal{K}\textbf{m}|^4 + |\lambda^{\textbf{m},1}|^2 ) + |\lambda^{\textbf{m},1}|^2(1-\Pra )|\mathcal{K}\textbf{m}|^2 }{|\lambda^{\textbf{m},1}|^2+4\Pra^2 \Sha^4 m_1^4} \Big ] \text{ , } \notag \\  T_3 & = (1+\Pra)\frac{(2\Pra + 1)|\mathcal{K}\textbf{m}|^2 \big ( \Pra |\mathcal{K}\textbf{m}|^4 + |\lambda^{\textbf{m},1}|^2 -  2\Pra \Sha^2 m_1^2 (1-\Pra )|\mathcal{K}\textbf{m}|^2 \big )}{2(|\lambda^{\textbf{m},1}|^2+4\Pra^2 \Sha^4 m_1^4)} \text{ , } \notag \\   T_4 & =  \big ( \Pra (1+2\Pra)\Sha^2 m_1^2 + (6\Pra^2 + \Pra - 2) m_3^2 \big ) |\mathcal{K}\textbf{m}|^2 + (2+\Pra + \frac{2m_3^2}{|\mathcal{K}\textbf{m}|^2} )|\lambda^{\textbf{m},1}|^2   \text{ . } \notag \end{align}
	Note that $T_2$ is a sum of strictly positive terms for $\Pra <1$.  One can easily check that if $m_3 \geq \sqrt{2}\Sha m_1$ then all of the negative terms in $T_1$, $T_3$ are more than compensated by the positive terms, hence these are strictly positive.  Finally, the prefactor for $T_4$ is clearly positive for $m_3 \geq \Sha m_1$, and furthermore one has $6 \Pra^2 + \Pra -2 \geq 0$ for $\Pra \geq 1/2$, hence $T_4$ is a sum of positive terms.  Therefore in this case, if one increases $\Ray$ through $\Ray^c$ a supercritical Hopf bifurcation occurs.  Note that for $m_3 \leq \Sha m_1$ it is possible that the integral of $C(\phi)$ is negative, in which case a sub-critical Hopf bifurcation occurs.  Theorem \ref{thm:OriginBifurcations} gives only a partial statement regarding the Hopf bifurcations at the origin only because of the complexity of determining the sign of the integral of $C(\phi)$, which is left open.
	
	\vspace{-.75 cm}
	\[ \textcolor{white!100}{ . } \]
\end{proof}

\section{Sketch of the proof regarding the upper bound on the attractor dimension}

\label{App:AttractDimUpperBound}

The crux of the proof studies how $d$-dimensional infinitesimal volumes are expanded or contracted by the flow, for arbitrary $d > 0$.  The goal is to find $d$ large enough such that all higher dimensional infinitesimal volumes are contracted by the flow, and the Hausdorff dimension of the attractor can be no larger than $d$ (see the general theorems on the relation between the Lyapunov exponents and the Hausdorff dimension \cite{temam_InfDimDynSys} Theorems V.3.1, V.3.3).  In order to study the distortion of infinitesimal volumes by the flow, one considers solutions of the linearization about the flow along the attractor.  Namely for an initial condition $\textbf{X}_0 \in \mathscr{A}$, let $\textbf{X}(t)$ denote the solution and for $j=1,...,d$, let $\textbf{Y}^j(t) = (\textbf{g}^j(t),\psi^j(t))$ be the solution of \eqref{LinearizedEquation} for some initial conditions $\textbf{Y}_0^j = (\textbf{g}^j_0,\psi^j_0) \in \textbf{H}^0$.	 The $d$-dimensional volume of the parallelepiped formed by these vectors is given by the determinant of the matrix of inner products, and it can be shown (for instance \cite{temam_InfDimDynSys} Lemma V.1.2) that this volume evolves via the formula:
\[ \mathsf{det} \Big ( \big \langle \textbf{Y}^i(t) \cdot \textbf{Y}^j(t) \big \rangle \Big )_{i,j=1}^{d} = \mathsf{det} \Big ( \big \langle \textbf{Y}^i_0 \cdot \textbf{Y}^j_0 \big \rangle \Big )_{i,j=1}^{d} \exp \Big [ \int_0^t   \mathsf{Tr} \Big [ \mathcal{L}^{ \textbf{X}_0 }(\tau)  \circ \mathcal{Q}^d(\tau) \Big ] d\tau \Big ] , \]
where $\mathcal{Q}^d(\tau)$ is the orthogonal projector of $\textbf{L}^2_{\sigma} \times L^2$ onto the subspace spanned by $\textbf{Y}^j(\tau)$, $j = 1,...,d$.  Letting $\uvec{Y}^j(\tau) = (\textbf{v}^j(\tau),\varphi^j(\tau))$, $j=1,...,d$ be orthonormal vectors spanning this same subspace, one has from the definition of $\mathcal{L}^{ \textbf{X}_0 }(\tau)$ 
\begin{equation} \notag
	\begin{split} \mathsf{Tr} & \Big [ \mathcal{L}^{ \textbf{X}_0 }(\tau)  \circ \mathcal{Q}^d(\tau) \Big ] = \sum_{j=1}^d \Big \langle \uvec{Y}^j(\tau) \cdot  \mathcal{L}^{ \textbf{X}_0 }(\tau) \uvec{Y}^j(\tau) \Big \rangle  \\ & = \sum_{j=1}^d - \Pra \| \textbf{v}^j \|_{\textbf{H}^1}^2 -  \| \varphi^j \|_{H^1}^2 + (\Pra \Ray + 1) \big \langle \varphi^j v^j_3 \big \rangle - \big \langle  \varphi^j \big ( \textbf{v}^j \cdot \nabla \theta \big ) + \textbf{v}^j \cdot \big [ \textbf{v}^j \cdot \nabla \textbf{u} \big ] \big \rangle . 
	\end{split} 
\end{equation}
The first two sign indefinite terms are easily bounded as follows using Cauchy-Schwarz, the maximum principle $|\theta(\tau)| \leq \pi$, Young's inequality, and the fact that $( \textbf{v}^j, \varphi^j)$ are normalized:
\begin{equation} \notag
	\begin{split} & \big | \langle \varphi^j \big ( \textbf{v}^j \cdot \nabla \theta \big ) \rangle \big | = \big | \langle \theta \big ( \textbf{v}^j \cdot \nabla \varphi^j \big ) \rangle \big | \leq \pi \| \textbf{v}^j \|_{\textbf{L}^2_{\sigma}} \| \varphi^j \|_{H^1} \leq \pi \| \varphi^j \|_{H^1} \leq \frac{1}{2} \big ( \pi^2 + \| \varphi^j \|_{H^1}^2 \big ) , \\ & \big \langle \varphi^j v^j_3 \big \rangle \leq \frac{1}{2} \big ( \| v^j_3 \|_{L^2}^2 + \| \varphi^j \|_{L^2}^2  \big ) \leq \frac{1}{2} .\end{split} 
\end{equation}
The last term sign indefinite term is more difficult to bound, and requires use of the Sobolev-Lieb-Thirring inequality, which in this case states
\begin{equation} \label{LiebThirringIneq} \| \sum_{j=1}^d ( \textbf{v}^j )^2 \|_{\textbf{L}^2_{\sigma}} \leq \kappa_1 \sum_{j=1}^d \| \textbf{v}^j \|_{\textbf{H}^1_{\sigma}} , \end{equation}
for a constant $\kappa_1$ depending only on the domain $\Omega$.  Hence one obtains
\[ \big | \sum_{j=1}^d \langle \textbf{v}^j \big ( \textbf{v}^j \cdot \nabla \textbf{u} \big ) \rangle \big | \leq \| \textbf{u} \|_{\textbf{H}^1_{\sigma}} \| \sum_{j=1}^d  (\textbf{v}^j)^2 \|_{\textbf{L}^2_{\sigma}} \leq  \frac{1}{2} \big ( \frac{\kappa_1^2}{\Pra} \| \textbf{u} \|_{\textbf{H}^1_{\sigma}}^2 + \Pra \sum_{j=1}^d \| \textbf{v}^j \|_{\textbf{H}^1_{\sigma}}^2 \big ) . \]
Putting all of these equations together, one obtains
\[ \mathsf{Tr} \Big [ \mathcal{L}^{ \textbf{X}_0 }(\tau)  \circ \mathcal{Q}^d(\tau) \Big ] \leq \frac{\pi^2+1}{2} d + \frac{\kappa_1^2}{2\Pra} \| \textbf{u} \|_{\textbf{H}^1_{\sigma}}^2 - \frac{1}{2} \Big ( \sum_{j=1}^d  \Pra \| \textbf{v}^j \|_{\textbf{H}^1}^2 +  \| \varphi^j \|_{H^1}^2 \Big ) . \]
The fact that $(\textbf{v}^j$, $\theta^j)$ are mutually orthogonal puts constraints on their Fourier expansions, and hence the following sum is bounded below by choosing $(\textbf{v}^j$, $\theta^j)$ to have the $d$ lowest allowed wavenumbers $(m_1,m_3)$.  The $H^1$ norms give a factor $m_1^2+m_3^2$, but since there are $n+1$ wavenumbers with $m_1+m_3=n$, the sum behaves as a sum over $j$:
\[ \sum_{j=1}^d  \Pra \| \textbf{v}^j \|_{\textbf{H}^1}^2 +  \| \varphi^j \|_{H^1}^2 \geq \sum_{j=1}^d  (1+\Pra)\mathsf{min}(1,\Sha^2) j \geq (1+\Pra)\mathsf{min}(1,\Sha^2) \frac{d(d+1)}{2} . \]
Finally, since \eqref{Attractor_IntegralBound} applies for all time for initial conditions on the attractor, one obtains
\begin{align} \int_0^t  & \mathsf{Tr} \Big [ \mathcal{L}^{ \textbf{X}_0 }(\tau)  \circ \mathcal{Q}^d(\tau) \Big ] d\tau \notag \\ & \leq \Big ( -\frac{(1+\Pra)\mathsf{min}(1,\Sha^2)}{4} d(d+1) +\frac{\pi^2+1}{2}d +\frac{\kappa_1^2}{2\Pra} \frac{\Ray}{2} ( \Ray +1+\frac{2}{\Pra} )(\sqrt{ \pi |\Omega| } + 2\epsilon )  \Big ) t . \notag \end{align}
This is quadratic in $d$, with leading coefficient negative, hence one can solve for the $d$ such that all higher dimensional volumes are contracted.  By gathering together the constants in this argument, one arrives at \eqref{DimensionUpperBound_Constants}.

\begin{Remark}
	Note that the constant $\kappa_1$ in the Lieb-Thirring inequality \eqref{LiebThirringIneq} can be taken to be 
	\[ \kappa_1 = \frac{4\pi^2}{\mathsf{min}(1,\Sha)} \inf_{1 < k < 2} 2^k \int_0^{\infty} \frac{ds}{(1+s)^k} \int_0^1 \rho^{1-k} (1- \rho)^{k} d\rho . \]
	This can be extracted by following the arguments in the Appendix of \cite{temam_InfDimDynSys}, although in these arguments there is a confusing notational choice, perhaps a typo.  Specifically, in the proof of the Birman-Schwinger inequality (Proposition 2.1) the exponent $k \geq 1$ is not the same $k$ which arises in the definition of the operator $\mathfrak{A}$ in (1.1).
\end{Remark}

\section{A survey of heat transport phenomena in the HKC-1 model}

\label{app:HeatTransportSurvey}

Consider the following remarks regarding the dynamics in the HKC-1 model, each of which appears as an aspect in the heat transport plots in Figure \ref{fig:BasicPhenomena}.  First consider $\Rot = 0$, where the HKC-1 model reduces to the Lorenz '63 model.

\begin{Remark}
	\label{rem:Lorenz}
	The following statements can be proven analytically:
	\begin{enumerate}[label=\textbf{(\alph*)}]
		\item \textbf{(Fixed points and maximal transport)} For small Rayleigh HKC-1 admits a Lyapunov function, so the origin is the global attractor for $\Ray < \Ray^{(1,1),c} = 6.75$, thus for any initial condition one has $\mathsf{Nu}^1 = 1$ (see \cite{Sparrow}).  Similar to Theorem \ref{thm:OriginBifurcations} (b), a pair of nontrivial fixed points emerge from the origin in the pitchfork bifurcation at $\Ray = 6.75$.  These fixed points exhibit the maximal heat transport for the Lorenz model \cite{SouzaDoering1}.  They are locally stable for $ 6.75 \leq \Ray \lesssim 166.97$, and at $\Ray \approx 166.97$ they lose stability in a subcritical Hopf bifurcation, whereupon solutions tend to the famous chaotic attractor.  
		\item \textbf{(Chaotic region)} The chaotic attractor is ergodic \cite{Tucker}, so the heat transport is well-defined.  This proof implies that in principle one could perform a rigorous numerical integration along the chaotic attractor to obtain its heat transport value within some desired error bounds, although the author is unaware of any work where this has been carried out.
		\item \textbf{(Large Rayleigh region)} As the Rayleigh number becomes very large, the heat transport of the fixed points tends toward a constant value.  Furthermore, the chaotic attractor collapses down to a stable periodic orbit \cite{OvsyRadeWelterLu_2023_LorenzLargeRayleigh}, and the heat transport realized on this periodic orbit tends to another constant which is strictly less than that of the non-trivial fixed points.
	\end{enumerate}
	Furthermore, the following statements appear to be true from numerical experiments and past studies, although the author is unaware whether they have been settled analytically:
	\begin{enumerate}[label=\textbf{(\alph*)}] \setcounter{enumi}{3}
		\item \textbf{(Fixed points)} For $ 6.75 \leq \Ray \lesssim 166.97$ the non-trivial fixed points appear to have a large basin of attraction and may be globally attractive.  Thus in Figure \ref{fig:BasicPhenomena} (a) solutions with random initial conditions appear to converge to the fixed point value, hence their Nusselt number does as well.  As $\Ray$ approaches $166.97$ orbits exhibit long, "psuedo-chaotic" behavior before settling down on the fixed points, thus a very long integration time is required and a finite time approximation may fail to capture the infinite Nusselt number.  Conversely, above the critical threshold trajectories may spend an arbitrarily long time near the fixed points before eventually eventually tracking the chaotic attractor.  As depicted in Figure \ref{fig:HeatTransportTimeConvergence} the heat transport may erroneously appear to converge to some value for a long time before slowly transitioning to another value. 
		
		\begin{figure}[H]
			\centering
			\includegraphics[height=50mm]{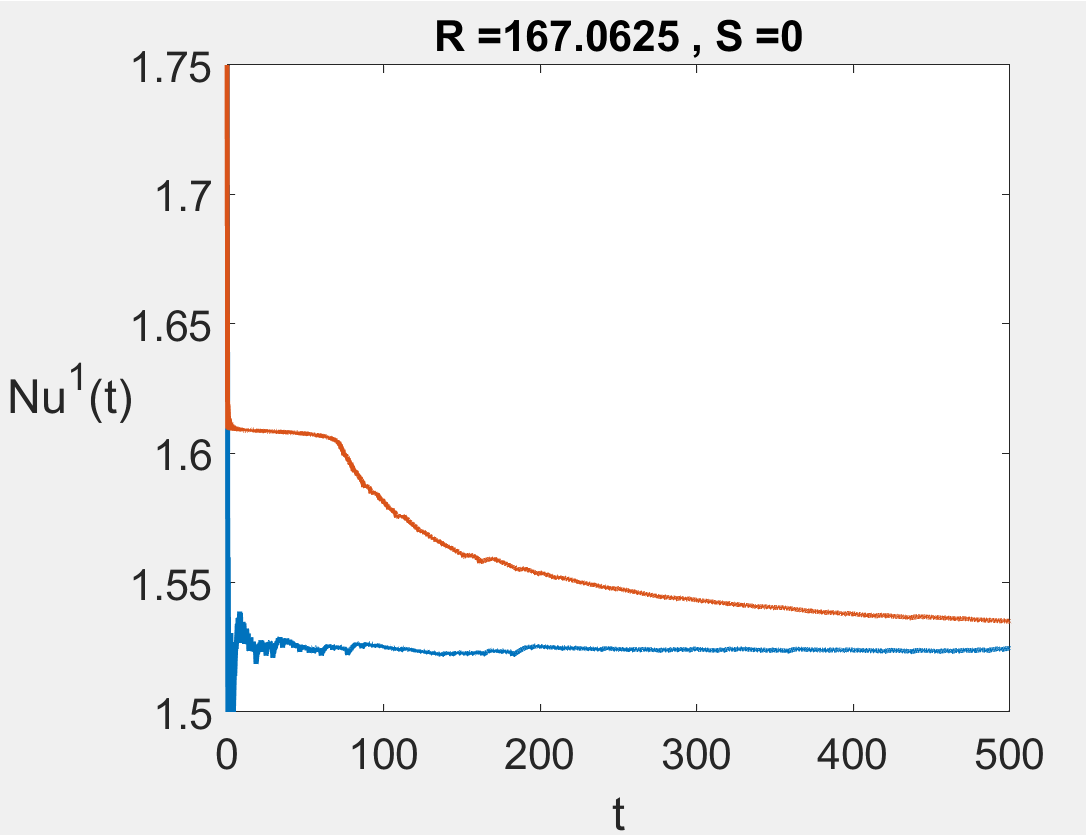}
			\caption{The apparent convergence of the finite time Nusselt number from HKC-1 towards an infinite time average for two different (random) initial conditions.  A time increment of $10^{-3}$ was used and the trajectories were computed for $5\cdot 10^{5}$ time steps.  The trajectory in orange was started closer to a non-trivial fixed point and requires a very long integration time for an accurate heat transport value.}
			\label{fig:HeatTransportTimeConvergence}
			\vspace{-.5 cm}
		\end{figure} 
		
		In this case both trajectories were determined from random perturbations from one of the non-trivial fixed points.  The trajectory in blue depicts the type of convergence scenario which was most commonly observed, namely after an initially large fluctuation the cumulative average quickly settles down to some value.  The trajectory in orange depicts another, more rare convergence scenario which was observed.  In this case the trajectory spends a fairly long time close to the fixed point before eventually converging to another value, hence the heat transport first appears to converge to one value before transitioning to another.  \verb|HeatTransport_Iterator.m| thus continues this trajectory until this threshold was met.
		
		\item \textbf{(Chaotic region)} Accurate numerical integration along the chaotic attractor is a difficult and computationally expensive task fraught with subtle sources of error.  For example, the persistent errors studied by Noethen \cite{noethen2023stepsize} do not average out over long integration times, but also require very small step sizes for accuracy.  In Figure \ref{fig:BasicPhenomena} (a) the trajectories were computed with step sizes of $10^{-4}$ for $4 \cdot 10^5$ time steps, and it was found that by the heat transport values changed by less than 5\% when decreasing the step size by a factor of $1/2$ and computing for $8 \cdot 10^{5}$ steps.  One observes here that trajectories appear to exhibit less heat transport than those of the unstable fixed points.  In this case, the fluid seems to be spending more energy moving around chaotically and hence transports less energy.
	\end{enumerate}
\end{Remark}

Introducing now $\Rot > 0$, consider the analogous statements and the effects of rotation:
\begin{Remark}
	\label{rem:HKC1rotation}
	
	The following statements can be proven analytically:
	\begin{enumerate}[label=\textbf{(\alph*)}]
		\item  \textbf{(Fixed points)} The origin is globally stable at least on the same region $\Ray < 6.75$, and locally stable on the larger parameter region $\Ray < \Ray^{(1,1),c} = 6.75+2\Rot^2$, and the non-trivial fixed points emerge in a pitchfork bifurcation only at this higher thermal forcing threshold.  After they emerge, the non-trivial fixed points are also locally stable on a larger parameter region.  This larger region has a boundary given by a somewhat complicated function in $\Rot$, but properly contains the region $ \Ray < 166.97 + 7.18 \Rot^2$ with simpler boundary. 
		\item \textbf{(Maximal transport)} At any fixed value of $\Ray$, the maximal heat transport among solutions with $\Rot > 0$ can never exceed the maximal heat transport among solutions with $\Rot = 0$.
	\end{enumerate}
	Furthermore, the following statements appear to be true from numerical experiments and past studies, although the author is unaware whether they have been settled analytically:
	\begin{enumerate}[label=\textbf{(\alph*)}]
		\item The origin appears to be the global attractor on the larger parameter range $0 \leq \Ray \leq \Ray^{(1,1),c}$.
		\item Fixing $\Rot > 0$ and increasing $\Ray$, one sees a similar story as for the $\Rot = 0$ case, although unfolding over a larger scale of Rayleigh numbers.  The non-trivial fixed points appear to be globally attractive for higher Rayleigh numbers, and the system enters an apparently chaotic region only after an increasingly large threshold.   
		\item On the other hand, by fixing $\Ray$ and increasing $\Rot$ one sees a reversal of the story above.  Beginning from a region with chaotic dynamics, the non-trivial fixed points stabilize as one increases the rotation, and eventually these non-trivial fixed points merge with the origin. 
	\end{enumerate}
\end{Remark}

\section*{Acknowledgements}

This paper is a contribution to the project M7 of the Collaborative Research Centre TRR 181 "Energy Transfers in Atmosphere and Ocean" funded by the Deutsche Forschungsgemeinschaft (DFG, German Research Foundation) - Projektnummer 274762653.  The author thanks Prof. Jens Rademacher for his leadership of the M7 project, Prof. Camilla Nobili and Prof. David Goluskin for sharing insights about Rayleigh-B\'enard convection, and Fabian Bleitner for discussions about PDE analysis.  Finally, many thanks to my lovely wife, Leigh Carroll.


\section*{Data Availability}

The code generated to create the figures in this work has been made available on GitHub at \url{https://github.com/rkwelter/HKC_CodeRepo}, while the data can be made available upon request to the author.

\section*{Compliance with Ethical Standards}

The author declares that he has no competing interests, and that he submits this article in accordance with \href{https://epubs.siam.org/journal/siads/editorial-policy}{SIAM editorial polices}.  

\printbibliography

@article{Chen1997Invariant,
title = {Invariant Foliations for $C^1$ Semigroups in Banach Spaces},
journal = {Journal of Differential Equations},
volume = {139},
number = {2},
pages = {283-318},
year = {1997},
issn = {0022-0396},
doi = {https://doi.org/10.1006/jdeq.1997.3255},
url = {https://www.sciencedirect.com/science/article/pii/S0022039697932551},
author = {Xu-Yan Chen and Jack K. Hale and Bin Tan}
}

@article{Chernyshenko2022Relationship,
author = {Chernyshenko, Sergei },
title = {Relationship between the methods of bounding time averages},
journal = {Philosophical Transactions of the Royal Society A: Mathematical, Physical and Engineering Sciences},
volume = {380},
number = {2225},
pages = {20210044},
year = {2022},
doi = {10.1098/rsta.2021.0044},
URL = {https://royalsocietypublishing.org/doi/abs/10.1098/rsta.2021.0044},
eprint = {https://royalsocietypublishing.org/doi/pdf/10.1098/rsta.2021.0044}
}

@article{ConstantinDoering_1996,
	title = {Variational bounds on energy dissipation in incompressible flows. III. Convection},
	author = {Doering, Charles R. and Constantin, Peter},
	journal = {Phys. Rev. E},
	volume = {53},
	issue = {6},
	pages = {5957--5981},
	numpages = {0},
	year = {1996},
	month = {June},
	publisher = {American Physical Society},
	doi = {10.1103/PhysRevE.53.5957},
	url = {https://link.aps.org/doi/10.1103/PhysRevE.53.5957}
}

@article{Constantin_1999,
	title = {Heat transport in rotating convection},
	journal = {Physica D: Nonlinear Phenomena},
	volume = {125},
	number = {3},
	pages = {275-284},
	year = {1999},
	issn = {0167-2789},
	doi = {https://doi.org/10.1016/S0167-2789(98)00252-8},
	url = {https://www.sciencedirect.com/science/article/pii/S0167278998002528},
	author = {Peter Constantin and Chris Hallstrom and Vachtang Putkaradze}
}

@article{Constantin_2001,
	title = {Logarithmic bounds for infinite Prandtl number rotating convection},
	journal = {Journal of Mathematical Physics},
	volume = {42},
	pages = {773–783},
	year = {2001},
	doi = {https://doi.org/10.1063/1.1336156},
	author = {Peter Constantin and Chris Hallstrom and Vachtang Putkaradze}
}

@article{FantuzziGoluskin,
	author = {Fantuzzi, G. and Goluskin, D. and Huang, D. and Chernyshenko, S. I.},
	title = {Bounds for Deterministic and Stochastic Dynamical Systems using Sum-of-Squares Optimization},
	journal = {SIAM Journal on Applied Dynamical Systems},
	volume = {15},
	number = {4},
	pages = {1962-1988},
	year = {2016},
	doi = {10.1137/15M1053347},
	URL = {https://doi.org/10.1137/15M1053347 },
	eprint = {https://doi.org/10.1137/15M1053347}
}

@article{FoiasManleyTemam_1987,
	title = {Attractors for the B\'enard problem: existence and physical bounds on their fractal dimension},
	journal = {Nonlinear Analysis: Theory, Methods and Applications},
	volume = {11},
	number = {8},
	pages = {939-967},
	year = {1987},
	issn = {0362-546X},
	doi = {https://doi.org/10.1016/0362-546X(87)90061-7},
	url = {https://www.sciencedirect.com/science/article/pii/0362546X87900617},
	author = {C. Foias and O. Manley and R. Temam}
}

@article{Giorgetta2018ICON,
	author = {Giorgetta, M. A. and Brokopf, R. and Crueger, T. and Esch, M. and Fiedler, S. and Helmert, J. and Hohenegger, C. and Kornblueh, L. and Köhler, M. and Manzini, E. and Mauritsen, T. and Nam, C. and Raddatz, T. and Rast, S. and Reinert, D. and Sakradzija, M. and Schmidt, H. and Schneck, R. and Schnur, R. and Silvers, L. and Wan, H. and Zängl, G. and Stevens, B.},
	title = {ICON-A, the Atmosphere Component of the ICON Earth System Model: I. Model Description},
	journal = {Journal of Advances in Modeling Earth Systems},
	volume = {10},
	number = {7},
	pages = {1613-1637},
	doi = {https://doi.org/10.1029/2017MS001242},
	url = {https://agupubs.onlinelibrary.wiley.com/doi/abs/10.1029/2017MS001242},
	eprint = {https://agupubs.onlinelibrary.wiley.com/doi/pdf/10.1029/2017MS001242},
	year = {2018}
}

@article{GluhovskyTongAgee_2002,
      author = "Alexander Gluhovsky and Christopher Tong and Ernest Agee",
      title = "Selection of Modes in Convective Low-Order Models",
      journal = "Journal of the Atmospheric Sciences",
      year = "2002",
      publisher = "American Meteorological Society",
      address = "Boston MA, USA",
      volume = "59",
      number = "8",
      doi = "https://doi.org/10.1175/1520-0469(2002)059<1383:SOMICL>2.0.CO;2",
      pages=      "1383 - 1393",
      url = "https://journals.ametsoc.org/view/journals/atsc/59/8/1520-0469_2002_059_1383_somicl_2.0.co_2.xml"
}

@article {Goluskin2018,
	AUTHOR = {Goluskin, David},
	TITLE = {Bounding averages rigorously using semidefinite programming: Mean moments of the {L}orenz system},
	JOURNAL = {J. Nonlinear Sci.},
	FJOURNAL = {Journal of Nonlinear Science},
	VOLUME = {28},
	YEAR = {2018},
	NUMBER = {2},
	PAGES = {621--651},
	ISSN = {0938-8974},
	MRCLASS = {34C29 (37A25 90C22)},
	MRNUMBER = {3770193},
	MRREVIEWER = {M. L. Blank},
	DOI = {10.1007/s00332-017-9421-2},
	URL = {https://doi.org/10.1007/s00332-017-9421-2},
}

@article{Goluskin2021,
	title = {Heat transport bounds for a truncated model of Rayleigh–Bénard convection via polynomial optimization},
	journal = {Physica D: Nonlinear Phenomena},
	volume = {415},
	pages = {132748},
	year = {2021},
	issn = {0167-2789},
	doi = {https://doi.org/10.1016/j.physd.2020.132748},
	url = {https://www.sciencedirect.com/science/article/pii/S0167278920302396},
	author = {Matthew L. Olson and David Goluskin and William W. Schultz and Charles R. Doering}
}

@article{HermizGuzdarFinn_1995,
  title = {Improved low-order model for shear flow driven by Rayleigh-B\'enard convection},
  author = {Hermiz, K. B. and Guzdar, P. N. and Finn, J. M.},
  journal = {Phys. Rev. E},
  volume = {51},
  issue = {1},
  pages = {325--331},
  numpages = {0},
  year = {1995},
  month = {January},
  publisher = {American Physical Society},
  doi = {10.1103/PhysRevE.51.325},
  url = {https://link.aps.org/doi/10.1103/PhysRevE.51.325}
}

@article{HowardKrishnamurti_1986, 
    title={Large-scale flow in turbulent convection: a mathematical model}, 
    volume={170}, 
    DOI={10.1017/S0022112086000940}, 
    journal={Journal of Fluid Mechanics}, 
    publisher={Cambridge University Press}, 
    author={Howard, L. N. and Krishnamurti, R.}, year={1986}, 
    pages={385–410}}

@article{iyer2020classical,
author = {Iyer, Kartik and Scheel, Janet and Schumacher, Joerg and Sreenivasan, Katepalli},
year = {2020},
month = {03},
pages = {201922794},
title = {Classical 1/3 scaling of convection holds up to Ra = 10 15},
volume = {117},
journal = {Proceedings of the National Academy of Sciences},
doi = {10.1073/pnas.1922794117}
}

@misc{noethen2023stepsize,
title = {Stepsize Variations for Lyapunov Exponents to Counter Persistent Errors},
year = {2023},
doi = {https://doi.org/10.48550/arXiv.2112.11388},
url = {https://arxiv.org/abs/2112.11388},
publisher = {arXiv},
author = {Florian Noethen}
}

@misc{OlsonDoering2022,
	doi = {10.48550/ARXIV.2203.02067},
	url = {https://arxiv.org/abs/2203.02067},
	author = {Olson, Matthew L. and Doering, Charles R.},
	title = {Heat transport in a hierarchy of reduced-order convection models},
	publisher = {arXiv},
	year = {2022},
	copyright = {arXiv.org perpetual, non-exclusive license}
}

@article{WenDingChiniKerswell2022,
	author = {Wen, Baole  and Ding, Zijing  and Chini, Gregory P.  and Kerswell, Rich R. },
	title = {Heat transport in Rayleigh–Bénard convection with linear marginality},
	journal = {Philosophical Transactions of the Royal Society A: Mathematical, Physical and Engineering Sciences},
	volume = {380},
	number = {2225},
	pages = {20210039},
	year = {2022},
	doi = {10.1098/rsta.2021.0039},
	
	URL = {https://royalsocietypublishing.org/doi/abs/10.1098/rsta.2021.0039},
	eprint = {https://royalsocietypublishing.org/doi/pdf/10.1098/rsta.2021.0039}
}

@article{wen_goluskin_doering_2022, title={Steady Rayleigh–Bénard convection between no-slip boundaries}, volume={933}, DOI={10.1017/jfm.2021.1042}, journal={Journal of Fluid Mechanics}, publisher={Cambridge University Press}, author={Wen, Baole and Goluskin, David and Doering, Charles R.}, year={2022}, pages={R4}}

@article{LStenflo_1996,
	doi = {10.1088/0031-8949/53/1/015},
	url = {https://dx.doi.org/10.1088/0031-8949/53/1/015},
	year = {1996},
	month = {01},
	publisher = {},
	volume = {53},
	number = {1},
	pages = {83},
	author = {L Stenflo},
	title = {Generalized Lorenz equations for acoustic-gravity waves in the atmosphere},
	journal = {Physica Scripta}
}

@article{Dhooge_2008_MatCont,
author = {A. Dhooge and W. Govaerts and Yu. A. Kuznetsov and H. G.E. Meijer and B. Sautois},
title = {New features of the software MatCont for bifurcation analysis of dynamical systems},
journal = {Mathematical and Computer Modelling of Dynamical Systems},
volume = {14},
number = {2},
pages = {147-175},
year = {2008},
publisher = {Taylor & Francis},
doi = {10.1080/13873950701742754},
URL = {https://doi.org/10.1080/13873950701742754},
eprint = {https://doi.org/10.1080/13873950701742754}
}

@article{OvsyRadeWelterLu_2023_LorenzLargeRayleigh,
	title = {Time Averages and Periodic Attractors at High Rayleigh Number for Lorenz-like Models},
	journal = {Journal of Nonlinear Science},
	volume = {33},
	year = {2023},
	doi = {10.1007/s00332-023-09933-x},
	url = {http://www.sciencedirect.com/science/article/pii/002203969290042L},
	author = {Ovsyannikov, Ivan and Rademacher, Jens D. M. and Welter, Roland and Lu, Bing-ying }
	
}

@article{Simon_1986,
	author = "Simon, Jacques",
	title = "Compact sets in the space {$ L^p(O,T; B) $}",
	journal = "Annali di Matematica Pura ed Applicata",
	year = "1986",
	volume = "146",
	doi = "10.1007/BF01762360",
	url = "https://doi.org/10.1007/BF01762360"
}

@article{SouzaDoering1,
	title = {Maximal transport in the Lorenz equations},
	journal = {Physics Letters A},
	volume = {379},
	number = {6},
	pages = {518-523},
	year = {2015},
	issn = {0375-9601},
	doi = {https://doi.org/10.1016/j.physleta.2014.10.050},
	url = {https://www.sciencedirect.com/science/article/pii/S0375960114012067},
	author = {Andre N. Souza and Charles R. Doering}
}

@book{Sparrow,
	title = {The Lorenz Equations: Bifurcations, Chaos, and Strange Attractors},
	author = {Sparrow, Colin},
	isbn = {9780387907758},
	series = {Applied Mathematical Sciences},
	year = {1982},
	DOI = {https://doi.org/10.1007/978-1-4612-5767-7},
	publisher = {Springer}
}

@article{StevensClercxLohse_2010,
	title = {Optimal Prandtl number for heat transfer in rotating Rayleigh–Bénard convection},
	journal = {New Journal of Physics},
	volume = {12},
	year = {2010},
	doi = {10.1088/1367-2630/12/7/075005},
	url = {https://iopscience.iop.org/article/10.1088/1367-2630/12/7/075005},
	author = {Richard J. A. M. Stevens and Herman J. H. Clercx and Detlef Lohse}
}

@article{ThiffeaultHorton_1996,
	title = {Energy‐conserving truncations for convection with shear flow},
	volume = {8},
	issn = {1070-6631},
	url = {https://doi.org/10.1063/1.868956},
	doi = {10.1063/1.868956},
	number = {7},
	journal = {Physics of Fluids},
	author = {Thiffeault, Jean‐Luc and Horton, Wendell},
	month = jul,
	year = {1996},
	pages = {1715--1719},
}

@article{Tucker,
	title = {The Lorenz attractor exists},
	journal = {Comptes Rendus de l'Académie des Sciences - Series I - Mathematics},
	volume = {328},
	number = {12},
	pages = {1197-1202},
	year = {1999},
	issn = {0764-4442},
	doi = {https://doi.org/10.1016/S0764-4442(99)80439-X},
	url = {https://www.sciencedirect.com/science/article/pii/S076444429980439X},
	author = {Warwick Tucker}
}

@article{WhiteheadDoering_2011,
	title = {Ultimate State of Two-Dimensional Rayleigh-B\'enard Convection between Free-Slip Fixed-Temperature Boundaries},
	author = {Whitehead, Jared P. and Doering, Charles R.},
	journal = {Phys. Rev. Lett.},
	volume = {106},
	issue = {24},
	pages = {244501},
	numpages = {4},
	year = {2011},
	month = {June},
	publisher = {American Physical Society},
	doi = {10.1103/PhysRevLett.106.244501},
	url = {https://link.aps.org/doi/10.1103/PhysRevLett.106.244501}
}

@article{Zelik_2014, 
    title={Inertial manifolds and finite-dimensional reduction for dissipative PDEs}, volume={144}, 
    DOI={10.1017/S0308210513000073}, 
    number={6}, 
    journal={Proceedings of the Royal Society of Edinburgh: Section A Mathematics}, author={Zelik, Sergey}, 
    year={2014},
    pages={1245–1327}}

@book{Chemin2006Mathematical,
	title = {Mathematical Geophysics},
	author = {Chemin, Jean-Yves and Desjardins, Benoit and Gallagher, Isabelle and Grenier, Emmanuel},
	isbn = {9780198571339},
	series = {Oxford Lecture Series in Mathematics and Its Applications},
	year = {2006},
	publisher = {Oxford University Press}
}

@book{eden2019energy,
  title={Energy Transfers in Atmosphere and Ocean},
  author={Eden, Carsten and Iske, Armin and Franzke, Christian L. E. and Oliver, Marcel and Rademacher, and Jens D. M. and Badin, Gualtiero and von Storch, Jin-Song and Olbers, Dirk and Pollmann, Friederike},
  isbn={978-3-030-05703-9},
  doi={https://doi.org/10.1007/978-3-030-05704-6},
  series={Mathematics of Planet Earth},
  url={https://link.springer.com/book/10.1007/978-3-030-05704-6},
  year={2019},
  publisher={Springer}
}

@book{haragus2011local,
author = {Haragus, Mariana and Iooss, Gerard},
year = {2011},
month = {01},
pages = {239-278},
title = {Local Bifurcations, Center Manifolds, and Normal Forms in Infinite-Dimensional Dynamical Systems},
isbn = {978-0-85729-111-0},
doi = {10.1007/978-0-85729-112-7_5}
}

@book{temam_InfDimDynSys,
	title = {Infinite-Dimensional Dynamical Systems in Mechanics and Physics},
	author = {Temam, Roger},
	series = {Applied Mathematical Sciences Series},
	isbn = {9781461268536},
	year = {2013},
	publisher = {Springer},
    DOI = {https://doi.org/10.1007/978-1-4612-0645-3}
}

\end{document}